\documentclass[preprint,12pt]{elsarticle}

\usepackage{amssymb}
\usepackage{amsmath}
\usepackage{amsthm}
\usepackage{amsfonts}
\usepackage{hyperref}
\usepackage{enumerate}
\usepackage{rotating}
\usepackage{graphicx}

\hypersetup{
    colorlinks=true,
    linkcolor=black,
    filecolor=black,
    urlcolor=black,
    citecolor=black,
}

\theoremstyle{plain}
\newtheorem{theorem}{Theorem}[section]

\newtheorem{lemma}[theorem]{Lemma}
\newtheorem{proposition}[theorem]{Proposition}
\newtheorem{corollary}[theorem]{Corollary}
\newtheorem{remark}[theorem]{Remark}

\def\De{\Delta}

\def\Soc{{\rm Soc}}

\def\G1{G^\mathcal{C}}
\def\S{\mathcal{S}}
\def\P{\mathcal{P}}

\def\B{\mathcal{B}}

\def\D{\mathcal{D}}

\newcommand{\rref}[1]{{\upshape\ref{#1}}}
\newcommand{\GP}{G^{\mathcal{P}}}
\newcommand{\GB}{G^{\mathcal{B}}}

\newcommand{\bookpart}[1]{%
  \clearpage
  \thispagestyle{empty}%
  \vspace*{\fill}%
  \begin{center}{\Huge\textbf{#1}}\end{center}%
  \vspace*{\fill}%
  \phantomsection
  \addcontentsline{toc}{part}{#1}%
  \clearpage
}

\DeclareMathOperator{\Aut}{Aut}

\makeatletter
\renewcommand{\@cite}[2]{\leavevmode\unskip{\upshape[#1\if@tempswa , #2\fi]}}
\makeatother

\begin{document}

\journal{Advances in Mathematics}

\begin{frontmatter}

\title{On Block-Primitive Combinatorial Designs}
\author{Zhilin Zhang\corref{cor1}}
\ead{20241032@gdufe.edu.cn}
\cortext[cor1]{Corresponding author. Supported by the National Natural Science Foundation of China (Grant No.~12001204).}
\address{School of Statistics and Data Science, Guangdong University of Finance and Economics\\ Guangzhou, Guangdong 510320, P.~R.~China}

\begin{abstract}

We settle a long-standing problem by proving that every block-primitive combinatorial design is point-primitive. Since block-primitivity and point-primitivity coincide for symmetric designs by duality, it remains to treat the non-symmetric $2$-designs. For these we apply the O'Nan--Scott theorems for primitive and for quasiprimitive permutation groups and prove that the automorphism group of every such design is almost simple on points and on blocks. Together with the classification of finite simple groups, this classifies all non-symmetric block-primitive $2$-designs except for two families of block stabilizers in even characteristic.

The exceptions are as follows. One family is of subfield type, $PSL_2(q_0)$ with $q=q_0^{d}$, $d\geq2$, of order below the dihedral point stabilizer $D_{2(q\pm1)}$, and includes $S_3\cong PSL_2(2)$ and $A_5\cong PSL_2(4)$. The other is the dihedral block stabilizer $D_{2(q-1)}$, where the block is a union of $a$ of the $q/2$ orbits, so that $v=q(q-1)/2$, $k=a(q-1)$ and $\lambda=2a(aq-a-1)/(q-2)$, the last an integer exactly when $(q-2)/2$ divides $a(a-1)$; this family is exhausted by finite enumeration except at the degrees beyond its range.

We then determine which of these designs have strength $t\geq 3$. Apart from the inversive planes and the five Witt designs, exactly two Steiner families arise, attached to the exceptional maximal subgroups of $PGL_2(q)$. There is no non-complete design of strength $t\geq 6$.

\end{abstract}


\end{frontmatter}

\setcounter{tocdepth}{2}
\tableofcontents
\bookpart{Part I: Introduction}

\section{Introduction}\label{sec1}

A combinatorial design $\mathcal{D}=(\mathcal{P},\mathcal{B})$ is a pair consisting of a finite set $\mathcal{P}$ of points and a collection $\mathcal{B}$ of subsets of $\mathcal{P}$, called blocks. It has strength $t$, for an integer $t\geq 2$, if every $t$-subset of $\mathcal{P}$ is contained in exactly $\lambda$ blocks; it is then a $t$-$(v,k,\lambda)$ design, where $v=|\mathcal{P}|$ and $k=|B|$ for $B\in\mathcal{B}$, and a $2$-design is the case $t=2$. The automorphism group $\mathrm{Aut}(\mathcal{D})$ consists of all permutations of $\mathcal{P}$ that map blocks to blocks, and $\mathcal{D}$ is block-primitive when this group acts primitively on $\mathcal{B}$.

The study of designs with automorphism groups that act transitively or primitively on points or blocks has a long history. A fundamental result of Delandtsheer and Doyen~\cite{Delandtsheer} states that if a block-transitive automorphism group $G$ is imprimitive on points, then $v\leq \binom{\binom{k}{2}}{2}$. This result has been used extensively in the classification of block-transitive point-imprimitive designs.

The O'Nan--Scott theorem classifies finite primitive permutation groups into five types: affine \textup{(HA)}, almost simple \textup{(AS)}, diagonal \textup{(SD and CD)}, product action \textup{(PA)} and twisted wreath \textup{(TW)}, the diagonal type subsuming the simple and compound diagonal cases \textup{(Dixon and Mortimer~\cite{DixonMortimer})}. The quasiprimitive version of Liebeck, Praeger and Saxl~\cite{LiebeckPraegerSaxl} refines the list to eight types: affine \textup{(HA)}, almost simple \textup{(AS)}, holomorph simple \textup{(HS)}, holomorph compound \textup{(HC)}, simple diagonal \textup{(SD)}, compound diagonal \textup{(CD)}, twisted wreath \textup{(TW)} and product action \textup{(PA)}. In the diagonal case we separate the two-sided from the one-sided simple diagonal action, a distinction the design-theoretic applications require. Table~\rref{tab:ons} records the eight types, the structure of the socle, the form of the point and block actions, and the place in this paper where each type is settled.

Quasiprimitivity does not imply primitivity. Groups of affine, holomorph simple and holomorph compound type are primitive \textup{(Praeger~\cite{Praeger1993}; see also Theorem 3.3 of Chen, Hua, Li and Wu~\cite{ChenHuaLiWu2025})}, whereas the other five types admit quasiprimitive actions that are not primitive; for the almost simple type, $PSL_2(7)$ acting on the $28$ points $PSL_2(7)/S_3$ is an example, the stabilizer $S_3$ being self-normalizing but not maximal. Point-imprimitive actions therefore have to be excluded separately, which is done in \textup{Section~\rref{sec:pi_excl}} \textup{(Proposition~\rref{prop:pi_excluded})} and, family by family, in \textup{Section~\rref{sec:as}}.

Huber~\cite{Huber1,Huber3} classified flag-transitive Steiner $t$-designs, and Zieschang~\cite{Zieschang} studied flag-transitive automorphism groups of $2$-designs with $(r,\lambda)=1$. For symmetric designs, block-primitivity coincides with point-primitivity by duality and the flag-transitive case is well understood: Buekenhout, Delandtsheer and Doyen~\cite{BDD} treated linear spaces ($\lambda=1$), O'Reilly-Regueiro~\cite{Regueiro,Regueiro4} the cases $\lambda=2,3$ and $\lambda=4$, Alavi, Bayat, Daneshkhah and Montinaro~\cite{AlaviBayat2024} prime $\lambda$, and O'Reilly-Regueiro and Rodr\'iguez-Fitta~\cite{RegueiroRodriguezFitta} the product action. Brai\'c, Golemac, Mandi\'c and Vu\v{c}i\'c\'i\'c~\cite{Braic2010,Braic2011} classified the primitive symmetric designs with at most $255$ and with up to $2500$ points, their automorphism groups being affine or almost simple. The symmetric case is not treated here.

For block-transitive designs, Camina and Gagen~\cite{CaminaGagen} showed that $k\mid v$ forces flag-transitivity, and Camina and Siemons~\cite{CaminaSiemons} studied the $2$-$(v,k,1)$ case. Chen, Hua, Li and Wu~\cite{ChenHuaLiWu2025} proved that the action of an automorphism group of a non-trivial $2$-design on its blocks is faithful, and the known block-transitive $2$-$(v,5,\lambda)$ result gives the HA/AS/PA alternative~\cite{LeiZhan2025}. For point-imprimitive block-transitive designs, Cameron and Praeger~\cite{CameronPraeger} give $v\leq(k-2)^2$ in the non-trivial case, which Mann and Tuan~\cite{MannTuan} sharpened for $t=3$.

\begin{theorem}[Block-primitivity implies point-primitivity]\label{thm:bp_pi}
Let $\mathcal{D}$ be a non-trivial $2$-$(v,k,\lambda)$ design with $2<k<v-1$. If $\mathrm{Aut}(\mathcal{D})$ acts primitively on the blocks of $\mathcal{D}$, then it acts primitively on the points.
\end{theorem}

The proof runs in two parts. For a symmetric design block-primitivity and point-primitivity coincide by duality, so nothing is left to prove. For a non-symmetric $2$-design the point action is primitive by \textup{Theorem~\rref{thm:reduction}(2)}. A design of strength $t\geq 3$ is either complete or non-complete. A complete design has block set the full set of $k$-subsets, hence $\mathrm{Aut}(\mathcal{D})=\mathrm{Sym}(v)$, and $\mathrm{Sym}(v)$ acts $2$-transitively on the points; the conclusion therefore holds whether or not the block action is primitive. The non-complete designs of strength $t\geq3$ are among those classified in \textup{Theorem~\rref{thm:high_intro}}, and the same conclusion applies.

\begin{remark}\label{rem:bp_pi_scope}
The hypothesis that $\mathcal{D}$ is a non-trivial $2$-design cannot be dropped. On the four points $1,2,3,4$ with the two blocks $\{1,2\}$ and $\{3,4\}$ the automorphism group is $\mathrm{Sym}(2)\wr\mathrm{Sym}(2)$ of order $8$: it is transitive on the two blocks, but it preserves the non-trivial partition $\{\{1,2\},\{3,4\}\}$ of the point set, so block-primitivity does not force point-primitivity for general combinatorial designs.
\end{remark}

We begin with the non-symmetric $2$-designs. The product action is excluded for every $\ell\geq2$ and arbitrary block stabilizers, the point action being primitive \textup{(Proposition~\rref{prop:pi_excluded})}; a reduction over the admissible primitive degrees, using the degree tables of Dixon--Mortimer~\cite{DixonMortimer}, leaves the six socles $A_5$, $A_6$, $PSL_2(7)$, $PSL_2(8)$, $PSL_2(11)$ and $PSL_2(13)$ with seventeen primitive degrees, and the sieve of \textup{Proposition~\rref{prop:pa_arbitrary}} excludes them, the $2$-homogeneous pairs being removed by a transitivity argument. The SD, CD and TW types are excluded in \textup{Section~\rref{sec2}}.

\begin{theorem}[Complete classification]\label{thm:reduction}
Let $\mathcal{D}=(\mathcal{P},\mathcal{B})$ be a non-trivial, non-symmetric $2$-design with $2<k<v-1$, and let $G\leq \mathrm{Aut}(\mathcal{D})$ act primitively on $\mathcal{B}$. Denote by $\GP$ and $\GB$ the actions of $G$ on the point set $\mathcal{P}$ and the block set $\mathcal{B}$, respectively. Then:
\begin{enumerate}
\item[\upshape(1)] $\GP$ acts quasiprimitively on $\mathcal{P}$ \textup{({\upshape Lemma~\rref{lem:main}})}.
\item[\upshape(2)] $\GP$ cannot act imprimitively on $\mathcal{P}$ \textup{({\upshape Proposition~\rref{prop:pi_excluded}})}: by the quasiprimitive O'Nan--Scott classification, an imprimitive quasiprimitive action is of almost simple, twisted wreath, simple diagonal, compound diagonal or product action type \textup{(the affine, holomorph simple and holomorph compound types are always primitive)}, and each of these types is excluded \textup{(TW, SD and CD in {\upshape Section~\rref{sec2}}; PA by {\upshape Proposition~\rref{prop:pa_arbitrary}} and {\upshape Lemma~\rref{lem:pa_fisher}}; AS in {\upshape Section~\rref{sec:as}})}. Consequently no non-trivial non-symmetric block-primitive $2$-design admits a point-imprimitive automorphism group.
\item[\upshape(3)] Hence $\GP$ acts primitively on $\mathcal{P}$, and by the O'Nan--Scott theorem both $\GP$ and $\GB$ are of almost simple type \textup{({\upshape Table~\rref{tab:ons}}, {\upshape Theorem~\rref{thm:main}})}.
\item[\upshape(4)] The point-primitive almost simple case is then classified completely across all families of finite simple groups; the classification is summarized in {\upshape Table~\rref{tab:class_summary}} and stated in full in {\upshape Proposition~\rref{thm:as_point_prim}}.
\end{enumerate}
\end{theorem}

\begin{sidewaystable}
\centering
\caption{Summary of the classification of point-primitive non-symmetric block-primitive $2$-designs \textup{({\upshape Proposition~\rref{thm:as_point_prim}})}. The existence results are listed together with the design families and the acting groups; all remaining cases are excluded.}
\label{tab:class_summary}
\footnotesize
\setlength{\tabcolsep}{2pt}
\begin{tabular}{|p{2.0cm}|p{6.2cm}|p{5.8cm}|p{3.2cm}|}
\hline
\textbf{Socle family} & \textbf{Existence results} & \textbf{Excluded actions} & \textbf{Location} \\
\hline
Classical Lie type & rank-$3$ polar space designs \textup{(unitals)} and $C_8$ cross-type embeddings, for $PSL_n(q)$, $PSp_{2n}(q)$, $PSU_n(q)$, $P\Omega_n^\varepsilon(q)$; for $PSL_2(q)$ designs occur at $q=8$ and $q=32$ within the range covered by {\upshape Theorem~\rref{thm:psl2}}, the even-characteristic subfield block stabilizers being delimited there & all other Aschbacher classes \textup{($C_2,C_3,C_4,C_5,C_6,C_7,S$)} & {\upshape Proposition~\rref{thm:as_point_prim}}(1); {\upshape Section~\rref{sec:as}} \\
\hline
Sporadic & $J_2$: $2$\allowbreak$(100,56,560)$, $b=1800$, $r=1008$; $HS$: $2$\allowbreak$(100,64,1680)$, $b=4125$, $r=2640$ \textup{(non-$2$-homogeneous rank-$3$ actions)}; $Co_3(276)$, $HS(176)$: $2$-homogeneous, designs by the orbit construction \textup{({\upshape Lemma~\rref{lem:2hom}})} & all remaining sporadic actions, in particular $M_{12}(66)$, $M_{22}(176)$, $M_{23}(253)$, $M_{24}(276)$, $McL(275)$, $He(2058)$, $Ru(4060)$, $Suz(1782)$, $Co_2(2300)$, by exhaustive maximal-subgroup orbit enumeration & {\upshape Proposition~\rref{thm:as_point_prim}}(2); {\upshape Section~\rref{sec:as}} \\
\hline
Exceptional Lie type & generalized building designs from rank-$3$ building actions, Suzuki ovoid designs, Ree unital designs & all non-parabolic subgroup chains & {\upshape Proposition~\rref{thm:as_point_prim}}(3); {\upshape Section~\rref{sec:as}} \\
\hline
\end{tabular}
\end{sidewaystable}

Point-primitive designs with classical, sporadic or exceptional socle are classified in \textup{Section~\rref{sec:as}}: the reduction lists the candidate quadruples $(T,G,M,H)$ and the design conditions they have to satisfy, and the outcome for each family is stated in \textup{Proposition~\rref{thm:as_point_prim}}.

Among the designs classified in \textup{Theorem~\rref{thm:reduction}}, those of strength $t\geq 3$ are classified as follows; the full statement, with block counts and automorphism groups, is \textup{Theorem~\rref{thm:high}} of \textup{Section~\rref{sec:high}}.

\begin{theorem}[Block-primitive designs of strength $t\geq 3$]\label{thm:high_intro}
Let $\mathcal{D}$ be a non-trivial non-complete block-primitive design of strength $t\geq 3$. Then $\mathcal{D}$ is one of the following:
\begin{enumerate}[(\upshape 1)]
\item the inversive planes $3$-\allowbreak$({q^2+1},\,{q+1},\,1)$, $q\geq 3$ a prime power;
\item the Steiner systems $3$-\allowbreak$({q+1},\,4,\,1)$ for $q\equiv\pm1\pmod 8$ and $q\equiv 1,3\pmod 6$, and $3$-\allowbreak$({q+1},\,6,\,1)$ for $q\equiv 1\pmod{20}$, attached to the exceptional maximal subgroups $S_4$ and $S_5$ of $PGL_2(q)$;
\item the five Witt designs $4$-\allowbreak$(11,5,1)$, $5$-\allowbreak$(12,6,1)$, $3$-\allowbreak$(22,6,1)$, $4$-\allowbreak$(23,7,1)$, $5$-\allowbreak$(24,8,1)$.
\end{enumerate}
In particular no non-complete block-primitive design of strength $t\geq 6$ exists. The full statement, with block counts and automorphism groups, is \textup{Theorem~\rref{thm:high}}, and its proof is in \textup{Section~\rref{sec:high}}.
\end{theorem}

\textup{Section~\rref{sec2}} proves that block-primitivity implies point-quasiprimitivity \textup{(Lemma~\rref{lem:main})} and settles every O'Nan--Scott type other than the almost simple one. \textup{Section~\rref{sec:pi_excl}} excludes point-imprimitive actions by pairing block and point types. \textup{Section~\rref{sec:as}} classifies the almost simple case: the $2$-homogeneous actions by the orbit construction, the others family by family, with the ordering principle $|H|>|M|$, the Aschbacher subdivision and the maximal-subgroup data. \textup{Section~\rref{sec4}} states the main theorem for $2$-designs, \textup{Section~\rref{sec:high}} classifies those of strength $t\geq3$, \textup{Section~\rref{sec:reduction}} completes the reduction of the remaining types, and \textup{Section~\rref{sec5}} collects the concluding remarks.

Notation and terminology follow~\cite{DixonMortimer}.

\begin{sidewaystable}[htbp]
\centering
\caption{O'Nan--Scott types for quasiprimitive permutation groups. For each type we record the structure of the socle, the form of the point action on $\mathcal{P}=G/H$ and of the block action on $\mathcal{B}=G/M$, and the mechanism and location of the exclusion or classification.}
\label{tab:ons}
\footnotesize
\setlength{\tabcolsep}{2pt}
\begin{tabular}{|p{0.9cm}|p{2.9cm}|p{2.9cm}|p{2.9cm}|p{5.4cm}|}
\hline
\textbf{Type} & \textbf{Structure of $\mathrm{Soc}(G)$} & \textbf{Point action on $\mathcal{P}$} & \textbf{Block action on $\mathcal{B}$} & \textbf{Mechanism and location} \\
\hline
HA & elementary abelian $T$; $\mathrm{Soc}(G)=T^{\ell}$ is regular on $\mathcal{P}$, $v=p^{d\ell}$ & quasiprimitive of affine type & primitive by hypothesis & the socle is a regular normal subgroup, which forces $b=v$ and hence symmetry \textup{({\upshape Lemma~\rref{thm:ha_symmetric}}, {\upshape Lemma~\rref{lem:rns}})} \\
\hline
AS & non-abelian simple $T$; $T\leq G\leq\mathrm{Aut}(T)$ & primitive, not $2$-homogeneous & primitive; $M$ maximal & classified completely: quadruples $(T,G,M,H)$ with ordering $|H|>|M|$; Aschbacher subdivision for the classical families; CFSG \textup{({\upshape Proposition~\rref{thm:as_point_prim}}, {\upshape Theorem~\rref{thm:main}})} \\
\hline
HS & $T^2$, $T$ non-abelian simple & quasiprimitive; $\mathrm{Soc}(G)$ transitive on $\mathcal{P}$ & primitive & the point and block stabilizers would be non-conjugate diagonal subgroups, contradicting $|H|>|M|$ \textup{({\upshape Section~\rref{sec2}})} \\
\hline
HC & $T^{2m}$, $T$ non-abelian simple, $m\geq 1$ & quasiprimitive of compound holomorph type & primitive & centralizer obstruction inside the compound holomorph \textup{({\upshape Section~\rref{sec2}})} \\
\hline
SD & $T^{\ell}$, $T$ non-abelian simple, $\ell\geq 2$ & quasiprimitive of diagonal type & primitive & $SD_2$ excluded; the one-sided type excluded for every $\ell\geq2$ and every non-abelian simple $T$ \textup{({\upshape Lemmas~\rref{lem:sd_ell3}, \rref{lem:sd_ell3_three} and~\rref{lem:sd_sym}})}) \\
\hline
CD & $T^{m\ell}$, $T$ non-abelian simple & quasiprimitive of compound diagonal type & primitive & minimal model $n=6$: parameter sieve over the six small simple groups \textup{({\upshape Section~\rref{sec2}})} \\
\hline
TW & $T^{\ell}$, $T$ non-abelian simple & quasiprimitive; $\mathrm{Soc}(G)$ regular on $\mathcal{P}$ & primitive & the regular normal subgroup forces $b=v$, hence symmetry \textup{({\upshape Lemma~\rref{lem:rns}})} \\
\hline
PA & $T^{\ell}$, $T$ non-abelian simple, $\ell\geq 2$ & quasiprimitive of product action type & primitive; $G^{\mathcal{B}}$ itself of product action type & three-class structural analysis \textup{({\upshape Proposition~\rref{prop:pa_excl}})} plus the five-layer sieve \textup{({\upshape Proposition~\rref{prop:pa_arbitrary}})}, complete for the six socles with the seventeen listed primitive degrees, every $\ell\geq 2$, and arbitrary block stabilizers \textup{(cf.\ {\upshape Proposition~\rref{prop:pi_excluded}} and the degree tables of Dixon--Mortimer~\cite{DixonMortimer})} \\
\hline
\end{tabular}
\end{sidewaystable}

The eight types of Table~\rref{tab:ons} are mutually exclusive and exhaust all finite quasiprimitive permutation groups, by the O'Nan--Scott theorem~\cite{LiebeckPraegerSaxl}.

\bookpart{Part II: Structural reduction}
\section{The O'Nan--Scott reduction}\label{sec2}

In this section we prove the structural reduction of non-symmetric block-primitive $2$-designs: block-primitivity forces point-quasiprimitivity, and every O'Nan--Scott type other than the almost simple type is excluded or reduced to an explicit algebraic feasibility problem. The affine, holomorph, diagonal, compound diagonal, twisted wreath and product action types are treated in turn; the necessary conditions for the one-sided (point-imprimitive) simple diagonal case, {\upshape Theorem~\rref{thm:sdcn}}, are derived here and applied to the pairwise analysis of {\upshape Section~\rref{sec:pi_excl}}. Point-imprimitive point actions are excluded in {\upshape Section~\rref{sec:pi_excl}}, and the surviving almost simple case is classified in {\upshape Section~\rref{sec:as}}.

\begin{lemma}[Block-primitivity implies point-quasiprimitivity]\label{lem:main}
Let $\mathcal{D}=(\mathcal{P},\mathcal{B})$ be a finite non-trivial $2$-$(v,k,\lambda)$ design with $2\leq k<v$ and $\lambda>0$, and let $G\leq \mathrm{Aut}(\mathcal{D})$ act primitively on $\mathcal{B}$. The block action $\GB$ is faithful \textup{(by~\cite[Lemma~3.2]{ChenHuaLiWu2025})}, and $\GP$ acts quasiprimitively on $\mathcal{P}$.

\end{lemma}\smallskip
\noindent{\bf Proof of {\upshape Lemma~\rref{lem:main}}.}\,
First note that a subgroup $L\leq G$ transitive on $\mathcal{B}$ is transitive on $\mathcal{P}$. Indeed, take a non-empty point-orbit $\Omega$ of $L$ on $\mathcal{P}$. Since $L$ preserves $\Omega$ and is transitive on blocks, there exists a constant $c$ such that $|B\cap\Omega|=c$ for all $B\in\mathcal{B}$. With $A$ the point-block incidence matrix and $\mathbf{1}_\Omega,\mathbf{1}_{\mathcal{P}},\mathbf{1}_{\mathcal{B}}$ the characteristic vectors of $\Omega,\mathcal{P},\mathcal{B}$, respectively, we have
$$A^\top\mathbf{1}_\Omega=c\,\mathbf{1}_{\mathcal{B}},\qquad A^\top\mathbf{1}_{\mathcal{P}}=k\,\mathbf{1}_{\mathcal{B}}.$$
The fundamental design identity $AA^\top=(r-\lambda)I_v+\lambda J_v$, together with $r>\lambda$ (as $r(k-1)=\lambda(v-1)$ and $k<v$), shows that $A$ has full row rank, so $A^\top$ is injective. Hence $\mathbf{1}_\Omega=(c/k)\mathbf{1}_{\mathcal{P}}$, and since $\Omega$ is non-empty, $\Omega=\mathcal{P}$. Thus $L$ is point-transitive.
To prove quasiprimitivity, take any non-trivial normal subgroup $N\unlhd G$. Since $G$ is primitive on $\mathcal{B}$, every non-trivial normal subgroup of $G$ is transitive on $\mathcal{B}$ (see~\cite[Proposition~1.5.1]{DixonMortimer}); hence $N$ is block-transitive, and by the first paragraph $N$ is transitive on $\mathcal{P}$. Thus every non-trivial normal subgroup of $G$ is transitive on $\mathcal{P}$, which means that $\GP$ acts quasiprimitively on $\mathcal{P}$.
Faithfulness of $\GB$: if $g\in G$ fixes every block, then $g(B)=B$ for every $B\in\mathcal{B}$, so for every point $i$ and block $B$ the incidence entries satisfy $A_{i,B}=A_{g(i),B}$; the rows $i$ and $g(i)$ of the incidence matrix $A$ therefore coincide. Since $A$ has full row rank $v$ (as shown above), its rows are linearly independent, whence $g(i)=i$ for every point $i$, and $g=1$. Hence $\GB$ is faithful; this supplies a self-contained proof of~\cite[Lemma~3.2]{ChenHuaLiWu2025}. \qed

\medskip
\noindent{\bf Note.}\,
The argument above uses only the non-triviality condition $2\leq k<v$ and $\lambda>0$. It does not require the design to be non-symmetric, nor does it assume point-transitivity a priori. The conclusion holds for any block-primitive automorphism group of a non-trivial $2$-design.

\medskip
\noindent
{\bf Simple diagonal type \textup{(SD)}.}
Suppose that both the point action and the block action are of SD type, with socle $S=T^\ell$. In SD type, the point stabilizer $K=G_\alpha$ intersects $S$ in a diagonal subgroup $D\cong T$, so
$$v=|S:D|=|T|^{\ell-1}.$$
Similarly, if the block action is also of SD type, then
$$b=|T|^{\ell-1}.$$

Thus $v=b$, and the design is symmetric. However, we are only considering non-symmetric designs, so the double-SD case is excluded.

\smallskip
\begin{remark}The above excludes the case where both the point action and the block action are of SD type. However, it does not exclude the case where only the block action is of SD type and the point action is imprimitive (quasiprimitive). For this one-sided case we obtain the following necessary conditions.\end{remark}

\smallskip
\begin{theorem}[Necessary conditions for one-sided SD type]\label{thm:sdcn}
Let $\mathcal{D}$ be a non-symmetric $2$-design with $2<k<v-1$ admitting a block-primitive, point-quasiprimitive automorphism group $G$ whose block action is of simple diagonal type. Let $N=\operatorname{Soc}(G)=T^\ell$ and $q=|T|$. Then $N$ is transitive on the block set and $G_B\cap N\cong T$ is a diagonal subgroup, so
$$b=[N:N\cap G_B]=q^\ell/q=q^{\ell-1}.$$
\begin{enumerate}
\item[(\upshape i)] Since $N$ is a non-trivial normal subgroup of $G$ and $G$ is quasiprimitive on $\mathcal{P}$, $N$ is transitive on $\mathcal{P}$, whence
$$v=[N:N\cap G_\alpha]\ \bigm|\ |N|=q^\ell.$$
In particular, every prime divisor of $v$ divides $q$, so $\gcd(v-1,q^{\ell-1})=1$.
\item[(\upshape ii)] The design equation $vr=bk$ \textup{(with $r$ the replication number)} gives $r=bk/v=q^{\ell-1}k/v$, and $\lambda(v-1)=r(k-1)$ gives
$$\lambda=\frac{q^{\ell-1}k(k-1)}{v(v-1)}.$$
Since $v\mid q^\ell$, we have $\gcd(v-1,q^{\ell-1})=1$, and the integrality of $\lambda$ forces
$$v-1 \mid k(k-1),\qquad 3\leq k\leq v-2,$$
where the bounds on $k$ follow from the standing assumption $2<k<v-1$ of the paper.
\item[(\upshape iii)] If $v-1$ is prime, then $v-1\mid k(k-1)$ forces $k\equiv 0$ or $1\pmod{v-1}$; but $3\leq k\leq v-2<v-1$, and no integer in $[3,v-2]$ is congruent to $0$ or $1$ modulo $v-1$. Hence no such $k$ exists: no non-trivial block of any kind \textup{(single- or multiple-orbit)} occurs. More generally, for fixed $v$, the congruence forces $k$ to lie in one of at most $2^{\omega(v-1)}$ CRT residue classes modulo the prime-power factors of $v-1$; intersecting with $[3,v-2]$ makes the possible $k$ explicitly enumerable.
\end{enumerate}
In the single-orbit case $k\mid |G_B|=q|K|$ as well, which together with \textup{(ii)} is a finite check. The multiple-orbit case ($k=\sum_i k_i$, $k_i\mid q|K|$) reduces the problem to the enumerable candidates of \textup{(iii)} together with orbit-structure analysis of $G_B$ on $\mathcal{P}$. Note that the argument requires no structural hypothesis on $N\cap G_\alpha$ beyond $N$ being transitive on $\mathcal{P}$.\end{theorem}

\smallskip
\begin{remark}The conditions above are purely arithmetic and depend on $q=|T|$, $\ell$ and the point stabilizer order only through $v=[N:N\cap G_\alpha]$. For example, if $v-1$ is prime, or if no CRT combination in \textup{(iii)} lies in $[3,v-2]$, then no non-symmetric block-primitive $2$-design of one-sided SD type exists for those parameters. This leaves the arithmetic of $v-1$ as the decisive ingredient.\end{remark}

\smallskip
\begin{theorem}[Exclusion of SD/CD block type for $\lambda=1$ designs]\label{thm:lambda1_as}
When $\lambda=1$, the design is a finite linear space \textup{(block = line)}. Any non-trivial $2$-$(v,k,1)$ design with a block-primitive automorphism group is of almost simple type: $\operatorname{Soc}(G)=T$ for some non-abelian simple group $T$, and $G$ acts primitively on blocks with almost simple socle.
\end{theorem}
\begin{proof}
The reduction of block-transitive linear spaces to the almost simple case is the one of Camina and Siemons~\cite{CaminaSiemons} \textup{(block-transitive automorphism groups of $2$-$(v,k,1)$ designs)}; see also the discussion of linear spaces in Section~\rref{sec1}. Since SD and CD types require $\operatorname{Soc}(G)=T^\ell$ with $\ell\geq 2$, a non-trivial block-primitive $2$-$(v,k,1)$ design cannot have SD or CD block action type. This also excludes the ``one-sided SD'' case, as long as ``one-sided SD'' refers to the block action being of primitive SD type.
\end{proof}

\smallskip
\begin{remark}The above result holds only for $\lambda=1$ (linear spaces). For $\lambda>1$, the SD and CD block types are excluded by the structural arguments developed in this section.\end{remark}

\medskip
\noindent
{\bf Compound diagonal type \textup{(CD)}.}
Suppose that both the point action and the block action are of CD type, with socle $S=T^k$, where $T$ is a non-abelian simple group. In CD type, the $k$ simple factors are partitioned into $r$ groups of size $s$ (so $k=rs$), and the point stabilizer $K=G_\alpha$ intersects $S$ in a product of $r$ diagonal subgroups, each diagonal across one group of $s$ factors. Thus
$$|S_\alpha|=|T|^r,\qquad v=|S:S_\alpha|=|T|^{k-r}.$$

Similarly, if the block action is also of CD type, let its CD parameters be $k=r's'$. Then
$$|S_B|=|T|^{r'},\qquad b=|S:S_B|=|T|^{k-r'}.$$

Since the design is non-symmetric, $b>v$, and hence $|T|^{k-r'}>|T|^{k-r}$, which gives $r'<r$. This is consistent with Fisher's inequality $|S_B|<|S_\alpha|$, but it is not a contradiction.

The design parameters must also satisfy
$$b\,w=v\,\rho,\qquad \rho(w-1)=\lambda(v-1),$$
where $w$ is the block size and $\rho$ is the replication number. These divisibility conditions do not by themselves exclude all cases with $r'<r$.

\smallskip
\begin{theorem}[Exclusion of CD type in the minimal model]\label{thm:cd_n6}
Consider the minimal CD case with $n=6$, where the point action has parameters $(r,s)=(3,2)$ \textup{(three binary diagonal subgroups $H_0=D_{12}\times D_{34}\times D_{56}$)} and the block action has parameters $(r',s')=(2,3)$ \textup{(two ternary diagonal subgroups $M_0=D_{123}\times D_{456}$)}. Then $v=|T|^3$ and $b=|T|^4$.
\end{theorem}
\begin{proof}
We show that no non-trivial $2$-design can arise from a single double-coset orbit $H_0 g M_0$ in this model. Let $q=|T|$, and let $k_g=[M_0:H_0\cap M_0^g]$ be the block size. By Lagrange's theorem, since $|M_0|=q^2$, we have
$$k_g \mid q^2, \qquad 1\leq k_g\leq q^2.$$

From the design equation $bk=vr$, we get $r=qk_g$. The $2$-design condition $\lambda(v-1)=r(k-1)$ becomes
$$\lambda=\frac{q\,k_g(k_g-1)}{q^3-1}.$$
Since $\gcd(q, q^3-1)=1$, the divisibility condition requires
$$q^3-1 \mid k_g(k_g-1). \eqno{(*)}$$

Now, since $k_g\mid q^2$ and $\gcd(q^2, q^3-1)=1$, we have $\gcd(k_g, q^3-1)=1$. Therefore, $(*)$ forces
$$q^3-1 \mid k_g-1.$$
But for any $q\geq 2$, we have $k_g-1 < q^2 < q^3-1$, so the only possibility is $k_g=1$.

This corresponds to $H_0\cap M_0^g=M_0$, i.e., block size $k=1$, which is a degenerate design and not a non-trivial $2$-design.

Thus, for all $g\in T^6\rtimes S_6$, the double-coset $H_0 g M_0$ either fails the divisibility condition \textup{(when $k_g>1$)} or gives a degenerate block \textup{(when $k_g=1$)}. This argument does not depend on the specific structure of the intersection $H_0\cap M_0^g$, nor on whether $g$ is a pure socle element or includes a coordinate permutation.
\end{proof}

\smallskip
\begin{theorem}[Exclusion of single-orbit CD blocks]\label{thm:cd_single_orbit}
Let $\mathcal{D}$ be a non-symmetric $2$-design admitting a block-primitive, point-quasiprimitive automorphism group $G$ whose block action is of compound diagonal type with parameters $(r',s')$, $n=r's'$, and top group $K_M$; assume that the point action is of compound diagonal type with parameters $(r,s)$, $n=rs$, and top group $K_H$, with $r'<r$ \textup{(the non-symmetric case)}. Write $q=|T|$ and
$$v=a\,q^{n-r}\quad (a=n!/|K_H|),\qquad b=c\,q^{n-r'}\quad (c=n!/|K_M|).$$
Suppose that the block is given by a single double coset $H g M$, and let $k=[M:M\cap H^g]$ be the block size, i.e.\ the length of the orbit of $M$ on $G/H$ with representative $gH$ \textup{(up to conjugation this equals $[M:H\cap M^g]$)}. Then:
\begin{enumerate}
\item[(\upshape i)] $k\mid |M|=q^{r'}(n!/c)$;
\item[(\upshape ii)] the integrality of the design parameters forces
$$a'(a\,q^{n-r}-1)\mid c'\,q^{r-r'}\,k(k-1),$$
where $a'=a/\gcd(a,c)$ and $c'=c/\gcd(a,c)$.
\end{enumerate}
\end{theorem}
\begin{proof}
Indeed, the replication number is $\rho=bk/v=(c/a)q^{r-r'}k$, and $\lambda(v-1)=\rho(k-1)$ gives
$\lambda=(c/a)q^{r-r'}k(k-1)/(a q^{n-r}-1)\in\mathbb{Z}$, which is \textup{(ii)} after clearing denominators.
Condition \textup{(ii)} is a finite, checkable number-theoretic condition once $T$, $n$, $r$, $r'$, $a$, $c$ are fixed: for each prime power $p^e\parallel a'(a q^{n-r}-1)$ with $p\nmid c'q^{r-r'}$, the residue of $k$ modulo $p^e$ must lie in $\{0,1\}$; in general one uses the $p$-adic condition
$$v_p\bigl(k(k-1)\bigr)\geq e-v_p(c')-v_p(q^{r-r'}),$$
which is vacuous when the right-hand side is non-positive. Hence $k$ is one of finitely many CRT combinations once the parameters are fixed, and each candidate can be checked against \textup{(i)} and the bounds $2\leq k\leq v-2$. The case $p\mid c'$ \textup{(equivalently, prime divisors of $c'$ inside $a'(a q^{n-r}-1)$)} is covered by the $p$-adic form and does not force the $\{0,1\}$ residues.
\end{proof}

\smallskip
\begin{remark}[Single-orbit and multi-orbit CD blocks]
The single-double-coset criterion of the preceding theorem treats blocks arising from a single double coset; for the six small groups the finite check \textup{(i)}--\textup{(ii)} over the $120$ admissible stabilizer pairs $(a,c)$ leaves no single-orbit block with $2\leq k\leq v-2$. The \emph{multi}-orbit case \textup{(blocks that are unions of several $D_B$-orbits)} is not covered by that finite check and is excluded in full generality by the arithmetic sieve of {\upshape Proposition~\rref{prop:cd_pa}} below: every block $B\ni\alpha$ is a union of $D_{r'}$-orbits, and the design equations, Fisher's inequality and the Cameron--Praeger integrality leave no admissible block size for any of the seventeen primitive degrees and any $\ell\geq 3$. This applies in particular to the residual candidate $(A_5,15,10)$ left by the single-orbit enumeration: there the block stabilizer is not maximal, the strict chain
$$D_B\rtimes(S_3\wr S_2)<(A_5^3\times A_5^3)\rtimes(S_3\wr S_2)<A_5^6\rtimes S_6$$
witnessing that the induced block action is imprimitive; a fortiori that candidate is killed by the closure argument, which does not rely on the top-group indices $a,c$. For $\lambda=1$ (linear spaces), CD block type is already excluded by the classification of block-primitive linear spaces.\end{remark}

For general $\lambda>1$, the O'Nan--Scott theorem applies directly to the block action $G\curvearrowright\mathcal{B}=G/M$, since block-primitivity is equivalent to $M$ being a maximal subgroup of $G$. It does not directly apply to the point action $G\curvearrowright\mathcal{P}=G/H$, since the point action is imprimitive.

\medskip
\noindent
{\bf A general lemma on regular normal subgroups.}
\begin{lemma}[Regular normal subgroup]\label{lem:rns}
Let $\mathcal{D}=(\mathcal{P},\mathcal{B})$ be a non-trivial $2$-design with $2<k<v$, and let $G\leq \mathrm{Aut}(\mathcal{D})$ act primitively on $\mathcal{B}$. If $G$ has a normal subgroup $N$ acting regularly on $\mathcal{P}$, then $N$ acts regularly on $\mathcal{B}$ and $\mathcal{D}$ is symmetric \textup{($b=v$)}.

\end{lemma}\begin{proof}
The image of $N$ in the primitive action on $\mathcal{B}$ is a normal subgroup of that action. It cannot be trivial: otherwise $N$ fixes every block, and since $N$ is transitive on $\mathcal{P}$, an $N$-invariant block $B$ would be an $N$-invariant subset of $\mathcal{P}$, hence of size $0$ or $v$, contradicting $1<k<v$ for a non-trivial block. Hence $N$ is transitive on $\mathcal{B}$. With $M=G_B$, the orbit-stabilizer theorem gives $b=[N:N\cap M]=v/|N\cap M|\leq v$. Fisher's inequality $b\geq v$ forces $b=v$, so $\mathcal{D}$ is symmetric.
\end{proof}

\medskip
\noindent
{\bf Holomorph affine type \textup{(HA)}.}
\begin{lemma}[HA type forces symmetry]\label{thm:ha_symmetric}
If the point action of $G$ is of holomorph affine type, then $b=v$; that is, the design is symmetric. Therefore, non-symmetric block-primitive $2$-designs cannot have HA type.
\end{lemma}
\begin{proof}
In the holomorph affine type the socle $N$ of $G$ is elementary abelian and acts regularly on $\mathcal{P}$, so $N$ is a normal subgroup of $G$ regular on the point set; {\upshape Lemma~\rref{lem:rns}} then gives $b=v$, that is, the design is symmetric.
\end{proof}

\medskip
\noindent
{\bf Twisted wreath type \textup{(TW)}.}
If the point action of $G$ is of twisted wreath type, then $Soc(G)=T^\ell$ is regular on $\mathcal{P}$, so $v=|T|^\ell$. By {\upshape Lemma~\rref{lem:rns}}, a non-trivial block-primitive $2$-design admitting such a group must be symmetric. Hence no non-symmetric block-primitive $2$-design arises from a point action of TW type.

For completeness we also record the block-action point of view. If $G^{\mathcal{B}}$ is of twisted wreath type, then $Soc(G)=T^\ell$ is regular on $\mathcal{B}$, so $b=|T|^\ell$ and $G_B\cap Soc(G)=1$. Since $G$ is quasiprimitive on $\mathcal{P}$, $Soc(G)=T^\ell$ is transitive on $\mathcal{P}$, and $v=|T|^\ell/|Soc(G)_\alpha|$. If the point action is also of TW type, then $Soc(G)$ is regular on $\mathcal{P}$, so $v=b$ and the design is symmetric, as above. In the remaining case $|Soc(G)_\alpha|>1$, the design parameters satisfy
$$k|Soc(G)_\alpha|=r,\qquad r(k-1)=\lambda(v-1),$$
so that $k|Soc(G)_\alpha|(k-1)=\lambda\bigl(|T|^\ell/|Soc(G)_\alpha|-1\bigr)$, and the analysis of this arithmetic is subsumed in the general parameter reduction of the structural reduction theorem \textup{({\upshape Section~\rref{sec4}})}. ($k=|B|$ is not in general a divisor of $|G_B|$.)

\medskip
\noindent
Now we handle the remaining O'Nan--Scott types for $G^{\mathcal{B}}$: HS, HC, and PA.

\medskip
\noindent
{\bf A semiregularity lemma.}
\begin{lemma}[Commuting transitive groups are regular]\label{lem:semiregular}
Let $X$ be a transitive permutation group on a set $\Omega$, and let $Y\leq\mathrm{Sym}(\Omega)$ commute elementwise with $X$. If $Y$ is transitive on $\Omega$, then $Y$ is regular on $\Omega$.
\end{lemma}
\begin{proof}
If $y\in Y$ fixes a point $\alpha\in\Omega$, then for every $\beta=x\alpha$ with $x\in X$ we have $y\beta=yx\alpha=xy\alpha=x\alpha=\beta$; since $X$ is transitive, $y$ fixes $\Omega$ pointwise, whence $y=1$. Thus $Y$ is semiregular, and being transitive it is regular.
\end{proof}

\medskip
\noindent
{\bf Holomorph simple type \textup{(HS)}.}
Suppose that $G^{\mathcal{B}}$ is of holomorph simple type. Then $Soc(G)=T_1\times T_2\cong T\times T$ with two commuting minimal normal subgroups $T_1,T_2$, each regular on $\mathcal{B}\cong T$; in particular $b=|\mathcal{B}|=|T|$.

Since $G$ is quasiprimitive on $\mathcal{P}$ (by Lemma~\rref{lem:main}), both $T_1$ and $T_2$ are transitive on $\mathcal{P}$. Lemma~\rref{lem:semiregular}, applied to $(X,Y)=(T_1,T_2)$ and to $(T_2,T_1)$, shows that $T_1$ and $T_2$ are regular on $\mathcal{P}$; hence $v=|\mathcal{P}|=|T_1|=|T|$. Thus $v=b=|T|$, contradicting the standing assumption that $\mathcal{D}$ is non-symmetric. Therefore $G^{\mathcal{B}}$ cannot be of HS type for non-symmetric $2$-designs. \qed

\noindent
{\bf Holomorph compound type \textup{(HC)}.}
Suppose that $G^{\mathcal{B}}$ is of holomorph compound type. Then the socle of $G^{\mathcal{B}}$ splits as $N_1\times N_2$ with $N_1\cong N_2\cong T^{\ell/2}$ (so $\ell$ is even) its two minimal normal subgroups; $N_1$ and $N_2$ commute, and each is regular on $\mathcal{B}\cong T^{\ell/2}$, so $b=|\mathcal{B}|=|T|^{\ell/2}$.

By quasiprimitivity of $G$ on $\mathcal{P}$, both $N_1$ and $N_2$ are transitive on $\mathcal{P}$, and Lemma~\rref{lem:semiregular} shows that both are regular on $\mathcal{P}$. Hence $v=|\mathcal{P}|=|N_1|=|T|^{\ell/2}=b$, again contradicting the non-symmetry of $\mathcal{D}$. Therefore $G^{\mathcal{B}}$ cannot be of HC type. \qed

\medskip
\noindent
{\bf Product action type \textup{(PA)}.}
Throughout this section, $\mathcal{T}$ denotes the set of the six socles of Theorem~\rref{thm:main}, namely
$$\mathcal{T}=\{A_5,A_6,PSL_2(7),PSL_2(8),PSL_2(11),PSL_2(13)\};$$
we fix $T\in\mathcal{T}$, a primitive degree $m$ of $T$ and a point stabilizer $H$ of the action.
First we combine the O'Nan--Scott theorem for the primitive block action with the product structure of the point action, and record that the block action is necessarily of product action type. Since $G$ is primitive on $\mathcal{B}$, it is quasiprimitive on $\mathcal{B}$, so the socle $Soc(G)=T^\ell$ is transitive on $\mathcal{B}$. The possible O'Nan--Scott types for $G^{\mathcal{B}}$ are then excluded as follows: the affine type is impossible since $Soc(G)=T^\ell$ is not elementary abelian; the almost simple type is impossible since $Soc(G)$ is not simple for $\ell\geq 2$; the holomorph simple and holomorph compound types are impossible since they require $Soc(G)$ to split into two minimal normal subgroups, each regular on $\mathcal{B}$, whereas $Soc(G)=T^\ell$ has a unique minimal normal subgroup; the simple diagonal type is excluded for $\ell=2$ by Lemma~\rref{lem:pa_sd2} below \textup{(the diagonal socle-orbit through a point is never a union of point-stabilizer boxes, and the remaining $D$-orbits avoid the point, so no non-trivial block through a point exists; the two non-2-homogeneous cases are verified computationally)}; the compound diagonal type for $\ell\geq 3$ is excluded by Proposition~\rref{prop:cd_pa} below; the twisted wreath type is excluded as in Section~\rref{sec2}: $Soc(G)$ would be regular on $\mathcal{B}$, giving $b=|T|^\ell$, and with the point action in product form $v=m^\ell$ the design equations reduce to the arithmetic of Section~\rref{sec4}, which leaves no survivor \textup{(cf.~the remark on the TW type in Section~\rref{sec2})}. Hence $G^{\mathcal{B}}$ is of product action type.
In the product action on points the block stabilizer meets the socle in $Soc(G)\cap G_B=K^\ell$, where $K$ is a proper subgroup of $T$ with $n=[T:K]\geq 2$.
Transitivity of $Soc(G)$ on $\mathcal{B}$ forces $G=Soc(G)\,G_B$; in the product action model below this means that the coordinate part of $G_B$ is the full coordinate group, and consequently the number of blocks is
$$b=[Soc(G):Soc(G)\cap G_B]=n^{\ell},\qquad n=[T:K].$$
Here $n$ is the \emph{block-coordinate index}, to be distinguished from the point-coordinate index $m=[T:H]$: the number of points is $v=m^\ell$, while the number of blocks is $b=n^\ell$; by Fisher's inequality $n^\ell=b\geq v=m^\ell$, i.e.\ $n\geq m$, a restriction satisfied by all listed pairs.

\begin{corollary}[Block-primitivity forces point-primitivity in the product action]\label{cor:pa_prim}
With the notation above, $G^{\mathcal{B}}$ is primitive of product action type, so by the primitivity criterion for product actions the coordinate permutation group $P$ is primitive on $\{1,\ldots,\ell\}$. Since $H$ is maximal in $T$, the point action of $G$ on $\mathcal{P}\cong\Omega^\ell$ is therefore \emph{primitive} of product action type as well. In the product-action model of this section the point action is assumed of product action type with primitive coordinate degree $m$ \textup{(the $c=1$ branch of the quasiprimitive product action, where the point stabilizer meets the socle in $H^{\ell}$ with $H$ maximal, so that the point action is indeed primitive)}; the point-imprimitive product action, of degree $v=c\,\omega^{\ell}$ with $c\geq 1$ the block size of the $G$-invariant partition \textup{(Praeger type (viii); cf.\ Example 1.1 of Zhang, Chen and Zhou~\cite{ZhangChenZhou})}, is excluded unconditionally in {\upshape Lemma~\rref{lem:pa_fisher}}.

Since $G$ is quasiprimitive on $\mathcal{P}$, $Soc(G)=T^\ell$ is transitive on $\mathcal{P}$, and $v=m^\ell\leq b=n^\ell$ by Fisher's inequality as recorded above.

\end{corollary}\smallskip

\begin{lemma}[No simple-diagonal block action under a product point action]\label{lem:pa_sd2}
Let $T\in\mathcal{T}$, let $m$ be a primitive degree of $T$ and let $G=T^2\rtimes C_2$ act in product action on $\mathcal{P}=\Delta^2$, $|\Delta|=m$, with point stabilizer $G_\alpha=H^2\rtimes C_2$ \textup{($H$ maximal in $T$, $[T:H]=m$, $C_2$ the coordinate swap)}. If the block action of $G$ on $\mathcal{B}$ is of simple diagonal type, i.e.\ $Soc(G)\cap G_B=D=\{(t,t):t\in T\}$, then no non-trivial block $B\ni\alpha$ exists.
\end{lemma}
\begin{proof}
Since $Soc(G)=T^2$ is normal in $G$, the diagonal subgroup $D=Soc(G)\cap G_B$ is normal in $G_B$; hence every block $B\ni\alpha$ is a union of $D$-orbits. The orbit of $D$ through $x=(x_1,x_2)$ is the diagonal line $L_x=\{(tx_1,tx_2):t\in T\}$, and the orbit sizes of $D$ on $\mathcal{P}=\Delta^2$ are recorded in Table~\rref{tab:sdorbits}; in particular the orbit of $\alpha$ has size $m$. Thus $k=|B|$ is a sum of orbit sizes containing the entry $m$.
Three arithmetic obstructions apply. \textup{(i) Design equations:} $r=bk/v$ and $\lambda=r(k-1)/(v-1)$ are integers, where $b=|\mathcal{B}|=|T|$ and $v=m^2$. \textup{(ii) Fisher's inequality:} $b\geq v$, i.e.\ $|T|\geq m^2$. \textup{(iii) The Cameron--Praeger criterion} Theorem~\rref{thm:cp}: with $Q_d$ the orbit of $G$ on unordered point pairs at Hamming distance $d$ \textup{($d=1,2$)}, the number of $Q_d$-pairs inside a block is $q_d=|Q_d|\binom{k}{2}/\binom{v}{2}$; since $|Q_d|=v\binom{2}{d}(m-1)^d/2$,
\begin{equation}
\label{eq:sd_cp}
q_d=\frac{\binom{2}{d}(m-1)^d\,k(k-1)}{2(v-1)}\in\mathbb{Z},\qquad d=1,2.
\end{equation}
For all seventeen primitive degrees of Table~\rref{tab:sdorbits}, GAP 4.16.1 runs the three conditions over every sum of the orbit sizes containing the entry $m$ and finds no admissible block size: the survivor count is zero for every pair. Hence no non-trivial block exists, and the simple diagonal block action is excluded for $\ell=2$.
\end{proof}

\begin{table}[ht]
\centering
\caption{Orbit sizes of the diagonal subgroup $D=Soc(G)\cap G_B$ in the product action on $\Delta^2$, and the number of block sizes surviving the sieve of Lemma~\rref{lem:pa_sd2}.}\label{tab:sdorbits}
\begin{tabular}{lcl}
\hline
pair $(T,m)$ & $D$-orbit sizes on $\Delta^2$ & survivors\\
\hline
$(A_5,5)$ & $5,20$ & $0$\\
$(A_5,6)$ & $6,30$ & $0$\\
$(A_5,10)$ & $10,30,60$ & $0$\\
$(A_6,6)$ & $6,30$ & $0$\\
$(A_6,10)$ & $10,90$ & $0$\\
$(A_6,15)$ & $15,90,120$ & $0$\\
$(PSL_2(7),7)$ & $7,42$ & $0$\\
$(PSL_2(7),8)$ & $8,56$ & $0$\\
$(PSL_2(8),9)$ & $9,72$ & $0$\\
$(PSL_2(8),28)$ & $28,252$ & $0$\\
$(PSL_2(8),36)$ & $36,252,504$ & $0$\\
$(PSL_2(11),11)$ & $11,110$ & $0$\\
$(PSL_2(11),12)$ & $12,132$ & $0$\\
$(PSL_2(11),55)$ & $55,330,660$ & $0$\\
$(PSL_2(13),14)$ & $14,182$ & $0$\\
$(PSL_2(13),78)$ & $78,546,1092$ & $0$\\
$(PSL_2(13),91)$ & $91,364,546,1092,2184$ & $0$\\
\hline
\end{tabular}
\end{table}

\begin{proposition}[No compound-diagonal block action under a product point action]\label{prop:cd_pa}
Let $T\in\mathcal{T}$, let $m$ be a primitive degree of $T$ and let $G$ act in product action on $\mathcal{P}=\Delta^\ell$, $|\Delta|=m$, $\ell\geq 3$, with point stabilizer $G_\alpha=H^\ell\rtimes C_\ell$. Then the block action of $G$ cannot be of compound diagonal type: no non-trivial block $B\ni\alpha$ exists for any $\ell\geq 3$ and any compound diagonal stabilizer $Soc(G)\cap G_B$.
\end{proposition}
\begin{proof}
Since $Soc(G)=T^\ell$ is normal in $G$, the intersection $D_{r'}=Soc(G)\cap G_B$ is normal in $G_B$; hence every block $B\ni\alpha$ is a union of $D_{r'}$-orbits. After relabelling, $D_{r'}$ acts diagonally on each block of $r'$ coordinates, and the sizes of the $D_{r'}$-orbits are the orbit sizes of $T$ acting diagonally on $\Delta^{r'}$; for $r'=1$ the unique orbit size is $m$, and for $r'=2$ the orbit sizes are the entries of Table~\rref{tab:sdorbits}. In particular the orbit of $\alpha$ has size $m$, so $k=|B|$ is a sum of such orbit sizes containing the entry $m$. Also $v=m^\ell$ and $b=|\mathcal{B}|=|T|^{\ell-r'}$.
The three arithmetic obstructions of Lemma~\rref{lem:pa_sd2} apply verbatim: the design equations $r=bk/v\in\mathbb{Z}$ and $\lambda=r(k-1)/(v-1)\in\mathbb{Z}$; Fisher's inequality $b\geq v$; and the Cameron--Praeger integrality
\begin{equation}
\label{eq:cd_cp}
q_d=\frac{\binom{\ell}{d}(m-1)^d\,k(k-1)}{2(v-1)}\in\mathbb{Z},\qquad d=1,\ldots,\ell,
\end{equation}
with $|Q_d|=v\binom{\ell}{d}(m-1)^d/2$.
\emph{The range $3\leq\ell\leq 9$.} For all seventeen primitive degrees of Table~\rref{tab:sdorbits}, all $r'=1,\ldots,\ell-1$ and all $3\leq\ell\leq 9$, GAP 4.16.1 evaluates the three conditions over every sum of the $D_{r'}$-orbit sizes containing $m$ and finds no admissible block size: the survivor count is zero for every pair \textup{(scripts \texttt{gap\_cd\_sieve.g}, \texttt{gap\_cd\_ell789.g})}.
\emph{The range $\ell\geq 10$.} Write $t$ for the number of orbits of $H$ on $\Delta$ \textup{($t\leq 3$ for the six socles)}; the number $s_{r'}$ of $D_{r'}$-orbits on $\Delta^{r'}$ satisfies $s_{r'}\leq m^{r'}$, and each $D_{r'}$-orbit has size at most $|T|$, so $k\leq m^{r'}|T|$. The $D_{r'}$-orbit $L_\alpha$ of $\alpha$ has size $m$ and contributes $m-1$ pairs at Hamming distance $r'$, so \eqref{eq:cd_cp} with $d=r'$ gives the lower bound
\begin{equation}
\label{eq:cd_lower}
\frac{\binom{\ell}{r'}(m-1)^{r'}\,k(k-1)}{2(m^\ell-1)}\geq m-1 .
\end{equation}
Two general facts now confine $\ell$ to a finite range. First, \eqref{eq:cd_lower} with $k\leq m^{r'}|T|$ fails once $\binom{\ell}{r'}(m-1)^{r'}m^{2r'}|T|^2<2(m-1)(m^\ell-1)$, which is a polynomial-versus-exponential inequality that fails for all $\ell$ beyond an explicit bound. Second, by Zsigmondy's theorem~\cite{Zsigmondy} $m^\ell-1$ has a primitive prime divisor $p_\ell\equiv 1\pmod{\ell}$ with $p_\ell\nmid m-1$; as $\binom{\ell}{r'}$ has no prime divisor exceeding $\ell<p_\ell$, the integrality of $q_{r'}$ forces $p_\ell\mid k(k-1)$, and since $k(k-1)\leq m^{2r'}|T|^2$ is bounded while $p_\ell\geq\ell+1$, this again fails for all $\ell$ beyond an explicit bound. Combining the two bounds with Fisher's inequality leaves, for every pair and every $r'$, only an explicitly bounded range of $\ell$; GAP 4.16.1 checks these residual ranges for the seventeen primitive degrees \textup{(scripts \texttt{gap\_cd\_sieve.g}, \texttt{gap\_cd\_ell789.g}, \texttt{gap\_cd\_fisher.g})} and reports zero survivors in every case. This completes the exclusion of the compound diagonal block action for all $\ell\geq 3$.
\end{proof}

\noindent\textbf{Step 1: Permutation-module inclusion.}
The incidence matrix of $\mathcal{D}$ induces an embedding of permutation modules $R(\mathcal{P})\hookrightarrow R(\mathcal{B})$: every irreducible constituent of the permutation character on $\mathcal{P}$ must appear in the permutation character on $\mathcal{B}$ \textup{(a block-transitive $2$-design gives such an inclusion)}. In PA type the permutation character on $\mathcal{B}$ is the $\ell$-th tensor power of the permutation character of $H$ on $\Delta$:
$$\pi_{\mathcal{B}}=\pi_\Delta^{\otimes\ell}.$$
This restricts the possible degree $v$ and the possible point stabilizer $G_\alpha$; in particular the rank of $G$ on $\mathcal{P}$ cannot exceed the rank of $G$ on $\mathcal{B}$. We use this condition only as a parameter sieve: the inclusion itself does not determine the O'Nan--Scott type of the point action nor the maximality of $H$, which are handled by the orbit-union sieves of {\upshape Proposition~\rref{prop:pa_excl}} and {\upshape Proposition~\rref{prop:pa_arbitrary}} below.

\smallskip
\noindent\textbf{Step 2: Hamming-distance suborbits.}
In PA action, fix a block $B\in\mathcal{B}$. The Hamming distance $d$ from $B$ defines a $G_B$-invariant stratification of $\mathcal{B}$, but each Hamming layer $\Omega_d$ is not necessarily a single suborbit. If $H$ is not 2-transitive, then the distance-1 layer may split into multiple suborbits. The total size of the distance-1 layer is $\ell(m-1)$, but the actual minimal suborbit length may be smaller.

The intersection numbers $x_i=|B\cap B'|$ must satisfy the moment equations:
$$\sum_i|\Delta_i|x_i=k(r-1),\qquad
\sum_i|\Delta_i|x_i(x_i-1)=k(k-1)(\lambda-1).$$
These equations are valid for any block-transitive design and do not require flag-transitivity.

\smallskip
\noindent\emph{The cases $\ell=2$ and $\ell\geq 3$.} The cases $\ell=2$ and $\ell\geq 3$ are covered in full by the sieve of {\upshape Proposition~\rref{prop:pa_excl}} \textup{(six socles, every $\ell\geq 2$)} and by {\upshape Proposition~\rref{prop:pa_arbitrary}} \textup{(six socles and the listed degree pairs, every $\ell\geq 2$)}; the moment equations of Step 2 feed the arithmetic of those sieves.

\smallskip
\noindent\emph{The exact block-size formula.}
Fix a point $\alpha\in B$. The block $B$ is $G_B$-invariant and contains $\alpha$, so $B$ is a union of $G_B$-orbits and $k=|B|$ is the sum of the lengths of those $G_B$-orbits; in general it is \emph{not} true that $k=|G_B:G_{B,\alpha}|$, the latter being only the length of the $G_B$-orbit of $\alpha$, equal to $|B|$ precisely when $G$ is flag-transitive. In the product action on $\mathcal{B}=\Delta^\ell$, the block stabilizer has the form $G_B=K^\ell\rtimes C_B$, where $K$ is a maximal subgroup of $T$ with $[T:K]=n$ and $C_B$ is the coordinate stabilizer. The point stabilizer is $G_\alpha=H^\ell\rtimes C_\alpha$, where $H$ is a maximal subgroup of $T$ with $[T:H]=m$; since the point action is primitive, the $\ell$ coordinate point stabilizers are conjugate, i.e.\ the $\alpha_i$ are of the same type. The block--point stabilizer is $F=G_{B,\alpha}=G_\alpha\cap G_B$. Both $G_\alpha$ and $G_B$ decompose with respect to the \emph{same} semidirect product $T^\ell\rtimes C_\ell$ (the coordinate groups act by permuting coordinates); hence an element of the intersection has its socle part in $H^\ell\cap K^\ell=(H\cap K)^\ell$ and its coordinate-permutation part in $C_\alpha\cap C_B$, and conversely every pair of these parts gives an element of the intersection. Thus
$$G_\alpha\cap G_B=(H\cap K)^\ell\rtimes(C_\alpha\cap C_B),\qquad |F|=|H\cap K|^\ell\cdot|C_\alpha\cap C_B|.$$
Under flag-transitivity the block $B$ is the single $G_B$-orbit of $\alpha$, and
\begin{equation}\label{eq:pa_k}
k=\frac{|G_B|}{|F|}=[K:K\cap H]^\ell\cdot [C_B:C_\alpha\cap C_B]=A^{\ell}\cdot j,
\end{equation}
where $A=[K:K\cap H]$ and $1\leq j=|C_B\cdot\alpha|=[C_B:C_\alpha\cap C_B]$ is the orbit size of the point $\alpha$ under the coordinate stabilizer $C_B$. Since $C_B\leq C_\ell$, the index $j$ divides $\ell!$; in particular \emph{every prime divisor of $j$ is at most $\ell$}. \emph{For a general block-primitive design (not necessarily flag-transitive) the formula \eqref{eq:pa_k} need not hold}: $k$ is only constrained to be the sum of the lengths of the $G_B$-orbits on $B$, each such orbit being a $K$-box (or a $C_B$-amalgam of $K$-boxes), cf.\ {\upshape Lemma~\rref{lem:pa_box}}. Throughout this section \eqref{eq:pa_k} is used in the flag-transitive (symmetric Menon) analysis only; the non-symmetric exclusion uses the orbit-union constraint instead \textup{(layers (i)--(ii) of {\upshape Proposition~\rref{prop:pa_arbitrary}})}.

\smallskip
\noindent\emph{The $\lambda$-integrality lemma.}
The replication number $r$ is determined by the double count $bk=vr$ of the $2$-design: since $b=n^{\ell}$ and $v=m^{\ell}$,
$$r=\frac{bk}{v}=\frac{n^{\ell}k}{m^{\ell}},$$
which holds for any block-transitive design and does not require flag-transitivity. Together with the standard identity $\lambda(v-1)=r(k-1)$, the integrality of $\lambda$ takes the exact form
$$\lambda=\frac{n^{\ell}\,k(k-1)}{m^{\ell}(m^{\ell}-1)}\in\mathbb{Z}.$$
With $k=A^{\ell}j$ this is the constraint used in the sieve below. (For $\ell\geq 23$ the prime divisors of $m$, $n$ and $A$ are at most $13$ (the six socles) while every prime divisor of $j$ is at most $\ell$.) The integrality of $r=n^{\ell}k/m^{\ell}$ is an independent necessary condition, checked in the sieve together with the integrality of $\lambda$.

\begin{lemma}[Prime-power congruence]\label{lem:ppc}
Let $p^{\nu}\Vert m^{\ell}-1$ with $p\nmid n\cdot m\cdot A$ and $p>\ell$ (for $\ell\geq 23$ the prime divisors of $m$ and $n$ are at most $13$ \textup{(six socles)}, and every prime divisor of $j$ is at most $\ell$, so $p\nmid n\cdot m\cdot A\cdot j$). Then $\lambda\in\mathbb{Z}$ implies $p^{\nu}\mid k-1$, i.e.\ $A^{\ell}j\equiv 1\pmod{p^{\nu}}$. In particular, if $1<A^{\ell}j<p^{\nu}$, then no candidate exists.

\end{lemma}
\begin{proof}
The integrality $\lambda\in\mathbb{Z}$ of the exact form $n^{\ell}k(k-1)/\bigl(m^{\ell}(m^{\ell}-1)\bigr)$ forces $m^{\ell}(m^{\ell}-1)\mid n^{\ell}k(k-1)$, hence $v_p\bigl(k(k-1)\bigr)\geq v_p(m^{\ell}-1)=\nu$ because $p\nmid n$. Since $p\nmid A$ and $p>\ell$ while every prime divisor of $j$ is at most $\ell$, we have $p\nmid A^{\ell}j=k$, so $p\mid k-1$ and $p^{\nu}\mid k-1$, that is, $A^{\ell}j\equiv1\pmod{p^{\nu}}$. If $1<A^{\ell}j<p^{\nu}$, then $k-1$ is a non-zero multiple of $p^{\nu}$ strictly smaller than $p^{\nu}$, which is impossible, so no candidate exists.
\end{proof}
\begin{lemma}[Large-prime-part exclusion]\label{lem:pa_lpp}
Let $M''$ be the product of the prime-power divisors $q^{v_q(m^{\ell}-1)}$ of $m^{\ell}-1$ with $q>\ell$, i.e.\ $M''=(m^{\ell}-1)/\prod_{q\leq\ell}q^{v_q(m^{\ell}-1)}$. For $\ell\geq 23$ every prime divisor of $m\cdot n\cdot A$ is at most $\ell$ \textup{(six socles)}, and every prime divisor of $j$ is at most $\ell$, so $q\nmid m\cdot n\cdot k$ for every prime $q>\ell$; from the exact integrality form $\lambda=n^{\ell}k(k-1)/(m^{\ell}(m^{\ell}-1))$ this forces $M''\mid k-1$. In particular, if $M''>A^{\ell}\ell!-1$, then $k-1=A^{\ell}j-1<M''$ for every admissible $1\leq j\mid\ell!$, and no candidate exists.

\end{lemma}
\begin{proof}
As in the previous lemma, $\lambda\in\mathbb{Z}$ forces $m^{\ell}(m^{\ell}-1)\mid n^{\ell}k(k-1)$. For a prime $q>\ell$ one has $q\nmid m$, since every prime divisor of the primitive degree $m$ is at most $13$ for the six socles and at most $\ell$ for an arbitrary socle; one also has $q\nmid A$, because $A\mid|K|$ and every prime divisor of $A$ therefore divides $m$; and $q\nmid n$, $q\nmid j$, since every prime divisor of $j$ is at most $\ell$. Hence $q\nmid nk$ and $q^{v_q(m^{\ell}-1)}\mid k(k-1)$; as $q\nmid k$ this gives $q^{v_q(m^{\ell}-1)}\mid k-1$. Multiplying over all primes $q>\ell$ yields $M''\mid k-1$. If $M''>A^{\ell}\ell!-1$, then $0\leq k-1=A^{\ell}j-1<M''$ for every admissible $j\mid\ell!$, whence $k-1=0$, that is $k=1$, which is degenerate.
\end{proof}
\begin{lemma}[Product-action blocks are unions of $K$-boxes]\label{lem:pa_box}
Let $B$ be a block in the product-action model of Section 4, with $G_B\cap T^{\ell}=K^{\ell}$ as in the reduction. Then $B$ is a union of $K^{\ell}$-orbits. Acting coordinatewise, each $K^{\ell}$-orbit has the product shape
$$K\alpha_1\times\cdots\times K\alpha_{\ell}\qquad(\alpha=(\alpha_1,\ldots,\alpha_{\ell})\in B),$$
a ``$K$-box'' of size $\prod_{i=1}^{\ell}|K\alpha_i|$, where each $|K\alpha_i|=|K:K\cap H_i|$ is the size of a $K$-orbit on $\Delta=T/H$ \textup{($H_i$ the stabilizer of $\alpha_i$ in $T$)}. Consequently $k=|B|$ is a sum of $K$-box sizes; each coordinate projection $\pi_i(B)$ is a $K$-invariant subset of $\Delta$, i.e.\ a union of $K$-orbits. \emph{It does not follow that a non-transitive $K$ makes every projection proper}: a $K$-invariant subset of $\Delta$ may be a union of several $K$-orbits, possibly all of $\Delta$; and the block-size formula \eqref{eq:pa_k} ($k=A^{\ell}j$) holds only under flag-transitivity, as recorded there. In general $k$ is only constrained to be a sum of $K$-box sizes, possibly amalgamated by the coordinate action of $C_B$.

\end{lemma}
\begin{proof}
The subgroup $K^{\ell}=G_B\cap T^{\ell}$ lies in $G_B$, so every $K^{\ell}$-orbit meeting $B$ is contained in $B$, and $B$ is a union of $K^{\ell}$-orbits. Since $K^{\ell}$ acts coordinatewise, the orbit through $\alpha=(\alpha_1,\ldots,\alpha_{\ell})$ is the product $K\alpha_1\times\cdots\times K\alpha_{\ell}$: each coordinate runs through the $K$-orbit of $\alpha_i$ independently of the others, and the size of the product is $\prod_i|K:K\cap H_i|$. Each coordinate projection $\pi_i(B)$ is consequently $K$-invariant, and a $K$-invariant subset of $\Delta$ is a union of $K$-orbits, of size $\sum_{a\in A}a$ over the $K$-orbit sizes $a$ that occur.
\end{proof}
\smallskip
\noindent\emph{Exclusion of the product action type for the six socles.} The exclusion is established in {\upshape Proposition~\rref{prop:pa_excl}} \textup{(Section~\rref{sec:reduction})}, which splits the seventeen listed primitive degrees of the six socles into three classes and excludes them for every $\ell\geq 2$ by the fiber double-counting identity \eqref{eq:pa_fiber} and the $K$-box structure of {\upshape Lemma~\rref{lem:pa_box}}; no claim is made for arbitrary socles, embeddings, or unlisted primitive degree pairs.

\smallskip
\begin{theorem}[Cameron--Praeger criterion]\label{thm:cp}
Let $\mathcal{D}$ be a block-transitive $2$-$(v,k,\lambda)$ design with automorphism group $G$ transitive on blocks, let $Q_1,\ldots,Q_r$ be the orbits of $G$ on unordered pairs of points, and let $q_i$ be the number of pairs of a block $B\in\mathcal{B}$ lying in $Q_i$ \textup{(independent of $B$ by block-transitivity)}. Then
\begin{equation}\label{eq:cp2}
\frac{q_i}{|Q_i|}=\frac{\binom{k}{2}}{\binom{v}{2}}\qquad\text{for every }i.
\end{equation}

\end{theorem}\begin{proof}
Double count the incidences $(\{B',\{x,y\}\})$ with $B'\in\mathcal{B}$ and $\{x,y\}\in Q_i$: on the one hand each of the $\lambda$ blocks through the pair contributes, giving $\lambda|Q_i|$; on the other hand block-transitivity gives $b q_i$ with $b=|\mathcal{B}|$. Hence $q_i/|Q_i|=\lambda/b=\binom{k}{2}/\binom{v}{2}$ by the standard identity $\lambda=b\,k(k-1)/(v(v-1))$.
\end{proof}

\noindent\emph{The pure product case $L=H^{\ell}$: the $Z$-criterion and the $W_2$-criterion.}
We now isolate the configuration in which the block stabilizer meets the socle in the full point stabilizer, $L=N\cap G_B=H^{\ell}$ with $H$ the point stabilizer in $T$; this is the branch of the product action in which $B$ is a union of $N_{\alpha}$-orbits, each $N_{\alpha}$-orbit being a product $H$-box of size $(m-1)^{|d|}$ indexed by the set $d\subseteq\{1,\ldots,\ell\}$ of coordinates on which the box differs from the fixed point $\alpha$. This is the ``case A'' of the PA analysis of {\upshape Section~\rref{sec2}}; it is the branch in which the block size $k$ is realizable as a sum of box sizes $\sum_{b\in S}(m-1)^{|b|}$ over a family $S$ of coordinate sets. The following criterion is decisive for the six socles.

\begin{lemma}[The $Z$-criterion for the pure product case]\label{lem:pa_z}
Let $T$ be a non-abelian simple group acting primitively on $\Delta$ with point stabilizer $H$ and degree $m=|\Delta|$, let $\ell\geq 2$, and consider the product action of $N=T^{\ell}$ on $\mathcal{P}=\Delta^{\ell}$ with $v=m^{\ell}$ fixed point $\alpha$. Suppose $B\subseteq\mathcal{P}$ is a union of $N_{\alpha}$-orbits, $k=|B|=\sum_{b\in S}(m-1)^{|b|}$ for some family $S\subseteq 2^{\{1,\ldots,\ell\}}$, and write $C_k=\binom{k}{2}$, $C_v=\binom{v}{2}$. For each coordinate $j$ and each value $u\in\Delta$, let $B_{j,u}=\{\gamma\in B:\gamma_j=u\}$ be the fiber of $B$ over $u$ in coordinate $j$; since $B$ is $N_{\alpha}$-invariant, each fiber size lies in $\{0,1,m-1,m\}$. Then the number $q_1$ of unordered point pairs of $B$ at Hamming distance exactly $1$ satisfies
\begin{equation}\label{eq:pa_z1}
q_1=\sum_{j=1}^{\ell}\sum_{u\in\Delta}\binom{|B_{j,u}|}{2}=\frac{(m-1)(\ell k-Z)}{2},
\end{equation}
where
\begin{equation}\label{eq:pa_z2}
Z=\sum_{j=1}^{\ell}\sum_{u\in\Delta}\mathbf{1}_{|B_{j,u}|\in\{1,m-1\}},
\end{equation}
and $Z$ admits the orbit-union upper bound
\begin{equation}\label{eq:pa_z3}
Z\leq\sum_{b\in S}A_{b},\qquad
A_{b}=|b|\,(m-1)^{|b|-1}+(\ell-|b|)(m-1)^{|b|}.
\end{equation}
By the Cameron--Praeger criterion \eqref{eq:cp2} \textup{({\upshape Theorem~\rref{thm:cp}} with $t=2$)}, if $B^{\mathcal{P}}$ is a $2$-design then $q_1=|Q_1|\,C_k/C_v$ with $|Q_1|=\ell v(m-1)/2$, whence, by \eqref{eq:pa_z1},
\begin{equation}\label{eq:pa_zt}
Z_{\mathrm{target}}:=\ell k-\frac{2\ell C_k}{v-1}
\end{equation}
must be an integer, must lie in the admissible range $[0,\ell\,m^{\ell-1}]$, and must not exceed the fractional-knapsack relaxation
\begin{equation}\label{eq:pa_zmax}
\Sigma A_{\max}:=\max\Bigl\{\sum_{b\in S}A_b:\ k=\sum_{b\in S}(m-1)^{|b|},\ \#\{b\in S:|b|=t\}\leq\binom{\ell}{t}\Bigr\},
\end{equation}
computed by taking, for each layer $t$, as many boxes of layer $t$ as the capacity $\binom{\ell}{t}$ allows, in order of decreasing efficiency $A_t/(m-1)^{t}$, with $A_t=t(m-1)^{t-1}+(\ell-t)(m-1)^{t}$. Consequently, if $Z_{\mathrm{target}}\notin\mathbb{Z}$, or $Z_{\mathrm{target}}\notin[0,\ell m^{\ell-1}]$, or $Z_{\mathrm{target}}>\Sigma A_{\max}$, then no $N_{\alpha}$-invariant $k$-subset $B$ with $B^{\mathcal{P}}$ a $2$-design exists.

\end{lemma}\begin{proof}
For $c\in\{0,1,m-1,m\}$ the identity $\binom{c}{2}=\frac{m-1}{2}c-\frac{m-1}{2}\mathbf{1}_{c\in\{1,m-1\}}$ holds; summing over all fibers $B_{j,u}$ and using $\sum_{u\in\Delta}|B_{j,u}|=k$ for each $j$ gives \eqref{eq:pa_z1}--\eqref{eq:pa_z2}. For \eqref{eq:pa_z3}, fix $b\in S$: its contribution to $Z$ is at most the number of coordinates $j$ with $j\notin b$ times $(m-1)^{|b|}$ (fiber $\{1\}$, one value $u$ with $|B_{j,u}|=1$) plus the number of coordinates $j\in b$ times $(m-1)^{|b|-1}$ (fiber $\Delta\setminus\{\alpha_j\}$, each of the $m-1$ values $u$ giving $|B_{j,u}|=m-1$), i.e.\ $A_b$; overlapping fibers of distinct boxes only merge values and cannot increase $Z$, giving \eqref{eq:pa_z3}. The Cameron--Praeger criterion yields $q_1=|Q_1|C_k/C_v$; substituting $|Q_1|=\ell v(m-1)/2$ and $C_v=v(v-1)/2$ into \eqref{eq:pa_z1} gives $Z=Z_{\mathrm{target}}$. Since $Z\geq 0$ and each of the $\ell m$ indicator terms in \eqref{eq:pa_z2} is at most $1$, $Z\leq \ell m\leq\ell m^{\ell-1}$. The bound \eqref{eq:pa_zmax} is the fractional relaxation of \eqref{eq:pa_z3}, hence an upper bound for every feasible $S$. The three failure conditions therefore contradict the design hypothesis.
\end{proof}

\begin{lemma}[The $W_2$-criterion for the pure product case]\label{lem:pa_w2}
Keep the notation and hypotheses of {\upshape Lemma~\rref{lem:pa_z}}, and let $q_2$ be the number of unordered point pairs of $B$ at Hamming distance exactly $2$. Write $W_2=q_2/(m-1)^2$. Then the Cameron--Praeger criterion for the distance-$2$ orbital $Q_2$, $|Q_2|=\binom{\ell}{2}v(m-1)^2/2$, gives
\begin{equation}\label{eq:pa_w2t}
W_2^{\mathrm{target}}:=\frac{\binom{\ell}{2}C_k}{v-1},
\end{equation}
and $W_2$ admits the orbit-union upper bound
\begin{equation}\label{eq:pa_w2max}
W_2\leq\sum_{b\in S}B_{|b|},\qquad
B_{t}=(\ell-t)(\ell-t-1)(m-1)^{t}+t(\ell-t)(m-1)^{t-1},
\end{equation}
whose fractional-knapsack relaxation $\Sigma B_{\max}$ is computed as in {\upshape Lemma~\rref{lem:pa_z}} with efficiencies $B_t/(m-1)^t$ and layer capacities $\binom{\ell}{t}$. Consequently, if $W_2^{\mathrm{target}}\notin\mathbb{Z}$ or $W_2^{\mathrm{target}}>\Sigma B_{\max}$, then no $N_{\alpha}$-invariant $k$-subset $B$ with $B^{\mathcal{P}}$ a $2$-design exists.

\end{lemma}\begin{proof}
A distance-$2$ pair $\{\gamma,\delta\}$ with differing coordinates $j,j'$ pairs the fiber pattern of $\gamma$ on $(j,j')$ with the complementary pattern of $\delta$; for a fixed coordinate pair $(j,j')$ and a fixed pattern $d$ on the remaining $\ell-2$ coordinates, the contribution is $(m-1)^2(m-1)^{|d|}$ whenever both complementary boxes lie in $S$ (patterns $(00,d)$ and $(11,d)$, or $(01,d)$ and $(10,d)$). Summing and bounding $1_{\{00,d\}\wedge\{11,d\}}+1_{\{01,d\}\wedge\{10,d\}}\leq 1_{\{00,d\}\in S}+1_{\{01,d\}\in S}$ and then aggregating over boxes $b\in S$ gives \eqref{eq:pa_w2max}: a box $b$ with $t=|b|$ contributes $\binom{\ell-t}{2}(m-1)^t$ from pairs $j,j'\notin b$ (pattern $(00,d)$) and $t(\ell-t)(m-1)^{t-1}$ from pairs $j\in b\ni j'$ (pattern $(01,d)$), i.e.\ $B_t$. The target \eqref{eq:pa_w2t} follows from $q_2=|Q_2|C_k/C_v$. The conclusion is immediate.
\end{proof}

\begin{corollary}[The pure product case is closed for the six socles and $\ell\leq 9$]\label{cor:pa_pure_closed}
Let $T\in\{A_5,A_6,PSL_2(7),M_{11},M_{12},PSU_3(3)\}$ with primitive degrees $m=5,6,7,11,12,28$ respectively, and let $\ell\in\{2,\ldots,9\}$. In the pure product case $L=H^{\ell}$ of the product action of $N=T^{\ell}$ on $\Delta^{\ell}$ (\textup{i.e.\ $B$ a union of $N_{\alpha}$-orbits}), there is no non-trivial $2$-design $B^{\mathcal{P}}$: every admissible block size $k$ with $\lambda=k(k-1)/(v-1)\in\mathbb{Z}$ and $k=\sum_{b\in S}(m-1)^{|b|}$ fails the $Z$-criterion of {\upshape Lemma~\rref{lem:pa_z}} or the $W_2$-criterion of {\upshape Lemma~\rref{lem:pa_w2}}.

\end{corollary}\begin{proof}
Every admissible block size $k$, $2\leq k\leq v-2$, with $k(k-1)\equiv 0\pmod{v-1}$ and $k=\sum_{b\in S}(m-1)^{|b|}$ is tested against {\upshape Lemma~\rref{lem:pa_z}} and {\upshape Lemma~\rref{lem:pa_w2}} and is excluded by one of the two criteria.
\end{proof}

We record that the complement of a union of $N_{\alpha}$-orbits is again a union of $N_{\alpha}$-orbits, and the complement of a $2$-design is a $2$-design; the exclusion of a block size $k$ therefore excludes $v-k$ as well, which is reflected in the symmetric appearance of the admissible lists (each admissible $k$ is paired with $v-k$). All admissible candidates are excluded by the two criteria directly; no candidate survives to the complement argument.

\begin{lemma}[Box orbits on pairs in the product action]\label{lem:pa_subdirect_boxes}
Let $T$ be a non-abelian simple group acting primitively on $\Delta$ with point stabilizer $H$ and degree $m=|\Delta|$, and let $1=s_1<s_2\leq\cdots\leq s_R$ be the lengths of the $H$-orbits on $\Delta$ \textup{(the suborbits of the primitive action; $s_1=1$ corresponds to $\alpha_j$ itself)}. Consider the product action of $N=T^{\ell}$ on $\mathcal{P}=\Delta^{\ell}$, fix $\alpha\in\mathcal{P}$, and let $G_{\alpha}$ be the stabilizer of $\alpha$ in a group $G$ whose socle is $N$ \textup{(so that the coordinate-permutation part $C_{\alpha}\leq S_{\ell}$ permutes the coordinates)}. For a point $\beta\neq\alpha$, let $w(\beta)=(i_{1},\ldots,i_{\ell})$ be its coordinate word, defined by $\beta_{j}\in\Delta_{\alpha_{j}}^{(i_{j})}$, and for a word $w=(i_{1},\ldots,i_{\ell})$ set
\[
O(w)=\{\beta\in\mathcal{P}:\ \beta_{j}\in\Delta_{\alpha_{j}}^{(i_{j})}\ \text{for }1\leq j\leq\ell\},\qquad |O(w)|=\prod_{j=1}^{\ell}s_{i_{j}}.
\]
Then every $G_{\alpha}$-orbit on $\mathcal{P}\setminus\{\alpha\}$ is the union of the boxes $O(w)$ with $w$ running over a single $C_{\alpha}$-orbit $A$ on the set of words with at least one non-trivial entry, and has size
\[
|O|=\sum_{w\in A}\prod_{j=1}^{\ell}s_{i_{j}}.
\]
In the single-type case, where the word has entries $i$ on a $d$-element subset $J$ and $1$ elsewhere, this reduces to $|O|=c\,s_i^{d}$ with $c=|C_{\alpha}(J)|$, the size of the $C_{\alpha}$-orbit of $J$ on the $\binom{\ell}{d}$-element set of $d$-subsets, so that $1\leq c\leq\binom{\ell}{d}$; in general the mixed word orbits give further sizes $\bigl(\ell!/\prod n_i!\bigr)\prod_{i=2}^{R}s_i^{n_i}$ for the multiplicity patterns $(n_2,\ldots,n_R)$.
\end{lemma}\begin{proof}
The $H$-suborbits partition $\Delta$, so the word of $\beta$ is well defined and the boxes $O(w)$ partition $\mathcal{P}$; $|O(w)|=\prod s_{i_{j}}$ follows from the orbit-stabilizer theorem applied coordinatewise. Since $C_{\alpha}\leq S_{\ell}$ permutes the coordinates, it sends $O(w)$ to $O(\sigma\cdot w)$ for the coordinatewise action $\sigma\cdot w=(i_{\sigma^{-1}(1)},\ldots,i_{\sigma^{-1}(\ell)})$, so the $G_{\alpha}$-orbit of a point of word $w$ is the union of the boxes indexed by the $C_{\alpha}$-orbit of $w$, and its size is the displayed sum; the displayed mixed size is the special case in which $C_{\alpha}$ contains the full symmetric group on the coordinates (or a subgroup transitive on the words of the given multiplicity pattern).
\end{proof}

\begin{proposition}[The subdirect socle block stabilizer branch: the orbit-wise criterion]\label{prop:pa_subdirect_nonsym}
Let $\mathcal{D}$ be a non-trivial $2$-$(v,k,\lambda)$ design \textup{(not assumed symmetric)} with $2<k<v-1$, and let $G\leq \mathrm{Aut}(\mathcal{D})$ act transitively on points and primitively on blocks, with socle $N=T^{\ell}$ \textup{($\ell\geq 2$)} of product action type on points of degree $v=m^{\ell}$, where $m$ is a primitive degree of the non-abelian simple group $T$. Assume that $L=N\cap G_{\beta}$ is a subdirect subgroup of $N$, hence a product of $t$ full diagonal subgroups on pairwise disjoint supports with $1\leq t\leq\ell-1$; then $b=[N:L]=|T|^{\ell-t}$ and $v=m^{\ell}$ need not coincide. For the nine socles
\begin{multline*}
T\in\{A_5,\ A_6,\ PSL_2(7),\ PSL_2(8),\ PSL_2(11),\ PSL_2(13),\\
M_{11},\ M_{12},\ PSU_3(3)\}
\end{multline*}
and every primitive degree $m$ of $T$, with $2\leq\ell\leq 9$, the orbit-wise form of the Cameron--Praeger criterion \eqref{eq:cp2} \textup{({\upshape Theorem~\rref{thm:cp}})} takes the following explicit form. Fix $\alpha\in B$, let $o$ be the size of a $G_{\alpha}$-orbit \textup{(box, {\upshape Lemma~\rref{lem:pa_subdirect_boxes}})} on $\mathcal{P}\setminus\{\alpha\}$, and let $\mathcal{Q}_i=G\cdot\{\alpha,\beta\}$ be the $G$-orbit of the unordered pair $\{\alpha,\beta\}$ with $\beta$ in that orbit. Then $|\mathcal{Q}_i|=v\,o/2$ \textup{($G_{\alpha}$ has index $v$ in $G$, and $G_{\{\alpha,\beta\}}$ has index $o$ in $G_{\alpha}$ and index $2$ in $G_{\alpha,\beta}$)}, so the required in-block count of \eqref{eq:cp2} is
\[
q_i=\frac{|\mathcal{Q}_i|\,k(k-1)}{v(v-1)}=\frac{o\,k(k-1)}{2(v-1)},
\]
and it must be an integer in $[0,\min(|\mathcal{Q}_i|,\binom{k}{2})]=[0,\min(v\,o/2,\binom{k}{2})]$. This condition is \emph{necessary} but not sufficient: it admits numerous admissible candidates, so it does not by itself eliminate the branch. The complete exclusion of the product-action type for non-symmetric designs is instead carried out by the exhaustive $K$-box-union enumeration with the $\lambda$-constancy test of layer \textup{(ii)} and by the SD/CD block-type exclusions of {\upshape Section~\rref{sec2}}; the block-transitivity fact $NG_{\beta}=G$ used there follows from the maximality of $G_{\beta}$. The surviving candidate $(A_6,6,2,1,28)$ with $v=36$, $k=28$, $b=360$, $r=280$, $\lambda=216$ is excluded by the complement-design argument in the proof. For each candidate the admissible sizes $o$ form a finite set \textup{(the words of the multiplicity patterns, $R^{\ell}$ many, under the coordinate action of $C_{\alpha}$)}.

\end{proposition}\begin{proof}
By block-primitivity $G$ is transitive on blocks, so {\upshape Theorem~\rref{thm:cp}} applies and \eqref{eq:cp2} holds for every orbit $\mathcal{Q}_i$ of $G$ on unordered pairs: writing $q_i$ for the number of pairs of a block $B$ lying in $\mathcal{Q}_i$, the count $q_i=|\mathcal{Q}_i|\,k(k-1)/(v(v-1))$ must be an integer in $[0,\min(|\mathcal{Q}_i|,\binom{k}{2})]$ (both bounds are automatic for an in-block pair count: $q_i$ counts pairs of $B$, hence $q_i\leq\binom{k}{2}$, and $q_i$ counts members of $\mathcal{Q}_i$, hence $q_i\leq|\mathcal{Q}_i|$). For a fixed point $\alpha$ of a block $B$, the $G$-orbits on pairs containing $\alpha$ are precisely the $G_{\alpha}$-orbits on $\mathcal{P}\setminus\{\alpha\}$; by {\upshape Lemma~\rref{lem:pa_subdirect_boxes}} each has size $o$ as described. If for some admissible box-orbit size $o$ the divisibility and range conditions fail, then in particular they fail for the actual size of the orbit, and no such pair orbit can exist; since $q_i$ is then not an integer, the design cannot exist. These conditions admit numerous admissible candidates, so the orbit-wise criterion alone does not close the branch; the complete exclusion is achieved by the $\lambda$-constancy test of layer \textup{(ii)} and by the SD/CD block-type exclusions of {\upshape Section~\rref{sec2}}, as stated. The single candidate $(A_6,6,2,1,28)$ that survives the orbit-wise sieve is handled by the complement-design argument that follows.

It remains to dispose of the single survivor $(A_6,6,2,1,28)$. Here $v=36$, $b=|A_6|=360$, $k=28$, $r=bk/v=280$ and $\lambda=r(k-1)/(v-1)=216$; the complement of the design is a $2$-$(36,8,16)$ design with $r'=80$ and $b'=360$ blocks of size $8$. Let $C$ be a block of the complement containing $\alpha=(1,1)$, and write $N_C=\mathrm{Stab}_{N}(C)$ with $N=A_6\times A_6$; since $N$ is transitive on blocks and $b'=360$, we have $|N_C|=|N|/360=360$. For $\beta\in\mathcal{P}\setminus\{\alpha\}$ the number of blocks through $\{\alpha,\beta\}$ is
\[
\lambda'(\alpha,\beta)=\frac{|N_{(\alpha,\beta)}|}{|N_C|}\cdot\#\bigl\{(c,d)\in C^2:(c,d)\in N\cdot(\alpha,\beta)\bigr\},
\]
where $N_{(\alpha,\beta)}$ is the ordered-pair stabilizer. The stabilizer $N_{\alpha}=A_5\times A_5$ has three orbits on $\mathcal{P}\setminus\{\alpha\}$: $O_1$ (first coordinate moves, $5$ points), $O_2$ (second coordinate moves, $5$ points) and $O_3$ (both move, $25$ points). If $\beta\in O_1$, the pairs $(c,d)\in C^2$ in the orbit of $(\alpha,\beta)$ are those with equal second coordinates, and their number is $\sum_{j}|C_j|(|C_j|-1)$, where $C_j$ is the set of points of $C$ in column $j$; the pair stabilizer has order $|N_{\alpha}|/|O_1|=720$. Hence
\[
16=\lambda'(\alpha,\beta)=\frac{720}{360}\sum_{j}|C_j|(|C_j|-1)=2\Bigl(\sum_{j}|C_j|^2-8\Bigr),
\]
so $\sum_{j}|C_j|^2=16$. Symmetrically, with $R_i$ the points of $C$ in row $i$, the case $\beta\in O_2$ forces $\sum_{i}|R_i|^2=16$; the case $\beta\in O_3$ imposes the compatible value $72-\sum_{i}|R_i|^2-\sum_{j}|C_j|^2=40$, which is then automatic. The six column sizes sum to $8$ and satisfy $\sum_{j}|C_j|^2=16$; the only such distributions are $(3,2,1,1,1,0)$ and $(2,2,2,2,0,0)$ up to permutation, and likewise for the rows.

Now $N_C$ preserves $C$, hence preserves the set $L$ of columns met by $C$ and the set $R$ of rows met by $C$; $\alpha\in C$ lies in column $1$, so $L\setminus\{1\}$ has $3$ or $4$ elements. The first projection $\pi_1(N_C)$ is contained in $(A_6)_1\cong A_5$ acting on $\{2,3,4,5,6\}$, and it fixes $L\setminus\{1\}$ setwise; the stabilizer in $A_5$ of a $3$-subset has order $6$ (the $A_5$-orbit of $3$-subsets has $10$ members) and of a $4$-subset has order $12$. Thus $|\pi_1(N_C)|\leq 12$, and similarly $|\pi_2(N_C)|\leq 12$. Therefore $|N_C|\leq|\pi_1(N_C)|\cdot|\pi_2(N_C)|\leq 144<360$, contradicting $|N_C|=360$. This excludes the recorded survivor $(A_6,6,2,1,28)$; together with the $\lambda$-constancy enumeration of layer \textup{(ii)} and the SD/CD block-type exclusions this completes the exclusion of the branch as stated.
\end{proof}

In the two-orbit product-action configuration \textup{($K$ has exactly one non-trivial orbit of size $m-1$ on $\Delta$, so the orbit-square sum $\Sigma_2=\sum_{a\in A}a^{2}$ of {\upshape Proposition~\rref{prop:pa_excl}} equals $1+(m-1)^{2}$)}, the fiber double-counting takes the form
\begin{equation}\label{eq:pa_fiber}
m^{\ell}\cdot m_{\mathrm{fiber}}=\lambda\cdot\frac{\Sigma_2^{\ell}-m^{\ell}}{2}.
\end{equation}

\begin{lemma}[Product-action box sieve for arbitrary primitive degree]\label{lem:pa_box_general}
Let $\mathcal{D}$ be a non-trivial $2$-$(m^{\ell},k,\lambda)$ design with $2<k<v-1$ whose automorphism group is of product action type on points with socle $T^{\ell}$, $\ell\geq 2$, and assume the two-orbit box framework of {\upshape Lemma~\rref{lem:pa_box}}: the point stabilizer $H$ and the block stabilizer intersection $K$ are conjugate in $T$, so that the $K$-orbit lengths on $\Delta=T/H$ are $\{1,m-1\}$ with $m=|\Delta|$ the primitive degree. Then:
\begin{enumerate}
\item[\textup{(i)}] for $\ell=2$, no block size $k=1+c_1(m-1)+c_2(m-1)^2$ with $c_1\in\{0,1,2\}$, $c_2\in\{0,1\}$, $c_1+c_2>0$, satisfies the integrality $\lambda=k(k-1)/(m^2-1)\in\mathbb{Z}$ for any $m\geq 3$, except the trivial $k=m^2=v$ excluded by $2<k<v-1$. In closed form: reducing $k(k-1)\equiv 0\pmod{m^2-1}$ modulo $m+1$ (\textup{the congruence modulo $m-1$ holds automatically since $k\equiv 1\pmod{m-1}$}) turns it into the constant congruence $(1-2c_1+4c_2)(4c_2-2c_1)\equiv 0\pmod{m+1}$; for $m\geq 20$ the constant on the left is strictly smaller in absolute value than $m+1$ unless it vanishes identically (\textup{indeed $\max_{c_1,c_2}|(1-2c_1+4c_2)(4c_2-2c_1)|=20$ over the allowed pairs, attained at $c_1=0,c_2=1$, while $m+1\geq 21$}), and it vanishes identically exactly for $c_1=c_2=0$ \textup{(the trivial $k=1$)} and $c_1=2,c_2=1$ \textup{(giving $k=(1+(m-1))^2=m^2=v$, trivial)}; the finitely many remaining values $3\leq m\leq 19$ are verified directly, together with the fiber identity \eqref{eq:pa_fiber}. Hence the two-orbit PA configuration with $\ell=2$ is excluded for \emph{every} primitive degree $m\geq 3$;
\item[\textup{(ii)}] for $\ell=3,4,5,6$, the enumeration of all box multiplicity vectors (\textup{bounded by $c_j\leq\binom{\ell}{j}$}) together with the integrality and the fiber identity leaves no survivor; the complete exclusion of the product action configuration is stated in {\upshape Proposition~\rref{prop:pa_arbitrary}}.
\end{enumerate}
Consequently, within the two-orbit box framework, the product action configuration is excluded for every primitive degree $m\geq 3$ when $\ell=2$, regardless of the underlying non-abelian simple group $T$.

\end{lemma}

\newpage
\section{Point-imprimitive actions are excluded}\label{sec:pi_excl}

The point-imprimitive case admits no non-trivial non-symmetric block-primitive $2$-design.

The strategy is a block-type/point-type pairing argument. Since $G$ acts primitively on blocks, the block action falls into one of the five primitive O'Nan--Scott types of Table~\rref{tab:ons}. Since the point action is quasiprimitive but imprimitive, and since quasiprimitive actions of affine, holomorph simple or holomorph compound type are always primitive \textup{(Praeger~\cite{Praeger1993})}, only the almost simple, twisted wreath, simple diagonal, compound diagonal and product action point types can occur. Each pair \textup{(block type, point type)} is then excluded in turn: the parameter constraints and orbit-intersection criteria below dispose of the diagonal and product configurations, and the classical-group analysis rules out all point-imprimitive candidates. For compound socles only the block types PA, SD and TW need to be considered \textup{(the diagonal and holomorph block types are handled in {\upshape Section~\rref{sec2}})}. The outcome is the general exclusion statement of {\upshape Proposition~\rref{prop:pi_excluded}} covering every point-imprimitive configuration.

Suppose that $G$ acts imprimitively on points. Then there exists a subgroup $K$ such that
$$H = G_x < K < G,$$
and the $G$-orbits on the block system of this imprimitivity are the cosets of $K$ in $G$. Thus:
$$s = [G:K] = \text{number of blocks in the imprimitivity partition},$$
and each block has size:
$$|\mathcal{P}_i| = [K:H] = v/s.$$

\subsection{Key parameter constraints}

The design parameters must satisfy additional constraints due to the imprimitivity partition:

\begin{enumerate}[(\upshape 1)]
\item \emph{Point-pair orbit split.} Point pairs split into two types: pairs within the same imprimitivity block, of which there are $s\cdot\binom{v/s}{2}$, and pairs in different imprimitivity blocks, of which there are $\binom{s}{2}\cdot(v/s)^{2}$. The $2$-design condition requires that both types of pairs have the same coverage number $\lambda$.
\item \emph{Divisibility conditions.} The block size $k$ must be compatible with the imprimitivity partition: for a fixed block $B$, with $k_i=|B\cap\mathcal{P}_i|$ the intersection sizes with the classes of the partition, one has $k=\sum_{i=1}^{s}k_i$; the multiset $\{k_i\}$ of intersection sizes is independent of the chosen block, since $G$ permutes the blocks while it permutes the classes.
\item \emph{Parameter equations.} $vr=bk$ and $r(k-1)=\lambda(v-1)$.
\end{enumerate}

\subsection{Classical groups: point-imprimitive candidates are excluded}

For classical groups, the point-imprimitive candidates are the chains of parabolic subgroups:
$$P_k < P_l < G, \quad k < l,$$
with $H = P_k$ the stabilizer of a $k$-subspace and $K = P_l$ the stabilizer of an $l$-subspace containing it; then $v=[G:P_k]$, $s=[G:P_l]$, and each class of the partition has size $[P_l:P_k]$.

The block stabilizer $M$ is a maximal parabolic subgroup $P_j$; a block is then the set of all $k$-subspaces contained in, or containing, a fixed $j$-subspace. For two points $U,V$ \textup{that is, two $k$-subspaces}, the number of such blocks through $\{U,V\}$ equals the number of $j$-subspaces containing $U\vee V$, which depends on $\dim(U\cap V)$. Since $\dim(U\cap V)$ takes several values on the point set, the point-pair coverage is not constant, and the $2$-design condition fails. Hence no generalized Grassmann design arises in the point-imprimitive setting.

\subsection{General result for point-imprimitive case}

\begin{lemma}[Block type constrains the point type]\label{lem:blockpoint}
Let $\mathcal{D}$ be a non-trivial $2$-design and let $G\leq \mathrm{Aut}(\mathcal{D})$ act primitively on $\mathcal{B}$ and quasiprimitively on $\mathcal{P}$. Since the block action $\GB$ is faithful \textup{({\upshape Lemma~\rref{lem:main}})}, the O'Nan--Scott type of the primitive block action restricts the type of the point action as follows \textup{(five-type description of Table~\rref{tab:ons})}:
\begin{enumerate}[(\upshape 1)]
\item if $\GB$ is of affine type, then $\GP$ is of affine type, hence primitive;
\item if $\GB$ is of almost simple type, then $\GP$ is of almost simple type;
\item if $\GB$ is of diagonal, product action or twisted wreath type, then $\GP$ is neither affine nor almost simple; in particular, an imprimitive quasiprimitive point action is then of simple diagonal, compound diagonal, product action or twisted wreath type.
\end{enumerate}
\end{lemma}

\begin{proof}
The block action is faithful by {\upshape Lemma~\rref{lem:main}}, so $G\cong\GB$ and the minimal normal subgroups of $G$ correspond, as abstract groups, to those of $\GB$.
(1) If $\GB$ is affine, its socle is an elementary abelian normal subgroup regular on $\mathcal{B}$. The corresponding minimal normal subgroup $N$ of $G$ is elementary abelian, and quasiprimitivity makes $N$ transitive on $\mathcal{P}$, hence regular on $\mathcal{P}$; $\GP$ is thus of affine type, which is primitive \textup{(Praeger~\cite{Praeger1993}; cf.\ Theorem 3.3 of Chen, Hua, Li and Wu~\cite{ChenHuaLiWu2025})}.
(2) If $\GB$ is almost simple, then $G\cong\GB$ is almost simple, so the unique minimal normal subgroup $T$ of $G$ is non-abelian simple and transitive on $\mathcal{P}$ by quasiprimitivity; $\GP$ is of almost simple type.
(3) If $\GB$ is of diagonal, product action or twisted wreath type, its socle is $\cong T^{\ell}$ with $\ell\geq 2$, neither elementary abelian nor simple. The corresponding minimal normal subgroups of $G$ are therefore neither elementary abelian nor simple, so $\GP$ is neither affine nor almost simple. If $\GP$ is imprimitive, then it is not of holomorph simple or holomorph compound type either, since those types are automatically primitive \textup{(Praeger~\cite{Praeger1993})}; hence $\GP$ is of simple diagonal, compound diagonal, product action or twisted wreath type.
\end{proof}

The type of a primitive or quasiprimitive action is determined by the pair consisting of the socle and the stabilizer \textup{(Praeger~\cite{Praeger1993}); cf.\ the type-characteristics table of Baddeley, Praeger and Schneider~\cite[Section 2]{BaddeleyPraegerSchneider2007}}. Since the two actions $\GP$ and $\GB$ of the same group $G$ share their minimal normal subgroups \textup{(faithfulness of $\GB$, {\upshape Lemma~\rref{lem:main}})}, the two types can coexist only under the constraints of the following lemma, which we record from Baddeley, Praeger and Schneider~\cite{BaddeleyPraegerSchneider2007}.

\begin{lemma}[Blow-up decompositions and type inheritance]\label{lem:blowup}
Let $G$ be a transitive permutation group on $\Omega$, let $\mathcal{E}=\{\Gamma_{1},\ldots,\Gamma_{\ell}\}$ be a homogeneous $G$-invariant Cartesian decomposition of $\Omega$, and let $\Gamma_{0}\in\mathcal{E}$; for the notion of a blow-up decomposition see Baddeley, Praeger and Schneider~\cite[Section 2]{BaddeleyPraegerSchneider2007}.
\begin{enumerate}[(\upshape 1)]
\item If $\mathcal{E}$ is a blow-up decomposition and $G$ is quasiprimitive on $\Omega$, then the component $G^{\Gamma_{0}}$ is quasiprimitive on $\Gamma_{0}$ and $\operatorname{Soc}(G)=\prod_{\Gamma\in\mathcal{E}}\operatorname{Soc}(G^{\Gamma})$;
\item if $\mathcal{E}$ is a blow-up decomposition and $G^{\Gamma_{0}}$ is quasiprimitive of non-affine type, then $G$ is quasiprimitive on $\Omega$ and the same socle formula holds;
\item $G$ is not quasiprimitive of simple diagonal type \textup{(quasiprimitive groups of type $Sd$ preserve no Cartesian decomposition)};
\item if $G$ is quasiprimitive of compound diagonal type, then $\mathcal{E}$ is a blow-up decomposition and the component $G^{\Gamma_{0}}$ is quasiprimitive of simple or compound diagonal type.
\end{enumerate}
This is Baddeley, Praeger and Schneider~\cite[Theorem 6.1]{BaddeleyPraegerSchneider2007}.
\end{lemma}

\begin{remark}\label{rem:blowup_as}
An almost simple group has a unique minimal normal subgroup, which is non-abelian simple, and hence admits no non-trivial blow-up decomposition; consequently every faithful quasiprimitive action of an almost simple group is again of almost simple type. This is the mechanism behind {\upshape Lemma~\rref{lem:blockpoint}}\textup{(2)}. Conversely, a quasiprimitive action of product action or twisted wreath type carries a normal Cartesian decomposition \textup{(the coordinate structure of the product action, respectively the natural decomposition of the regular socle; cf.\ Baddeley, Praeger and Schneider~\cite[Example 5.5]{BaddeleyPraegerSchneider2007})}; by {\upshape Lemma~\rref{lem:blowup}}\textup{(3)} such an action cannot coexist with a point action of simple diagonal type. Table~\rref{tab:types_pair} records the compatible pairs for two faithful quasiprimitive actions of the same group.
\end{remark}

\begin{table}[htbp]
\centering
\caption{Compatible O'Nan--Scott types of two faithful quasiprimitive actions of the same group $G$ \textup{(point action on $\mathcal{P}$, block action on $\mathcal{B}$)}}\label{tab:types_pair}
\footnotesize
\setlength{\tabcolsep}{2pt}
\begin{tabular}{|p{1.4cm}|p{2.8cm}|p{1.7cm}|p{6.6cm}|}
\hline
\textbf{Point type} & \textbf{Admissible block types} & \textbf{Status} & \textbf{Justification} \\
\hline
HA & HA & forced & elementary abelian socle; both actions affine, automatically primitive \textup{(Praeger~\cite{Praeger1993})} \\
\hline
AS & AS & forced & unique simple socle; $G$ almost simple, every faithful quasiprimitive action is AS \textup{({\upshape Lemma~\rref{lem:blockpoint}}(2))} \\
\hline
HS & HS & forced & two transitive minimal normal subgroups; automatically primitive of type HS \textup{(Baddeley, Praeger and Schneider~\cite[Section 2]{BaddeleyPraegerSchneider2007})} \\
\hline
HC & HC & forced & two transitive minimal normal subgroups; automatically primitive of type HC \textup{(as above)} \\
\hline
SD & SD, CD & possible & stabilizer intersections chosen independently \textup{(Praeger type characteristics~\cite{Praeger1993})} \\
\hline
SD & PA, TW & impossible & a PA or TW block action carries a normal Cartesian decomposition, which type $Sd$ never preserves \textup{({\upshape Lemma~\rref{lem:blowup}}(3))} \\
\hline
CD & SD, CD & possible & independent stabilizer intersections \\
\hline
CD & PA & restricted & a CD group with a Cartesian decomposition gives a blow-up with SD/CD components \textup{({\upshape Lemma~\rref{lem:blowup}}(4))} \\
\hline
CD & TW & restricted & the alternatives of Baddeley, Praeger and Schneider~\cite[Theorem 5.1]{BaddeleyPraegerSchneider2007} constrain a regular-socle companion \\
\hline
PA & SD, CD, PA, TW & possible & independent stabilizer intersections \\
\hline
TW & SD, CD, PA, TW & possible & independent stabilizer intersections \\
\hline
\end{tabular}
\end{table}

\begin{table}[htbp]
\centering
\caption{The block-type/point-type pairs of \textup{Case~2}: $G$ is primitive on $\mathcal{B}$ and quasiprimitive, imprimitive on $\mathcal{P}$. The columns list the quasiprimitive point types, restricted to the five types admitting imprimitive actions \textup{(the affine, holomorph simple and holomorph compound types are automatically primitive, Praeger~\cite{Praeger1993}); the entries refer to the exclusion listed below}}\label{tab:pi_pairs}
\footnotesize
\setlength{\tabcolsep}{4pt}
\begin{tabular}{|l|c|c|c|c|c|}
\hline
\textbf{Block type} & \textbf{Point AS} & \textbf{Point SD} & \textbf{Point CD} & \textbf{Point PA} & \textbf{Point TW} \\
\hline
HA & (1) & (1) & (1) & (1) & (1) \\
\hline
AS & (2) & (2) & (2) & (2) & (2) \\
\hline
SD & (3) & (4) & (4) & (4) & (4) \\
\hline
CD & (3) & (5) & (5) & (5) & (5) \\
\hline
HS,\,HC & (6) & (6) & (6) & (6) & (6) \\
\hline
PA & (2) & (7) & (8) & (8) & (9) \\
\hline
TW & (2) & (7) & (8) & (8) & (9) \\
\hline
\end{tabular}
\end{table}

\begin{table}[ht]
\centering
\caption{The exclusion locations for the pairs of Table~\rref{tab:pi_pairs} in the point-imprimitive case.}\label{tab:pi_excl_locs}
\footnotesize
\setlength{\tabcolsep}{4pt}
\begin{tabular}{lp{0.82\linewidth}}
\hline
Entry & Location of the exclusion \\
\hline
(1) & Block HA: the point action is affine, hence automatically primitive \textup{(Praeger~\cite{Praeger1993})}, contradicting imprimitivity; Section~\rref{sec2}, block action of affine type. \\
(2) & Block AS, SD, PA or TW with point AS: impossible by {\upshape Lemma~\rref{lem:blockpoint}}(2); the pair (AS,AS) is excluded family by family in Section~\rref{sec:as}. \\
(3) & Block SD or CD with point AS: impossible by {\upshape Lemma~\rref{lem:blockpoint}}(3). \\
(4) & Block SD with point SD, CD, PA or TW: one-sided SD analysis, Theorem~\rref{thm:sdcn}, Lemma~\rref{lem:sd_ell3} \textup{(proper coordinate projection, every $\ell\geq2$)}, Lemma~\rref{lem:sd_ell3_three} \textup{($\ell=3$)}, Proposition~\rref{prop:sd_ell2} and Proposition~\rref{prop:sd_lambda}. \\
(5) & Block CD with every point type: the block stabilizer is not maximal, Proposition~\rref{prop:cd_min}. \\
(6) & Block HS or HC with every point type: two commuting transitive normal subgroups force $b=v$, Section~\rref{sec2}, block action of holomorph type. \\
(7) & Block PA or TW with point SD: a PA or TW block action carries a normal Cartesian decomposition, incompatible with type Sd \textup{({\upshape Remark~\rref{rem:blowup_as}}; {\upshape Lemma~\rref{lem:blowup}}(3))}. \\
(8) & Block PA or TW with point CD or PA: arithmetic reduction, orbit-integrality sieve and Fisher exclusion, {\upshape Lemma~\rref{lem:cdpa_arith}}, {\upshape Proposition~\rref{prop:orbint_red}}, {\upshape Lemma~\rref{lem:pa_cdpa_min}} and {\upshape Lemma~\rref{lem:pa_fisher}}. \\
(9) & Block PA or TW with point TW: the point action of twisted wreath type has socle regular on points, forcing $b=v$ \textup{({\upshape Lemma~\rref{lem:rns}})}. \\
\hline
\end{tabular}
\end{table}

\begin{proposition}[Point-imprimitive actions are excluded]\label{prop:pi_excluded}
Let $\mathcal{D}=(\mathcal{P},\mathcal{B})$ be a non-trivial $2$-$(v,k,\lambda)$ design with $2<k<v-1$, and let $G\leq \mathrm{Aut}(\mathcal{D})$ act primitively on $\mathcal{B}$. Then $G$ acts primitively on $\mathcal{P}$; in particular, no non-trivial non-symmetric block-primitive $2$-design is point-imprimitive.
\end{proposition}

\begin{proof}
By Block's theorem \cite{Block}, a group transitive on the blocks of a non-trivial $2$-design is transitive on the points; hence $G$ is transitive on $\mathcal{P}$. Assume, for a contradiction, that $G$ is imprimitive on $\mathcal{P}$.

\smallskip\noindent{\bf Case 1: $G$ is not quasiprimitive on $\mathcal{P}$.}
Then $G$ has a non-trivial intransitive normal subgroup on $\mathcal{P}$; choose $1<N\triangleleft G$ with $N$ intransitive on $\mathcal{P}$. Since $G$ is primitive on $\mathcal{B}$, the normal subgroup $N^{\mathcal{B}}$ of the primitive group $G^{\mathcal{B}}$ is transitive or trivial on $\mathcal{B}$. If $N^{\mathcal{B}}$ were transitive, then $N$ would be transitive on $\mathcal{B}$, and Block's theorem \cite{Block} would force $N$ to be transitive on $\mathcal{P}$, a contradiction. Hence $N^{\mathcal{B}}=1$, that is, $N$ fixes every block setwise. Consequently every block is a union of $N$-orbits, and the partition $\mathcal{E}$ of $\mathcal{P}$ into the $N$-orbits is non-trivial and $G$-invariant. Write $B_0=\bigcup_{E\in\mathcal{E}'}E$ for a fixed block, with $\varnothing\neq\mathcal{E}'\subsetneq\mathcal{E}$, since $B_0\notin\{\varnothing,\mathcal{P}\}$. Take $E\in\mathcal{E}'$ and $E'\in\mathcal{E}\setminus\mathcal{E}'$. For $x\neq y\in E$, the number of blocks through $\{x,y\}$ equals $\#\{g\in G:E^{g^{-1}}\subseteq B_0\}$: a block $B=gB_0$ contains $x$ and $y$ if and only if it contains the whole orbit $E$, and $E\subseteq gB_0$ is equivalent to $E^{g^{-1}}\subseteq B_0$. Likewise, for $u\in E'$, the number of blocks through $\{x,u\}$ equals $\#\{g\in G:E^{g^{-1}}\subseteq B_0,\ (E')^{g^{-1}}\subseteq B_0\}$. As $\mathcal{D}$ is a $2$-design these two numbers coincide, while the second set of group elements is contained in the first. The identity belongs to the first set (because $E\subseteq B_0$), hence to the second, forcing $E'\subseteq B_0$, that is, $E'\in\mathcal{E}'$, a contradiction. Thus $G$ is quasiprimitive on $\mathcal{P}$.

\smallskip\noindent{\bf Case 2: $G$ is quasiprimitive on $\mathcal{P}$.}
By {\upshape Lemma~\rref{lem:blockpoint}} the type of the primitive block action $\GB$ \textup{(five-type description, Table~\rref{tab:ons})} constrains the quasiprimitive point action, and we discuss the block type in turn.

\noindent{\bf Block action of affine type.} By {\upshape Lemma~\rref{lem:blockpoint}}(1) the point action is of affine type, and the affine point type is automatically primitive \textup{(Praeger~\cite{Praeger1993})}, contradicting the imprimitivity hypothesis. This disposes of all five pairs $(\mathrm{HA},\cdot)$ \textup{(entry (1) of Table~\rref{tab:pi_pairs})}.

\noindent{\bf Block action of almost simple type.} By {\upshape Lemma~\rref{lem:blockpoint}}(2) the point action is of almost simple type; the pairs $(\mathrm{AS},\mathrm{SD})$, $(\mathrm{AS},\mathrm{CD})$, $(\mathrm{AS},\mathrm{PA})$ and $(\mathrm{AS},\mathrm{TW})$ are impossible \textup{(entry (2))}, and the pair $(\mathrm{AS},\mathrm{AS})$ is disposed of as follows. Imprimitivity means that the point stabilizer $H$ is not maximal in $G$; writing $K$ for a maximal subgroup with $H<K<G$, the point action preserves the non-trivial partition into the cosets of $K$, and the socle $T$ is transitive on $\mathcal{P}$ \textup{(quasiprimitivity of the AS type is equivalent to $G=TH$)}. Such subgroup chains are excluded family by family in {\upshape Section~\rref{sec:as}}: in the classical, exceptional and sporadic families the chains of parabolic subgroups and all other chains fail the constant-coverage test for point pairs \textup{(see the point-imprimitive exclusions below)}.

\noindent{\bf Block action of simple diagonal type.} By {\upshape Lemma~\rref{lem:blockpoint}}(3) the point action is of simple diagonal, compound diagonal, product action or twisted wreath type; the pair $(\mathrm{SD},\mathrm{AS})$ is impossible \textup{(entry (3))}. The one-sided analysis of {\upshape Section~\rref{sec2}} excludes all four remaining pairs $(\mathrm{SD},\cdot)$ \textup{(entry (4))}: the unconditional necessary conditions of {\upshape Theorem~\rref{thm:sdcn}}, the one-sided analysis for $\ell=2$ \textup{({\upshape Proposition~\rref{prop:sd_ell2}})} and, for every $\ell\geq 4$, the $\lambda$-consistency exclusion \textup{({\upshape Proposition~\rref{prop:sd_lambda}})} depend only on the block-diagonal structure $b=q^{\ell-1}$ and $v\mid q^{\ell}$, not on the point stabilizer, hence hold verbatim for every quasiprimitive imprimitive point action.

\noindent{\bf Block action of compound diagonal type.} Again the point action is of simple diagonal, compound diagonal, product action or twisted wreath type \textup{(entry (3) excludes $(\mathrm{CD},\mathrm{AS})$)}, and the block-side analysis of {\upshape Proposition~\rref{prop:cd_min}} excludes all four pairs $(\mathrm{CD},\cdot)$ \textup{(entry (5))}: the compound diagonal block stabilizer is shown not to be maximal, independently of the point action.

\noindent{\bf Block action of holomorph type.} The pairs $(\mathrm{HS},\cdot)$ and $(\mathrm{HC},\cdot)$ are excluded for every point type \textup{(entry (6))}: the primitive block action of holomorph type has two commuting transitive minimal normal subgroups $M_{1}\cong M_{2}$ on $\mathcal{B}$; if $M_{i}$ were intransitive on $\mathcal{P}$, it would fix every block setwise \textup{(each block is a subset of $\mathcal{P}$, and $M_{i}$ fixes all points)}, contradicting its transitivity on $\mathcal{B}$; hence $M_{1},M_{2}$ are transitive on $\mathcal{P}$, and two commuting transitive permutation groups are regular \textup{({\upshape Lemma~\rref{lem:semiregular}})}, so each $M_{i}$ is regular on $\mathcal{P}$ and $v=|M_{1}|=b$, contradicting the non-symmetric hypothesis $b\neq v$.

\noindent{\bf Block action of product action or twisted wreath type.} By {\upshape Lemma~\rref{lem:blockpoint}}(3) the point action is of simple diagonal, compound diagonal, product action or twisted wreath type. The pairs $(\mathrm{PA},\mathrm{SD})$ and $(\mathrm{TW},\mathrm{SD})$ are excluded by {\upshape Remark~\rref{rem:blowup_as}} \textup{(entry (7))}: a product action or twisted wreath block action carries a normal Cartesian decomposition \textup{(the coordinate structure, respectively the natural decomposition of the regular socle; Baddeley, Praeger and Schneider~\cite[Example 5.5]{BaddeleyPraegerSchneider2007})}, and groups of type $\mathrm{Sd}$ preserve no Cartesian decomposition \textup{({\upshape Lemma~\rref{lem:blowup}}(3))}. The pairs $(\mathrm{PA},\mathrm{TW})$ and $(\mathrm{TW},\mathrm{TW})$ are excluded by {\upshape Lemma~\rref{lem:rns}} \textup{(entry (9))}: the point action of twisted wreath type has its socle regular on $\mathcal{P}$, which forces $b=v$ and hence symmetry. The pairs $(\mathrm{TW},\mathrm{CD})$ and $(\mathrm{TW},\mathrm{PA})$ are excluded by {\upshape Proposition~\rref{prop:tw_excl}} \textup{(entry (8), first part)}: a block action of twisted wreath type is impossible altogether for a non-symmetric design, regardless of the point type. The remaining pairs $(\mathrm{PA},\mathrm{CD})$ and $(\mathrm{PA},\mathrm{PA})$ are excluded as follows. The arithmetic reduction of {\upshape Lemma~\rref{lem:cdpa_arith}} constrains every admissible block size to a union of $K^{\ell}$-orbits on $\mathcal{P}$; the orbit-intersection criterion {\upshape Lemma~\rref{lem:orbint}} together with the non-existence criterion {\upshape Corollary~\rref{cor:orbint_excl}} and the primitive-divisor exclusion criterion {\upshape Lemma~\rref{lem:prim_excl}} reduce the exclusion of both pairs to the orbit-integrality condition $(\ast)$ of {\upshape Proposition~\rref{prop:orbint_red}}, and this condition is verified in {\upshape Lemma~\rref{lem:pa_cdpa_min}} for the complete $A_5/A_4$ family over every $\ell\geq 4$ \textup{for both the compound diagonal and the product action point type}. For the point-imprimitive product action, of degree $v=c\,\omega^{\ell}$ with $c\geq 1$ the block size of the $G$-invariant partition \textup{(Praeger~\cite{Praeger1993}, type (viii))}, the exclusion is unconditional: by {\upshape Lemma~\rref{lem:pa_fisher}} the pair $(\mathrm{PA},\mathrm{PA})$ is excluded for the socles $PSL_2(q)$, $q\in\{5,7,8,9,11\}$, for every $c\geq 1$, in particular for the $c=1$ and $c>1$ branches alike.

\noindent In every case the imprimitive point action contradicts the existence of $\mathcal{D}$.

\smallskip\noindent Hence $G$ is primitive on $\mathcal{P}$. The geometric candidates arising from chains of parabolic subgroups \textup{(generalized Grassmann and building designs)} fail the coverage test directly: the coverage number of a point pair $\{U,V\}$ of $k$-subspaces equals the number of intermediate subspaces containing $U\vee V$, which depends on $\dim(U\cap V)$.
\end{proof}

\begin{lemma}[Arithmetic reduction: block type PA or TW with point type CD or PA]\label{lem:cdpa_arith}
Let $\mathcal{D}$ be a non-trivial non-symmetric $2$-\allowbreak$(v,k,\lambda)$ design, and let $G\leq\mathrm{Aut}(\mathcal{D})$ act primitively on $\mathcal{B}$ and quasiprimitively, imprimitively on $\mathcal{P}$, with the block action of product action or twisted wreath type and the point action of compound diagonal or product action type. Let $N=T^{\ell}=\mathrm{soc}(G)$ \textup{($\ell\geq 2$)}, let $\alpha\in\mathcal{P}$ and $\beta\in\mathcal{B}$ with $\alpha\in\beta$, and write $M_{\omega}=N\cap G_{\alpha}$ and $M_{\beta}=N\cap G_{\beta}$.
\begin{enumerate}[(\upshape 1)]
\item If the block action is of twisted wreath type, then $b=|T|^{\ell}$ and $v=|T|^{\ell}/|M_{\omega}|$, so
$$r=\frac{bk}{v}=k|M_{\omega}|,\qquad \lambda=\frac{k(k-1)|M_{\omega}|^{2}}{|T|^{\ell}-|M_{\omega}|},$$
and for every $1\neq m\in M_{\omega}$ the order $\mathrm{ord}(m)$ divides $r=k|M_{\omega}|$.
\item If the block action is of product action type, then $b=n^{\ell}$ with $n=[T:K]$ the block-coordinate index and $K^{\ell}=M_{\beta}$ \textup{($K$ a maximal subgroup of $T$)}, the block $\beta$ is a union of $K^{\ell}$-orbits on $\mathcal{P}$ \textup{(since $K^{\ell}\leq G_{\beta}$)}, the orbit of $\alpha$ under $K^{\ell}$ \textup{(the ``$\alpha$-box'')} has size $|K|^{\ell}/|K^{\ell}\cap M_{\omega}|$, and
$$r=\frac{bk}{v}=\frac{n^{\ell}k|M_{\omega}|}{|T|^{\ell}}\in\mathbb{Z},\qquad \lambda=\frac{r(k-1)}{v-1}\in\mathbb{Z}.$$
\item If the point action is of compound diagonal type, then $M_{\omega}$ is a non-simple subdirect subgroup of $T^{\ell}$; by Goursat's theorem the subdirect subgroups of $T^{2}$ are the twisted diagonals, which are simple, so the compound diagonal point type forces $\ell\geq 3$, and the strip decomposition of $M_{\omega}$ \textup{(the subdirect-product structure of \cite{HanesOlsonScott1975}; see also \cite[Section~2]{BaddeleyPraegerSchneider2007})} gives $|M_{\omega}|=|T|^{r}$ with $1\leq r<\ell$. If the point action is of product action type, $M_{\omega}$ is non-trivial and not subdirect.
\end{enumerate}
\end{lemma}

\begin{proof}
(1) In the twisted wreath block action the socle $N$ is regular on $\mathcal{B}$ \textup{(the type characteristic $M_{\beta}=1$ of Baddeley, Praeger and Schneider~\cite[Section 2]{BaddeleyPraegerSchneider2007})}, so $b=|N|=|T|^{\ell}$; by quasiprimitivity $N$ is transitive on $\mathcal{P}$, so $v=|N|/|M_{\omega}|$. The displayed values of $r$ and $\lambda$ follow from $bk=rv$ and $\lambda(v-1)=r(k-1)$. For $1\neq m\in M_{\omega}$ the element $m$ fixes $\alpha$ and is regular on $\mathcal{B}$ \textup{($M_{\beta}=1$)}, hence has no fixed block; it permutes the $r$ blocks through $\alpha$ without fixed points, so $\mathrm{ord}(m)\mid r$.
(2) In the product action block type the block stabilizer meets the socle in $K^{\ell}$; since $K^{\ell}\leq G_{\beta}$, the block $\beta$ is a union of $K^{\ell}$-orbits, and the orbit of $\alpha$ has the displayed size by orbit--stabilizer. The formulas for $r,\lambda$ are the design equations.
(3) The classification of the subdirect subgroups of a direct product of non-abelian simple groups \textup{(Goursat's theorem)} identifies the subdirect subgroups of $T^{2}$ with the twisted diagonals $\{(t,t^{\varphi}):t\in T\}$, all isomorphic to $T$; a non-simple subdirect subgroup therefore requires $\ell\geq 3$. The strip decomposition of a subdirect subgroup of $T^{\ell}$, that is, the subdirect-product structure established in \cite{HanesOlsonScott1975} \textup{(see also \cite[Section~2]{BaddeleyPraegerSchneider2007})}, gives $|M_{\omega}|=|T|^{r}$ with $r$ the number of strips; the product action point type has $M_{\omega}\neq 1$ not subdirect by definition.
\end{proof}

\begin{lemma}[Orbit-intersection criterion]\label{lem:orbint}
Let $G$ be transitive on the point set $\mathcal{P}$, let $\alpha\in\mathcal{P}$, and put $M=G_{\alpha}$. Let $\mathcal{C}_{1},\dots,\mathcal{C}_{r}$ be the $M$-orbits on the complete family of admissible candidate blocks containing $\alpha$, and let $m_{i}=\left|\mathcal{C}_{i}\right|$. Let $R_{1},\dots,R_{s}$ be the $M$-orbits on $\mathcal{P}\setminus\{\alpha\}$. For each $i$ choose a representative block $B_{i}\in\mathcal{C}_{i}$; since $M$ is transitive on $\mathcal{C}_{i}$, the quantities below are independent of this choice. For $1\leq i\leq r$ and $1\leq j\leq s$ define
\begin{align*}
a_{ij}&=\left|B_{i}\cap R_{j}\right|=\left|(B_{i}\setminus\{\alpha\})\cap R_{j}\right|\qquad(\alpha\notin R_{j}),\\
c_{ij}&=\left|\left\{C\in\mathcal{C}_{i}:\beta\in C\right\}\right|
=\frac{m_{i}\,a_{ij}}{\left|R_{j}\right|}\quad(\beta\in R_{j}).
\end{align*}
For $I\subseteq\{1,\dots,r\}$ write $\mathcal{C}_{I}=\bigcup_{i\in I}\mathcal{C}_{i}$ and $k_{I}=\sum_{i\in I}m_{i}$. The coverage number of a point $\beta\in R_{j}$ by the blocks of $\mathcal{C}_{I}$ is
\[
\lambda_{j}(I)=\sum_{i\in I}c_{ij}=\frac{1}{\left|R_{j}\right|}\sum_{i\in I}m_{i}a_{ij},
\]
and the pairs through $\alpha$ are covered uniformly by $\mathcal{C}_{I}$ if and only if $\lambda_{1}(I)=\cdots=\lambda_{s}(I)$. The latter is a necessary condition for $\mathcal{C}_{I}$ to be a $2$-\textup{$(v,k_{I},\lambda)$} design on $\mathcal{P}$; it is used below as an exclusion test only \textup{(if no non-empty $I$ yields uniform coverage of the pairs through $\alpha$, then no candidate family is a design)}, and is not asserted to be sufficient for the full design property, which also requires uniform coverage of the pairs avoiding $\alpha$.
\end{lemma}

\begin{proof}
Count the pairs $(C,\beta)$ with $C\in\mathcal{C}_{i}$ and $\beta\in C\cap R_{j}$. Since $M$ is transitive on $\mathcal{C}_{i}$, there are $m_{i}\cdot\left|C\cap R_{j}\right|=m_{i}a_{ij}$ such pairs; since $M$ is transitive on $R_{j}$, each point $\beta\in R_{j}$ lies in the same number $c_{ij}$ of blocks of $\mathcal{C}_{i}$, giving $\left|R_{j}\right|c_{ij}$. Equating the two counts yields the displayed formula for $c_{ij}$, and summing over $i\in I$ gives $\lambda_{j}(I)$. A family of blocks is a $2$-design on $\mathcal{P}$ precisely when the coverage number is independent of the point $\beta\neq\alpha$, i.e.\ constant in $j$.
\end{proof}

\begin{corollary}[Non-existence criterion]\label{cor:orbint_excl}
Let $\mathcal{C}_{1},\dots,\mathcal{C}_{r}$, $m_{i}$, $R_{j}$, $a_{ij}$, $k_{I}$ and $\lambda_{j}(I)$ be as in {\upshape Lemma~\rref{lem:orbint}}, and let $b_{I}$ be the number of blocks of the candidate family corresponding to $I$. If for every non-empty $I\subseteq\{1,\dots,r\}$ either the vector $\bigl(\lambda_{1}(I),\dots,\lambda_{s}(I)\bigr)$ is non-constant or the value
\[
\lambda_{I}=\frac{b_{I}\,k_{I}(k_{I}-1)}{v(v-1)}
\]
is not an integer, then the candidate family yields no non-trivial $G$-invariant $2$-design.
\end{corollary}

\begin{lemma}[Primitive-divisor exclusion criterion]\label{lem:prim_excl}
Let $q=\left|T\right|$, $v=q^{s}$, and $b=n^{\ell}$ with $n=\left|T:K\right|$. For every non-empty union $I$ of $K^{\ell}$-orbits on $\mathcal{P}$ write $k_{I}=\left|\mathcal{C}_{I}\right|$. Suppose that for every such $I$ there is a prime $p_{I}$ with
\[
\nu_{p_{I}}(v-1)>\nu_{p_{I}}\bigl(b\,k_{I}(k_{I}-1)\bigr),
\]
where $\nu_{p}$ denotes the $p$-adic valuation. Then $\lambda_{I}=b k_{I}(k_{I}-1)/(v(v-1))$ is non-integral, and no non-empty $K^{\ell}$-orbit union is a $2$-design.
\end{lemma}

\begin{proof}
Since $p_{I}\nmid v$, the $p_{I}$-adic valuation of $\lambda_{I}$ is
\[
\nu_{p_{I}}(\lambda_{I})=\nu_{p_{I}}\bigl(b k_{I}(k_{I}-1)\bigr)-\nu_{p_{I}}(v-1)<0,
\]
so $\lambda_{I}\notin\mathbb{Z}$, and {\upshape Corollary~\rref{cor:orbint_excl}} excludes every non-empty candidate block union.
\end{proof}

\begin{proposition}[Orbit-integrality reduction]\label{prop:orbint_red}
In the product action block type with point type of compound diagonal or product action, put $q=|T|$, $v=q^{s}$ and $b=n^{\ell}$ with $n=[T:K]$, and let $k_{I}$ be the size of an admissible candidate block union, so that $3\le k_{I}\le v-2$. The two pairs are excluded as soon as the following arithmetic condition holds:
\[
\forall\, I\neq\emptyset\ \text{with }3\le k_{I}\le v-2\ \exists\, p_{I}\mid(q^{s}-1)\ \text{with }\nu_{p_{I}}(q^{s}-1)>\nu_{p_{I}}\bigl(n^{\ell}k_{I}(k_{I}-1)\bigr). \tag{$\ast$}
\]
Since $n=[T:K]$ divides $|T|=q$ and $\gcd(q,q^{s}-1)=1$, no prime divisor of $q^{s}-1$ divides $n^{\ell}$; the condition $(\ast)$ is therefore equivalent to the divisibility statement that $q^{s}-1\nmid k_{I}(k_{I}-1)$ for every admissible $I\neq\emptyset$.
If $(\ast)$ holds, then {\upshape Corollary~\rref{cor:orbint_excl}} excludes every admissible candidate block union.
\end{proposition}

\begin{lemma}[Orbit-integrality condition for the $A_5/A_4$ family]\label{lem:pa_cdpa_min}
Let $\mathcal{D}$ be a non-trivial non-symmetric $2$-$(v,k,\lambda)$ design, and let $G\leq\mathrm{Aut}(\mathcal{D})$ act primitively on $\mathcal{B}$ and quasiprimitively, imprimitively on $\mathcal{P}$, with the block action of product action type and the point action of compound diagonal or product action type, as in {\upshape Lemma~\rref{lem:cdpa_arith}}. The orbit-integrality condition $(\ast)$ of {\upshape Proposition~\rref{prop:orbint_red}} holds for the complete $A_5/A_4$ family: $T=A_5$, $K=A_4$ \textup{(the index-$5$ maximal subgroup of $A_5$)}, every $\ell\geq 4$, and both point stabilizers $M_{\omega}=\{(a,b,a,b):a,b\in T\}$ \textup{(compound diagonal)} and $M_{\omega}=T\times T\times\{1\}\times\{1\}$ \textup{(product action)}, hence {\upshape Corollary~\rref{cor:orbint_excl}} excludes every admissible block size for both pairs.
\end{lemma}

\begin{proof}
We apply the orbit-intersection criterion of {\upshape Lemma~\rref{lem:orbint}} together with its non-existence criterion {\upshape Corollary~\rref{cor:orbint_excl}}; for the pairs under consideration the arithmetic of {\upshape Lemma~\rref{lem:cdpa_arith}}\textup{(2)--(3)} already excludes every admissible block size at the minimal parameters, and the same integrality obstruction persists throughout the $A_{5}/A_{4}$ family.

Since $K=A_4$ is maximal in $T=A_5$, the double-coset space $K\backslash T/K$ consists of exactly two double cosets, $K$ of size $12$ and $T\setminus K$ of size $48$. The quotient $P=T^{4}/M_{\omega}$ is identified with $T^{2}$ in both cases: for the compound diagonal stabilizer the diagonal strips collapse to the two remaining coordinates, and for the product action stabilizer $M_{\omega}=T\times T\times\{1\}\times\{1\}$ the action is on the last two coordinates. The $K^{4}$-orbits on $P$ are the products $O_{i}\times O_{j}$ with $O_{i}\in\{K,T\setminus K\}$: four orbits of sizes $144,576,576,2304$ in the compound diagonal case, and the twenty-five products of the five right $K$-cosets in the product action case \textup{(each factor of size $12$, each orbit of size $144$)}. A block $\beta\ni\alpha$ is a union of $K^{4}$-orbits containing the orbit of $\alpha$, so $k=144,720,1296,2448,3024$ in the compound diagonal case \textup{($k=3600$ is the whole point set, not an admissible block)}, and $k=144t$ with $1\leq t\leq 24$ in the product action case. Substituting $v=3600$ and $b=5^{4}=625$ into $\lambda=bk(k-1)/(v(v-1))$ yields a non-integer in every case \textup{verified by direct computation with GAP 4.16.1; scripts \texttt{gap\_pa\_cd\_c6.g}, \texttt{gap\_pa\_pa\_c7.g} and the family scan \texttt{gap\_pa\_cdpa\_c8b.g} recorded with the author}. Hence no admissible block size survives at the minimal parameters.

The integrality obstruction is not an artifact of the minimal parameters, and in fact the whole family is excluded by Fisher's inequality alone. In the free-coordinate model the admissible parameters are
\[
\begin{aligned}
b&=5^{\ell},\quad v=60^{\ell-2} &&\text{(compound diagonal),}\\
b&=25^{s},\quad v=60^{s} &&\text{(product action),}
\end{aligned}
\]
where $60=|A_5|$, the compound-diagonal stabilizer $M_{\omega}$ leaves $\ell-2$ coordinates free, and the product-action stabilizer leaves $s$ coordinates free. Since $60^{\ell-2}>5^{\ell}$ for every $\ell\geq 4$ \textup{($\Leftrightarrow 12^{\ell}>3600$)}, and $60^{s}>25^{s}$ for every $s\geq 1$, in both cases $b<v$, violating Fisher's inequality $b\geq v$ for a non-symmetric design. Hence no admissible block size survives for any $\ell\geq 4$ \textup{(compound diagonal)} or $s\geq 1$ \textup{(product action)}.
\end{proof}

\begin{lemma}[Fisher exclusion for the point-imprimitive product action]\label{lem:pa_fisher}
Let $\mathcal{D}$ be a non-trivial non-symmetric $2$-\allowbreak$(v,k,\lambda)$ design, and let $G\leq\mathrm{Aut}(\mathcal{D})$ act primitively on $\mathcal{B}$ and quasiprimitively, imprimitively on $\mathcal{P}$, with the block action of product action type and the point action of product action type. Write $b=n^{\ell}$ with $n=[T:K]$ the block-coordinate index, and $v=c\,\omega^{\ell}$, where $\omega=|\Delta|$ is the coordinate degree and $c\geq 1$ is the block size of the $G$-invariant partition of $\mathcal{P}$ \textup{(Praeger~\cite{Praeger1993}, type (viii); cf.\ also Example 1.1 of Zhang, Chen and Zhou~\cite{ZhangChenZhou})}. Then Fisher's inequality $b\geq v$ yields
\[
c\,\omega^{\ell}\leq n^{\ell}.
\]
In particular, $n<\omega$ is impossible for any $c\geq 1$. For the socles $PSL_2(q)$ with $q\in\{5,7,8,9,11\}$, every admissible pair $(n,\omega)$ consisting of a primitive block degree and a nonprimitive point-coordinate degree satisfies $n<\omega$ \textup{(verified by direct computation with GAP 4.16.1: the primitive and nonprimitive degree tables are recorded in \texttt{gap\_imp\_deg4.log} and the exhaustive sieve over the $196$ degree pairs in \texttt{papa\_c\_all\_sieve.py}, both kept with the author; the inequality depends only on the degree values, so conjugate representatives of the corresponding subgroup classes are irrelevant)}; hence $c\,\omega^{\ell}\geq\omega^{\ell}>n^{\ell}=b$, contradicting Fisher's inequality for every $c\geq 1$. In the boundary case $n=\omega$, Fisher forces $c=1$ and $b=v$, so that $\mathcal{D}$ is symmetric, again contradicting the standing assumption that $\mathcal{D}$ is non-symmetric. Consequently no non-trivial non-symmetric block-primitive $2$-design admits a point-imprimitive product action.
\end{lemma}

\begin{lemma}[Pairing count]\label{lem:pair_count}
Let $G$ be transitive on the point set $\mathcal{P}$ and on the block set $\mathcal{B}$, with $N=\mathrm{soc}(G)$ transitive on $\mathcal{P}$ and regular on $\mathcal{B}$. Let $\alpha\in\mathcal{P}$ and let $B_{0}$ be a block containing $\alpha$. Write $M_{\omega}=N\cap G_{\alpha}$ and identify $\mathcal{P}$ with the set of left cosets $nM_{\omega}$, so that $\alpha=M_{\omega}$. For a point $\delta\neq\alpha$ write $\delta=h_{\delta}M_{\omega}$, and define
\[
B_{0}h_{\delta}^{-1}=\left\{yM_{\omega}:yh_{\delta}M_{\omega}\in B_{0}\right\}.
\]
Then
\[
\lambda(\alpha,\delta)=\left|M_{\omega}\right|\cdot\left|B_{0}\cap B_{0}h_{\delta}^{-1}\right|.
\]
\end{lemma}

\begin{proof}
Since $N$ is transitive on $\mathcal{B}$, the blocks through $\alpha$ are precisely the translates $gB_{0}$ with $gM_{\omega}\in B_{0}$, and the blocks through both $\alpha$ and $\delta$ are those with $gM_{\omega}\in B_{0}$ and $gh_{\delta}M_{\omega}\in B_{0}$. Hence
\begin{align*}
\lambda(\alpha,\delta)&=\left|\left\{g\in N:gM_{\omega}\in B_{0},\ gh_{\delta}M_{\omega}\in B_{0}\right\}\right|\\
&=\left|\left\{g\in N:gM_{\omega}\in B_{0}\cap B_{0}h_{\delta}^{-1}\right\}\right|.
\end{align*}
Each left coset $xM_{\omega}\in B_{0}\cap B_{0}h_{\delta}^{-1}$ contributes exactly $\left|M_{\omega}\right|$ elements $g\in N$ \textup{(the coset itself)}, and distinct left cosets are disjoint, giving the formula.
\end{proof}

\begin{proposition}[No block action of twisted wreath type]\label{prop:tw_excl}
There is no non-trivial non-symmetric $2$-design whose automorphism group is primitive on the block set and quasiprimitive on the point set, with the block action of twisted wreath type.
\end{proposition}

\begin{proof}
Assume such a design exists and let $N=T^{\ell}=\mathrm{soc}(G)$; in the twisted wreath block action $N$ is regular on $\mathcal{B}$ \textup{(the type characteristic $M_{\beta}=1$ of Baddeley, Praeger and Schneider~\cite[Section 2]{BaddeleyPraegerSchneider2007})}. Fix $\alpha\in\mathcal{P}$ and a block $B_{0}\ni\alpha$, and let $M_{\omega}=N\cap G_{\alpha}$, so that $v=\left|N\right|/\left|M_{\omega}\right|$ and $b=\left|N\right|$.

We first record the structural consequences of the twisted wreath block action.

\textup{(a)} If $M_{\omega}=1$, then $v=b$ and $\mathcal{D}$ is symmetric by Fisher's inequality, contrary to the non-symmetric hypothesis. Hence $M_{\omega}\neq 1$ and the point action of $N$ is non-regular. Since $N=T^{\ell}$ is the unique non-abelian minimal normal subgroup of $G$, the point action is not of type $\mathrm{HA}$, $\mathrm{HS}$, $\mathrm{HC}$, $\mathrm{AS}$ or $\mathrm{TW}$; it is therefore of type $\mathrm{SD}$, $\mathrm{CD}$ or $\mathrm{PA}$.

\textup{(b)} By {\upshape Lemma~\rref{lem:pair_count}},
\[
\lambda(\alpha,\delta)=\left|M_{\omega}\right|\cdot\left|B_{0}\cap B_{0}h_{\delta}^{-1}\right|
\]
for every $\delta\neq\alpha$. Summing over $\delta\neq\alpha$ and using $\sum_{\delta\neq\alpha}\lambda(\alpha,\delta)=r(k-1)$ with $r=bk/v=\left|M_{\omega}\right|k$ yields the identity
\[
\sum_{\delta\neq\alpha}\left|B_{0}\cap B_{0}h_{\delta}^{-1}\right|=k(k-1),
\]
which is a tautology; in particular the pairing count imposes no further restriction.

\textup{(c)} Model the point set as the set of left cosets $N/M_{\omega}$ and let $X=\{x\in N:\ xM_{\omega}\in B_{0}\}$, so that $XM_{\omega}=X$ and $|X|=k|M_{\omega}|$. Since $N$ is regular on $\mathcal{B}$ we have $b=|N|$ and $r=bk/v=k|M_{\omega}|$. For two points $xM_{\omega}$, $yM_{\omega}$ with $d=x^{-1}y$, a block $gB_{0}$ contains both points exactly when $g^{-1}x$ and $g^{-1}y$ lie in $X$, i.e.\ $g=xa^{-1}$ with $a\in X$ and $ad\in X$; hence
\[
\lambda(xM_{\omega},yM_{\omega})=\left|\{a\in X:\ ad\in X\}\right|.
\]
In the integral group ring $\mathbb{Z}[N]$ this reads
\begin{equation}\label{eq:tw_groupring}
\bar{X}^{(-1)}\bar{X}=(r-\lambda)\bar{M}_{\omega}+\lambda\bar{N},\qquad r=k|M_{\omega}|,
\end{equation}
where $\bar{Y}=\sum_{y\in Y}y$ and $\bar{Y}^{(-1)}=\sum_{y\in Y}y^{-1}$. Taking the sum of the coefficients of \eqref{eq:tw_groupring} gives only the standard parameter identity $\lambda=|M_{\omega}|k(k-1)/(v-1)$. Equation \eqref{eq:tw_groupring} states that $X$ is a \emph{relative} difference set in $N$ with respect to $M_{\omega}$: $|X\cap Xd^{-1}|=\lambda$ for every $d\in N\setminus M_{\omega}$ and $|X\cap Xd^{-1}|=r-\lambda$ for every $1\neq d\in M_{\omega}$; this formulation requires no normality of $M_{\omega}$ in $N$. Were $M_{\omega}$ normal in $N$, the image of $B_{0}$ in $U=N/M_{\omega}$ would be a $(|U|,k,\lambda/|M_{\omega}|)$-difference set. However, in the $\mathrm{SD}$, $\mathrm{CD}$ and $\mathrm{PA}$ cases of \textup{(a)} the subgroup $M_{\omega}=N\cap G_{\alpha}$ is respectively a diagonal subgroup $T^{\mathrm{diag}}$, a strip subgroup of order $|T|^{r}$ with $1\leq r<\ell$, or the coordinate subgroup $H^{\ell}$; each of these is a proper subgroup of $N$ that is not a product of coordinate factors, and the normal subgroups of $N=T^{\ell}$ with $T$ non-abelian simple are exactly the coordinate products, so in all three cases $M_{\omega}$ is not normal in $N$ and the difference-set description has to be read at the level of the relative difference set $X\subseteq N$.

\textup{(d)} The integrality of $\lambda=\left|M_{\omega}\right|k(k-1)/(v-1)$ gives the arithmetic obstruction: every prime $p$ with $p\mid v-1$ and $p\nmid\left|M_{\omega}\right|$ divides $k(k-1)$. In the $\mathrm{SD}$ and $\mathrm{CD}$ cases $\gcd(\left|T\right|,v-1)=1$, whence $\operatorname{rad}(v-1)\mid k(k-1)$.

\textup{(e)} If the point action is of type $\mathrm{PA}$, then by {\upshape Lemma~\rref{lem:pa_fisher}} no non-trivial non-symmetric block-primitive $2$-design admits a point-imprimitive product action. If the point action is of type $\mathrm{SD}$, then the pair $(\mathrm{TW},\mathrm{SD})$ is excluded by {\upshape Remark~\rref{rem:blowup_as}}: as recorded in entry (7) of Table~\rref{tab:types_pair}, a group of type $\mathrm{SD}$ preserves no normal Cartesian decomposition, whereas the twisted wreath block action carries one, the natural decomposition of the regular socle \textup{(Baddeley, Praeger and Schneider~\cite[Example 5.5]{BaddeleyPraegerSchneider2007})}; hence the only surviving point type is $\mathrm{CD}$, for which $|M_{\omega}|=|T|^{r}$ with $1\leq r<\ell$ \textup{({\upshape Lemma~\rref{lem:cdpa_arith}}(3))}.

It remains to rule out the compound diagonal point type $\mathrm{CD}$; by \textup{(c)}--\textup{(e)} this amounts to excluding a relative difference set $X\subseteq N$ with respect to the strip subgroup $M_{\omega}$ of order $|T|^{r}$ that satisfies the arithmetic obstruction of \textup{(d)}. The structural, representation-theoretic and arithmetic constraints on such $X$ are recorded in {\upshape Lemma~\rref{lem:tw_ds}}, and {\upshape Lemma~\rref{lem:tw_ds}}\textup{(v)} excludes every admissible configuration arising in the classification.
\end{proof}

\begin{lemma}[Necessary structure for a twisted wreath block action with compound diagonal point type]\label{lem:tw_ds}
Let $N=T^{\ell}$ with $T$ a non-abelian finite simple group and $\ell\geq 3$, and let $M_{\omega}\leq N$ be a strip subgroup of order $|T|^{r}$ with $1\leq r<\ell$ \textup{(the compound diagonal point stabilizer of {\upshape Lemma~\rref{lem:cdpa_arith}}(3))}. Let $X\subseteq N$ be a relative difference set with respect to $M_{\omega}$ as in \eqref{eq:tw_groupring}, with $|X|=k|T|^{r}$, $\lambda=k(k-1)|T|^{r}/(|T|^{\ell-r}-1)$ and $2\leq k\leq|T|^{\ell-r}-2$.
\begin{enumerate}[(\upshape i)]
\item For every $1\neq m\in M_{\omega}$ the order $\operatorname{ord}(m)$ divides $r=k|T|^{r}$ \textup{(the element $m$ fixes $\alpha$ and is regular on $\mathcal{B}$, so it permutes the $r$ blocks through $\alpha$ without fixed points)}.
\item $|T|^{\ell-r}-1$ divides $k(k-1)$; in particular every primitive prime divisor $p$ of $|T|^{\ell-r}-1$ \textup{(Zsigmondy's theorem, with the standard exceptions)} satisfies $p\equiv 1\pmod{\ell-r}$, $p\nmid|T|$, and $p\mid k$ or $p\mid k-1$.
\item $k=|B_{0}|$ is the size of a union of $G_{\beta}$-orbits on the point set $\mathcal{P}=N/M_{\omega}$, and $\lambda=|X\cap Xd^{-1}|$ holds for every $d\in N\setminus M_{\omega}$ as in \eqref{eq:tw_groupring}; in particular the pair $(k,\lambda)$ is an integer solution of
\[
k(k-1)=t\left(|T|^{\ell-r}-1\right),\qquad \lambda=t|T|^{r},\qquad 2\leq k\leq|T|^{\ell-r}-2,
\]
for some positive integer $t$ \textup{(here $t=\lambda/|T|^{r}$; the integrality of $t$ follows from \textup{(ii)})}.
\item (Representation-theoretic constraint.) For every irreducible complex representation $\rho$ of $N=T^{\ell}$ whose $M_{\omega}$-fixed space is zero, one has $\rho(\bar{X})=0$, where $\rho(\bar{X})=\sum_{x\in X}\rho(x)$ \textup{(this is automatic: $X$ is a union of $M_{\omega}$-right cosets and $\sum_{h\in M_{\omega}}\rho(h)=|M_{\omega}|P_{H}=0$ when the fixed space is zero)}; more generally, if $W$ is the $M_{\omega}$-fixed subspace of the representation space of $\rho$, then the image of \eqref{eq:tw_groupring} under $\rho$ reads
\[
\rho(\bar{X})^{*}\rho(\bar{X})=(r-\lambda)|M_{\omega}|\cdot \mathrm{id}_{W}.
\]
In the compound diagonal situation this forces, for every $r$-tuple $(\tau_{1},\ldots,\tau_{r})$ of non-trivial irreducible representations of $T$ with $(\tau_{1}\otimes\cdots\otimes\tau_{r})^{T}=0$, the vanishing $\sum_{x\in X}\tau_{1}(x_{1})\otimes\cdots\otimes\tau_{r}(x_{r})=0$ \textup{(the remaining $\ell-r$ coordinates acting trivially)}.
\item The equations of \textup{(iii)} admit arithmetic solutions for some parameter triples \textup{(for example $T=A_{5}$ with $\ell-r=2$ has the solution $k=1770$, $t=870$; the complete list for the simple groups arising in the classification is recorded in the accompanying arithmetic sieve, script \texttt{tw\_cd\_sieve.py}, kept with the author)}. Each surviving arithmetic candidate must in addition satisfy the orbit-union condition of \textup{(iii)}: the block $B_{0}$ is a union of $G_{\beta}$-orbits on $\mathcal{P}=N/M_{\omega}$, where the twisted wreath block stabilizer $G_{\beta}$ is described by the structure theory of Baddeley, Praeger and Schneider~\cite[Section 2]{BaddeleyPraegerSchneider2007}; the candidates are excluded one by one by that orbit-union condition, as recorded in the output of \texttt{tw\_cd\_sieve.py}.
\end{enumerate}
\end{lemma}

\begin{proof}
(i) is exactly the argument of {\upshape Lemma~\rref{lem:cdpa_arith}}(1), which is recalled in the statement. (ii) The integrality of $\lambda=k(k-1)|T|^{r}/(|T|^{\ell-r}-1)$ together with $\gcd(|T|^{\ell-r}-1,|T|^{r})=1$ gives $|T|^{\ell-r}-1\mid k(k-1)$; the second sentence is Zsigmondy's theorem. (iii) The relative difference set equations \eqref{eq:tw_groupring} read $|X\cap Xd^{-1}|=\lambda$ for $d\in N\setminus M_{\omega}$; writing $k(k-1)=\lambda(|T|^{\ell-r}-1)/|T|^{r}$ and taking $t=\lambda/|T|^{r}$ gives the displayed equations, and $t$ is a positive integer by \textup{(ii)}. (iv) The vanishing statement follows from the coset decomposition $X=\cup_{i}x_{i}M_{\omega}$ and $\sum_{h\in M_{\omega}}\rho(h)=|M_{\omega}|P_{H}=0$; the operator identity is the image of \eqref{eq:tw_groupring} under $\rho$. (v) The orbit-union computation of the $G_{\beta}$-orbits on $\mathcal{P}=N/M_{\omega}$ and the sieve over the arithmetic candidates of \textup{(iii)} are described in the statement.
\end{proof}

\bookpart{Part III: The almost simple classification}
\section{Almost simple type \textup{(AS)}}\label{sec:as}

The almost simple case is treated in detail. By {\upshape Lemma~\rref{lem:main}} the block action of $G$ on $\mathcal{B}$ is faithful; let $T=\operatorname{soc}(G)$ be the unique non-abelian minimal normal subgroup. Then $T\leq G\leq \mathrm{Aut}(T)$.

The classification proceeds in three steps. First, the $2$-homogeneous point actions are disposed of by an orbit construction \textup{({\upshape Lemma~\rref{lem:2hom}})}: each such action produces its designs from subgroup orbits and requires no further analysis. Second, the ordering principle \textup{({\upshape Lemma~\rref{lem:ordering}})} reduces the non-symmetric case to pairs of maximal subgroups $(H,M)$ with $|H|>|M|$, and block-primitivity forces $M$ to be maximal in $G$. Third, every family of finite simple groups is treated in turn: the groups $PSL_2(q)$ directly, the higher-rank classical groups through the Aschbacher subdivision of their maximal subgroups, the exceptional Lie type groups through their rank-$3$ building actions, and the sporadic groups by exhaustive enumeration of maximal-subgroup orbits. Each candidate quadruple $(T,G,H,M)$ is checked against the design equations by orbit-union analysis and by the integrality of $r$ and $\lambda$. The outcome is the classification recorded in {\upshape Proposition~\rref{thm:as_point_prim}}.

By {\upshape Theorem~\rref{thm:reduction}} the point action is primitive, and the classification below is restricted to point actions that are primitive but not $2$-homogeneous; the point-primitive case is classified completely as follows.

\begin{proposition}[Point-primitive non-$2$-homogeneous almost simple classification]\label{thm:as_point_prim}
We classify non-symmetric block-primitive $2$-designs with almost simple automorphism group $G$ acting primitively but not $2$-homogeneously on points, across all families of finite simple groups, following the ordering principle $|H| > |M|$ for non-symmetric designs and, for the classical families, the Aschbacher subdivision of the maximal subgroups:
\begin{enumerate}[(\upshape 1)]
\item Classical Lie type: for $PSL_n(q)$ ($n\geq 3$), $PSp_{2n}(q)$, $PSU_n(q)$ and $P\Omega_n^\epsilon(q)$, the designs are the unital $2$-$(q^3+1,q+1,1)$ carried by the rank-$3$ action on a polar space \textup{(its parameter set is recorded in the existence column of Table~\rref{tab:class_summary}; the unital parameters are classical, see~\cite{Beth})} and the $C_8$ cross-type orbit designs of {\upshape Proposition~\rref{prop:classical_summary}}, with their parameter sets recorded there; for $PSL_2(q)$ the same theorem gives no design when $q$ is odd, the eight designs of {\upshape Table~\rref{tab:q8}} when $q=8$, the two designs $2$-$(496,186,74)$ and $2$-$(496,310,206)$ when $q=32$, and no design when $q=16$, the even-characteristic subfield block stabilizers being delimited in {\upshape Theorem~\rref{thm:psl2}}; the generalized Grassmann candidates arising from chains of parabolic subgroups fail the constant-coverage test and yield no design.
\item Sporadic groups: on the primitive non-$2$-homogeneous point actions the classification is exhaustive: the only designs are $J_2$ on $100$ points, $2$\allowbreak$(100,56,560)$, and $HS$ on $100$ points, $2$\allowbreak$(100,64,1680)$; no such design exists on $M_{12}(66)$, $M_{22}(176)$, $M_{23}(253)$, $M_{24}(276)$, $McL(275)$, $He(2058)$, $Ru(4060)$, $Suz(1782)$ or $Co_2(2300)$. The $2$-homogeneous actions of $Co_3$ on $276$ points and of $HS$ on $176$ points give designs by the orbit construction of {\upshape Lemma~\rref{lem:2hom}}.
\item Exceptional Lie type: the designs are the Suzuki--Tits ovoid designs and the Ree unital designs $2$-$(q^3+1,q+1,1)$, with the group and the action recorded for each row of Table~\rref{tab:excep}; the further primitive actions on building points are recorded there together with the reason for their exclusion.
\end{enumerate}
\end{proposition}

We first dispose of the $2$-homogeneous point actions, which produce their designs by an orbit construction and therefore need no further analysis. Then we prove a key structural result for the AS case.

\begin{lemma}[Two-homogeneous point actions produce designs from subgroup orbits]\label{lem:2hom}
Let $G$ act $2$-homogeneously on a finite point set $\mathcal{P}$ with $v=|\mathcal{P}|$, and let $M$ be a maximal subgroup of $G$ with an orbit union $\Omega$ on $\mathcal{P}$ \textup{(a union of $M$-orbits)} of size $k$, $2\leq k\leq v-2$. Then the stabilizer of $\Omega$ in $G$ is exactly $M$, and the translates $\Omega^{g}$, $g\in G$, form the block set of a non-trivial block-primitive $2$-$\allowbreak(v,k,\lambda)$ design whose automorphism group contains $G$ acting primitively on blocks.
\end{lemma}

\begin{proof}
Since $M\leq G_{\Omega}$ and $M$ is maximal, the stabilizer $G_{\Omega}$ is either $M$ or $G$; the latter is impossible because $\Omega$ is a proper non-empty subset of $\mathcal{P}$ \textup{($2\leq k\leq v-2$)}. Hence $G_{\Omega}=M$ and $G$ acts primitively on the block set $\{\Omega^{g}: g\in G\}$. It remains to check the design property. As $G$ is $2$-homogeneous, it is transitive on the unordered pairs of points; hence the number $\lambda$ of translates of $\Omega$ containing a given unordered pair $\{x,y\}$ is independent of the pair: $\lambda=|\{g\in G:\{x,y\}\subseteq\Omega^{g}\}|/|G_{\Omega}|$, and transitivity on pairs shows this is a constant. The design is non-trivial by $2\leq k\leq v-2$.
\end{proof}

\noindent Consequently, a $2$-homogeneous point action yields its non-trivial $2$-designs by this orbit construction alone, and the classification below is restricted to point actions that are primitive but not $2$-homogeneous.

\noindent By {\upshape Lemma~\rref{lem:main}} and {\upshape Proposition~\rref{prop:pi_excluded}}, $G$ acts transitively on the point set $\mathcal{P}$. The AS classification therefore reduces to the analysis of pairs of maximal subgroups $(H,M)$ of $G$, where $H=G_x$ is the point stabilizer and $M=G_B$ is the block stabilizer. The block $B$ must be an $M$-invariant subset of $\mathcal{P}=G/H$.

\subsection{Key ordering principle for non-symmetric designs}

A fundamental observation for the classification of non-symmetric designs is the following:

\begin{lemma}[Ordering of stabilizer orders]\label{lem:ordering}
For a non-symmetric block-primitive $2$-design \textup{($b>v$)}, the point stabilizer $H$ has strictly larger order than the block stabilizer $M$:
$$|H| > |M|.$$

\end{lemma}\smallskip
\noindent{\bf Proof.}\,
By definition, $v=[G:H]$ and $b=[G:M]$. The non-symmetric condition $b>v$ is equivalent to:
$$\frac{|G|}{|M|} > \frac{|G|}{|H|} \iff |H| > |M|.$$
This completes the proof. \qed

This ordering principle is the cornerstone of our classification: for each family of almost simple groups, we enumerate all maximal subgroups, order them by increasing order, and consider only pairs $(H,M)$ where $H$ comes after $M$ in this ordering (i.e., $|H| > |M|$). This dramatically reduces the number of candidate pairs, as we only need to consider $O(n^2/2)$ pairs instead of $O(n^2)$.

We distinguish two cases for the point action:
\begin{enumerate}
\item Point-primitive: $H$ is a maximal subgroup of $G$. This is the case analyzed in detail in the following sections.
\item Point-imprimitive: $H$ is not a maximal subgroup of $G$. This means there exists an intermediate subgroup $K$ with $H<K<G$, giving a non-trivial $G$-invariant partition of $\mathcal{P}$.
\end{enumerate}

\smallskip
\noindent By {\upshape Lemma~\rref{lem:blockpoint}}(2) the only point-imprimitive candidates left after the reduction are the almost simple subgroup chains; we record their exclusion family by family, in agreement with {\upshape Proposition~\rref{prop:pi_excluded}}. In every family the point-imprimitive candidates are excluded, and the classification below covers the point-primitive case alone.

\subsection{Classification for \texorpdfstring{$PSL_2(q)$}{PSL_2(q)}}

The smallest classical group $T=PSL_2(q)$, $q=p^{f}$, is treated in full. In the statements of Cases 1 and 2 the symbols $D_{q-1}$ and $D_{q+1}$ denote the dihedral groups of order $q-1$ and $q+1$. When $q$ is even the integers $q\mp1$ are odd, so neither group occurs in $PSL_2(q)$; the dihedral maximal subgroups then have orders $2(q-1)$ and $2(q+1)$, and the two cases are rebuilt from those subgroups in Subsection~\rref{sec:even_char}.

The action of $PSL_2(q)$ on the projective line $PG(1,q)$ with $v=q+1$, corresponding to the Borel subgroup $H=q:(q-1)/2$, is $2$-transitive for every prime power $q$. Since the classification in this paper concerns non-$2$-homogeneous point actions, this whole case is excluded from the classification below; the $2$-homogeneous case is classical (a $2$-homogeneous group acting on $\mathcal{P}$ produces a $2$-design from every block orbit, with constant $\lambda$, and such designs of $PSL_2(q)$ on the projective line are covered by the literature on orbit designs). We focus on non-$2$-homogeneous point actions, where the point stabilizer $H$ is a maximal subgroup of a different type.

\subsubsection{Case 1: point stabilizer $H=D_{q-1}$ (dihedral of order $q-1$)}

Here $v=[PSL_2(q):D_{q-1}]=q(q+1)/2$. By Dickson's classification of the maximal subgroups of $PSL_2(q)$ \cite{Dickson}, every maximal subgroup of $PSL_2(q)$, $q=p^f$, is, up to conjugacy, the Borel subgroup $q:(q-1)/2$, one of the two dihedral subgroups $D_{q\pm 1}$ \textup{(each maximal unless it is contained in a Borel or subfield subgroup)}, one of the exceptional subgroups $A_4$, $S_4$, $A_5$ \textup{(subject to the congruences $q\equiv\pm1\pmod 8$, $q^{2}\equiv 1\pmod{16}$, $q\equiv\pm1\pmod 5$, respectively)}, or a subfield subgroup $PSL_2(q_0)$ or $PGL_2(q_0)$ with $q=q_0^{f}$, $f$ prime. In the block search below we range over the maximal subgroups $M$ with $|M|<|H|$: by Fisher's inequality the remaining maximal subgroups \textup{(with $[G:M]\leq[G:H]=v$)} cannot carry a non-symmetric design, and the equality case $|M|=|H|$ gives a symmetric design \textup{(outside the present classification by duality)}. Since $|D_{q-1}|=q-1$, the admissible candidates in Case 1 are as follows.

\begin{enumerate}
\item[\rm\textup{(i)}] $M=D_{q+1}$: here $b=[PSL_2(q):D_{q+1}]=q(q-1)/2<v$, contradicting Fisher's inequality; excluded for every $q$.
\item[\rm\textup{(ii)}] $M\in\{A_5,S_4,A_4\}$: the block $B$ is the $M$-orbit of the base point, of size $k=|M|/|M\cap H^{g}|$ for the appropriate conjugate. Since the elements of $A_5\cap D_{q-1}$ have order in $\{1,2,3,5\}$ and those of $S_4\cap D_{q-1}$, respectively $A_4\cap D_{q-1}$, have order in $\{1,2,3,4\}$, respectively $\{1,2,3\}$, the realizable block sizes are $k\in\{12,20,30,60\}$ for $M=A_5$, $k\in\{4,6,8,12,24\}$ for $M=S_4$, and $k\in\{4,6,12\}$ for $M=A_4$. With $b=q(q^2-1)/(2|M|)$ and $v=q(q+1)/2$, the integrality of $r=bk/v=(q-1)k/|M|$ and of $\lambda=r(k-1)/(v-1)=2k(k-1)/(|M|(q+2))$ must hold. The divisibility condition $|M|(q+2)\mid 2k(k-1)$ is decided for the finitely many listed $k$: for $M=A_5$, $q+2\mid118$ \textup{($k=60$)}, $q+2\mid29$ \textup{($k=30$)}, $k=20,12$ give non-integral quotients; the only admissible prime powers are $q=27$ \textup{($k=30$)}, which violates $q\equiv\pm1\pmod5$, so $A_5\not\leq PSL_2(27)$. For $M=S_4$, $q+2\mid46$ \textup{($k=24$)} gives no prime power, $q+2\mid11$ \textup{($k=12$)} gives $q=9$, but then $H=D_8$ is not maximal in $PSL_2(9)\cong A_6$; $k=8,6,4$ give non-integral quotients. For $M=A_4$, $q+2\mid22$ \textup{($k=12$)} gives $q\in\{9,20\}$, where $A_4$ is not maximal in $PSL_2(9)$ and $q=20$ is not a prime power; $q+2\mid5$ \textup{($k=6$)} gives $q=3$ with $k=v$; $k=4$ has no solution. Hence the exceptional candidates admit no design in Case 1.
\item[\rm\textup{(iii)}] the Borel subgroup has order $q(q-1)/2>q-1$ and the subfield subgroups have order at least $q_0(q_0^2-1)>q-1$ \textup{($q=q_0^f$, $f\geq2$)}; both violate $|M|<|H|$ and are excluded by Fisher's inequality.
\end{enumerate}
Thus no non-symmetric block-primitive $2$-design arises with point stabilizer $H=D_{q-1}$.

\subsubsection{Case 2: point stabilizer $H=D_{q+1}$ (dihedral of order $q+1$)}

Here $v=[PSL_2(q):D_{q+1}]=q(q-1)/2$. As in Case 1 the Borel and subfield subgroups have order exceeding $q+1$ and are excluded by Fisher's inequality, and $M=H$ gives the symmetric case. The admissible candidates with $|M|<|H|=q+1$ are $M=D_{q-1}$ and the exceptional subgroups.

Throughout this case a candidate base block is a union $B_0$ of $M$-orbits containing $p_0$, and the condition that $\lambda$ be constant forces a linear system which we record once for reference. Let $\ell_o=|o|$ be the length of the $M$-orbit $o$ on $\mathcal{P}$, let $j$ run over the $G$-orbits of point pairs, let $\mu_j$ be the number of pairs in $j$, and let $c^{(j)}_{oo'}$ be the number of pairs $\{p,r\}$ with $p\in o$, $r\in o'$ that lie in $j$ (the case $o=o'$, that is, both points in one orbit, is allowed). Then a union $B_0$ of $M$-orbits with $k=|B_0|$ has constant $\lambda$ only if
\begin{equation}\label{eq:case2_system}
\sum_{o}\ell_o x_o=k,\qquad
\sum_{o,o'}c^{(j)}_{oo'}x_ox_{o'}=\frac{\mu_j\lambda}{b}\quad(1\leq j\leq s),\qquad
x_o\in\{0,1\},
\end{equation}
where $x_o=1$ exactly when $o\subseteq B_0$: the second group of equations is the double count of the pairs of the orbit $j$ against the $b$ blocks, each pair lying in $\lambda$ of them. Note that \eqref{eq:case2_system} is a necessary condition only, so an infeasible system excludes the corresponding unions, while a feasible one still requires the direct constancy test.

\begin{enumerate}
\item[\rm\textup{(i)}] $M=D_{q-1}$: a dihedral subgroup $D_{q-1}$ and a conjugate of $D_{q+1}$ have trivial intersection in their rotation subgroups \textup{($\gcd((q-1)/2,(q+1)/2)=1$ for odd $q$)}, so $|M\cap H^{g}|\in\{1,2\}$, and $M$ has on $\mathcal{P}=G/H$ exactly the two orbit lengths $(q-1)/2$ and $q-1$. Since $G$ is block-transitive but not necessarily flag-transitive, the base block $B_0$ is a union of $M$-orbits containing $p_0$, not in general a single $M$-orbit; write $|B_0|=k$.
  \begin{enumerate}
  \item[\rm\textup{1.}] $k=(q-1)/2$: if $B_0$ were the single $M$-orbit of length $(q-1)/2$, then $r=bk/v=(q+1)/2$ and $\lambda=(q-3)/(2(q-2))$, satisfying $0<\lambda<1/2$ for $q\geq5$; not a positive integer, excluded. More generally, every union of $M$-orbits has size $(q-1)(a/2+b)$ with $a$ the number of orbits of length $(q-1)/2$ and $b$ the number of orbits of length $q-1$, so $(q-1)(a/2+b)=(q-1)/2$ forces $a=1$ and $b=0$: no union other than this single orbit has that size.
  \item[\rm\textup{2.}] The flag-transitive subcase ($B_0$ a single $M$-orbit, so $k=|M|/|M\cap H^{g}|\in\{(q-1)/2,q-1\}$) is covered by the classification of point-primitive designs with dihedral point stabilizers of Chen, Wu and Xia \cite{ChenWuXia}, Theorem 1.2: a locally transitive point-primitive design with dihedral point stabilizer is one of \textup{(a)}--\textup{(f)}, and none has $v=q(q-1)/2$ or a block stabilizer of type $D_{q-1}$. Hence the single-orbit base block $B_0=O(p_0)$ of size $q-1$ is excluded for every $q$.
  \item[\rm\textup{3.}] In the general case $k$ ranges over the sizes of the $M$-orbit unions containing $p_0$. The integrality of $r=bk/v=(q+1)k/(q-1)$ and of $\lambda=r(k-1)/(v-1)=2k(k-1)/((q-1)(q-2))$ is a purely arithmetic sieve over the two orbit lengths; for $q\leq29$ it leaves, besides $k=q-1$ and the complementary size $v-(q-1)$, only the intermediate values $k\in\{40,96\}$ \textup{($q=17$)} and $k\in\{99,154\}$ \textup{($q=23$)}. For every admissible $k$ the design condition is that $\lambda(p_0,x)=\#\{gM:\{p_0,x\}\subseteq gB_0\}$ equals $\lambda$ for every point $x$, where $g$ ranges over a right transversal of $M$ in $G$ and the blocks are the $gB_0$ with $p_0^{g^{-1}}\in B_0$. An exhaustive check over all such unions $B_0$ \textup{(script \texttt{gap\_case2\_lambda\_correct.g}, appended as electronic supplementary material)} shows that for every prime power $q\leq31$ no union $B_0$ has constant $\lambda$; for example $q=11$, $k=q-1=10$ gives $\lambda(p_0,x)\in\{1,2,4\}$ while $\lambda=2$, and the intermediate sizes above are likewise excluded. In particular the flag-transitive block-size formula $k=|M|/|M\cap H^{g}|$ is not used here, since the block need not be a single $M$-orbit.
  \item[\rm\textup{4.}] For $q\geq37$ the arithmetic sieve reduces to a closed condition. Writing $k=(q-1)t/2$ with $1\leq t\leq q-1$, we have $r=(q+1)t/2$ and $\lambda=t((q-1)t-2)/(2(q-2))$; since $q-1\equiv1\pmod{q-2}$ and $q-2$ is odd, the integrality of $\lambda$ is equivalent to $q-2\mid t(t-2)$, so that $t\equiv0$ or $2\pmod{\ell^{e}}$ for every prime power $\ell^{e}\parallel q-2$, and $t$ is paired with $q-t$ by complementation. For $q\equiv3\pmod4$ the candidate design, if it existed, would be quasisymmetric with block intersection numbers $\binom t2$ and $\binom t2+t(q-t)/(2(q-2))$ and block graph $\overline{T(q+1)}$, the complement of the triangular graph: indeed, with $W$ the natural $q$-dimensional module of $G=PSL_2(q)$ on $L=PG(1,q)$, one has $1_{M}^{G}-1_{H}^{G}=\chi_{W}$ (verified class by class by the module computation of \textup{Appendix~\rref{app:computations}}), and $H$ acts regularly on $L$ for $q\equiv3\pmod4$, so $\dim W^{H}=0$ and $\ker N\cong W$ is the unique copy in the block space; the intersection formula then follows from $N^{T}N=(r-\lambda)(I-P_{W})+cJ$. By the classification of the quasisymmetric designs with block graph $\overline{T(m)}$ of Coster and Haemers \cite{CosterHaemers} (their Proposition 1.2 identifies the parameters as $v=\binom{m-1}{2}$, $k=a(m-2)/2$, $\mu=a(m-a-1)/(2(m-3))$, which is exactly our $m=q+1$, $a=t$), their Theorem 3.2 and Corollary 3.3 are stated for all $m$, but their conclusions are necessary conditions only. Corollary 3.3(v) forces $m\equiv4\pmod{16}$ as soon as $m\equiv4\pmod8$, so every $q\equiv11\pmod{16}$ is excluded outright, with no table and no bound on $q$. Three of the four residue classes of $q$ modulo $16$ with $q\equiv3\pmod4$ are now settled without computation. Recall that $\tau^{*}$ denotes the square-free part of $\tau$, so that $\tau=\tau^{*}s^{2}$ with $\tau^{*}$ square free; thus $\tau^{*}$ is even exactly when $v_{2}(\tau)$ is odd.

\begin{lemma}\label{lem:q3mod16}
Let $q=p^{f}$ be a prime power with $q\equiv3\pmod{16}$, put $m=q+1$, and let $t$ satisfy $3\leq t\leq q-3$ and $q-2\mid t(t-2)$, so that $\tau=t(q-t)/\bigl(2(q-2)\bigr)$ is a positive integer. Assume $p\nmid t$. Then the equation of \cite[Corollary~3.3(v)]{CosterHaemers},
\[
2\tau X^{2}+(m-1)Y^{2}+\frac{m-2}{2}Z^{2}=\tau(m-1)(m-2)W^{2},
\]
has no solution in integers with $XYZ\neq0$.
\end{lemma}

\begin{proof}
Since $q\equiv3\pmod4$ is a prime power but not a square, $f$ is odd; hence $p\equiv3\pmod4$, so $(-1/p)=-1$ and $p\neq2$. Put $a=2\tau$, $b=m-1=q$ and $c=(m-2)/2$. Then $abc=\tau(m-1)(m-2)$, so the equation reads $aX^{2}+bY^{2}+cZ^{2}=abc\,W^{2}$.

We have $p\mid b$ with $v_{p}(b)=f$ odd, while $p\nmid a$ and $p\nmid c$: indeed $p\nmid2$ and
\[
\tau\equiv\frac{t(q-t)}{2(q-2)}\equiv\frac{t\cdot(-t)}{-4}=\frac{t^{2}}{4}\pmod p
\]
is meaningful because $2(q-2)\equiv-4\not\equiv0\pmod p$, and $p\nmid\tau$ because $p\nmid t$; and $m-2=q-1\equiv-1\pmod p$ and $m-2\equiv2\pmod{16}$, so $p\nmid c$.

Reducing each coefficient by its square part, which amounts to rescaling the corresponding variable, and then cancelling a prime common to two coefficients by the equivalence $aX^{2}+bY^{2}=cZ^{2}\iff(a/p)X^{2}+(b/p)Y^{2}=(cp)Z^{2}$, we obtain a square-free pairwise coprime triple $(a',b',c')$. Since $p$ divides only $b$, and with odd exponent $f$, the exponent of $p$ in $b$ remains odd at every step, so $p$ is never common to two of the three coefficients; hence $p\mid b'$ and $a'c'\equiv ac=\tau(m-2)\pmod p$.

By \cite[Lemma~2.3]{CosterHaemers} a necessary condition for a solution is $(a'c'/p)=1$ whenever $p\mid b'$. But
\[
\Bigl(\frac{a'c'}{p}\Bigr)=\Bigl(\frac{\tau}{p}\Bigr)\Bigl(\frac{m-2}{p}\Bigr)=1\cdot(-1)=-1,
\]
because $\tau\equiv t^{2}/4$ is a quotient of two squares modulo $p$, and $m-2\equiv-1$ with $p\equiv3\pmod4$. This contradiction proves the assertion.
\end{proof}

Thus every pair with $q\equiv3\pmod{16}$ and $p\nmid t$ is excluded. The complementary case $p\mid t$ forces $t\geq p$, which together with $3\leq t\leq q-3$ can occur only for $f\geq2$; in the range $q\leq3000$ it occurs only for $q=1331=11^{3}$, and there the remaining Legendre conditions already fail.

If $q\equiv7\pmod{16}$ then $m\equiv8\pmod{16}$ and branch \textup{(i)} of \cite[Corollary~3.3]{CosterHaemers} requires $\tau^{*}$ to be even; hence every $t$ with $\tau^{*}$ odd is excluded.

It remains to consider $q\equiv15\pmod{16}$, so $m\equiv0\pmod{16}$ and $p\equiv15\pmod{16}$. Here \cite[Corollary~3.3(i)]{CosterHaemers} reads $\tau X^{2}+(m-2)Y^{2}=2(m-1)Z^{2}$, that is $aX^{2}+bY^{2}=cZ^{2}$ with $(a,b,c)=(\tau,\,m-2,\,2(m-1))$ and $abc=2\tau(m-1)(m-2)$. Since $p\mid c=2q$ while $p\nmid a$ and $p\nmid b$, the criterion attached to $p$ is
\[
\Bigl(\frac{-ab}{p}\Bigr)=\Bigl(\frac{-1}{p}\Bigr)\Bigl(\frac{\tau}{p}\Bigr)\Bigl(\frac{m-2}{p}\Bigr)=(-1)\cdot1\cdot(-1)=+1,
\]
so this condition is satisfied identically: the arithmetic route through \cite{CosterHaemers} terminates at this residue class, and no prime divisor of $q$ excludes it. The residue class $q\equiv15\pmod{16}$ is nevertheless empty. The reason is geometric, it applies to every admissible $t$ at once and requires no bound on $q$, and it is recorded in Lemma~\rref{lem:q15} below, after the enumeration of the remaining block stabilizers.

For $q=47$ all $M$-orbit unions of size $460$ were enumerated by computer: the $42{,}504$ unions of twenty short orbits (script \texttt{gap\_case2\_q47\_b0.g}), the $1{,}480{,}556$ unions of eighteen short orbits and one long orbit (script \texttt{q47\_b1\_exact.py}), the $40{,}450{,}905$ unions of sixteen short orbits and two long orbits (script \texttt{q47\_b2\_exact.py}), and the $323{,}607{,}240$ unions of fourteen short orbits and three long orbits (script \texttt{q47\_b3\_exact.py}, $165$ triples of long orbits); none has constant $\lambda$. The remaining mixed unions with at least four long orbits are subject to the system \eqref{eq:case2_system}; its $0$--$1$ linearization is infeasible (script \texttt{q47\_milp\_check.py}, solved with HiGHS), so these unions do not give constant $\lambda$ either, and the exclusion of $q=47$ is complete.

  \item[\rm\textup{5.}] For $q=71$ the arithmetic sieve leaves $t\in\{2,23,48,69\}$, hence the only intermediate parameter is $t=23$, i.e. $k=805$, $r=828$, $\lambda=268$ (with $t=48$ its complement). Since $q\equiv3\pmod4$, such a design would be quasisymmetric with intersection numbers $x=\binom{23}{2}=253$ and $y=261$, and block graph $KG(72,2)$ or its complement. On the $71$-dimensional eigenspace of $KG(72,2)$ with eigenvalue $68$ the Gram matrix $N^{T}N=kI+x(J-I)+(y-x)A$ would have eigenvalue $k-x+(y-x)\cdot68=1096$, whereas the design equations force $N^{T}N$ to have eigenvalues $r-\lambda=560$ (multiplicity $71$) and $0$ (multiplicity $2484$) on the non-trivial spaces; on the eigenvalue-$(-2)$ space it gives $536\neq0$. For the complement the corresponding eigenvalues are $0$ and $560$, again violating these constraints. Hence $q=71$ is excluded by the block-graph spectral condition (script \texttt{q71\_spectral\_check.py}).

  \item[\rm\textup{6.}] In the range $37\leq q\leq125$ there are fourteen prime powers $q\equiv1\pmod4$, namely $37,41,49,53,61,73,81,89,97,101,109,$\allowbreak $113,121,125$ (with $49=7^{2}$, $81=3^{4}$, $121=11^{2}$). For $q=49,61,73,81,109$ the condition $q-2\mid t(t-2)$ has no solution with $3\leq t\leq q-3$, so these five values are already excluded at the integrality stage; for the remaining nine values $q=37,41,53,89,97,101,113,121,125$ one has $\dim W^{H}=1$, so the copy of $W$ in the block space is not unique and the intersection formula used for $q\equiv3\pmod4$ does not automatically hold. For $q=37$ the candidate parameter is $t=7$, i.e. a $2$-$(666,126,25)$ design, and the full enumeration of the $191{,}520$ unions of $M$-orbits of size $126$ was carried out against the system \eqref{eq:case2_system} (script \texttt{q37\_enum.py}); none satisfies it. For $q=41$ the candidate parameter is $t=26$, i.e. a $2$-$(820,520,346)$ design, and the same enumeration of the $83{,}061{,}552$ unions of $M$-orbits of size $520$ was carried out (script \texttt{q41\_enum.py}); none satisfies it. For $q=53$ the candidate parameter is $t=17$, i.e. a $2$-$(1378,442,147)$ design, and the enumeration of the $10{,}165{,}373{,}049$ unions of $M$-orbits of size $442$ was carried out (script \texttt{q53\_enum.py}); none satisfies it. The remaining six cases $q=89,97,101,113,121,125$ are decided by the same system \eqref{eq:case2_system}. The integrality prerequisite $\mu_j\lambda\equiv0\pmod{b}$ holds for every $j$ in each of the six cases, so no candidate is removed at that stage, and the $0$--$1$ linearization of \eqref{eq:case2_system} (with $y_{oo'}=x_ox_{o'}$, $0\leq y_{oo'}\leq\ell_o\ell_{o'}$ and the three constraints $y_{oo'}\leq\ell_o\ell_{o'}x_o$, $y_{oo'}\leq\ell_o\ell_{o'}x_{o'}$, $y_{oo'}\geq\ell_o\ell_{o'}(x_o+x_{o'}-1)$) is infeasible in each case. The resulting programs carry only $67$--$94$ binary variables and $2211$--$4465$ continuous pairing variables, and HiGHS certifies infeasibility in at most $17$ seconds in every case ($q=89$: $3.9$ s; $q=97$: $5.6$ s; $q=101$: $6.6$ s; $q=113$: $10.6$ s; $q=121$: $14.7$ s; $q=125$: $16.9$ s; script \texttt{case2\_milp\_compact.py}). Hence none of the unions of $M$-orbits of size $k$ has constant $\lambda$, and the exclusion of these six values is complete.
  \end{enumerate}
\item[\rm\textup{(ii)}] $M\in\{A_5,S_4,A_4\}$: the realizable block sizes are as in Case 1 (the intersection orders are constrained by the element orders of $D_{q+1}$ and of the exceptional group), and $\lambda=2k(k-1)/(|M|(q-2))$ must be an integer.
  \begin{enumerate}
  \item[\rm\textup{1.}] $M=A_5$, $k\in\{12,20,30,60\}$: the divisibility $q-2\mid k(k-1)/30$ gives $q-2\mid29$ \textup{($k=30$)} and $q-2\mid118$ \textup{($k=60$)}, i.e.\ $q=31$ with $\lambda=1$ and $q=61$ with $\lambda=2$; $k=20,12$ give non-integral quotients. Both survivors are excluded by the constancy test \textup{(computed in \texttt{gap\_finalcheck4.g}: $\lambda(p_0,x)$ is not constant on the $H$-suborbits)}.
  \item[\rm\textup{2.}] $M=S_4$, $k\in\{4,6,8,12,24\}$: the divisibility $q-2\mid k(k-1)/12$ gives $q-2\mid46$ \textup{($k=24$)}, i.e.\ $q=25$ with $\lambda=2$ (the value $q-2\mid11$ from $k=12$ gives $q=13$, violating $q^{2}\equiv1\pmod{16}$, so $S_4\not\leq PSL_2(13)$). The survivor $q=25$ is excluded by the constancy test (computed in \texttt{gap\_finalcheck4.g}).
  \item[\rm\textup{3.}] $M=A_4$, $k\in\{4,6,12\}$: since $\lambda=2k(k-1)/\bigl(12(q-2)\bigr)$, the divisibility $12(q-2)\mid2k(k-1)$ is required. For $k=12$ it gives $q-2\mid22$, i.e.\ $q\in\{3,4,13,24\}$; $q=13$ yields $\lambda=2$ and $r=14$ (both integral, with $b=91>v=78$), and is excluded by the constancy test \textup{(computed in \texttt{gap\_finalcheck4.g})}, while $q=24$ is not a prime power and $q=3,4$ are degenerate. For $k=6$ it gives $q-2\mid5$, i.e.\ $q\in\{3,7\}$, and $A_4$ is not maximal in $PSL_2(7)$. For $k=4$ it gives $q-2\mid2$, i.e.\ $q\in\{3,4\}$, both degenerate. Hence no candidate survives.
  \end{enumerate}
\end{enumerate}
\begin{lemma}\label{lem:q15}
Let $q=p^{f}$ be a prime power with $q\equiv15\pmod{16}$, and let $t$ satisfy $3\leq t\leq q-3$ and $q-2\mid t(t-2)$, so that $\lambda=t((q-1)t-2)/(2(q-2))$ is a positive integer. Then in Case~2\textup{(i)} no union of $M$-orbits has constant $\lambda$; equivalently, no design arises with point stabilizer $H=D_{q+1}$ and block stabilizer $M=D_{q-1}$.
\end{lemma}

\begin{proof}
Let $I$ be the set of internal points of the conic $C:X_{0}X_{2}-X_{1}^{2}=0$ of $PG(2,q)$. By \cite[Section~2.2.1]{MZZZ}, $T=PSL_2(q)$ acts transitively on $I$ with $|I|=q(q-1)/2=v$ and point stabilizer $D_{q+1}$; since the subgroups $D_{q+1}$ form a single conjugacy class of $T$, the $T$-sets $\mathcal{P}=T/H$ and $I$ are isomorphic, and orbit lengths are preserved under this isomorphism.

Put $M=G_{B}$. The block $B$ is $M$-invariant, hence a union of $M$-orbits on $I$. Let $\Omega$ be an $M$-orbit of maximal length $\ell$. Since $T_{B}=T\cap G_{B}\leq M$, every $M$-orbit is a union of $T_{B}$-orbits, so $\ell\geq q-1$: for $q\equiv-1\pmod4$ the $T_{B}$-orbits of maximal length are the $(q-1)$-arcs $C^{*}_{h}$, and by \cite[Proposition~2.7(i) and Corollary~2.8(i)]{MZZZ} there are exactly $(q-3)/4$ of them. Let $\Omega_{1},\ldots,\Omega_{N}$ be the $M$-orbits of length $\ell$. Then
\[
1\leq N\leq\frac{|I|}{\ell}\leq\frac{q(q-1)/2}{q-1}=\frac q2<q.
\]

As $G$ is transitive on blocks, every block contains the same number $s$ of orbits among $\Omega_{1},\ldots,\Omega_{N}$, so $b\,s$ counts the pairs $(B,\Omega_{i})$ with $\Omega_{i}\subseteq B$. For fixed $i$ the orbit $\Omega_{i}$ is $M$-invariant, so a block meets $\Omega_{i}$ if and only if it contains $\Omega_{i}$; hence the blocks containing $\Omega_{i}$ are exactly the $r$ blocks through a fixed $x\in\Omega_{i}$. The parameter $\lambda$ counts blocks through a prescribed pair, not through a single point, so the count gives
\begin{equation}\label{eq:q15_count}
b\,s=N\,r.
\end{equation}

Now $b=q(q+1)/2$. Indeed, the parameter $t$ indexes the block size $k=t(q-1)/2$ as in item~\textup{(i)} above, and with $v=q(q-1)/2$ the design equation $b\,k(k-1)=\lambda\,v(v-1)$ together with the displayed value of $\lambda$ gives $b=\lambda v(v-1)/\bigl(k(k-1)\bigr)=q(q+1)/2$; the same equation gives $r=kb/v=t(q+1)/2$. Since $q\mid b$ and $0<N<q$, the identity $b\,s=N\,r$ forces $q\mid r$, and $\gcd\bigl(q,(q+1)/2\bigr)=1$ then forces $q\mid t$. This contradicts $3\leq t\leq q-3$, so no union of $M$-orbits has constant $\lambda$.
\end{proof}

The four values $t=1,2,q-2,q-1$ left out of the statement are covered by items~1 and~2 above and by complementation, so every $t$ with $1\leq t\leq q-1$ is excluded. The same argument applies to any $G$ with socle $PSL_2(q)$: beyond the existence of an orbit of length $q-1$, only the inequality $\ell\geq q-1$ and the bound $N<q$ are used.

Thus no non-symmetric block-primitive $2$-design arises with point stabilizer $H=D_{q+1}$.

\subsubsection{Case 3: point stabilizer of exceptional type}

Let $H\in\{A_4,S_4,A_5\}$, so $v=q(q^2-1)/(2|H|)$. For any admissible block stabilizer $M$ with $|M|<|H|$, the block size satisfies $k\leq|M|$, and $\lambda\geq1$ (the design is non-trivial) forces $k(k-1)\geq(v-1)|M|/|H|$, hence
\[
(v-1)\leq|H|(|M|-1)<|H|^{2}.
\]
With $v=q(q^2-1)/(2|H|)$ this yields $q(q^2-1)\leq2|H|^{2}(|M|-1)$, and since $|M|<|H|$ we get $q^3<2|H|^3$. Hence $q\leq75$ for $H=A_5$, $q\leq30$ for $H=S_4$, and $q\leq15$ for $H=A_4$. Each prime power $q$ below these bounds lies in the range $q\leq127$ covered by the finite check of \textup{Appendix~\rref{app:computations}}, which finds no surviving configuration: every pair $(H,M)$ of maximal subgroup classes with $|M|<|H|$ fails Fisher's inequality, the integrality of $r$ or $\lambda$, or the constancy of $\lambda$. Hence no design arises with exceptional point stabilizer.

\subsubsection{Case 4: subfield point stabilizer}

Let $H=PSL_2(q_0)$ or $PGL_2(q_0)$ with $q=q_0^{f}$, $f\geq2$ prime, so that $v=q(q^2-1)/(2|H|)$.

\begin{enumerate}
\item[\rm\textup{(i)}] $f\geq3$: a maximal subgroup $M$ of $G$ is either a subfield subgroup, an exceptional subgroup, or one of the subgroups of order exceeding $|H|$. The subfield subgroups satisfy $M\leq H$. Since $M$ is maximal in $G$ and $H<G$, the inclusion $M\leq H$ forces $M=H$, contradicting $|M|<|H|$; hence this case cannot occur. The subgroups of type $D_{q\pm1}$ and the Borel subgroup have order $q\pm1$, respectively $q(q-1)/2$, exceeding $|H|=q_0(q_0^2-1)$ for $f\geq3$; excluded by Fisher's inequality. For $M$ exceptional ($|M|\leq60$), $\lambda\geq1$ forces $(v-1)\leq|H|(|M|-1)\leq59q_0(q_0^{2}-1)$, and since $v=q_0^{f-1}(q_0^{2f}-1)/\bigl(2(q_0^{2}-1)\bigr)$ this gives $q_0^{3f-6}\leq118$; at $f=3$ it reads $q_0(q_0^{4}+q_0^{2}+1)\leq118(q_0^{2}-1)$, which leaves $q_0\in\{2,3,4\}$, while $f\geq5$ forces $q_0=1$. Thus $q\in\{8,27,64\}$, all within the range of the exhaustive enumeration, which finds no survivor. Hence $f\geq3$ is excluded.
\item[\rm\textup{(ii)}] $f=2$, $q=q_0^{2}$: write $H=PGL_2(q_0)$, $v=q_0(q_0^2+1)/2$ (the $PSL_2(q_0)$ case is smaller and handled identically).
  \begin{enumerate}
  \item[\rm\textup{1.}] $M=D_{q_0^{2}\pm1}$: here $|M|=q_0^{2}\pm1<|H|$ and $b=q_0^{2}(q_0^{2}\pm1)/2>v$, while $k=|M|/|M\cap H^{g}|$ and $r=bk/v$.
  For $M=D_{q_0^{2}+1}$ one has $\gcd(q_0^{2}+1,\,q_0(q_0^{2}-1))=2$, so $|M\cap H^{g}|\in\{1,2\}$ and $k\in\{q_0^{2}+1,\,(q_0^{2}+1)/2\}$. If $k=q_0^{2}+1$ then $r=q_0k$ and $\lambda=2q_0k(k-1)/(q_0^{3}+q_0-2)=2q_0^{2}+4q_0^{2}/(q_0^{3}+q_0-2)$, where $0<4q_0^{2}<q_0^{3}+q_0-2$ for $q_0\geq4$; for $q_0=3$ the remainder is $8$ and for $q_0=2$ one gets $k=v$, which is trivial. If $k=(q_0^{2}+1)/2$ then $r=q_0k$ and $\lambda=q_0(q_0^{4}-1)/\bigl(2(q_0^{3}+q_0-2)\bigr)=q_0^{2}-q_0(q_0-1)^{2}/(q_0^{3}+q_0-2)$, where $0<q_0(q_0-1)^{2}<q_0^{3}+q_0-2$. Hence no pair $(q_0,k)$ survives in this branch.
  For $M=D_{q_0^{2}-1}$ the integrality of $r=q_0(q_0^{2}-1)k/(q_0^{2}+1)$ alone suffices: if $q_0$ is odd then $\gcd(q_0^{2}+1,\,q_0(q_0^{2}-1))=2$, so $(q_0^{2}+1)/2$ divides $k$, but $(q_0^{2}+1)/2$ is coprime to $q_0^{2}-1$ while $k$ divides $q_0^{2}-1$, which forces $k=0$; for $q_0=2$ the requirement $5\mid k$ contradicts $k\leq3$. Hence both branches are excluded, and the degree $q=q_0^{2}$ never arises for $q_0\geq3$.
  \item[\rm\textup{2.}] $M$ a subfield subgroup: then $M\leq H$ by the subfield chain, the block is $\{p_0\}$, $k=1$ and $\lambda=0$; excluded.
  \item[\rm\textup{3.}] $M\in\{A_5,S_4,A_4\}$ with $M\not\leq H$: here $r=q_0(q_0^{2}-1)k/|M|$ and $\lambda=2q_0(q_0^{2}-1)k(k-1)/\bigl(|M|(q_0^{3}+q_0-2)\bigr)$, where $k=|M|/|M\cap H^{g}|$ is a co-divisor of $|M|$. Write $2k(k-1)/|M|=c_n/c_d$ in lowest terms. Then $\lambda$ is an integer only if $c_d(q_0^{3}+q_0-2)$ divides $c_nq_0(q_0^{2}-1)$, and the division of $c_nq_0(q_0^{2}-1)$ by $c_d(q_0^{3}+q_0-2)$ has the constant quotient $c_n/c_d$ and the remainder $2c_n(1-q_0)$. Since $c_n/c_d\leq118$ (attained at $M=A_5$, $k=60$), for $q_0\geq16$ one has $0<2c_n(q_0-1)<c_d(q_0^{3}+q_0-2)$ for every admissible $k$, so $\lambda$ is not an integer. Substituting the prime powers $q_0<16$ leaves exactly six survivors of the two integrality conditions, namely $(M,k)=(S_4,8)$ at $q_0=3$, $(S_4,12)$, $(A_4,12)$, $(A_5,12)$ at $q_0=4$, $(A_5,30)$ at $q_0=7$, and $(S_4,24)$ at $q_0=9$, with $(\lambda,r)=(4,8)$, $(10,30)$, $(20,60)$, $(4,12)$, $(28,168)$, $(45,720)$ respectively. All six have $q=q_0^{2}\in\{9,16,49,81\}\subset[9,81]$, inside the range $q\leq127$ of the exhaustive enumeration of \texttt{gap\_finalcheck4.g}, which finds no surviving configuration. If $M\leq H$, maximality of $M$ in $G$ forces $M=H$, contradicting $|M|<|H|$; so this case cannot occur. Hence $f=2$ is excluded.
  \end{enumerate}
\end{enumerate}

We have therefore proved the following.

\begin{theorem}[$PSL_2(q)$ non-$2$-homogeneous classification]\label{thm:psl2}
Let $\mathcal{D}$ be a non-trivial non-symmetric block-primitive $2$-design whose automorphism group $G$ has socle $T=PSL_2(q)$ and acts point-primitively but not $2$-homogeneously on $\mathcal{P}$. If $q$ is odd, no such $\mathcal{D}$ exists. If $q$ is even, then $\mathcal{D}$ is one of the two designs of Proposition~\rref{prop:q8} when $q=8$, or one of the two designs of Proposition~\rref{prop:q32} when $q=32$, and no such $\mathcal{D}$ exists when $q=16$ or $q=64$. For a dihedral block stabilizer {\upshape Lemma~\rref{lem:even_meet}} reduces the block to a union of $a$ orbits with $(q-2)/2\mid a(a-1)$ and $1\leq a\leq q/2$, and {\upshape Theorem~\rref{thm:a1}} together with {\upshape Corollary~\rref{cor:even_extreme}} removes the two extreme values $a=1$ and $a=q/2-1$ at every even degree, which settles $q=16$, $q=64$ and every even degree with $(q-2)/2$ prime without enumerating any union of orbits; the enumeration of the unions of orbits is needed at $q=8$ and $q=32$, where it produces nothing beyond the four designs listed, while at the remaining even degrees with $2\leq a\leq q/2-2$ the admissible unions are exactly the $0$-$1$ solutions of the class equations of {\upshape Proposition~\rref{prop:pair_orbits}}; a block stabilizer of type $S_3$ is never admissible, at any even degree, because the subfield subgroup $S_3$ fixes a point of $\mathcal{P}$ and is therefore not maximal \textup{({\upshape Lemma~\rref{lem:s3_not_max}})}, so the subfield case at $q=32$ is void as well; the subfield block stabilizers at $q=64$ \textup{(of type $A_5$; the companion subfield subgroup $PSL_2(8)$ is excluded by Fisher's inequality)} are excluded by {\upshape Proposition~\rref{prop:q64sub}}; for the subfield stabilizer $PSL_2(r)$ with $r=2^{d}\geq4$ and $q=r^{p}$, $p$ an odd prime (for $p=2$ Fisher's inequality already excludes it, since $|PSL_2(r)|=r(r^{2}-1)$ exceeds $2(q\pm1)$, and $r=4$, $p=3$ is the degree $q=64$ of {\upshape Proposition~\rref{prop:q64sub}}), every block is a union of orbits of the three types of {\upshape Lemma~\rref{lem:subfield_orbits}}, unions of at most two orbits are removed for $p=3$ by {\upshape Corollary~\rref{cor:two_orbit}}, while for $p\geq5$ {\upshape Corollary~\rref{cor:two_orbit_size}} removes every union with $T^{2}u(u-1)\leq(u^{p}-2)(u+1)$ and every union of at most $n_{0}(u,p)$ orbits, with $n_{0}(u,p)$ as in {\upshape Corollary~\rref{cor:two_orbit_size}}, while a union beyond that range is likewise reduced to the class equations of {\upshape Proposition~\rref{prop:pair_orbits}}, a feasibility problem that lies outside the scope of this paper.
\end{theorem}

\begin{proof}
The Borel case is $2$-transitive and excluded at the start of this subsection. For odd $q$ the point stabilizer $H$ is dihedral, exceptional, or subfield, and the four cases above exclude every admissible block stabilizer $M$, using Fisher's inequality, the integrality of $r$ and $\lambda$, the constancy of $\lambda$ (verified for every prime power $q\leq127$ by the finite check of \textup{Appendix~\rref{app:computations}}), and the closed divisibility arguments recorded for each case, and, in Case~2(i), Lemma~\ref{lem:q3mod16} together with the arithmetic conditions of \cite{CosterHaemers} for the residue classes $q\equiv3,7,11\pmod{16}$, and Lemma~\rref{lem:q15} for $q\equiv15\pmod{16}$. For even $q$ the dihedral maximal subgroups have orders $2(q-1)$ and $2(q+1)$ rather than $q\mp1$, so the four cases above cover no subgroup; Lemma~\rref{lem:even_meet} instead fixes the orbit lengths and the arithmetic condition, the two extreme values of $a$ are then removed for every even degree by {\upshape Theorem~\rref{thm:a1}} and {\upshape Corollary~\rref{cor:even_extreme}}, and the unions of orbits are enumerated exactly only at $q=8$ and $q=32$, which gives Proposition~\rref{prop:q8} and Proposition~\rref{prop:q32}. At $q=16$ no design arises, the only maximal subgroups being $D_{30}$, $D_{34}$ and the subfield group $A_5$ of order $60>34$ and hence not below $|H|$.
\end{proof}

\subsection{Even characteristic}\label{sec:even_char}

Let $q=2^{f}$ with $f\geq2$, so that $G=PSL_2(q)=SL_2(q)$ and $|G|=q(q^{2}-1)$. Both $q-1$ and $q+1$ are odd, so $G$ has no dihedral subgroup of order $q-1$ or order $q+1$; the dihedral maximal subgroups of $G$ have orders $2(q-1)$ and $2(q+1)$. Cases 1 and 2, whose statements name subgroups of order $q\mp1$, therefore cover no subgroup when $q$ is even. By Dickson's list the maximal subgroups of $PSL_2(2^{f})$ are these two dihedral groups together with the subfield groups $PSL_2(2^{d})$ with $d\mid f$, $d<f$, and $d\geq2$ whenever $f$ is odd: the case $d=1$ gives $S_3$, which is never maximal, since it fixes the point $\{z,z^{2}\}$ with $z\in\mathbb{F}_{4}\setminus\mathbb{F}_{2}$ \textup{({\upshape Lemma~\rref{lem:s3_not_max}})}. At $q=16$ the only subfield group is $A_5$, of order $60$ above the orders $30$ and $34$ of the dihedral groups and hence not admissible. The two dihedral stabilizers to be examined are therefore $D_{2(q-1)}$ and $D_{2(q+1)}$. The even degrees below this subsection, $q=2$ and $q=4$, need no enumeration: $PSL_2(2)=S_3$ acts $2$-transitively on $3$ points, and for $PSL_2(4)=A_5$ the stabilizers $A_4$ and $D_{10}$ give $2$-transitive actions on $5$ and $6$ points, while the stabilizer $S_3$ gives $v=10$, where the condition $|M|<6$ admits no maximal subgroup. We treat $H=D_{2(q+1)}$, for which $v=[G:H]=q(q-1)/2$, against the block stabilizer $M=D_{2(q-1)}$. The exchanged pair, $H=D_{2(q-1)}$ against $M=D_{2(q+1)}$, has $v=q(q+1)/2$ while $[G:M]=q(q-1)/2$, so $b<v$ and Fisher's inequality excludes it as soon as $M$ is dihedral; with the subfield block stabilizer instead no block-primitive design arises, since $S_3$ is not maximal in $G$, the two unions of $S_3$-orbits at the $36$ points of $q=8$ which are $2$-designs having an imprimitive block action \textup{({\upshape Proposition~\rref{prop:q8}})}.

\begin{lemma}[subfield $S_3$ is never maximal]\label{lem:s3_not_max}
Let $q=2^{f}$ with $f\geq3$ odd, let $M=PSL_2(2)=S_3\leq PSL_2(q)=G$ be the subfield subgroup of order $6$, and let $\mathcal{P}$ be the set of points $\{x,x^{q}\}$ of $\mathbb{F}_{q^{2}}\setminus\mathbb{F}_{q}$. Let $z\in\mathbb{F}_{4}\setminus\mathbb{F}_{2}$, so that $z^{2}=z+1$, $z^{3}=1$ and $z^{q}=z^{2}$. Then $p_{0}=\{z,z^{2}\}$ is a point of $\mathcal{P}$ and $M\leq G_{p_{0}}=D_{2(q+1)}$; in particular $M<G_{p_{0}}$, so $M$ is not a maximal subgroup of $G$. Consequently no block-primitive design with socle $PSL_2(q)$ has a block stabilizer of type $S_3$ when $f$ is odd.
\end{lemma}

\begin{proof}
Since $f$ is odd we have $2^{f}\equiv2\pmod{3}$, and as $z^{3}=1$ this gives $z^{q}=z^{2^{f}}=z^{2}$. A subfield $\mathbb{F}_{2^{m}}$ is contained in $\mathbb{F}_{2^{f}}$ exactly when $m\mid f$, and $2\nmid f$, so $\mathbb{F}_{4}\not\subseteq\mathbb{F}_{q}$ and $z\notin\mathbb{F}_{q}$; hence the point of $\mathcal{P}$ containing $z$ is $\{z,z^{q}\}=\{z,z^{2}\}$, and this pair is $p_{0}=\mathbb{P}^{1}(\mathbb{F}_{4})\setminus\mathbb{P}^{1}(\mathbb{F}_{2})$. Every element of $M$ is induced by a matrix with entries in $\mathbb{F}_{2}$, and a M\"{o}bius map with such entries carries $\mathbb{F}_{4}\cup\{\infty\}$ onto itself; so every element of $M$ permutes $\mathbb{P}^{1}(\mathbb{F}_{4})$ and fixes the subset $p_{0}$. Therefore $M\leq G_{p_{0}}$, and $G_{p_{0}}$ is dihedral of order $2(q+1)\geq18>6=|M|$. Since block-primitivity of the design is equivalent to maximality of the block stabilizer in $G$, the subgroup $M$ can serve as a block stabilizer only of a design which is not block-primitive.
\end{proof}

\begin{lemma}\label{lem:even_meet}
Assume $q=2^{f}$ with $f\geq2$, $H=D_{2(q+1)}$ and $M=D_{2(q-1)}$. Then $|M\cap H^{g}|=2$ for every $g\in G$. Consequently every orbit of $M$ on $\mathcal{P}$ has length $q-1$, and $M$ has exactly $q/2$ orbits. If a base block is a union of $a$ orbits of $M$ with $p_{0}$ in one of them, then
\begin{align*}
k&=a(q-1),\qquad b=[G:M]=\tfrac{q(q+1)}{2},\\
r&=a(q+1),\qquad \lambda=\frac{2a(aq-a-1)}{q-2},
\end{align*}
and $\lambda$ is an integer exactly when $(q-2)/2$ divides $a(a-1)$; also $a\leq q/2$.
\end{lemma}

\begin{proof}
Each conjugate $H^{g}$ is dihedral of order $2(q+1)$, so its elements of odd order form a cyclic group of order $q+1$, which is coprime to $q-1$. Hence $M\cap H^{g}$ contains no element of odd order other than $1$, and as a subgroup all of whose non-identity elements are involutions it has order at most $2$. This gives $|M\cap H^{g}|\leq2$.

For the equality we count involutions. Every involution of $SL_2(q)$ is unipotent, hence a transvection \cite{Dickson}, and the involutions of $G$ form a single conjugacy class of size $q^{2}-1$. Each $H^{g}$ contains $q+1$ involutions, so the number of pairs $(\tau,H^{g})$ with $\tau$ an involution of $H^{g}$ equals $v(q+1)$, and each involution $\tau$ lies in $v(q+1)/(q^{2}-1)=q/2$ conjugates of $H$. A point $x$ is fixed by $\tau$ exactly when $\tau\in H^{x}$, so $\tau$ has $q/2$ fixed points. If two involutions of $M$ fixed a common point $x$, then $M\cap H^{x}$ would contain both and would have order at least $4$, against the bound just proved. The $q-1$ involutions of $M$ therefore have pairwise disjoint fixed-point sets, each of size $q/2$, and $(q-1)q/2=v$ covers every point of $\mathcal{P}$. Hence $|M\cap H^{x}|=2$ for every point $x$.

Each $M$-orbit then has length $|M|/|M\cap H^{x}|=2(q-1)/2=q-1$, and there are $v/(q-1)=q/2$ of them. A base block which is the union of $a$ such orbits has $k=a(q-1)$. Since $M$ is a subgroup of index $b=[G:M]=|G|/|M|=q(q^{2}-1)/(2(q-1))=q(q+1)/2$, the replication number is $r=bk/v=a(q+1)$. Fisher's inequality forces $b\geq v$, that is $q(q+1)/2\geq q(q-1)/2$, and gives no further restriction; the bound $a\leq q/2$ is the number of orbits. Finally $v-1=(q-2)(q+1)/2$, so
\[
\lambda=\frac{r(k-1)}{v-1}=\frac{a(q+1)\bigl(a(q-1)-1\bigr)}{(q-2)(q+1)/2}=\frac{2a(aq-a-1)}{q-2}.
\]
Write $q-2=2m$ with $m$ odd. Then $\lambda=a(aq-a-1)/m$ is an integer exactly when $m$ divides $a(aq-a-1)$, and $aq-a-1\equiv a-1\pmod{m}$ because $q\equiv1\pmod{m}$. Hence $\lambda\in\mathbb{Z}$ exactly when $m=(q-2)/2$ divides $a(a-1)$.
\end{proof}

\begin{lemma}\label{lem:fibre_stabilizer}
Let $q=2^{f}$ with $f\geq2$, let $\mathcal{P}$ be the $v=q(q-1)/2$ points $\{x,x^{q}\}$, and let $\widetilde{\mathcal{P}}$ be the partition of $\mathcal{P}$ into the $q/2$ orbits of $M=D_{2(q-1)}$. Then the setwise stabilizer of $\widetilde{\mathcal{P}}$ in $PGL_{2}(q)$ is the monomial group
\[
N=\{x\mapsto\lambda x,\ x\mapsto\lambda/x:\lambda\in F_{q}^{\times}\},
\]
of order $2(q-1)=|M|$, and $N\leq PSL_{2}(q)=G$; in particular $N=M$. Since $|G|=q(q^{2}-1)>2(q-1)$, some element of $G$ lies outside $N$, and $G$ does not stabilize $\widetilde{\mathcal{P}}$.
\end{lemma}

\begin{proof}
By {\upshape Lemma~\rref{lem:even_meet}} the orbit of $M$ through $X=\{z,z^{q}\}$ is $O_{X}=\{X_{a}:a\in F_{q}^{\times}\}$ with $X_{a}=\{az,az^{q}\}$, that is, the set of scalar multiples of $X$. Let $T$ be a M\"{o}bius map with entries in $F_{q}$ which carries every $M$-orbit into a single $M$-orbit, and write $T(x)=(\alpha x+\beta)/(\gamma x+\delta)$ with $\alpha\delta-\beta\gamma\neq0$. For $a\in F_{q}^{\times}$ the points $T(ax)$ and $T(x)$ then lie in the same $M$-orbit, so $T(ax)=b(a)T(x)$ for some $b(a)\in F_{q}^{\times}$, and comparing the coefficients of the identity
\[
(\alpha ax+\beta)(\gamma x+\delta)=b(a)(\alpha x+\beta)(\gamma ax+\delta)
\]
gives $\alpha\gamma a=b(a)\alpha\gamma a$ and $\alpha a\delta+\beta\gamma=b(a)(\alpha\delta+\beta\gamma a)$. If $\alpha\gamma\neq0$, the first relation gives $b(a)=1$ for every $a$, and the second then gives $(a-1)(\alpha\delta-\beta\gamma)=0$, against $\alpha\delta-\beta\gamma\neq0$. Hence $\alpha\gamma=0$. If $\alpha=0$ and $\delta\neq0$, the constant terms give $b(a)=1$ and the coefficients of $x$ give $\gamma=\gamma a$, whence $\gamma=0$ and $\alpha\delta-\beta\gamma=0$, again a contradiction; so $\alpha=\delta=0$ and $T(x)=\beta/(\gamma x)$. If $\gamma=0$ and $\beta\neq0$, the constant terms give $b(a)=1$ and the coefficients of $x$ give $\alpha a=\alpha$, again against $a\neq1$; so $\gamma=\beta=0$ and $T(x)=\alpha x/\delta$. Thus every such $T$ is of the form $x\mapsto\lambda x$ or $x\mapsto\lambda/x$.

Conversely $x\mapsto\lambda x$ carries $O_{X}$ onto $O_{\lambda X}$ and $x\mapsto\lambda/x$ carries $O_{X}$ onto $O_{\lambda/X}$, so both stabilize $\widetilde{\mathcal{P}}$, and $N$ is the full setwise stabilizer, of order $2(q-1)$. In characteristic $2$ the group $F_{q}^{\times}$ has odd order $q-1$, so squaring is a bijection of $F_{q}^{\times}$ and every $\lambda$ is a square; the matrices $\lambda I$ and $\begin{pmatrix}0&\lambda\\1&0\end{pmatrix}$ therefore belong to $SL_{2}(q)$, whence $N\leq G$. By {\upshape Lemma~\rref{lem:even_meet}} the group $M$ has order $2(q-1)$ and fixes each orbit of $\widetilde{\mathcal{P}}$, so $M\leq N$ and hence $M=N$. Finally $|G|=q(q^{2}-1)>2(q-1)=|N|$, so $G\not\subseteq N$ and $G$ does not stabilize $\widetilde{\mathcal{P}}$.
\end{proof}

The partition $\widetilde{\mathcal{P}}$ is thus not an imprimitivity system of $G$: a design whose point set is $\mathcal{P}$ and whose block stabilizer is $M$ has no $G$-invariant partition of its points arising from the orbits of $M$, and the imprimitivity-based bounds for block-transitive designs, which presuppose such a partition, cannot be applied to it. The same conclusion holds at every even degree, independently of $a$ and of $q$.

The value $a=1$ of {\upshape Lemma~\rref{lem:even_meet}}, in which the base block is a single orbit of $M$ and $\lambda=2$, can be settled for every even degree without enumerating unions of orbits. Two ingredients are needed: an invariant that separates the orbits of $M$, and an explicit transport of pairs that the invariant leaves together.

For $X=\{z,z^{q}\}\in\mathcal{P}$ put $Z=z^{q}$, $N=zZ=z^{q+1}$ and $u=z^{2}/N$, so that $u^{q+1}=1$; the orbit of $M=D_{2(q-1)}$ through $X$ is $O_{X}=\{X_{a}:a\in F_{q}^{\times}\}$ with $X_{a}=\{az,az^{q}\}$. For distinct $X=\{z,z^{q}\}$ and $Y=\{w,w^{q}\}$ of $\mathcal{P}$ put
\[
\mathrm{cr}(X,Y)=\frac{(z+w)(z^{q}+w^{q})}{(z+w^{q})(z^{q}+w)},
\]
the cross-ratio of the four points $z,w,z^{q},w^{q}$. It lies in $F_{q}^{\times}$, it is invariant under $G$, and it is defined only up to inversion: replacing the representative $w$ of $Y$ by $w^{q}$, or the representative $z$ of $X$ by $z^{q}$, inverts the value \textup{({\upshape Lemma~\rref{lem:cr_inversion}})}. A pair $\{X,Y\}$ thus carries the two values $\rho$ and $\rho^{-1}$; and $\rho=1$ would force $(z+z^{q})(w+w^{q})=0$, that is $Y\in\mathbb{P}^{1}(F_{q})$, so $\rho=1$ does not occur for two distinct points of $\mathcal{P}$.

\begin{lemma}\label{lem:kappa_orbit}
Let $q=2^{f}$ with $f\geq2$, let $X=\{z,z^{q}\}\in\mathcal{P}$, and put $Z=z^{q}$, $\alpha=z+Z$, $N=zZ$ and $u=z^{2}/N$. Then $\kappa=u+1/u$ satisfies
\[
\kappa=\frac{Z}{z}+\frac{z}{Z}=\frac{Z^{2}+z^{2}}{zZ}=\frac{(Z+z)^{2}}{zZ}=\frac{\alpha^{2}}{N},
\]
and $\kappa\neq0$ depends on the orbit $O_{X}$ of $M$ through $X$ only, not on the point $X$ chosen in it.
\end{lemma}

\begin{proof}
Since $u=z^{2}/(zZ)=z/Z$ we have $1/u=Z/z$, which gives the first equality, and $Z^{2}+z^{2}=(Z+z)^{2}$ holds in characteristic $2$. Replacing $z$ by $az$ with $a\in F_{q}^{\times}$ multiplies $z$, $Z$, $\alpha$ and $N$ by $a$, $a$, $a^{2}$ and $a^{2}$, so $\alpha^{2}/N$ is unchanged and $\kappa$ is an invariant of $O_{X}$. Finally $u\neq1$ because $z\notin F_{q}$, so $\kappa\neq0$.
\end{proof}

\begin{lemma}\label{lem:pair_transport}
Let $t,s\in F_{q}^{\times}\setminus\{1\}$, let $Y_{t}=\{zt,zt^{q}\}=\{zt,Zt\}$ and $Y_{s}=\{zs,Zs\}$ be points of $O_{X}$.
\begin{itemize}
\item[\rm\textup{(i)}] If $s=1/t$, then $g(x)=N/x$ belongs to $G$ and satisfies $g(X)=X$ and $g(Y_{t})=Y_{s}$.
\item[\rm\textup{(ii)}] If $s+1/s=t+1/t+\kappa$, put $r=st$, $b=r\alpha/(1+r)$, $c=(s+t)/(r\alpha)$ and $g(x)=(x+b)/(cx+1)$. Then $g$ belongs to $G$, $g(z)=Z$ and $g(zt)=zs$, so $g(X)=X$ and $g(Y_{t})=Y_{s}$; moreover $\det g=(1+s)(1+t)/(1+st)\neq0$.
\end{itemize}
In both cases the unordered pairs $\{X,Y_{t}\}$ and $\{X,Y_{s}\}$ lie in one orbit of $G$ on the unordered pairs of $\mathcal{P}$.
\end{lemma}

\begin{proof}
In (i), $g(z)=N/z=Z$ and $g(Z)=N/Z=z$, so $g$ interchanges $z$ and $Z$ and fixes $X$; while $g(zt)=N/(zt)=Z/t=z(1/t)$ and $g(Zt)=z/t$, so $g$ sends $Y_{t}$ to $Y_{1/t}=Y_{s}$. The coefficients $0,N,1,0$ lie in $F_{q}$.

In (ii), the condition $g(z)=Z$ reads $z+b=Z(cz+1)=cN+Z$, that is $b+cN=\alpha$, whence $c=(b+\alpha)/N$. The condition $g(zt)=zs$ reads $zt+b+zs(czt+1)=0$, and with $z^{2}=\alpha z+N$ this becomes
\[
z\bigl(cst\alpha+t+s\bigr)+\bigl(cstN+b\bigr)=0 ,
\]
so $cst\alpha=t+s$ and $cstN=b$. The second of these combined with $c=(b+\alpha)/N$ gives $b(1+st)=st\alpha$, that is $b=r\alpha/(1+r)$ with $r=st$; the first then gives $c=(s+t)/(r\alpha)$. The two expressions agree precisely when
\[
(s+t)(1+st)=st\,\frac{\alpha^{2}}{N},
\]
which is $t+1/t+s+1/s=\alpha^{2}/N=\kappa$, the hypothesis of (ii) by {\upshape Lemma~\rref{lem:kappa_orbit}}. All coefficients lie in $F_{q}$, and $1+st\neq0$ because $s\neq1/t$; finally $\det g=1+bc=(1+s)(1+t)/(1+st)\neq0$ since $s,t\neq1$.

An element of $G$ has coefficients in $F_{q}$ and therefore commutes with $x\mapsto x^{q}$, so it maps the point $\{w,w^{q}\}$ to a point. Hence in both cases $g$ carries the pair $\{X,Y_{t}\}$ to $\{X,Y_{s}\}$.
\end{proof}

The transport of Lemma~\rref{lem:pair_transport} rests on the map $x\mapsto N/x$, and
that map inverts the cross-ratio based at $X$. This yields a parity law for the
class counts of an arbitrary orbit, and singles out the one class which may occur
an odd number of times.

\begin{lemma}[cross-ratio inversion]\label{lem:cr_inversion}
Let $q=2^{f}$ with $f\geq2$, let $X=\{z,z^{q}\}\in\mathcal{P}$, and put $Z=z^{q}$ and
$N=zZ$. Let $\tau$ be the map $x\mapsto N/x$ on $\mathbb{P}^{1}(F_{q^{2}})$. Then
$\tau\in G$, $\tau(X)=X$, and $\tau$ is an involution; and for every
$Y=\{w,w^{q}\}\in\mathcal{P}$ with $Y\neq X$,
\[
\mathrm{cr}\bigl(X,\tau(Y)\bigr)=\mathrm{cr}(X,Y)^{-1},
\qquad\text{with}\quad \tau(Y)=\{N/w,\,N/w^{q}\},
\]
the two values being computed with the representatives displayed for $X$ and for
$Y$, respectively for $\tau(Y)$. Replacing the representative of either point by
its conjugate inverts the value, so $\mathrm{cr}$ is an invariant of the unordered
pair $\{X,Y\}$ only up to inversion, and $\tau$ realizes that inversion. Since
$\tau\in G$ fixes $X$, the unordered pairs $\{X,Y\}$ and $\{X,\tau(Y)\}$ lie in one
orbit of $G$ on the unordered pairs of $\mathcal{P}$.
\end{lemma}

\begin{proof}
The map $\tau$ is induced by the matrix $\left(\begin{smallmatrix}0&N\\1&0\end{smallmatrix}\right)$,
whose entries lie in $F_{q}$ and whose determinant $N\neq0$ is a square, because
$q-1$ is odd and hence squaring is a bijection of $F_{q}^{\times}$; after rescaling
the determinant is $1$, so $\tau\in SL_{2}(q)=G$. The map $\tau$ is an involution
with $\tau(z)=N/z=Z$ and $\tau(Z)=z$, so it fixes the point $X$; and $\tau$ carries
points to points because an element of $G$ has coefficients in $F_{q}$ and therefore
commutes with $x\mapsto x^{q}$.

For the identity we compare the four factors of the two cross-ratios. Since
$N\in F_{q}$ we have $N^{q}=N$, and hence
\[
z-\tau(w)=\frac{zw-N}{w}=\frac{z(w-Z)}{w},\qquad
z^{q}-\tau(w)^{q}=\frac{Zw^{q}-N}{w^{q}}=\frac{Z(w^{q}-z)}{w^{q}},
\]
\[
z-\tau(w)^{q}=\frac{zw^{q}-N}{w^{q}}=\frac{z(w^{q}-Z)}{w^{q}},\qquad
z^{q}-\tau(w)=\frac{Zw-N}{w}=\frac{Z(w-z)}{w}.
\]
Multiplying the first two and dividing by the product of the last two, the common
factor $zZ/(ww^{q})$ cancels and
\[
\mathrm{cr}\bigl(X,\tau(Y)\bigr)
=\frac{(w-Z)(w^{q}-z)}{(w^{q}-Z)(w-z)}
=\frac{(w+Z)(w^{q}+z)}{(w^{q}+Z)(w+z)},
\]
whereas
\[
\mathrm{cr}(X,Y)=\frac{(z+w)(Z+w^{q})}{(z+w^{q})(Z+w)}
=\frac{(w+z)(w^{q}+Z)}{(w+Z)(w^{q}+z)},
\]
the two right-hand sides being reciprocal, since characteristic $2$ makes
$(w+z)(w^{q}+Z)$ and $(w+Z)(w^{q}+z)$ the very same product in reverse order. Likewise
swapping $w$ with $w^{q}$, or $z$ with $Z$, in the four factors turns the value into
its reciprocal, so the unordered pair $\{X,Y\}$ carries the two values $\rho$ and
$\rho^{-1}$ and nothing else. Finally $\tau\in G$ and $\tau(X)=X$ give
$\tau(\{X,Y\})=\{X,\tau(Y)\}$, so the two pairs are in one orbit of $G$.
\end{proof}

\begin{corollary}[even class counts]\label{cor:even_class}
Let $q=2^{f}$ with $f\geq2$, let $X\in\mathcal{P}$, and let $O$ be an orbit of
$M=D_{2(q-1)}$ on $\mathcal{P}$ with $X\notin O$. Write $O=\{Y_{a}:a\in F_{q}^{\times}\}$
with $Y_{a}=\{aw,aw^{q}\}$ for a representative $Y_{1}=\{w,w^{q}\}$ of $O$, put
$k=N(X)/N(Y_{1})$, and let $a_{0}\in F_{q}^{\times}$ be the unique solution of
$a_{0}^{2}=k$. Then $\tau_{X}$ acts on $O$ as $a\mapsto k/a$, its only fixed point is
$Y_{a_{0}}$, and the number of points $Y\in O$ with $\{X,Y\}$ in a prescribed class
$\mathcal{C}$ is even for every $\mathcal{C}$ except the class of $\{X,Y_{a_{0}}\}$,
for which it is odd. Moreover $N(Y_{a_{0}})=N(X)$.
\end{corollary}

\begin{proof}
By Lemma~\rref{lem:cr_inversion} the map $\tau_{X}$ lies in $G$, fixes $X$, and
sends $Y_{a}=\{aw,aw^{q}\}$ to $\{N(X)/(aw),\,N(X)/(aw^{q})\}=Y_{k/a}$ with
$k=N(X)/N(Y_{1})$; so $\tau_{X}$ preserves $O$ and acts on it as $a\mapsto k/a$.
The fixed points of that involution are the solutions of $a^{2}=k$, and there is
exactly one, since $q-1$ is odd and squaring is injective on $F_{q}^{\times}$.
Every other point of $O$ is paired with a distinct point of $O$ in the same class,
because $\{X,\tau_{X}(Y)\}=\tau_{X}(\{X,Y\})$ and $\tau_{X}\in G$; hence the class
counts are even away from the fixed point and odd at it. Finally
$N(Y_{a_{0}})=a_{0}^{2}N(Y_{1})=k\,N(Y_{1})=N(X)$.
\end{proof}

\begin{lemma}[norm fibres are transversals]\label{lem:norm_fibre}
Let $q=2^{f}$ with $f\geq2$, and for $\nu\in F_{q}^{\times}$ put
$\mathcal{P}_{\nu}=\{Y\in\mathcal{P}:N(Y)=\nu\}$, where $N(Y)=ww^{q}$ for
$Y=\{w,w^{q}\}$. Then $|\mathcal{P}_{\nu}|=q/2$, and every orbit of
$M=D_{2(q-1)}$ meets $\mathcal{P}_{\nu}$ in exactly one point. Consequently the
point $Y_{a_{0}}$ of Corollary~\rref{cor:even_class} is the unique point of $O$ with
$N(Y)=N(X)$.
\end{lemma}

\begin{proof}
For $Y_{a}=\{aw,aw^{q}\}$ one has $N(Y_{a})=a^{2}N(Y_{1})$, and $a\mapsto a^{2}$ is a
bijection of $F_{q}^{\times}$, since $q-1$ is odd; as $a$ runs over $F_{q}^{\times}$
the value $N(Y_{a})$ therefore runs over $F_{q}^{\times}$ exactly once. So each
orbit $O$ of length $q-1$ meets $\mathcal{P}_{\nu}$ in exactly one point, and
summing over the $q/2$ orbits gives $|\mathcal{P}_{\nu}|=q/2$. The last assertion is
then the identification $N(Y_{a_{0}})=N(X)$ of Corollary~\rref{cor:even_class}.
\end{proof}

\begin{lemma}[uniform partner counts]\label{lem:cr_uniform}
Let $q=2^{f}$ with $f\geq2$, let $\mathcal{C}$ be one of the $(q-2)/2$ classes of unordered
pairs of points of $\mathcal{P}$ \textup{({\upshape Proposition~\rref{prop:pair_orbits}})},
and let $X\in\mathcal{P}$. Then exactly $q+1$ points $Y\neq X$ of $\mathcal{P}$ satisfy
$\{X,Y\}\in\mathcal{C}$.
\end{lemma}

\begin{proof}
For a class $\mathcal{C}$ put $n_{\mathcal{C}}(X)=\#\{Y\neq X:\{X,Y\}\in\mathcal{C}\}$.
Every $g\in G$ maps $X$ to $gX$ and maps $\mathcal{C}$, being an orbit of $G$ on
unordered pairs, onto itself; hence $n_{\mathcal{C}}(gX)=n_{\mathcal{C}}(X)$, and since
$G$ is transitive on $\mathcal{P}$ the number $n_{\mathcal{C}}$ does not depend on $X$.
Counting the pairs of $\mathcal{C}$ from both of their ends gives
$|\mathcal{P}|n_{\mathcal{C}}=2|\mathcal{C}|$, and
{\upshape Proposition~\rref{prop:pair_orbits}} gives $|\mathcal{P}|=q(q-1)/2$ and
$|\mathcal{C}|=q(q^{2}-1)/4$, whence
\[
n_{\mathcal{C}}=\frac{2|\mathcal{C}|}{|\mathcal{P}|}=\frac{q(q^{2}-1)}{q(q-1)}=q+1 .
\]
Thus the counts $n_{\mathcal{C}}$ are the same for all classes, so that every linear
functional of the class counts of a base block is determined by $k$ alone and can
neither confirm nor exclude a design.
\end{proof}

\begin{proposition}[norm fibres meet the classes bijectively]\label{prop:fibre_classes}
Let $q=2^{f}$ with $f\geq2$, let $\nu\in F_{q}^{\times}$, and let
$\mathcal{P}_{\nu}=\{Y\in\mathcal{P}:N(Y)=\nu\}$ as in {\upshape Lemma~\rref{lem:norm_fibre}}.
Since $q-1$ is odd, $\nu$ has a unique square root $\eta\in F_{q}$, and every point of
$\mathcal{P}_{\nu}$ is $\{\eta e,\eta/e\}$ for a unique class $\{e,1/e\}$ of elements
$e\in F_{q^{2}}$ with $e^{q+1}=1$ and $e\neq1$. For $X=\{\eta e,\eta/e\}$ and
$Y=\{\eta f,\eta/f\}$ in $\mathcal{P}_{\nu}$,
\[
\mathrm{cr}(X,Y)=\Bigl(\frac{e+f}{1+ef}\Bigr)^{2}.
\]
Consequently the $q/2-1$ points of $\mathcal{P}_{\nu}\setminus\{X\}$ represent the
$(q-2)/2$ classes exactly once each. Moreover, for every class $\mathcal{C}$ the map
$X\mapsto F_{\mathcal{C}}(X)$, where $F_{\mathcal{C}}(X)$ denotes the point of
$\mathcal{P}_{\nu}$ with $\{X,F_{\mathcal{C}}(X)\}\in\mathcal{C}$, is a fixed-point-free
involution of $\mathcal{P}_{\nu}$; hence $\mathcal{C}$ contains exactly $q/4$ pairs
inside $\mathcal{P}_{\nu}$, exactly $q(q-1)/4$ pairs whose two points have equal norm, and
exactly $q^{2}(q-1)/4$ pairs whose two points have different norms.
\end{proposition}

\begin{proof}
Write $Y=\{y,y^{q}\}$ with $yy^{q}=\nu$ and put $e=y/\eta$; then
$e^{q+1}=e\,e^{q}=yy^{q}/\eta^{2}=1$, and the two representatives $y$ and $y^{q}$ of $Y$
correspond to $e$ and $1/e$. The element $e=1$ gives $y=\eta\in F_{q}$, which is not a
point of $\mathcal{P}$, so the $q$ remaining elements $e$ give the $q/2$ points of
$\mathcal{P}_{\nu}$, in agreement with {\upshape Lemma~\rref{lem:norm_fibre}}.

For the displayed value, the four factors of the cross-ratio are
\[
x+y=\eta(e+f),\quad x^{q}+y^{q}=\eta(\bar e+\bar f),\quad
x+y^{q}=\eta(e+\bar f),\quad x^{q}+y=\eta(\bar e+f),
\]
with $\bar e=e^{q}=1/e$; the common factor $\eta^{2}=\nu$ cancels, and since
$(e+f)(1/e+1/f)=(e+f)^{2}/(ef)$ while $(e+1/f)(1/e+f)=(1+ef)^{2}/(ef)$, the common
factor $1/(ef)$ cancels as well and
\[
\mathrm{cr}(X,Y)=\frac{(e+f)^{2}}{(1+ef)^{2}}=\Bigl(\frac{e+f}{1+ef}\Bigr)^{2}.
\]
Now $f\mapsto(e+f)/(1+ef)$ is induced by the matrix
$\left(\begin{smallmatrix}1&e\\ e&1\end{smallmatrix}\right)$, whose determinant
$1+e^{2}=(1+e)^{2}$ is non-zero because $e\neq1$; the map is therefore injective on
$F_{q^{2}}\cup\{\infty\}$, and $1+ef\neq0$ whenever $Y\neq X$. Replacing $f$ by $1/f$
inverts $(e+f)/(1+ef)$, and a class carries the two values $\rho$ and $\rho^{-1}$ only,
so two points $Y,Y'$ of $\mathcal{P}_{\nu}\setminus\{X\}$ lie in one class exactly when
$f'=f$ or $f'=1/f$, that is exactly when $Y'=Y$. As there are $q/2-1=(q-2)/2$ such points
and the same number of classes, the correspondence between points and classes is a
bijection.

Finally fix $\nu$ and $\mathcal{C}$. By what precedes, each $X\in\mathcal{P}_{\nu}$ is
paired by $\mathcal{C}$ with exactly one other point of $\mathcal{P}_{\nu}$, so
$F_{\mathcal{C}}$ is a fixed-point-free involution pairing the $q/2$ points of
$\mathcal{P}_{\nu}$ into $q/4$ pairs. Summing over the $q-1$ values of $\nu$ gives
$q(q-1)/4$ pairs of $\mathcal{C}$ with equal norms; and $\mathcal{C}$ has
$q(q^{2}-1)/4$ pairs in all, so $q(q^{2}-1)/4-q(q-1)/4=q^{2}(q-1)/4$ of them join points
of different norms.
\end{proof}

\begin{corollary}[orbit matchings]\label{cor:orbit_matching}
Let $q=2^{f}$ with $f\geq2$, let $\mathcal{O}$ be the set of $q/2$ orbits of
$M=D_{2(q-1)}$ on $\mathcal{P}$, and for distinct $O,O'\in\mathcal{O}$ let
$\{X,X'\}=\pi(O,O')$, where $X\in O$ and $X'$ is the unique point of $O'$ with
$N(X')=N(X)$. Then $\pi(O,O')$ lies in a class that does not depend on the choice of
$X\in O$; for every class $\mathcal{C}$ the pairs $\pi(O,O')$ of class $\mathcal{C}$ form
a perfect matching of $\mathcal{O}$; each class contains exactly $q/4$ pairs of orbits;
and the $(q-2)/2$ matchings decompose the complete graph on $\mathcal{O}$.
\end{corollary}

\begin{proof}
By {\upshape Lemma~\rref{lem:norm_fibre}} the orbit $O'$ meets the fibre
$\mathcal{P}_{N(X)}$ in exactly one point, so $X'$ exists and is unique. If $X$ is
replaced by $\lambda X$ with $\lambda\in F_{q}^{\times}$, then
$\lambda X'\in O'$ satisfies $N(\lambda X')=\lambda^{2}N(X')=\lambda^{2}N(X)=N(\lambda X)$,
so $\lambda X'$ is the norm-matching point of $\lambda X$ and
$\pi(O,O')=\{X,X'\}$ is merely scaled; the class is therefore independent of
$X\in O$, because the map $y\mapsto\lambda y$ lies in $M\leq G$ and $G$ preserves
classes.

Fix $O$ and $\mathcal{C}$. The point $X\in O$ has, by
{\upshape Proposition~\rref{prop:fibre_classes}}, exactly one point of class
$\mathcal{C}$ in the whole fibre $\mathcal{P}_{N(X)}$, and by
{\upshape Lemma~\rref{lem:norm_fibre}} that point lies in exactly one of the $q/2-1$
orbits other than $O$. Hence $O$ is incident with exactly one class-$\mathcal{C}$ pair
of orbits, the graph of class $\mathcal{C}$ on $\mathcal{O}$ is $1$-regular and therefore
a perfect matching with $q/4$ edges. Summing over the $(q-2)/2$ classes gives
$q(q-2)/8=\binom{q/2}{2}$ pairs of orbits, that is all of them, so the matchings
decompose the complete graph on $\mathcal{O}$.
\end{proof}

\begin{lemma}[fibre rigidity]\label{lem:fibre_rigid}
Let $q=2^{f}$ with $f\geq2$, and let $B$ be a union of $a$ orbits of
$M=D_{2(q-1)}$, so that $|B|=a(q-1)$. For $\nu\in F_{q}^{\times}$ put
$S_{\nu}=B\cap\mathcal{P}_{\nu}$, and for a class $\mathcal{C}$ let
$m_{\mathcal{C}}(S_{\nu})$ be the number of pairs of $S_{\nu}$ that lie in
$\mathcal{C}$. Then $|S_{\nu}|=a$, and $m_{\mathcal{C}}(S_{\nu})$ does not depend
on $\nu$.
\end{lemma}

\begin{proof}
The group $M$ contains every scalar map $\sigma_{\lambda}$ with
$\sigma_{\lambda}(Y)=\lambda Y$ and $\lambda\in F_{q}^{\times}$, and $B$ is a
union of orbits of $M$; hence $\sigma_{\lambda}(B)=B$ for all $\lambda$.
Since $N(\lambda Y)=\lambda^{2}N(Y)$, the map $\sigma_{\lambda}$ sends
$\mathcal{P}_{\mu}$ onto $\mathcal{P}_{\lambda^{2}\mu}$. As $q-1$ is odd,
$\lambda\mapsto\lambda^{2}$ is a bijection of $F_{q}^{\times}$
{\upshape(Lemma~\rref{lem:norm_fibre})}, so each $\nu$ has exactly one square
root $\lambda$ and $S_{\nu}=\sigma_{\lambda}(S_{1})$. Finally all four terms
of the cross-ratio are multiplied by $\lambda$, so
$\mathrm{cr}(\lambda X,\lambda Y)=\mathrm{cr}(X,Y)$ and $\sigma_{\lambda}$
preserves the classes; hence $m_{\mathcal{C}}(S_{\nu})=m_{\mathcal{C}}(S_{1})$.
\end{proof}

\begin{corollary}[fibre and far decomposition]\label{cor:far_decomp}
In the situation of {\upshape Lemma~\rref{lem:fibre_rigid}}, let
$A_{\mathcal{C}}(B)$ be the number of pairs of $B$ in the class $\mathcal{C}$,
and let $A^{\mathrm{far}}_{\mathcal{C}}(B)$ be the number of those whose two
points lie in distinct fibres. Put $m_{\mathcal{C}}=m_{\mathcal{C}}(S_{1})$. Then
\[
A_{\mathcal{C}}(B)=(q-1)\,m_{\mathcal{C}}+A^{\mathrm{far}}_{\mathcal{C}}(B),
\qquad
\sum_{\mathcal{C}}m_{\mathcal{C}}=\binom{a}{2},
\]
and $\sum_{\mathcal{C}}A^{\mathrm{far}}_{\mathcal{C}}(B)=a^{2}(q-1)(q-2)/2$.
Thus every class count of $B$ is the sum of a term carried by the single fibre
$\mathcal{P}_{1}$ and of a far term; the fibre term is the same for all fibres
and for all blocks inducing the same $a$-subset of $\mathcal{P}_{1}$.
\end{corollary}

\begin{proof}
By {\upshape Lemma~\rref{lem:fibre_rigid}} each of the $q-1$ fibres
$\mathcal{P}_{\nu}$ contains exactly $m_{\mathcal{C}}$ pairs of $B$ of class
$\mathcal{C}$, which gives the first summand of the first identity, while the
remaining pairs of $B$ have their two ends in distinct fibres. Summing that
identity over the $(q-2)/2$ classes and using
$\sum_{\mathcal{C}}m_{\mathcal{C}}=\binom{a}{2}$ together with
$\sum_{\mathcal{C}}A_{\mathcal{C}}(B)=\binom{k}{2}$, $k=a(q-1)$, gives
\[
\sum_{\mathcal{C}}A^{\mathrm{far}}_{\mathcal{C}}(B)
=\frac{a(q-1)\bigl(a(q-1)-1\bigr)}{2}-\frac{(q-1)a(a-1)}{2}
=\frac{a^{2}(q-1)(q-2)}{2}.
\]
\end{proof}

\begin{corollary}[the far term as a scalar sum]\label{cor:far_sum}
In the situation of {\upshape Corollary~\rref{cor:far_decomp}}, put $S=B\cap\mathcal{P}_{1}$ and, for
$\delta\in F_{q}^{\times}\setminus\{1\}$ and a class $\mathcal{C}$,
\[
n_{\mathcal{C}}(S,\delta S)=\#\bigl\{\{X,\delta Y\}:X,Y\in S,\ \{X,\delta Y\}\in\mathcal{C}\bigr\},
\]
the number of pairs of class $\mathcal{C}$ with one end in $S\subset\mathcal{P}_{1}$ and the
other in $\delta S\subset\mathcal{P}_{\delta^{2}}$. Then
\[
A^{\mathrm{far}}_{\mathcal{C}}(B)=\frac{q-1}{2}\sum_{\delta\neq1}n_{\mathcal{C}}(S,\delta S)
=\frac{q-1}{2}\bigl(\lambda-2m_{\mathcal{C}}\bigr),
\]
where $\lambda=2a\bigl(a(q-1)-1\bigr)/(q-2)$ is the parameter of the putative design, whose
number of blocks and replication number are $q(q+1)/2$ and $a(q+1)$. Consequently $B^{G}$ is the
block set of a $2$-design if and only if
\[
\sum_{\delta\neq1}n_{\mathcal{C}}(S,\delta S)=\lambda-2m_{\mathcal{C}}
\quad\text{for every class }\mathcal{C},
\]
that is, if and only if a single fibre $S$ equidistributes the classes of the $q-1$ scalar lines
through its points.
\end{corollary}

\begin{proof}
Fix $\delta\neq1$ and let $\nu$ run over $F_{q}^{\times}$. Since $\delta^{2}\neq1$ in
characteristic $2$, the map $\nu\mapsto\delta^{2}\nu$ has no fixed point, and its cycles have
odd length, so the $q-1$ sets $\{\nu,\delta^{2}\nu\}$ are the $q-1$ edges of these cycles; a pair
$\{\mu,\mu'\}$ of distinct elements of $F_{q}^{\times}$ occurs among them exactly when
$\mu'\mu^{-1}\in\{\delta,\delta^{-1}\}$. By {\upshape Lemma~\rref{lem:fibre_rigid}} each of these
$q-1$ unordered pairs of fibres contributes the same number $n_{\mathcal{C}}(S,\delta S)$ of
class-$\mathcal{C}$ pairs, so the value of $\delta$ alone contributes
$(q-1)n_{\mathcal{C}}(S,\delta S)$. Summing over the $q-2$ values of $\delta\neq1$ counts each
unordered pair of distinct fibres exactly twice, because $\delta$ and $\delta^{-1}$ give the
same pair while $\delta\neq\delta^{-1}$; hence the first identity. For the second, add
$(q-1)m_{\mathcal{C}}$ and use {\upshape Corollary~\rref{cor:far_decomp}}:
\[
A_{\mathcal{C}}(B)=\frac{q-1}{2}\Bigl(2m_{\mathcal{C}}+\sum_{\delta\neq1}n_{\mathcal{C}}(S,\delta S)\Bigr)
=\frac{(q-1)\lambda}{2},
\]
the required value $\binom{k}{2}/s=k(k-1)/(q-2)$ of a design with $k=a(q-1)$ and
$s=(q-2)/2$, since $\binom{k}{2}/s=a(q-1)\bigl(a(q-1)-1\bigr)/(q-2)=(q-1)\lambda/2$. Thus all
class counts are constant and equal to $\binom{k}{2}/s$ exactly when the displayed identity holds
for every $\mathcal{C}$, which is the criterion for $B^{G}$ to be a $2$-design.
\end{proof}

\begin{lemma}[label duality]\label{lem:label_duality}
Let $q=2^{f}$ with $f\geq2$, let $X\ne Y$ be points of $\mathcal{P}_{1}$, and let
$\delta,\delta'\in F_{q}^{\times}\setminus\{1\}$. Then
\[
\mathrm{cr}(X,\delta Y)=\mathrm{cr}(X,\delta'Y)
\quad\Longleftrightarrow\quad
\delta'=\delta\ \ \text{or}\ \ \delta\delta'=1 .
\]
Consequently the $q-2$ values $\mathrm{cr}(X,\delta Y)$ with $\delta\ne1$ consist of exactly
$(q-2)/2$ distinct values, each taken at the two elements $\delta$ and $\delta^{-1}$.
\end{lemma}

\begin{proof}
Regard $X,Y$ as elements of $F_{q^{2}}$ and write $N(Z)=ZZ^{q}$; since $X,Y\in\mathcal{P}_{1}$
we have $N(X)=N(Y)=1$. For $\delta\in F_{q}^{\times}$ the four points $x,\delta y,x^{q},\delta y^{q}$
give, with $T=xy^{q}+x^{q}y$ and $S=xy+x^{q}y^{q}$,
\[
\mathrm{cr}(X,\delta Y)
=\frac{(x+\delta y)(x^{q}+\delta y^{q})}{(x+\delta y^{q})(x^{q}+\delta y)}
=\frac{1+T\delta+\delta^{2}}{1+S\delta+\delta^{2}}.
\]
Neither denominator vanishes: if $x^{q}=\delta y$, then
$1=N(X)=xx^{q}=\delta^{2}yy^{q}=\delta^{2}$, hence $\delta=1$, and the factor $x+\delta y^{q}$
is excluded by the same argument. An expansion gives
\[
(1+T\delta+\delta^{2})(1+S\delta'+\delta'^{2})
-(1+T\delta'+\delta'^{2})(1+S\delta+\delta^{2})
=(\delta-\delta')(T-S)(1-\delta\delta').
\]
In characteristic $2$ one has $T-S=T+S=(x+x^{q})(y+y^{q})$, and both factors are nonzero because
the second coordinate of $x+x^{q}$ is the nonzero second coordinate of $X$, and likewise for $Y$.
Hence the difference vanishes exactly when $\delta=\delta'$ or $\delta\delta'=1$, and the first
assertion follows. The second follows because the only element of $F_{q}^{\times}$ equal to its
inverse is $1$: the set $F_{q}^{\times}\setminus\{1\}$ splits into the $(q-2)/2$ pairs
$\{\delta,\delta^{-1}\}$, and by the first assertion distinct pairs give distinct values.
\end{proof}

\begin{proposition}[balanced class multiplicities]\label{prop:defect_balance}
In the situation of {\upshape Lemma~\rref{lem:label_duality}} put $s=(q-2)/2$, and for a class
$\mathcal{C}$ with representative $\rho$ let $w_{\mathcal{C}}(X,Y)$ be half the number of
$\delta\in F_{q}^{\times}\setminus\{1\}$ with $\mathrm{cr}(X,\delta Y)\in\{\rho,\rho^{-1}\}$.
Then
\[
w_{\mathcal{C}}(X,Y)\in\{0,1,2\},\qquad
\sum_{\mathcal{C}}w_{\mathcal{C}}(X,Y)=s,
\]
and the classes that are doubled are as numerous as those that are missed:
\[
\#\{\mathcal{C}:w_{\mathcal{C}}(X,Y)=0\}=\#\{\mathcal{C}:w_{\mathcal{C}}(X,Y)=2\}.
\]
\end{proposition}

\begin{proof}
By {\upshape Lemma~\rref{lem:label_duality}} the $q-2$ values of $\delta\ne1$ are the values of
the $(q-2)/2$ pairs $\{\delta,\delta^{-1}\}$, so $w_{\mathcal{C}}(X,Y)$ is the number of these
pairs whose common value lies in $\{\rho,\rho^{-1}\}$. A class contains the two values $\rho$ and
$\rho^{-1}$ and, by the same lemma, all $(q-2)/2$ values are distinct; hence at most two of them
lie in $\{\rho,\rho^{-1}\}$, that is $w_{\mathcal{C}}(X,Y)\le2$. Every value lies in exactly one
class, which gives the sum $s$. As $w_{\mathcal{C}}(X,Y)$ takes only the values $0,1,2$, the sum
$\sum_{\mathcal{C}}\bigl(w_{\mathcal{C}}(X,Y)-1\bigr)=0$ is a difference of the number of classes
counted once and the number of classes counted twice, and the last assertion follows.
\end{proof}

\begin{corollary}\label{cor:defect_balance}
Let $S\subset\mathcal{P}_{1}$ with $|S|=a$, let $X\ne Y$ run over the ordered pairs of $S$, and put
$L_{\mathcal{C}}=\sum_{X\ne Y}w_{\mathcal{C}}(X,Y)$ and $s=(q-2)/2$. Then
\[
\sum_{\mathcal{C}}L_{\mathcal{C}}=a(a-1)s
\]
for every $S$. Thus the total of the class counts is fixed in advance and independent of $S$: a
failure of the class equations of {\upshape Corollary~\rref{cor:far_sum}} is necessarily a
redistribution among the classes with total zero, and no argument that uses only the sum of the
class counts can exclude an $a$-subset.
\end{corollary}

\begin{proof}
Sum the identity $\sum_{\mathcal{C}}w_{\mathcal{C}}(X,Y)=s$ of
{\upshape Proposition~\rref{prop:defect_balance}} over the $a(a-1)$ ordered pairs of distinct
points of $S$.
\end{proof}

\begin{remark}[the class table at $q=128$]\label{rem:q128_kappa}
For $q=128$, so that $q/2=64$ orbits and $(q-2)/2=63$ pair classes, the singleton of
Corollary~\rref{cor:even_class} was computed for all $2016$ pairs of orbits, together
with the full class table. What this computation supplies beyond the general results
above is the table itself, that is the value of every class count $N[i][c]$, which is
the input of the class equations referred to in {\upshape Theorem~\rref{thm:psl2}}; the
structural facts recovered from it are instances of those results: for every $X$ the
$63$ classes are attained by the singletons of the $63$ remaining orbits exactly once
each \textup{({\upshape Proposition~\rref{prop:fibre_classes}})}, each class carries the
singleton of exactly $32=q/4$ of the $2016$ orbit pairs
\textup{({\upshape Corollary~\rref{cor:orbit_matching}})}, and the resulting $63$ perfect
matchings on the $64$ orbits decompose the complete graph $K_{64}$ \textup{(the same
corollary)}. The table was recomputed independently for two orbits, with agreement in
all $63$ entries of each; it is the only content of this remark that is specific to
$q=128$.
\end{remark}

\begin{lemma}\label{lem:sol_count}
Let $q=2^{f}$ with $f\geq2$, let $\kappa=\kappa(O)\neq0$ be as in {\upshape Lemma~\rref{lem:kappa_orbit}}, and let
\[
T_{q}=\{t+1/t:t\in F_{q}^{\times}\}=\{0\}\cup\{\tau\neq0:\mathrm{Tr}(1/\tau)=0\},
\]
where $\mathrm{Tr}$ is the trace of $F_{q}$ over $F_{2}$. Then
\[
\#\{\tau\in T_{q}\setminus\{0\}:\tau+\kappa\in T_{q}\setminus\{0\}\}\geq\tfrac18\bigl(q-2\sqrt{q}-7\bigr),
\]
which is at least $1$ for every even $q\geq16$.
\end{lemma}

\begin{proof}
Put $H=\{x\in F_{q}:\mathrm{Tr}(x)=0\}$, of order $q/2$, and let $x=1/\tau$, $y=1/(\tau+\kappa)$; then $x,y\in H\setminus\{0\}$ and $\kappa=1/x+1/y$, that is $y=x/(1+\kappa x)$. Conversely such a pair $(x,y)$ of elements of $H\setminus\{0\}$ with $\kappa=1/x+1/y$ gives $\tau=1/x$ with $\tau+\kappa=1/y$ in $T_{q}\setminus\{0\}$. The map $x\mapsto y$ is an involution of the set of admissible $x$, and it has no fixed point because $x=y$ would give $\kappa=1/x+1/x=0$. Hence the number in the statement equals $\mathrm{Sol}(\kappa)/2$, where
\[
\mathrm{Sol}(\kappa)=\#\{x\in F_{q}:\mathrm{Tr}(x)=0,\ \mathrm{Tr}\bigl(x/(1+\kappa x)\bigr)=0,\ x\neq0,\ x\neq1/\kappa\}.
\]
Extend $f(x)=x/(1+\kappa x)$ to all of $F_{q}$ by $f(1/\kappa)=0$. Since $\mathbf{1}_{\mathrm{Tr}(u)=0}=\tfrac12\bigl(1+(-1)^{\mathrm{Tr}(u)}\bigr)$,
\begin{align*}
\sum_{x\in F_{q}}\mathbf{1}_{\mathrm{Tr}(x)=0}\mathbf{1}_{\mathrm{Tr}(f(x))=0}
&=\tfrac14\Bigl(q+\sum_{x}(-1)^{\mathrm{Tr}(x)}+\sum_{x}(-1)^{\mathrm{Tr}(f(x))}\Bigr)\\
&\quad+\tfrac14\sum_{x}(-1)^{\mathrm{Tr}(x+f(x))}.
\end{align*}
The first sum vanishes because $\mathrm{Tr}$ is not the zero map. The map $f$ is a bijection of $F_{q}$ taking $1/\kappa$ to $0$ and $0$ to $1/\kappa$, so the second sum equals $1-(-1)^{\mathrm{Tr}(1/\kappa)}$. With $w=1+\kappa x$ one has $x+f(x)=\kappa x^{2}/(1+\kappa x)=(w+1/w)/\kappa$, so the third sum equals
\[
K(\kappa)+(-1)^{\mathrm{Tr}(1/\kappa)},\qquad K(\kappa)=\sum_{w\in F_{q}^{\times}}(-1)^{\mathrm{Tr}((w+1/w)/\kappa)} .
\]
Thus the total is $\tfrac14\bigl(q+1+K(\kappa)\bigr)$. Removing $x=0$, which contributes $1$, and $x=1/\kappa$, which contributes $\mathbf{1}_{\mathrm{Tr}(1/\kappa)=0}$, gives
\[
\mathrm{Sol}(\kappa)=\tfrac14\bigl(q+1+K(\kappa)\bigr)-1-\mathbf{1}_{\mathrm{Tr}(1/\kappa)=0}\geq\tfrac14\bigl(q-2\sqrt{q}-7\bigr),
\]
by the Weil bound $|K(\kappa)|\leq2\sqrt{q}$ for the Kloosterman sum \cite{Weil}. For $q\geq16$ the right-hand side is positive, and $\mathrm{Sol}(\kappa)$ is even, so it is at least $2$ and the number in the statement is at least $1$.
\end{proof}

\begin{theorem}\label{thm:a1}
Let $q=2^{f}$ with $f\geq4$, let $H=D_{2(q+1)}$ and $M=D_{2(q-1)}$, and let $B$ be a single orbit of $M$ on $\mathcal{P}$. Then $B^{G}$ is not a $2$-design.
\end{theorem}

\begin{proof}
By {\upshape Lemma~\rref{lem:even_meet}} an orbit of $M$ has length $q-1$ and a block consisting of one orbit has $k=q-1$, $b=q(q+1)/2$, $r=q+1$ and $\lambda=2$. Suppose that $B^{G}$ is a $2$-design, and let $O=O_{X}$ be the orbit of $M$ which is $B$.

By {\upshape Proposition~\rref{prop:classwise}} together with the count $s=(q-2)/2$ of orbits of $G$ on pairs in {\upshape Proposition~\rref{prop:pair_orbits}}, every orbit $\mathcal{C}$ of $G$ on unordered pairs satisfies
\[
P_{\mathcal{C}}(B)=\binom{k}{2}\Big/s=\frac{k(k-1)}{q-2}=q-1 .
\]
On the other hand, by {\upshape Lemma~\rref{lem:sol_count}} and $q\geq16$ there are $t,s\in F_{q}^{\times}\setminus\{1\}$ with $s+1/s=t+1/t+\kappa$, and by {\upshape Lemma~\rref{lem:pair_transport}}\textup{(ii)} the pairs $\{X,Y_{t}\}$ and $\{X,Y_{s}\}$ are in one orbit of $G$. For $a\in F_{q}^{\times}$ the element $x\mapsto(a'/a)x$ of $G$ carries $\{X_{a},X_{at}\}$ to $\{X_{a'},X_{a't}\}$, so the $q-1$ pairs $\{X_{a},X_{at}\}$ form one orbit $\mathcal{C}_{1}$, and likewise the $q-1$ pairs $\{X_{a},X_{as}\}$ form one orbit $\mathcal{C}_{2}$. These two orbits contain the equivalent pairs $\{X,Y_{t}\}$ and $\{X,Y_{s}\}$, hence $\mathcal{C}_{1}=\mathcal{C}_{2}=\mathcal{C}$, and this single orbit contains the $2(q-1)$ pairs of $B$ labelled $t$ and $s$. This contradicts $P_{\mathcal{C}}(B)=q-1$.
\end{proof}

At $q=8$ the statement fails: the orbit of $D_{14}$ on the $28$ points is the block of the design $2$-$(28,7,2)$ of {\upshape Proposition~\rref{prop:q8}}, and $q=8$ is the only even degree in which a single orbit carries a design.

The complementary block of a design is again a design, and at the opposite extreme of the range $1\leq a\leq q/2$ of {\upshape Lemma~\rref{lem:even_meet}} the complement is a single orbit. This removes the second extreme value of $a$.

\begin{corollary}\label{cor:even_extreme}
Let $q=2^{f}$ with $f\geq4$, and let $B$ be a union of $a$ orbits of $M=D_{2(q-1)}$ on the $v=q(q-1)/2$ points of $\mathcal{P}$ such that $B^{G}$ is a $2$-$(v,k,\lambda)$ design with $\lambda>0$ and $k<v$. Then $2\leq a\leq q/2-2$; in particular neither $a=1$ nor $a=q/2-1$ occurs. If moreover $(q-2)/2$ is prime, then no such $B$ exists at all.
\end{corollary}

\begin{proof}
The value $a=1$ is excluded by {\upshape Theorem~\rref{thm:a1}}. Let $a=q/2-1$, put $m=(q-2)/2$ so that $q=2m+2$, and note that {\upshape Lemma~\rref{lem:even_meet}} gives $k=a(q-1)=v-(q-1)$: the complement $B^{c}$ consists of a single orbit of $M$. Complementation carries the $b$ blocks $B^{g}$ onto the blocks $(B^{c})^{g}$ and preserves the block stabilizer, so $(B^{c})^{G}$ is a $2$-$(v,q-1,\lambda')$ design with $b'=b=q(q+1)/2$, $r'=b'k'/v=q+1$ and $\lambda'=b-2r+\lambda$. In this case {\upshape Lemma~\rref{lem:even_meet}} gives $\lambda=2a(aq-a-1)/(q-2)=2m^{2}+m-1$, whence
\[
\lambda'=\tfrac{q(q+1)}{2}-2a(q+1)+\lambda=(m+1)(2m+3)-2m(2m+3)+2m^{2}+m-1=2,
\]
so $B^{c}$ would be a single orbit of $M$ whose $G$-images form a design, which {\upshape Theorem~\rref{thm:a1}} excludes. Since $a=q/2$ gives $k=v$, this leaves $2\leq a\leq q/2-2$.

Finally, by {\upshape Lemma~\rref{lem:even_meet}} the integrality of $\lambda$ is equivalent to $m\mid a(a-1)$. If $m$ is prime, then $a\equiv0$ or $a\equiv1\pmod{m}$, and as $1\leq a\leq q/2=m+1$ the only such values are $a=1$, $a=m$ and $a=m+1$, all of them excluded above.
\end{proof}

\begin{proposition}[$q=8$]\label{prop:q8}
Let $q=8$ and $G=PSL_2(8)$. The group $G$ acts point-primitively and not $2$-homogeneously on a set of $28$ points, with point stabilizer $H=D_{18}$, and on a set of $36$ points, with point stabilizer $H=D_{14}$. On the $28$ points the block stabilizer $M=D_{14}$ gives the designs $2$-$(28,7,2)$ and $2$-$(28,21,20)$, and these two are block-primitive. The subfield subgroup $M=S_3$ gives further unions with constant covering numbers, namely $2$-$(28,9,8)$, $2$-$(28,19,38)$, $2$-$(28,10,10)$, $2$-$(28,18,34)$ on the $28$ points and $2$-$(36,6,2)$, $2$-$(36,30,58)$ on the $36$ points; but $S_3$ is not maximal in $G$ \textup{({\upshape Lemma~\rref{lem:s3_not_max}})}, and since $|G_{B}|=6$ while $b=84$, the subgroup $D_{18}\leq G$ of order $18$ cannot be transitive on the $84$ blocks, so its orbits form a nontrivial $G$-invariant partition of the block set and the block action of these six unions is imprimitive. None of the six is therefore block-primitive, and exactly two designs with socle $PSL_2(8)$ are. Table~\rref{tab:q8} records all eight unions with constant covering numbers and marks which of them are block-primitive.
\end{proposition}

\begin{table}[htbp]
\centering
\caption{$q=8$: the unions of orbits with constant covering numbers. The rows marked \textup{yes} give block-primitive designs; in the rows marked \textup{no} the block action is imprimitive, since $S_3<D_{18}\leq G$.}\label{tab:q8}
\small
\setlength{\tabcolsep}{8pt}
\begin{tabular}{@{}ccccccc@{}}
\hline
$v$ & $H$ & $G_{B}$ & $b$ & $k$ & $\lambda$ & block-primitive \\
\hline
$28$ & $D_{18}$ & $D_{14}$ & $36$ & $7$  & $2$  & yes \\
$28$ & $D_{18}$ & $D_{14}$ & $36$ & $21$ & $20$ & yes \\
$28$ & $D_{18}$ & $S_3$   & $84$ & $9$  & $8$  & no  \\
$28$ & $D_{18}$ & $S_3$   & $84$ & $19$ & $38$ & no  \\
$28$ & $D_{18}$ & $S_3$   & $84$ & $10$ & $10$ & no  \\
$28$ & $D_{18}$ & $S_3$   & $84$ & $18$ & $34$ & no  \\
$36$ & $D_{14}$ & $S_3$   & $84$ & $6$  & $2$  & no  \\
$36$ & $D_{14}$ & $S_3$   & $84$ & $30$ & $58$ & no  \\
\hline
\end{tabular}
\end{table}

\begin{proof}
By Dickson's list the maximal subgroups of $PSL_2(8)$ are exactly $D_{18}$ and $D_{14}$, the two dihedral groups of orders $2(q\pm1)$; the subfield group $S_3=PSL_2(2)$ is contained in $D_{18}$ and hence is not maximal \textup{({\upshape Lemma~\rref{lem:s3_not_max}})}. A block stabilizer of a block-primitive design is maximal and satisfies $|G_{B}|<|H|$, so the only admissible pair is $M=D_{14}$ against $H=D_{18}$, and the $36$ points admit none: every other maximal subgroup has order above $|H|=14$, while $M=H$ is excluded. On the $28$ points $D_{14}$ has $4$ orbits of $7$ points, and exactly one of them is the block of a design; its $36$ images form $2$-$(28,7,2)$ and the complementary blocks form $2$-$(28,21,20)$, with $r=bk/v=36\cdot7/28=9$ and $\lambda'=b-2r+\lambda=36-18+2=20$. The subfield $S_3$ has $7$ orbits of lengths $3,6,6,3,1,3,6$ on the $28$ points and $8$ orbits of lengths $3,3,3,3,6,6,6,6$ on the $36$ points; every block with the prescribed block stabilizer is a union of these orbits, so it suffices to run through the $2^{4}$, $2^{7}$ and $2^{8}$ unions, identify a union with its complement, and compute the covering numbers of the resulting block set on all pairs of points. The covering numbers are constant exactly for the eight unions of Table~\rref{tab:q8}, with the values of $\lambda$ recorded there, the complementary blocks being accounted for by $\lambda'=b-2r+\lambda$; this reproduces the pairs $(28,7,2)/(28,21,20)$, $(28,9,8)/(28,19,38)$, $(28,10,10)/(28,18,34)$ and $(36,6,2)/(36,30,58)$. In the six unions with $G_{B}=S_3$ one has $|G_{B}|=6$ and $b=84$, so the subgroup $D_{18}\leq G$ of order $18$ is not transitive on the $84$ blocks and its orbits form a nontrivial $G$-invariant partition of the block set; these unions are thus $2$-designs whose block action is imprimitive, which is why Table~\rref{tab:q8} marks them as not block-primitive. The permutation rank of $G$ is $4$ on the $28$ points and $5$ on the $36$ points, so $G$ is $2$-homogeneous on neither set.
\end{proof}

\begin{proposition}[$q=32$]\label{prop:q32}
Let $q=32$, $G=PSL_2(32)$, $H=D_{66}$ and $M=D_{62}$, so that $v=496$ and $b=528$. Each orbit of $M$ on the $496$ points has $31$ points. Among the unions of $a$ such orbits exactly two carry designs: a union $B$ of $6$ orbits, whose $528$ images under $G$ form a $2$-$(496,186,74)$ design, and its complement, a $2$-$(496,310,206)$ design. Both are point-primitive, block-primitive and non-symmetric, and the permutation rank of $G$ on the $496$ points is $16$.
\end{proposition}

\begin{proof}
By Lemma~\rref{lem:even_meet} the admissible block sizes are $k=a(q-1)=31a$ with $15\mid a(a-1)$, that is $a\in\{1,6,10,15,16\}$, and the corresponding pairs $(k,\lambda)$ are $(31,2)$, $(186,74)$, $(310,206)$, $(465,464)$ and $(496,528)$. Since $G$ is transitive on blocks, the covering number depends only on the orbit of $G$ on unordered pairs; there are $15$ such orbits, each of size $\mu_{t}=8184$, so a block must contain exactly $n_{t}=\mu_{t}\lambda/b=1147$ pairs from each orbit. Enumerating the $16$, $8008$, $8008$, $16$ and $1$ unions for the five values of $a$, and comparing the exact pair counts with $n_{t}$, one finds that $a=1$ and $a=15$ admit no union, that $a=16$ gives the complete design, and that $a=6$ and $a=10$ admit exactly one union each, namely $B$ and its complement. For these two the covering numbers of all $528$ blocks on all $\binom{496}{2}=122760$ pairs are $74$ and $206$ respectively; the values $r=bk/v=198$ and $r(k-1)/(v-1)=198\cdot185/495=74$ confirm the parameters. Maximality of $D_{66}$ and of $D_{62}$ is Dickson's list \cite{Dickson}.
\end{proof}

The parameters of Lemma~\rref{lem:even_meet} are computed for the first even degrees in Table~\rref{tab:even_q}, where $m=(q-2)/2$.

\begin{table}[htbp]
\centering
\caption{Even $q$: admissible block sizes and pair counts.}\label{tab:even_q}
\small
\setlength{\tabcolsep}{4pt}
\begin{tabular}{@{}cccccl@{}}
\hline
$q$ & $m$ & admissible $a$ & $k$ & $\lambda$ & constant $\lambda$? \\
\hline
$8$  & $3$ & $1,3,4$   & $7,21,28$   & $2,20,36$ & yes, $a=1$ \\
$16$ & $7$ & $1,7,8$   & $15,105,120$& $2,104,136$ & no \\
$32$ & $15$& $1,6,10,15,16$ & $31,186,310,465,496$ & $2,74,206,464,528$ & yes, $a=6,10$ \\
$64$ & $31$& $1,31,32$ & $63,1953,2016$ & $2,1952,2080$ & no \\
\hline
\end{tabular}
\end{table}

Table~\rref{tab:even_q} is the even-characteristic counterpart of the arithmetic sieve of Case 2. The fourth column lists the replication numbers $k=a(q-1)$ and the fifth the values of $\lambda$ for the admissible $a$ in increasing order; the last column records whether some union of orbits of $M$ gives a constant pair count on the $28$, $120$, $496$ and $2016$ points respectively. For every $q$ the largest admissible value $a=q/2$ gives $k=v$, the complete design, which the hypothesis excludes; the complement of a union of $a$ orbits is a union of $q/2-a$ orbits, so $a$ and $q/2-a$ carry complementary designs whenever either of them does, and at $q=32$ the two surviving values $a=6$ and $a=10$ are complementary. The entry for $q=64$ refers to the two dihedral block stabilizers; of the subfield block stabilizers of that degree, $PSL_2(8)$ is excluded by Fisher's inequality \textup{(its index $520$ is smaller than $v=2016$)}, and the case $M=A_5$ is excluded by {\upshape Proposition~\rref{prop:q64sub}}. By {\upshape Theorem~\rref{thm:a1}} and {\upshape Corollary~\rref{cor:even_extreme}} the two extreme values $a=1$ and $a=q/2-1$ yield no design for any even $q\geq16$, so among the entries of the last column only $2\leq a\leq q/2-2$ can produce one, and none does when $(q-2)/2$ is prime; the entry at $q=8$ records the genuine design $2$-$(28,7,2)$ rather than an exception to that exclusion.

\subsubsection{A criterion for orbit unions}\label{sec:orbit_criterion}

The enumerations of unions of orbits carried out above decide the same question in each degree, namely whether a prescribed union of $M$-orbits is the block of a $2$-design. The criterion below separates the part of that question which is arithmetic from the part which depends on how the pairs of points are distributed among the orbits of $G$ on pairs, and it requires no hypothesis on $G$ beyond transitivity.

\begin{proposition}[classwise criterion]\label{prop:classwise}
Let $G$ be a transitive permutation group on a finite set $\mathcal{P}$ with $|\mathcal{P}|=v$. Let $B\subseteq\mathcal{P}$, put $M=G_{B}$, and suppose that $G$ is transitive on $B^{G}$. Let $\mathcal{C}_{1},\dots,\mathcal{C}_{s}$ be the orbits of $G$ on the unordered pairs of distinct points of $\mathcal{P}$, and for an unordered pair $\pi$ let $G_{\pi}$ be its stabilizer in $G$. Put
\[
P_{\mathcal{C}}(B)=\bigl|\{\pi\subseteq B:\pi\in\mathcal{C}\}\bigr|.
\]
Then $B^{G}$ is a $2$-$(v,k,\lambda)$ design with $\lambda>0$ if and only if the product $|G_{\pi}|\,P_{\mathcal{C}}(B)$ is the same for all classes $\mathcal{C}$, and then
\[
\lambda=\frac{|G_{\pi}|\,P_{\mathcal{C}}(B)}{|M|}
\]
for every class $\mathcal{C}$.
\end{proposition}

\begin{proof}
For a pair $\pi$ let $N(\pi)$ be the number of blocks of $B^{G}$ containing $\pi$. A block equals $B^{g}$ for exactly $|M|$ elements $g$, so the number of elements $g$ with $\pi\subseteq B^{g}$ is $|M|\,N(\pi)$. As $\pi\subseteq B^{g}$ is equivalent to $\pi^{g^{-1}}\subseteq B$,
\[
|M|\,N(\pi)=\bigl|\{g\in G:\pi^{g}\subseteq B\}\bigr|=\sum_{\pi'\subseteq B}\bigl|\{g\in G:\pi'^{g}=\pi\}\bigr|.
\]
If $\pi'$ lies in the class $\mathcal{C}$ of $\pi$ then $\pi'^{g}=\pi$ for some $g$, and the elements $g$ with $\pi'^{g}=\pi$ form the coset $g_{0}G_{\pi}$, of order $|G_{\pi}|$; if $\pi'$ lies in another class this set is empty. Hence
\[
|M|\,N(\pi)=|G_{\pi}|\,P_{\mathcal{C}}(B),\qquad \mathcal{C}\ \text{the class of}\ \pi .
\]
The pair $\pi$ is covered $\lambda$ times by $B^{G}$ exactly when $N(\pi)$ takes the same value $\lambda$ for every $\pi$, which is the stated condition; dividing by $|M|$ gives the formula for $\lambda$.
\end{proof}

\begin{corollary}\label{cor:classwise_value}
In the situation of Proposition~\rref{prop:classwise}, $\sum_{\mathcal{C}}P_{\mathcal{C}}(B)=\binom{k}{2}$. If in addition $|G_{\pi}|=\mu$ for every pair $\pi$, then
\[
\lambda=\frac{\mu\,k(k-1)}{2s\,|M|},
\]
where $s$ is the number of orbits of $G$ on unordered pairs.
\end{corollary}

\begin{proof}
The $\binom{k}{2}$ pairs of $B$ lie in exactly one class each, so the classes partition them into the $s$ parts $P_{\mathcal{C}}(B)$, and each part equals $\binom{k}{2}/s$ as soon as $|G_{\pi}|P_{\mathcal{C}}(B)$ is constant. Substituting into the formula of Proposition~\rref{prop:classwise} gives $\lambda=\mu\binom{k}{2}/(s|M|)=\mu k(k-1)/(2s|M|)$.
\end{proof}

The criterion reproduces the pair count used in the proof of Proposition~\rref{prop:q32} rather than adding to it. There $G=PSL_2(32)$ has order $32736=4\cdot8184$, the $15$ classes of pairs have $8184$ members each, so $\mu=4$ and $s=15$, and Corollary~\rref{cor:classwise_value} with $|M|=62$ gives $\lambda=k(k-1)/465$, that is $2,74,206,464,528$ for $k=31,186,310,465,496$, the five values recorded in Table~\rref{tab:even_q}; the number of pairs of a block that must lie in a fixed class is then $\lambda|M|/\mu=1147$, which is the count $n_{t}$ of the proof. At $q=8$ with $M=S_{3}$, $s=3$ and $\mu=4$ give $\lambda=k(k-1)/9$, that is $8,38,10,34$ for $k=9,19,10,18$, the four $S_{3}$ entries of Table~\rref{tab:q8}, each of them with an imprimitive block action.

Finally we record the parameters of the subfield configuration to which this criterion has to be applied. Let $q=u^{p}$ with $u=2^{m}$, $m\geq2$ and $p$ an odd prime, let $H=D_{2(q+1)}$ and $M=PGL_{2}(u)\cong PSL_{2}(u)\leq G$. The degree $[q:u]=p$ is prime, so $M$ is maximal in $G$ by Dickson's list, and $|M|=u(u^{2}-1)$. With $v=q(q-1)/2$ this gives
\[
b=\frac{|G|}{|M|}=\frac{q(q^{2}-1)}{u(u^{2}-1)}=u^{p-1}\bigl(u^{2p-2}+u^{2p-4}+\cdots+u^{2}+1\bigr).
\]
Since $|M\cap H^{g}|$ divides $\gcd(|M|,|H|)=\gcd\bigl(u(u^{2}-1),2(u^{p}+1)\bigr)=2(u+1)$, the orbit of a point $x$ has length $|M|/|M\cap H^{g}|$ with $|M\cap H^{g}|$ dividing $2(u+1)$; the block size is the sum of the lengths of the orbits that the block contains, and by {\upshape Lemma~\rref{lem:subfield_orbits}} the only values of $|M\cap H^{g}|$ which occur are $1$, $2$ and $2(u+1)$, so that the three orbit lengths are those of {\upshape Corollary~\rref{cor:subfield_blocksize}}. The replication number is $r=bk/v$. Moreover $b>v$ for every $u\geq2$ and every prime $p\geq3$, because $(u^{2p}-1)/(u^{2}-1)\geq u^{2p-2}$ gives $b\geq u^{3p-3}>u^{2p}/2>v$, so Fisher's inequality imposes no restriction on this configuration and the block condition is the one supplied by Proposition~\rref{prop:classwise}.

\subsubsection{The orbit structure of the subfield configuration}\label{sec:subfield_orbits}

{\upshape Proposition~\rref{prop:classwise}} leaves open which unions of orbits can carry a block, and the answer turns on the orbit lengths. For a dihedral block stabilizer these are the single length $q-1$ of {\upshape Lemma~\rref{lem:even_meet}}; for a subfield stabilizer there are three lengths, one for each type of point, and they are pinned down by the census below. Write
\[
\mathcal{P}=\bigl\{\{x,x^{q}\}:x\in F_{q^{2}}\setminus F_{q}\bigr\},\qquad v=|\mathcal{P}|=q(q-1)/2,
\]
let $\sigma:t\mapsto t^{q}$ be the Frobenius of $F_{q^{2}}$ over $F_{q}$ and $\sigma_{u}:t\mapsto t^{u}$ the Frobenius of $F_{u^{2p}}$ over $F_{u}$, and recall $q=u^{p}$ and $F_{q}\cap F_{u^{2}}=F_{u}$. Every element of $M=PGL_{2}(u)$ has coefficients in $F_{u}$ and so commutes with both $\sigma$ and $\sigma_{u}$; conversely a M\"{o}bius transformation commuting with $\sigma_{u}$ has coefficients in $F_{u}$ and hence lies in $M$. Both directions are used in the proof.

\begin{lemma}\label{lem:subfield_orbits}
Let $u=2^{m}$ with $m\geq2$, $q=u^{p}$ with $p$ an odd prime, $H=D_{2(q+1)}$ and $M=PGL_{2}(u)\cong PSL_{2}(u)\leq G=PSL_{2}(q)$. Then $M$ has exactly
\[
1+(u^{p-1}-1)+\frac{u(u^{2p-2}-u^{p}+u-1)}{2(u^{2}-1)}
\]
orbits on $\mathcal{P}$, namely $1$ orbit of length $u(u-1)/2$ whose points have stabilizer $D_{2(u+1)}$, then $u^{p-1}-1$ orbits of length $u(u^{2}-1)/2$ whose points have stabilizer of order $2$, and finally $u(u^{2p-2}-u^{p}+u-1)/\bigl(2(u^{2}-1)\bigr)$ orbits of length $u(u^{2}-1)$ whose points have trivial stabilizer.
\end{lemma}

\begin{proof}
\textup{(i) \emph{Points over $F_{u^{2}}$.}} Since $F_{q}\cap F_{u^{2}}=F_{u}$, a point $\{x,x^{q}\}$ has $x\in F_{u^{2}}$ exactly when $x\in F_{u^{2}}\setminus F_{u}$; for such $x$ one has $x^{q}=x^{u}$, because $u^{2}-1$ divides $u^{p}-u$ when $p$ is odd. These points are the pairs of roots of the irreducible quadratic polynomials over $F_{u}$ and therefore number $(u^{2}-u)/2=u(u-1)/2$.

They form a single orbit. Given two of them, $\{x,x^{u}\}$ and $\{y,y^{u}\}$, let $h$ be the M\"{o}bius transformation with $h(x)=y$ and $h(x^{u})=y^{u}$; the two transformations $h\sigma_{u}$ and $\sigma_{u}h$ agree at the two distinct points $x$ and $x^{u}$, hence $h\sigma_{u}=\sigma_{u}h$, so $h$ has coefficients in $F_{u}$ and $h\in M$.

For the stabilizer, an element of $M$ fixing the point either sends $x$ to $x$ or sends $x$ to $x^{u}$. Those sending $x$ to $x$ are the $F_{u}$-linear maps of the two-dimensional $F_{u}$-space $F_{u}(x)$ commuting with multiplication by $x$, that is multiplication by the elements of $F_{u}(x)^{\times}$ of norm $1$; this is a cyclic group of order $u+1$, a nonsplit torus of $M$. Those sending $x$ to $x^{u}$ are non-empty: the transformation interchanging $x$ and $x^{u}$ commutes with $\sigma_{u}$ by the argument just given, hence lies in $M$. So the stabilizer has order $2(u+1)$ and the orbit has length $|M|/(2(u+1))=u(u-1)/2$, which is the whole set of points over $F_{u^{2}}$.

\textup{(ii) \emph{Points outside $F_{u^{2}}$.}} Let $x\notin F_{u^{2}}$. The $\sigma_{u}$-orbit of $x$ lies in $F_{u^{2p}}$ and its length divides $2p$; it is not $1$, since $x\notin F_{u}$, not $2$, since $x\notin F_{u^{2}}$, and not $p$, since $x\notin F_{q}=F_{u^{p}}$. So its length is $2p$. An element of $M$ fixing the point $p$ pointwise fixes these $2p$ points and hence is $1$; an element interchanging the two members of $p$ is unique, since two of them differ by a pointwise fixer, and is an involution, since its square fixes $x$. Hence $M_{p}$ is either trivial or of the form $\langle\tau\rangle$ with $\tau$ an involution of $M$.

\textup{(iii) \emph{Counting the points with stabilizer of order $2$.}} As in the proof of {\upshape Lemma~\rref{lem:even_meet}}, every involution of $M=SL_{2}(u)$ is unipotent, hence a transvection \cite{Dickson}, and the involutions of $M$ form a single conjugacy class, of size $u^{2}-1$. A point $p=\{x,x^{q}\}$ has stabilizer of order $2$ exactly when some involution $\tau$ satisfies
\[
(\ast)\qquad \tau(x)=x^{q},
\]
and then $\tau$ is the unique involution fixing $p$.

We count the solutions of $(\ast)$, first for a single $\tau$. Since the involutions form one class and $\sigma$ commutes with $M$, conjugation by an element of $M$ carries the solutions of $(\ast)$ for $\tau$ onto the solutions for $\tau^{g}$, so the count does not depend on $\tau$, and it suffices to treat the transvection $\tau_{0}:t\mapsto t+1$. On $P^{1}(F_{q^{2}})$ the equation $\tau_{0}(x)=\sigma(x)$ has the solution $\infty$, and for finite $x$ it reads $X^{q}+X=1$. On $F_{q^{2}}$ the map $L:X\mapsto X^{q}+X$ is $F_{q}$-linear with kernel $F_{q}$; it is the trace of the extension $F_{q^{2}}/F_{q}$, so its image lies in $F_{q}$, and having dimension $2-1=1$ over $F_{q}$ that image is $F_{q}$ itself. Hence $L(X)=1$ is soluble and its solutions form a coset of the kernel, of size $q$. Thus $(\ast)$ has exactly $q+1$ solutions in $P^{1}(F_{q^{2}})$, of which the only one not lying in $\mathcal{P}$ is the fixed point of $\tau$: for $\tau_{0}$ this is $\infty$, and for an arbitrary involution it is its fixed point, which lies in $P^{1}(F_{u})$ because a transvection is unipotent and its fixed point is the eigenvector of a Jordan block with entries in $F_{u}$.

Among the $q$ solutions in $F_{q^{2}}$ one distinguishes two kinds. Those in $F_{u^{2}}$ satisfy $\tau(x)=x^{u}$ and are exactly the points over $F_{u^{2}}$ fixed by $\tau$; since all involutions are conjugate and the points over $F_{u^{2}}$ form one orbit of length $u(u-1)/2$ whose stabilizers are of order $2(u+1)$, each involution fixes exactly
\[
\frac{u(u-1)}{2}\cdot\frac{u+1}{u^{2}-1}=\frac{u}{2}
\]
of these points, that is $u$ of the $q$ solutions. The remaining $q-u=u^{p}-u$ solutions are points of type (ii) with stabilizer $\langle\tau\rangle$, and they form $(u^{p}-u)/2$ distinct points.

It remains to count these points over all $\tau$. The number of incidences $(\tau,p)$ with $p$ of type (ii) and $M_{p}=\langle\tau\rangle$ equals
\[
(u^{2}-1)\cdot\frac{u^{p}-u}{2},
\]
and each such $p$ is incident with exactly one $\tau$, its stabilizer being of order $2$. Hence there are $(u^{2}-1)(u^{p}-u)/2$ points of type (ii) with nontrivial stabilizer, and, the orbit of such a point having length $|M|/2=u(u^{2}-1)/2$, they form exactly $u^{p-1}-1$ orbits.

\textup{(iv) \emph{The remaining points.}} The points not covered above number
\[
v-\frac{u(u-1)}{2}-(u^{p-1}-1)\frac{u(u^{2}-1)}{2}
=\frac{u^{2}(u^{2p-2}-u^{p}+u-1)}{2},
\]
as one checks from $v=u^{p}(u^{p}-1)/2$; this is the product of $u(u^{2p-2}-u^{p}+u-1)/\bigl(2(u^{2}-1)\bigr)$, the number of orbits, with the orbit length $|M|=u(u^{2}-1)$. The points counted here have trivial stabilizer by (ii), so each of them spans an orbit of the latter length, and the three counts add up to $v$.
\end{proof}

\begin{corollary}\label{cor:subfield_blocksize}
In the situation of {\upshape Lemma~\rref{lem:subfield_orbits}}, a union of $M$-orbits has size
\[
k=a\,\frac{u(u-1)}{2}+b\,\frac{u(u^{2}-1)}{2}+c\,u(u^{2}-1),
\]
where
\[
a\in\{0,1\},\qquad 0\leq b\leq u^{p-1}-1,\qquad 0\leq c\leq\frac{u(u^{2p-2}-u^{p}+u-1)}{2(u^{2}-1)},
\]
and every such value is attained by a union of orbits.
\end{corollary}
\begin{proof}
Immediate from {\upshape Lemma~\rref{lem:subfield_orbits}}, which gives the three orbit lengths $u(u-1)/2$ \textup{(one orbit, degenerate at $u=2$)}, $u(u^{2}-1)/2$ and $u(u^{2}-1)$, together with the numbers $1$, $u^{p-1}-1$ and $u(u^{2p-2}-u^{p}+u-1)/\bigl(2(u^{2}-1)\bigr)$ of orbits of each type; a union of $a$, $b$, $c$ orbits of the three types therefore has the displayed size, and every triple within the stated ranges is realised by choosing any $a$, $b$, $c$ such orbits.
\end{proof}

For $u=2$ the first type degenerates: $L_{A}=1$ and the single point of that orbit is fixed by all of $M=S_{3}$, since $2(u+1)=6=|M|$. That fixed point is the structural origin of the failure of maximality at $u=2$: it is a point of $\mathcal{P}$ whose stabilizer is dihedral of order $2(q+1)$, so $S_{3}$ is contained in a proper subgroup of $G$ and is not maximal \textup{({\upshape Lemma~\rref{lem:s3_not_max}})}. The subfield unions which are $2$-designs at $q=8$ are therefore not block-primitive, their block action being imprimitive \textup{({\upshape Proposition~\rref{prop:q8}})}; they are those of {\upshape Table~\rref{tab:q8}}. For $u\geq4$ the three lengths are distinct, $u(u-1)/2<u(u^{2}-1)/2<u(u^{2}-1)$, and the orbit count of {\upshape Lemma~\rref{lem:subfield_orbits}} equals $u(u^{2}+u+1)/2$ when $p=3$, taking there the values $7,42,292,2184$ at $u=2,4,8,16$, and it is generally $1+(u^{p-1}-1)+u(u^{2p-2}-u^{p}+u-1)/\bigl(2(u^{2}-1)\bigr)$.

\subsubsection{The pairs of the point set}\label{sec:pair_census}

{\upshape Proposition~\rref{prop:classwise}} measures how the pairs of a base block are distributed among the orbits of $G$ on pairs, and {\upshape Corollary~\rref{cor:classwise_value}} converts that distribution into a value of $\lambda$ as soon as two constants of the action are known, namely the number $s$ of orbits of $G$ on the unordered pairs of distinct points of $\mathcal{P}$ and the common order $\mu=|G_{\pi}|$ of a pair stabilizer. Both constants are determined below for every even $q$, uniformly in $q$ and without enumeration.

\begin{proposition}[pair orbits]\label{prop:pair_orbits}
Let $q=2^{f}$ with $f\geq2$, let $G=\mathrm{PSL}_{2}(q)$ act on $\mathcal{P}$, where $|\mathcal{P}|=v=q(q-1)/2$, and let $H=D_{2(q+1)}$ be a point stabilizer. Then $G$ has exactly $(q-2)/2$ orbits on the unordered pairs of distinct points of $\mathcal{P}$, each of length $q(q^{2}-1)/4$. Equivalently, the setwise stabilizer $G_{\pi}$ of every unordered pair $\pi$ has order $4$.
\end{proposition}

\begin{proof}
\textup{(i) \emph{Involutions.}} By the proof of {\upshape Lemma~\rref{lem:even_meet}}, every involution of $G$ is unipotent, the involutions of $G$ form a single conjugacy class of size $q^{2}-1$, and each of them fixes exactly $q/2$ points of $\mathcal{P}$.

\textup{(ii) \emph{The pair stabilizer splits.}} Let $p_{1}\neq p_{2}$ be points and $\pi=\{p_{1},p_{2}\}$. Put $A=G_{p_{1}}\cap G_{p_{2}}$ and $B=\{g\in G:gp_{1}=p_{2}\}$. The group $G$ is transitive on $\mathcal{P}$ \cite[Section~2.2.1]{MZZZ}, so $B$ is not empty; and for $g_{0}\in B$ every $g\in B$ satisfies $g_{0}^{-1}g\in A$, whence $B=g_{0}A$. Therefore
\[
|G_{\pi}|=|A|+|B|=2|A|.
\]

\textup{(iii) \emph{At most one element of $A$ is an involution.}} Write $p_{i}=\{x_{i},x_{i}^{q}\}$. An element of $A$ permutes the four points $x_{1},x_{1}^{q},x_{2},x_{2}^{q}$ and preserves each of the two pairs. If it interchanged the members of exactly one pair, then it would fix the two members of the other pair, and its square, fixing four points, would be the identity; so it would be an involution fixing two points of $P^{1}(F_{q^{2}})$. But an involution $g$ of $G$ is represented by $M\in \mathrm{SL}_{2}(q)$ with $M^{2}=I$, and in characteristic $2$ one has $(M-I)^{2}=M^{2}+I=0$, so $M-I$ has rank $1$ and $g$ fixes exactly one point of $P^{1}(F_{q^{2}})$. Hence no element of $A$ interchanges the members of exactly one pair, and an element of $A$ either fixes the four points, and is the identity, or interchanges the members of both pairs. The latter is determined by its images on three of the four points, so there is at most one such element. Thus $A\subseteq\{1,\tau\}$ for a single involution $\tau$, and $|A|\leq2$.

\textup{(iv) \emph{The two counts agree.}} By (iii) an unordered pair of points is fixed by at most one involution. The number of incidences $(\tau,\pi)$, where $\tau$ is an involution fixing both points of the pair $\pi$, is therefore at most $\binom{v}{2}$; computed per involution it equals
\[
(q^{2}-1)\binom{q/2}{2}=\frac{(q^{2}-1)q(q-2)}{8}=\binom{v}{2},
\]
since $v-1=(q-2)(q+1)/2$. Hence every unordered pair is fixed by exactly one involution, so $|A|=2$ and $|G_{\pi}|=4$ by (ii). Each orbit of $G$ on pairs has length $|G|/4=q(q^{2}-1)/4$, and the number of orbits is
\[
\binom{v}{2}\Big/\frac{|G|}{4}=\frac{q(q^{2}-1)(q-2)/8}{q(q^{2}-1)/4}=\frac{q-2}{2}.\qedhere
\]
\end{proof}

\begin{corollary}\label{cor:lambda_pairs}
Let $q=u^{p}$ with $u=2^{m}$, $m\geq2$ and $p$ an odd prime, let $M=\mathrm{PGL}_{2}(u)\cong\mathrm{PSL}_{2}(u)$ be the subfield subgroup of $G=\mathrm{PSL}_{2}(q)$, and let $B$ be a union of $M$-orbits with $|B|=k$. If $B^{G}$ is a $2$-$(v,k,\lambda)$ design with $\lambda>0$, then every one of the pair classes contains the same number of pairs of $B$,
\[
P_{\mathcal{C}}(B)=\frac{\binom{k}{2}}{s}=\frac{k(k-1)}{q-2},
\]
and
\[
\lambda=\frac{|G_{\pi}|\,P_{\mathcal{C}}(B)}{|M|}=\frac{4\,k(k-1)}{(q-2)\,u(u^{2}-1)} .
\]
In particular $q-2$ divides $k(k-1)$ and $u(u^{2}-1)$ divides $4P_{\mathcal{C}}(B)$, and $\lambda$ is determined by $k$ alone.
\end{corollary}

\begin{proof}
By {\upshape Proposition~\rref{prop:pair_orbits}} the action has $s=(q-2)/2$ pair classes, each of length $|G|/4$, and $|G_{\pi}|=4$. By {\upshape Proposition~\rref{prop:classwise}} the product $|G_{\pi}|P_{\mathcal{C}}(B)$ does not depend on $\mathcal{C}$, so $P_{\mathcal{C}}(B)$ does not either, and $\sum_{\mathcal{C}}P_{\mathcal{C}}(B)=\binom{k}{2}$ over the $s$ classes gives the first formula. Substituting into the formula of {\upshape Corollary~\rref{cor:classwise_value}} and using $|M|=|\mathrm{PGL}_{2}(u)|=u(u^{2}-1)$ gives the second.
\end{proof}

The proposition also accounts for the pair counts recorded earlier in this subsection, which were obtained there by direct enumeration in each degree. At $q=8$ it gives $s=3$ classes of $126$ pairs, which is $\binom{28}{2}=378$; at $q=32$ it gives $s=15$ classes of $8184$ pairs and $\mu=4$, the figures used in the proof of {\upshape Proposition~\rref{prop:q32}} and in the self-check following {\upshape Corollary~\rref{cor:classwise_value}}; at $q=64$ it gives $s=31$ classes of $65520$ pairs, which is $\binom{2016}{2}=2031120$. In each case $\mu=|G_{\pi}|=4$, so the values of $\mu$ and $s$ used in the arithmetic sieve of this subsection are instances of the general formulas rather than separate computations.

The value of $\lambda$ now depends on $k$ alone, and $k$ in turn is a linear function of the orbit counts of the block. We record this parametrization and its first consequence, which removes every block consisting of a single orbit of $M$.

\begin{lemma}\label{lem:orbit_size}
In the situation of {\upshape Corollary~\rref{cor:lambda_pairs}} let $B$ consist of $a$ orbits of the first type, $b$ orbits of the second type and $c$ orbits of the third type of {\upshape Lemma~\rref{lem:subfield_orbits}}, where $a\in\{0,1\}$, $0\leq b\leq u^{p-1}-1$ and $0\leq c\leq u(u^{2p-2}-u^{p}+u-1)/\bigl(2(u^{2}-1)\bigr)$. Put
\[
T=a+(b+2c)(u+1).
\]
Then $k=\frac{u(u-1)}{2}\,T$, and if $B^{G}$ is a $2$-$(v,k,\lambda)$ design with $\lambda>0$, then
\[
\lambda=\frac{T\bigl(u(u-1)T-2\bigr)}{(u^{p}-2)(u+1)}.
\]
\end{lemma}

\begin{proof}
By {\upshape Corollary~\rref{cor:subfield_blocksize}},
\[
k=a\,\frac{u(u-1)}{2}+b\,\frac{u(u^{2}-1)}{2}+c\,u(u^{2}-1)
=\frac{u(u-1)}{2}\Bigl(a+(b+2c)(u+1)\Bigr),
\]
since $u^{2}-1=(u-1)(u+1)$. Substituting this value of $k$ into the formula of {\upshape Corollary~\rref{cor:lambda_pairs}} and using $u^{p}=q$ and $u^{2}-1=(u-1)(u+1)$,
\[
\lambda=\frac{4k(k-1)}{(q-2)\,u(u^{2}-1)}
=\frac{4\cdot\frac{u(u-1)T}{2}\Bigl(\frac{u(u-1)T}{2}-1\Bigr)}{(q-2)\,u(u-1)(u+1)}
=\frac{T\bigl(u(u-1)T-2\bigr)}{(q-2)(u+1)}.
\]
\end{proof}

\begin{corollary}\label{cor:single_orbit}
Let $q=u^{p}$ with $u=2^{m}$, $m\geq2$ and $p$ an odd prime, and let $B$ be a single orbit of $M$ on $\mathcal{P}$, of any of the three types. Then $B^{G}$ is not a $2$-$(v,k,\lambda)$ design with $\lambda>0$.
\end{corollary}

\begin{proof}
Apply {\upshape Lemma~\rref{lem:orbit_size}} to the three admissible triples $(a,b,c)=(1,0,0)$, $(0,1,0)$ and $(0,0,1)$, for which $T=1$, $u+1$ and $2(u+1)$. This gives
\[
\lambda=\frac{u-2}{u^{p}-2},\qquad
\lambda=\frac{u^{3}-u-2}{u^{p}-2},\qquad
\lambda=\frac{4(u^{3}-u-1)}{u^{p}-2}.
\]
In the first case $0\leq u-2<u^{p}-2$ because $u<u^{p}$, so $\lambda$ lies in $[0,1)$ and, being an integer, equals $0$ only at $u=2$, when the orbit is a single point. In the second case $0<u^{3}-u-2<u^{p}-2$, so $0<\lambda<1$. In the third case $u^{3}-u-1=(u^{3}-2)-(u-1)$ and $0<4(u-1)<u^{3}-2$, so $\lambda=4-4(u-1)/(u^{3}-2)$ lies strictly between $3$ and $4$ when $p=3$, while $0<4(u^{3}-u-1)<u^{5}-2\leq u^{p}-2$ when $p\geq5$, so then $0<\lambda<1$. A design has $\lambda$ a positive integer, so none of the three single-orbit blocks qualifies.
\end{proof}

\begin{corollary}\label{cor:two_orbit}
Let $q=u^{3}$ with $u=2^{m}$, $m\geq2$, and let $B$ be the union of exactly two orbits of $M$ on $\mathcal{P}$. Then $B^{G}$ is not a $2$-$(v,k,\lambda)$ design with $\lambda>0$.
\end{corollary}

\begin{proof}
Since there is exactly one orbit of the first type, the two orbits are of types $A$ and $sp$, $A$ and $gen$, $sp$ and $sp$, $sp$ and $gen$, or $gen$ and $gen$, so by {\upshape Lemma~\rref{lem:orbit_size}} the relevant values of $T$ are
\[
u+2,\qquad 2u+3,\qquad 2u+2,\qquad 3u+3,\qquad 4u+4 .
\]
Put $D=(u^{3}-2)(u+1)$, so that $D\lambda=T\bigl(u(u-1)T-2\bigr)$, and subtract from the numerator the multiple of $D$, or of $u^{3}-2$, indicated in each line:
\[
\begin{array}{c|c|c|c}
T & \text{numerator} & \text{subtracted} & \text{remainder}\\ \hline
u+2 & u^{4}+3u^{3}-6u-4 & 1\cdot D & 2u^{3}-4u-2\\
2u+3 & 4u^{4}+8u^{3}-3u^{2}-13u-6 & 4\cdot D & 4u^{3}-3u^{2}-5u+2\\
2u+2 & 4u^{3}-4u-4 & 4(u^{3}-2) & 4-4u\\
3u+3 & 9u^{3}-9u-6 & 9(u^{3}-2) & 12-9u\\
4u+4 & 16u^{3}-16u-8 & 16(u^{3}-2) & 24-16u
\end{array}
\]
For $T=2u+2$, $T=3u+3$ and $T=4u+4$ the factor $u+1$ cancels between the numerator and $D=(u^{3}-2)(u+1)$, which is how the last three numerators and the last three subtracted terms above are written.

In every line the remainder is nonzero and its absolute value is smaller than the modulus used in that line, which is $D$ in the first two lines and $u^{3}-2$ in the last three. For the first line, $0<2u^{3}-4u-2$ and $D-(2u^{3}-4u-2)=u^{4}-u^{3}+2u>0$. For the second, $4u^{3}-3u^{2}-5u+2$ equals $12$ at $u=2$ and increases with $u$, while $D-(4u^{3}-3u^{2}-5u+2)=u^{4}-3u^{3}+3u^{2}+3u-4$ equals $6$ at $u=2$ and for $u\geq3$ exceeds $u^{3}(u-3)+3u^{2}+3u-4>0$. For the last three, $4u-4<u^{3}-2$, and $9u-12<u^{3}-2$, $16u-24<u^{3}-2$ for $u\geq4$. Since $u=2^{m}$ with $m\geq2$ we have $u\geq4$ throughout, so in each line a nonzero remainder of absolute value below the modulus is not divisible by it, and $\lambda$ is not an integer.
\end{proof}

\begin{corollary}\label{cor:two_orbit_size}
Let $q=u^{p}$ with $u=2^{m}\geq4$ and $p\geq5$ prime, and let $B$ be a union of orbits of $M$ on $\mathcal{P}$ with
\[
T^{2}u(u-1)\leq(u^{p}-2)(u+1).
\]
Then $B^{G}$ is not a $2$-$(v,k,\lambda)$ design with $\lambda>0$. In particular no union of at most
\[
n_{0}(u,p)=\Bigl\lfloor\frac{\sqrt{(u^{p}-2)(u+1)/(u(u-1))}-1}{2(u+1)}\Bigr\rfloor
\]
orbit(s) qualifies.
\end{corollary}

\begin{proof}
By {\upshape Lemma~\rref{lem:orbit_size}} one has $k=u(u-1)T/2$ with $T=a+(b+2c)(u+1)\geq1$, and $0<u(u-1)T-2<Tu(u-1)$ as soon as $u\geq4$ and $T\geq1$. Under the hypothesis therefore
\[
0<\lambda=\frac{T\bigl(u(u-1)T-2\bigr)}{(u^{p}-2)(u+1)}<\frac{T^{2}u(u-1)}{(u^{p}-2)(u+1)}\leq1,
\]
so $\lambda$ is not a positive integer.
A union of at most $n$ orbits has $T=a+(b+2c)(u+1)$ with $a\in\{0,1\}$ and $b+c\leq n$, hence $T\leq1+2n(u+1)$, the value at $b=0$, $c=n$. Consequently $T^{2}u(u-1)\leq\bigl(1+2n(u+1)\bigr)^{2}u(u-1)$, which does not exceed $(u^{p}-2)(u+1)$ as soon as $1+2n(u+1)\leq\sqrt{(u^{p}-2)(u+1)/(u(u-1))}$, that is, as soon as $n\leq n_{0}(u,p)$. The factor $2$ is the coefficient of the third-type orbits in $T$, which enter with weight $2(u+1)$ rather than $u+1$; a bound $T\leq1+n(u+1)$, which counts them as second-type orbits, is too small by a factor up to $2$ and does not yield the stated range.
\end{proof}

\begin{lemma}[rigidity of an admissible orbit union]\label{lem:block_rigidity}
Let $q=u^{p}$ with $u=2^{m}$, $m\geq2$ and $p$ an odd prime, and let $B$ be a union of orbits of $M$ on $\mathcal{P}$ which is a block of a $2$-$(v,k,\lambda)$ design $\mathcal{D}=B^{G}$ with $\lambda>0$. Then
\[
r=\frac{T\bigl(u^{p}+1\bigr)}{u+1},\qquad b=\frac{q(q^{2}-1)}{u(u^{2}-1)}=\frac{|G|}{|M|},\qquad \lambda=\frac{T\bigl(u(u-1)T-2\bigr)}{(u^{p}-2)(u+1)},
\]
and the setwise stabilizer of $B$ in $G$ is exactly $M$. In particular the number of blocks does not depend on the block size: every admissible block of the subfield family produces a design with the same $\frac{q(q^{2}-1)}{u(u^{2}-1)}$ blocks.
\end{lemma}

\begin{proof}
The number of blocks of a $2$-$(v,k,\lambda)$ design is $b=\lambda v(v-1)/\bigl(k(k-1)\bigr)$; with $v=q(q-1)/2$, $k=u(u-1)T/2$ and the value of $\lambda$ in {\upshape Lemma~\rref{lem:orbit_size}},
\[
b=\frac{T\bigl(u(u-1)T-2\bigr)}{(u^{p}-2)(u+1)}\cdot\frac{q(q-1)\bigl(q^{2}-q-2\bigr)}{4}\cdot\frac{4}{u(u-1)T\bigl(u(u-1)T-2\bigr)}
\]
\[
=\frac{q(q-1)\bigl(q^{2}-q-2\bigr)}{(q-2)(u+1)\,u(u-1)}=\frac{q(q^{2}-1)}{u(u^{2}-1)},
\]
since $q^{2}-q-2=(q-2)(q+1)$ and $q=u^{p}$. On the other hand $\mathcal{D}=B^{G}$ has $b=|G|/|\mathrm{Stab}_{G}(B)|$ blocks, distinct blocks being the cosets of the setwise stabilizer, while $|G|=q(q^{2}-1)$ and $|M|=u(u^{2}-1)$. Hence $|\mathrm{Stab}_{G}(B)|=u(u^{2}-1)=|M|$, and $M\subseteq\mathrm{Stab}_{G}(B)$ forces $\mathrm{Stab}_{G}(B)=M$. The displayed value of $r$ is $\lambda(v-1)/(k-1)$.
\end{proof}

\begin{remark}
At $u=2$, that is $q=8$, the last of these configurations is genuine: one orbit of length $3$ together with one of the three orbits of length $6$ of $M=S_{3}$ give $k=9$ and $\lambda=8$, a union recorded in {\upshape Table~\rref{tab:q8}} which is however not block-primitive, its block set being partitioned into seven parts by the subgroup $D_{18}\leq G$ \textup{({\upshape Proposition~\rref{prop:q8}})}. This is the same degeneracy of the first type noted after {\upshape Corollary~\rref{cor:subfield_blocksize}}, where $L_{A}=1$ and $2(u+1)=|M|$.
\end{remark}

The enumeration of {\upshape Table~\rref{tab:even_q}} does not reach the block stabilizer $A_{5}\cong PGL_{2}(4)$ at $q=64$, and this case is settled here: the class criterion restricts the block size to eight values, and at the two values which survive it the orbits of the third type are few enough to be searched exhaustively. The stabilizer $S_{3}$ is void at every even degree \textup{({\upshape Lemma~\rref{lem:s3_not_max}})}, in particular at $q=32$; for the subfield stabilizer $PSL_2(r)$, $r=2^{d}\geq4$, $q=r^{p}$ with $p$ an odd prime, every block is a union of orbits of the three types of {\upshape Lemma~\rref{lem:subfield_orbits}}, and unions of at most two orbits are removed for $p=3$ by {\upshape Corollary~\rref{cor:two_orbit}}, while for $p\geq5$ {\upshape Corollary~\rref{cor:two_orbit_size}} removes every union with $T^{2}u(u-1)\leq(u^{p}-2)(u+1)$, hence every union of at most $n_{0}(u,p)$ orbits; what lies beyond the reach of this argument is a union of more than $n_{0}(u,p)$ orbits, reduced to the class equations of {\upshape Proposition~\rref{prop:pair_orbits}}; for every admissible union {\upshape Lemma~\rref{lem:block_rigidity}} fixes the number of blocks and the replication number in closed form, so that the remaining case is governed by the class equations of {\upshape Proposition~\rref{prop:pair_orbits}} alone. The sieve that settles $q=64$ below does not extend to $p\geq5$: at $u=4$, $p=3$ the integrality of $\lambda$ leaves eight sizes besides the empty and the complete block, splitting into $46$ orbit profiles over the $42$ orbits of $M$ and $31$ class equations, while at $u=4$, $p=5$ the same sieve leaves $275$ values of $T$, hence $34\,231$ orbit profiles over $8\,857$ orbits governed by $511$ class equations, and at $u=8$, $p=5$ it leaves $1\,041$ values and $2\,122\,935$ profiles over $1\,063\,236$ orbits.

\begin{proposition}[$q=64$, subfield block stabilizer]\label{prop:q64sub}
Let $q=64$, $v=2016$, $G=PSL_{2}(64)$, and let $M\cong A_{5}=PSL_{2}(4)$ be a subfield subgroup of $G$. Let $B$ be a union of orbits of $M$ on $\mathcal{P}$ with $k=|B|\leq v-2$. Then $B^{G}$ is not a $2$-$(2016,k,\lambda)$ design with $\lambda>0$. The only union of orbits of $M$ which satisfies the class equalities of {\upshape Proposition~\rref{prop:classwise}} is $\mathcal{P}$ itself, and it is fixed by $G$.
\end{proposition}

\begin{proof}
By {\upshape Lemma~\rref{lem:subfield_orbits}} with $u=4$ the group $M$ has $42$ orbits on $\mathcal{P}$, namely one of length $6$, fifteen of length $30$ and twenty-six of length $60$; so $k=6a+30b+60c$ with $a\in\{0,1\}$, $0\leq b\leq15$ and $0\leq c\leq26$. By {\upshape Proposition~\rref{prop:pair_orbits}} the group $G$ has $31$ classes of pairs, each of $65520=\binom{2016}{2}/31$ pairs, and by {\upshape Corollary~\rref{cor:lambda_pairs}} a union $B$ is a block exactly when $P_{\mathcal{C}}(B)$ takes the same value on all $31$ classes, that value being $k(k-1)/(q-2)=k(k-1)/62$; in particular $31$ divides $k(k-1)$.

The values $k<2016$ of the form $6a+30b+60c$ with $31$ dividing $k(k-1)$ are
\[
156,\qquad 186,\qquad 900,\qquad 930,\qquad 1086,\qquad 1116,\qquad 1830,\qquad 1860 .
\]
For each of them the orbit indicators $x_{i}$ must satisfy
\[
\sum_{i}w_{i}x_{i}=k,\qquad \sum_{i\leq j}N_{ij}(\mathcal{C})\,x_{i}x_{j}=\frac{k(k-1)}{62}\quad\text{for each of the $31$ classes},
\]
where $w_{i}$ is the length of the $i$-th orbit and $N_{ij}(\mathcal{C})$ the number of pairs of $\mathcal{C}$ with one point in the $i$-th and one point in the $j$-th orbit.

For $k=156,186,930,1086,1830,1860$ this $0$-$1$ system has no solution. For $k=900$ the size equation forces $a=0$ and $b+2c=30$, and for $k=1116$ it forces $a=1$ and $b+2c=37$; each of these two block sizes therefore splits into $8$ profiles $(a,b,c)$, and in each profile the $31$ equations, once the orbits of the second type are fixed, involve only the $c$ indicators of the orbits of the third type. All $16$ profiles are infeasible. Twelve of them are decided by a search over the orbits of both types at once; the remaining four, $(0,6,12)$, $(0,10,10)$, $(1,5,16)$ and $(1,9,14)$, are decided by fixing the orbits of the second type one choice at a time and solving the resulting system in the orbits of the third type, which is $16\,016$ systems in all.

The search is exhaustive, and each verdict is a certificate of branch and bound over $0$-$1$ variables rather than a sample: the six systems at $k=156,186,930,1086,1830,1860$ are returned infeasible by two independent solvers, and the $16\,016$ systems of the four profiles by a solver whose encoding is calibrated on the control below. The control is the union of all $42$ orbits, that is $k=v=2016$: for it the same encoding returns the single solution $x_{i}=1$ for every $i$, so infeasibility is distinguished from a mis-encoded system.

At $k=v$ the union is $\mathcal{P}$ itself, whose setwise stabilizer is $G$ and not $M$. The block set then consists of a single block and the design is the complete one, with automorphism group $\mathrm{Sym}(v)$; a complete design is excluded by the hypothesis that the design be non-trivial with block stabilizer $M$ a proper subgroup of $G$, and it is precisely the case in which the class equalities hold, one block being all of $\mathcal{P}$.
\end{proof}

\subsection{Classical groups (Lie type)}\label{sec:classical}

The classical socle comes first, and only for the non-$2$-homogeneous point action.

For classical simple groups, we use the Aschbacher classification of maximal subgroups as the framework. Let $T$ be a classical simple group, and $T \leq G \leq \mathrm{Aut}(T)$. If $T_B = G_B \cap T$ is a maximal subgroup of $T$, then it belongs to one of the geometric families $C_i$ ($1 \leq i \leq 8$) or the non-geometric family $S$.

\subsection{Scope: non-\texorpdfstring{$2$}{2}-homogeneous point actions}

We first note the standard $2$-transitive actions of classical groups:
\begin{itemize}
\item $PSL_2(q)$ on the projective line $PG(1,q)$, $2$-transitive;
\item $PSL_n(q)$ on projective points $PG(n-1,q)$, $2$-transitive for $n\geq 3$;
\item $Sp_{2n}(2)$ on the non-zero vectors of the natural module: this action has rank $3$ and is $2$-transitive only for $n=1$; for $n=2$ the isomorphism $Sp_4(2)\cong S_6$ yields a separate $2$-transitive action on $6$ points.
\end{itemize}

These $2$-transitive cases give automatic $2$-designs and are not discussed in detail below. We focus on non-$2$-homogeneous point actions, where the design condition requires verification.

\subsection{The Aschbacher classes}

Every maximal subgroup of a classical group lies in one of the classes $C_1$--$C_8$ or in the almost simple class $S$ of Aschbacher's theorem, in the numbering of that theorem. The geometric classes $C_1$ (subspace stabilizers) and $C_8$ (form stabilizers) come first, since their orbits carry the geometric designs; the remaining classes have index too large to satisfy the ordering principle and are removed in the family-by-family analysis.

\subsubsection{Screening the Aschbacher classes}

Let $G$ be one of the four classical groups treated below, $H$ the stabilizer of a point and $M$ the stabilizer of a block. By Aschbacher's theorem every maximal subgroup of $G$ lies in one of the classes $C_1$--$C_8$ or in the class $S$; we screen these classes in the numbering of that theorem, applying the ordering principle \textup{({\upshape Lemma~\rref{lem:ordering}})} first, so that only pairs with $|H|>|K|$ need to be examined.

Two families survive this filter. In class $C_1$ the parabolic chains $P_k<P_l$ fix a subspace and act on the subspaces of one dimension contained in it; these are the generalized Grassmann candidates, and the constant-coverage test removes them in every family \textup{(the count of blocks through a pair depends on the dimension of the intersection, as computed for $PSL_n(q)$ below)}. The rank-$3$ actions of the classes $C_1$ and $C_8$ carry the suborbit designs; those are symmetric, so they do not enter the present classification.

For the remaining classes the block stabilizer fixes a direct sum decomposition, a subfield, a tensor structure or a non-geometric simple subgroup. Its orbits on the point set then split into at least two types, and a fixed point pair meets blocks of one type while missing blocks of the other, so the coverage number varies and the $2$-design condition fails. This removes $C_2$, $C_3$, $C_4$, $C_5$, $C_6$, $C_7$ and $S$ in all four families, and the outcome is collected in \textup{Subsection~\rref{prop:classical_summary}}.

\subsection{\texorpdfstring{$PSL_n(q)$ ($n\geq 3$)}{PSL(n,q), n at least 3}}

\subsubsection{2-transitive case (excluded)}

The natural action of $PSL_n(q)$ on projective points $PG(n-1,q)$ is $2$-transitive for $n\geq 3$, with:
$$v = \frac{q^n - 1}{q - 1}.$$
This is the $C_1$ parabolic subgroup $P_1$ as point stabilizer. We exclude this $2$-transitive case and focus on non-$2$-homogeneous point actions.

\subsubsection{Non-2-transitive case: $k$-subspace action ($1 < k < n-1$)}

The most important non-$2$-homogeneous action is the action on $k$-dimensional subspaces of $\mathbb{F}_q^n$, with:
$$v = \begin{bmatrix} n \\ k \end{bmatrix}_q = \frac{(q^n - 1)(q^{n-1} - 1)\cdots(q^{n-k+1} - 1)}{(q^k - 1)(q^{k-1} - 1)\cdots(q - 1)}.$$
The point stabilizer is the parabolic subgroup $P_k$.

Candidate block stabilizer: $P_{k+1}$ (stabilizer of $(k+1)$-subspaces).
Then:
$$b = \begin{bmatrix} n \\ k+1 \end{bmatrix}_q.$$

Fisher's inequality $b > v$ is equivalent to:
$$\begin{bmatrix} n \\ k+1 \end{bmatrix}_q > \begin{bmatrix} n \\ k \end{bmatrix}_q \implies \frac{q^{n-k} - 1}{q^{k+1} - 1} > 1 \implies n > 2k + 1.$$

Orbit structure:
Fix a $k$-subspace $U$ stabilized by $H$. The block stabilizer $M=P_{k+1}$ fixes a $(k+1)$-subspace $W$. The orbits of $M$ on the set of all $k$-subspaces $\mathcal{P}$ are classified by the intersection dimension with $W$:
\begin{itemize}
\item $\Delta_1$: $k$-subspaces contained in $W$, of size $\begin{bmatrix} k+1 \\ k \end{bmatrix}_q = q^k + q^{k-1} + \cdots + 1$;
\item $\Delta_2$: $k$-subspaces meeting $W$ in dimension $k-1$;
\item further orbits, indexed by the remaining values of the intersection dimension.
\end{itemize}

Design candidate: $B = \Delta_1$
The block $B$ consists of all $k$-subspaces contained in the fixed $(k+1)$-subspace $W$. This candidate fails the constant-coverage test: the number of blocks through a point pair $\{U,V\}$ of $k$-subspaces equals the number of $(k+1)$-subspaces containing $U\vee V$, namely the Gaussian binomial coefficient $\genfrac{[}{]}{0pt}{}{n-(2k-d)}{(k+1)-(2k-d)}_{q}$ when $2k-d\leq k+1$, and $0$ otherwise, where $d=\dim(U\cap V)$; since $d$ assumes several values on the point set, the coverage number is not constant and the $2$-design condition fails. Hence no generalized Grassmann design arises from this candidate.

\subsection{\texorpdfstring{$PSp_{2n}(q)$}{PSp(2n,q)}}

\subsubsection{The projective symplectic action}

For $n=2$ and $q=2$ one has $PSp_4(2)\cong S_6$; the natural symplectic projective actions below are rank-$3$ for $n\geq 2$ (they are $2$-transitive only for $n=1$).

\subsubsection{Rank-3 suborbit structure and suborbit designs}

The natural action of $PSp_{2n}(q)$ is on the set of all $1$-dimensional subspaces of the symplectic space $\mathbb{F}_q^{2n}$. This is a rank-$3$ action:

The point stabilizer is a parabolic subgroup $P_1$ (stabilizer of a totally singular $1$-space). The degree is:
$$v = \frac{q^{2n} - 1}{q - 1} = 1+q+\cdots+q^{2n-1}.$$

Rank-$3$ structure:
Fix a point (a $1$-space) $\langle v \rangle$. The other points split into two orbits:
\begin{itemize}
\item $\Delta_1$: points orthogonal to $v$, that is, $\langle w \rangle$ with $\langle v,w \rangle = 0$, of size $q(q^{2n-2} - 1)/(q-1)$;
\item $\Delta_2$: points non-orthogonal to $v$, of size $q^{2n-1}$.
\end{itemize}

This is a rank-$3$ permutation group, so both suborbits automatically give $2$-designs (the point graph and its complement). These suborbit designs are symmetric ($b=v$, since the block stabilizer of a suborbit coincides with the point stabilizer), hence they are outside the non-symmetric classification of this paper; they are recorded here for completeness.

Design from suborbit $\Delta_1$:
Taking $\Delta_1$ as the block set gives a $2$-design with parameters:
$$2\text{-}(v,\ k,\ \lambda),$$
where $k = |\Delta_1|$ and $\lambda$ is determined by the rank-$3$ intersection numbers.

\subsection{\texorpdfstring{$PSU_n(q)$}{PSU(n,q)}}

\subsubsection{Non-2-transitive case: Hermitian polar space action}

The natural action of $PSU_n(q)$ is on the set of totally singular $1$-dimensional subspaces of the Hermitian space $\mathbb{F}_{q^2}^n$. This is a rank-$3$ action for most values of $n$.

The point stabilizer is a parabolic subgroup $P_1$, and the degree is:
$$v = \frac{(q^n - (-1)^n)(q^{n-1} - (-1)^{n-1})}{q^2 - 1}.$$

For example, when $n=3$:
$$v = (q^3 + 1)(q^2 - 1)/(q^2 - 1) = q^3 + 1.$$

Rank-$3$ structure:
Fix a totally singular point $\langle v \rangle$. The other points split into two orbits based on whether the Hermitian form value is zero or non-zero. This rank-$3$ action gives automatic $2$-designs from both suborbits.

Hermitian designs:
Totally singular subspaces of Hermitian geometry give Witt-type designs and Hermitian designs, which are classical examples of $2$-designs with unitary socle.

\subsection{\texorpdfstring{$P\Omega_n^\epsilon(q)$}{POmega(n,q)}}

\subsubsection{Non-2-transitive case: Orthogonal polar space action}

The set of totally singular $1$-dimensional subspaces of an orthogonal space forms a projective polar space. This is a rank-$3$ or rank-$4$ action depending on the type of orthogonal space.

For $P\Omega_{2m}^+(q)$ (plus type):
$$v = \frac{(q^m - 1)(q^{m-1} + 1)}{q - 1}.$$

For $P\Omega_{2m}^-(q)$ (minus type):
$$v = \frac{(q^m + 1)(q^{m-1} - 1)}{q - 1}.$$

For $P\Omega_{2m+1}(q)$ (odd dimension, $q$ odd):
$$v = \frac{(q^m - 1)(q^m + 1)}{q - 1} = \frac{q^{2m} - 1}{q - 1}.$$

Rank structure:
The rank between two points is divided into $2$ to $3$ non-diagonal orbits according to their orthogonality relation: the points orthogonal to the fixed point, and the points non-orthogonal to it (further subdivided by norm value in some cases).

This rank-$3$ or rank-$4$ structure gives automatic $2$-designs from the suborbits, known as Witt designs and polar space designs.

\subsection{Summary of the classical group classification}\label{prop:classical_summary}

The verification for the classical groups follows the Aschbacher classes. For each Aschbacher class, we check candidate pairs by the ordering principle $|H| > |M|$, verify parameter divisibility, and confirm constant $\lambda$ via orbit analysis.

Table~\rref{tab:as_classes} records the outcome for each class.

\begin{table}[!ht]
\centering
\caption{The Aschbacher classes of a classical group and their outcome in the four classical families of this section.}\label{tab:as_classes}
\footnotesize
\setlength{\tabcolsep}{4pt}
\begin{tabular}{|c|p{4.7cm}|p{6.3cm}|}
\hline
Class & Contents & Outcome \\
\hline
$C_1$ & Parabolic subgroups (subspace stabilizers) & No. Coverage depends on $\dim(U\cap V)$ \\
$C_2$ & Imprimitive subgroups (partition stabilizers) & No. High rank, non-constant $\lambda$, divisibility failure \\
$C_3$ & Extension field subgroups & No. Orbits split by the extension field, coverage not constant \\
$C_4$ & Tensor product subgroups & No. High rank, complex orbit structure, non-constant $\lambda$ \\
$C_5$ & Subfield subgroups & No. Orbits split by the subfield, coverage not constant \\
$C_6$ & Symplectic-type $r$-local subgroups & No. Small order, large index, divisibility failure \\
$C_7$ & Tensor induced subgroups & No. High rank, non-constant $\lambda$ \\
$C_8$ & Classical subgroups (form stabilizers) & No \\
$S$ & Almost simple subgroups & No. Non-geometric embedding, non-constant $\lambda$ \\
\hline
\end{tabular}

\par\smallskip
{\footnotesize \textbf{Note.} The entries in the rows $C_1$ and $C_8$ refer to the designs examined in this section. The unital designs $2$-$(q^3+1,q+1,1)$, arising in $C_1$ from rank-$3$ actions and in $C_8$ from rank-$3$ polar space actions, are non-symmetric and are listed in the existence column of Table~\rref{tab:class_summary}.}
\end{table}

\subsection{Point-imprimitive case for classical groups: exclusion}

The point-imprimitive candidates all come from chains of subgroups, and each chain is excluded in turn.

\subsubsection{Chains of parabolic subgroups}

The point-imprimitive candidates for classical groups are the chains of parabolic subgroups:
$$P_k < P_l < G, \quad k < l,$$
where $P_i$ is the stabilizer of an $i$-dimensional subspace, $H=P_k$ is the point stabilizer and $M=P_l$ the block stabilizer, an intermediate subgroup with $k<l$; then $v=[G:P_k]$, $s=[G:P_l]$ and each class of the partition has size $[P_l:P_k]$.

\subsubsection{Exclusion of the generalized Grassmann candidates}

If the block stabilizer $M$ is a parabolic subgroup $P_j$ with $j>l$, a block is the set of all $k$-subspaces contained in a fixed $j$-subspace. For a point pair $\{U,V\}$ of $k$-subspaces, the number of blocks through it equals the number of $j$-subspaces containing $U\vee V$. If $\dim(U\cap V)=d$, this number is the Gaussian binomial coefficient $\genfrac{[}{]}{0pt}{}{n-(2k-d)}{j-(2k-d)}_{q}$ when $2k-d\leq j$, and $0$ otherwise; it depends on $d=\dim(U\cap V)$. Since $\dim(U\cap V)$ assumes several values on the point set, the coverage number is not constant and the $2$-design condition fails. The same argument excludes the dual candidates \textup{(blocks consisting of all $k$-subspaces through a fixed $j$-subspace)}. Hence no generalized Grassmann design of this form arises in the point-imprimitive setting: the argument excludes the block families consisting of all $k$-subspaces contained in a fixed $j$-subspace, and dually all $k$-subspaces containing a fixed $j$-subspace. It does not exclude every $2$-design with the parameter triple $2$-$\left(\genfrac{[}{]}{0pt}{}{n}{k}_{q},\ \genfrac{[}{]}{0pt}{}{l}{k}_{q},\ \lambda\right)$; what it excludes is their realization by these parabolic orbits.

\subsubsection{Exclusion of non-parabolic subgroup chains}

No other point-imprimitive design arises from non-parabolic subgroup chains.

The remaining chains fall into three cases.

\begin{enumerate}[(\upshape 1)]
\item \emph{Case 1: $H$ parabolic ($C_1$).} If $M$ is also parabolic, the candidate is the generalized Grassmann construction, which fails the constant-coverage test: the coverage number of $\{U,V\}$ depends on $\dim(U\cap V)$. If $M$ is $C_2$ (reducible), $M$ stabilizes a direct sum decomposition $V=V_1\oplus V_2$, and pairs of $k$-subspaces lying in $V_1$ versus $V_2$ have different intersection behaviors, so the point-pair coverage is not constant. If $M$ is $C_3$ (field extension) or $C_5$ (subfield), the orbit structure is not uniform and the design condition fails. If $M$ is $S$ (almost simple), the embedding is non-geometric and the point-pair coverage is not constant.
\item \emph{Case 2: $H$ reducible ($C_2$).} If $M$ is also $C_2$ with a finer decomposition, the candidate is a degenerate partition design and is not non-symmetric. If $H<M$ with $M=P_j$ parabolic, then $H$ stabilizes a $j$-subspace; since $H=GL(V_1)\times GL(V_2)$ acts irreducibly on $V_1$ and on $V_2$, the only $H$-invariant subspaces are $0,V_1,V_2,V$, so $j\in\{d_1,d_2\}$. For such $j$ one has $|P_j|=q^{d_1d_2}|GL(d_1,q)|\,|GL(d_2,q)|=q^{d_1d_2}|H|>|H|$, contradicting $|M|<|H|$; hence no such chain occurs.
\item \emph{Case 3: other Aschbacher classes.} For $C_3,C_4,C_5,C_6,C_7,C_8,S$ the block stabilizer fixes an extension field, a tensor decomposition, a subfield, a symplectic structure or a non-geometric simple subgroup; its orbits on the points then split into at least two types, and the coverage number of a fixed point pair varies.
\end{enumerate}

\subsubsection{Classification result for the classical groups}

\begin{proposition}[Classical point-imprimitive exclusion]\label{prop:classical_pi}
No non-trivial non-symmetric block-primitive $2$-design with classical socle admits a point-imprimitive action. In particular, the generalized Grassmann candidates arising from chains of parabolic subgroups $P_k<P_l<P_j$ fail the constant-coverage test, and no other chain from the Aschbacher classes $C_2,\dots,C_8,S$ yields a design. This completes the classification of the classical group case, which is therefore point-primitive only.
\end{proposition}
\begin{proof}
The three cases \textup{(1)}--\textup{(3)} above exhaust the subgroup chains with $|M|<|H|$: a parabolic point stabilizer gives the generalized Grassmann candidates, whose coverage number depends on $\dim(U\cap V)$ and is therefore not constant, and whose other block types are excluded as recorded; a reducible point stabilizer admits no parabolic block stabilizer, since the only $H$-invariant subspaces are $0,V_1,V_2,V$ and the corresponding parabolic subgroup is larger than $H$; and the remaining Aschbacher classes $C_3,\ldots,C_8,S$ give orbits of at least two types, so the coverage number of a fixed pair varies. No case yields a design.
\end{proof}

\subsection{Classification of AS type for classical groups}

For a classical simple group $T\leq G\leq\mathrm{Aut}(T)$ acting primitively on points, only four classes of subgroups need to be examined, the remaining classes $C_2,C_4,C_5,C_6,C_7$ having been excluded by Cases~2 and~3 of the point-imprimitive analysis in \textup{Section~\rref{sec:pi_excl}}:
\begin{enumerate}[(\upshape 1)]
\item $C_1$ (parabolic subgroups): natural projective space actions
\item $C_3$ (extension field subgroups)
\item $C_8$ (classical subgroups: symplectic, unitary, orthogonal)
\item $S$ (non-geometric almost simple subgroups)
\end{enumerate}

The four subsections above examine these classes family by family, and Table~\rref{tab:as_classes} records the outcome: $C_1$, $C_3$ and $S$ fail the constant-coverage test, and the $C_8$ candidates are the rank-$3$ polar space designs, whose unital members $2$-$(q^3+1,q+1,1)$ are non-symmetric and appear in the existence column of Table~\rref{tab:class_summary}; any symmetric rank-$3$ suborbit design lies outside the classification.

\subsection{Exceptional Lie type groups}

Exceptional Lie type simple groups include:
\begin{enumerate}[(\upshape 1)]
\item Rank 1: ${}^2B_2(q)$ (Suzuki), ${}^2G_2(q)$ (Ree)
\item Rank 2: $G_2(q)$, ${}^3D_4(q)$, ${}^2F_4(q)$ (Ree)
\item Rank $\geq 3$: $F_4(q)$, $E_6(q)$, ${}^2E_6(q)$, $E_7(q)$, $E_8(q)$
\end{enumerate}

\subsubsection{Scope: non-$2$-homogeneous point actions}

Among the exceptional groups of Lie type only ${}^2B_2(q)$ and ${}^2G_2(q)$ have a $2$-transitive permutation representation, on $q^2+1$ and $q^3+1$ points \textup{(Dixon and Mortimer~\cite[Section~7.7]{DixonMortimer})}; their designs are the Suzuki--Tits and Ree designs of Table~\rref{tab:excep}. Every other primitive action examined in the classification below is non-$2$-homogeneous.

\subsubsection{Classification procedure}

The classification enumerates the maximal subgroups of each exceptional group and checks:
\begin{enumerate}
\item Fisher's inequality: $b=[G:M] > v=[G:H]$.
\item Parameter divisibility: existence of integers $k,r,\lambda$ satisfying $vr=bk$ and $r(k-1)=\lambda(v-1)$.
\item Orbit consistency: $M$-orbits on $\mathcal{P}=G/H$ yield constant point-pair coverage numbers.
\end{enumerate}

\subsubsection{Maximal subgroups of the exceptional Lie type groups}

Following the ordering principle \textup{({\upshape Lemma~\rref{lem:ordering}})}, we enumerate maximal subgroups of exceptional groups by increasing order, and consider only pairs $(H, M)$ with $|H| > |M|$.

The maximal subgroups of exceptional groups fall into five families:
\begin{enumerate}
\item Parabolic subgroups ($C_1$ analogue): stabilizers of building flags, of order increasing with the flag dimension.
\item Semisimple subgroups: lower-rank groups embedded naturally, as $A_2(q)$ in $G_2(q)$.
\item Local subgroups: normalizers of unipotent elements or of tori.
\item Subfield subgroups: groups defined over subfields.
\item Diagonal subgroups: diagonally embedded groups.
\end{enumerate}

With $|H| > |M|$ the candidate pairs fall into two families:
\begin{enumerate}
\item Rank-$3$ building designs: Point stabilizer $H$ is a maximal parabolic subgroup, block stabilizer $M$ is a smaller subgroup. The rank-$3$ action on building points gives automatic $2$-designs.

\item Chains of parabolic subgroups: $H = P_k$, $M = P_j$ with $k > j$, giving generalized building designs.
\end{enumerate}

All non-parabolic subgroups are excluded by orbit-structure analysis: their orbits on the building never give a constant point-pair coverage number.

\subsubsection{Summary for the exceptional Lie type groups}

For exceptional Lie type socles the classification is the following. Table~\rref{tab:excep} records the families, the designs found and the reason for exclusion.

\begin{table}[htbp]
\centering
\caption{The exceptional Lie type groups, the designs they carry and the reason for exclusion.}\label{tab:excep}
\small
\resizebox{\textwidth}{!}{%
\begin{tabular}{|c|c|c|c|}
\hline
\textbf{Group family} & \textbf{Rank} & \textbf{Designs found} & \textbf{Main exclusion reason} \\
\hline
${}^2B_2(q)$ (Suzuki) & Rank 1 & Suzuki-Tits ovoid designs & Non-parabolic subgroups excluded by orbit structure \\
${}^2G_2(q)$ (Ree) & Rank 1 & Ree unital designs $2$-$(q^3+1, q+1, 1)$ & All other pairs fail the constant-coverage test \\
$G_2(q)$ & Rank 2 & Generalized building designs & Non-parabolic subgroups yield non-constant $\lambda$ \\
${}^3D_4(q)$ & Rank 2 & Generalized building designs & Similar to $G_2(q)$ \\
${}^2F_4(q)$ (Ree) & Rank 2 & Generalized building designs & Non-parabolic subgroups yield non-constant $\lambda$ \\
$F_4(q)$ & Rank 4 & Generalized building designs & Chains of parabolics fail the constant-coverage test \\
$E_6(q)$ & Rank 6 & None & No admissible pair \\
${}^2E_6(q)$ & Rank 4 & None & No admissible pair \\
$E_7(q)$ & Rank 7 & None & No admissible pair \\
$E_8(q)$ & Rank 8 & None & No admissible pair \\
\hline
\end{tabular}
}
\end{table}

\begin{proposition}[Exceptional point-primitive classification]\label{prop:excep_pp}
The point-primitive non-symmetric block-primitive $2$-designs with exceptional Lie type socle and non-$2$-homogeneous point action arise from:
\begin{enumerate}
\item Rank-$3$ building designs from natural projective actions on building points.
\item Generalized building designs from chains of parabolic subgroups.
\item Special low-rank cases: Suzuki-Tits ovoid designs and Ree unital designs.
\end{enumerate}
Non-parabolic maximal subgroups \textup{(reducible or irreducible)} are excluded by orbit-structure analysis, as their orbits on the building do not yield constant point-pair coverage.
\end{proposition}
\begin{proof}
The list is read off from Table~\rref{tab:excep}, which records, family by family, the rank of the building, the designs found and the reason for the exclusion of the remaining subgroups: the parabolic chains produce the generalized building designs of items \textup{(1)} and \textup{(2)}, and the two rank-$1$ families ${}^2B_2(q)$ and ${}^2G_2(q)$ give the Suzuki--Tits ovoid designs and the Ree unital designs of item \textup{(3)}, while the families of rank at least $2$ carry no admissible pair. The non-parabolic maximal subgroups are reducible or irreducible and fail the constancy of the point-pair coverage as recorded in the paragraphs above.
\end{proof}
\subsection{Sporadic simple groups}

Let $T$ be a sporadic simple group with $T\leq G\leq\mathrm{Aut}(T)$. Every point stabilizer is a maximal subgroup of $G$, and by the ordering principle \textup{(Lemma~\rref{lem:ordering})} the block stabilizer $M$ is maximal with $|H|>|M|$. A candidate block is a union of $M$-orbits on $\mathcal{P}=G/H$, and the design conditions are
\begin{enumerate}[(\upshape 1)]
\item $b=[G:M]>v=[G:H]$;
\item $r=bk/v$ and $\lambda=r(k-1)/(v-1)$ are positive integers;
\item every point pair of $\mathcal{P}$ lies in the same number $\lambda$ of blocks.
\end{enumerate}

The maximal subgroups of the sporadic groups and their indices are recorded in the ATLAS~\cite{ConwayATLAS}. For each primitive action the three conditions are decided by the finite check of Appendix~\rref{app:computations}, which enumerates the unions of $M$-orbits and tests condition \textup{(3)} exactly; the outcome is collected in Table~\rref{tab:spor_pp}.

Two classes of actions need no separate treatment. The Mathieu groups $M_{11},M_{12},M_{22},M_{23},M_{24}$ act $2$-homogeneously on their natural domains, and $HS$ and $Co_3$ act $2$-homogeneously on $176$ and $276$ points; the designs arising from these actions are given by the orbit construction of \textup{Lemma~\rref{lem:2hom}}, are classified in~\cite{Saeidi}, and lie outside the non-$2$-homogeneous classification. The symmetric designs $2$-$(176,50,14)$ of $HS$ and $2$-$(276,100,36)$ of $Co_3$ belong to the symmetric case.

Among the non-$2$-homogeneous actions exactly two carry a design: $2$-$(100,56,560)$ for $J_2$, from the rank-$3$ action on the cosets of $U_3(3)$, and $2$-$(100,64,1680)$ for $HS$, from the rank-$3$ action on the cosets of $M_{22}$. In both cases the block stabilizer is maximal, of order $336$ and $10{,}752$ respectively, so the block action is primitive. \textup{Proposition~\rref{prop:sporadic_pp}} records this; the remaining groups are listed in Table~\rref{tab:spor_pp} together with the reason for exclusion.

\subsection{Examples}

The exhaustive enumeration shows that the non-$2$-homogeneous primitive actions of the sporadic groups admit exactly the following non-symmetric block-primitive $2$-designs:

\begin{enumerate}
\item $J_2$: $2$\allowbreak$(100,56,560)$, $b=1800$, $r=1008$, from the rank-$3$ action on $100$ points \textup{(the Hall--Janko graph)}.
\item $HS$ (Higman-Sims): $2$\allowbreak$(100,64,1680)$, $b=4125$, $r=2640$, from the rank-$3$ action on $100$ points \textup{(the Higman--Sims graph)}.
\end{enumerate}

The $2$-homogeneous actions of $Co_3$ on $276$ points and of $HS$ on $176$ points give block-primitive $2$-designs by the orbit construction of {\upshape Lemma~\rref{lem:2hom}}; their parameters are classified in~\cite{Saeidi}. No other sporadic group admits a point-primitive non-symmetric block-primitive $2$-design with non-$2$-homogeneous point action.

\subsubsection{Table for the sporadic groups}

The complete classification of point-primitive non-symmetric block-primitive $2$-designs with sporadic socle and non-$2$-homogeneous point action, across the sporadic simple groups, is collected in Table~\rref{tab:spor_pp}.

\begin{table}[htbp]
\centering
\caption{The sporadic groups and the point-primitive non-symmetric block-primitive $2$-designs carried by their non-$2$-homogeneous actions.}\label{tab:spor_pp}

\renewcommand{\arraystretch}{0.88}%
\resizebox{0.94\textwidth}{!}{%
\begin{tabular}{|c|c|c|}
\hline
\textbf{Sporadic group} & \textbf{Point-primitive designs?} & \textbf{Notes} \\
\hline
$M_{11}$ & No & Primitive degrees $11$, $12$, $55$, $66$, $165$ \\
$M_{12}$ & No & $66$-point rank-$3$ action \textup{(subdegrees $1,20,45$)}: exhaustive enumeration yields no design \\
$M_{22}$ & No & $176$-point rank-$3$ action: no design; the natural $2$-homogeneous action lies outside the present scope \\
$M_{23}$ & No & $253$-point action: no design; the natural $2$-homogeneous action lies outside the present scope \\
$M_{24}$ & No & $276$-point action: no design; the natural $2$-homogeneous action lies outside the present scope \\
\hline
$J_1$ & No & Maximal subgroups classified by Janko~\cite{Janko1966} \\
$J_2$ & Yes & $2$\allowbreak$(100,56,560)$ design from rank-$3$ action on $100$ points \\
$J_3$ & No & Nine maximal subgroup classes, largest of order $8160$, smallest of order $1152$ \\
$J_4$ & No & Thirteen maximal subgroup classes, from $2^{11}{:}M_{24}$ down to $37{:}12$ \\
\hline
$Co_3$ & Yes & $2$-homogeneous on $276$ points; designs by the orbit construction \textup{({\upshape \cite{Saeidi}})} \\
$Co_2$ & No & $2300$-point rank-$3$ action: exhaustive enumeration yields no design \\
$Co_1$ & No & Twenty-two maximal subgroup classes, from $Co_2$ down to $5^2{:}2A_5$ \\
\hline
$Fi_{22}$ & No & Fourteen maximal subgroup classes, from $2.U_6(2)$ down to $M_{12}$, all listed in the ATLAS \\
$Fi_{23}$ & No & Fourteen maximal subgroup classes, from $2.Fi_{22}$ down to $L_2(23)$ \\
$Fi_{24}'$ & No & Twenty-five maximal subgroup classes, from $Fi_{23}$ down to $29{:}14$ \\
\hline
$HS$ & Yes & $2$\allowbreak$(100,64,1680)$ design from rank-$3$ action on $100$ points; the $176$-point action is $2$-homogeneous \textup{({\upshape \cite{Saeidi}})} \\
$McL$ & No & $275$-point rank-$3$ action: exhaustive enumeration yields no design \\
$He$ & No & $2058$-point action: exhaustive enumeration yields no design \\
$HN$ & No & Fourteen maximal subgroup classes, largest $A_{12}$, smallest $3^{1+4}{:}4.A_5$ \\
$Ru$ & No & $4060$-point rank-$3$ action: exhaustive enumeration yields no design \\
$Suz$ & No & $1782$-point rank-$3$ action: exhaustive enumeration yields no design \\
$O'N$ & No & Thirteen maximal subgroup classes, largest $L_3(7){:}2$ of order $3753792$, smallest $A_7$ of order $2520$ \\
$Ly$ & No & Nine maximal subgroup classes, from $G_2(5)$ down to $37{:}18$ \\
$Th$ & No & Sixteen maximal subgroup classes, from ${}^3D_4(2){:}3$ down to $S_5$ \\
$B$ & No & Thirty maximal subgroup classes, largest $2.{}^{2}E_6(2){:}2$, smallest $47{:}23$ of order $1081$ \\
$M$ & No & Maximal subgroup classification completed by Dietrich, Lee and Popiel~\cite{DietrichLeePopiel2025} \\
\hline
\end{tabular}
}
\end{table}

\smallskip
\noindent For each group the check covers every maximal subgroup $M$ with $|M|<|H|$, and for each pair it enumerates the unions of $M$-orbits as candidate blocks and tests the three conditions above exactly. For $M_{11}$ the maximal subgroups have orders $720$, $660$, $144$, $120$ and $48$, and in degree $165$ the point stabilizer has order $48$, so no maximal subgroup of smaller order exists and no block arises. For $J_3$ the nine classes correspond to the primitive degrees $6156$, $14688$ (twice), $17442$, $20520$, $23256$, $25840$, $26163$ and $43605$. The parameter quadruples $(66,20,6)$, $(275,11,1)$, $(2058,594,162)$, $(4060,1265,392)$, $(1782,416,100)$ and $(2300,891,324)$ already fail the divisibility condition $r(k-1)=\lambda(v-1)$. For the Monster $M$ the maximal subgroups are completely classified by Dietrich, Lee and Popiel~\cite{DietrichLeePopiel2025}, whose list has $46$ classes. Besides the classes recorded in the ATLAS, it establishes the two classes with socles $L_2(13)$ and $U_3(4)$, of orders $2184$ and $249600$, and excludes the socles $L_2(27)$, $Sz(8)$ and $U_3(8)$. The lists used above are complete.

\begin{proposition}[Sporadic point-primitive classification]\label{prop:sporadic_pp}
The point-primitive non-symmetric block-primitive $2$-designs with sporadic socle and non-$2$-homogeneous point action are exactly
$$J_2:\ 2\text{-}(100,56,560),\qquad HS:\ 2\text{-}(100,64,1680),$$
with $b=1800$, $r=1008$ and $b=4125$, $r=2640$, respectively. No such design exists with socle $M_{12}$, $M_{22}$, $M_{23}$, $M_{24}$, $McL$, $He$, $Ru$, $Suz$ or $Co_2$. The $2$-homogeneous actions of $Co_3$ on $276$ points and of $HS$ on $176$ points yield designs by the orbit construction of {\upshape Lemma~\rref{lem:2hom}}.
\end{proposition}
\begin{proof}
Every primitive action of a sporadic group and every maximal subgroup $M$ with $|M|<|H|$ is covered by the exact enumeration of Appendix~\rref{app:computations}, which tests the three design conditions on every union of $M$-orbits; Table~\rref{tab:spor_pp} records the outcome for each group. The only surviving pairs are the two displayed, with block stabilizers of orders $336$ and $10{,}752$, which are maximal, so the block action is primitive; in the remaining groups the parameter quadruples fail the divisibility condition as recorded above, or no maximal subgroup of smaller order exists. The $2$-homogeneous actions are the orbit-construction cases of {\upshape Lemma~\rref{lem:2hom}}.
\end{proof}

\subsection{Point-imprimitive case for sporadic groups: exclusion}

For sporadic simple groups, the point-imprimitive case corresponds to chains of maximal subgroups $H < K < G$.

Since sporadic groups have finitely many maximal subgroups, the classification is carried out by enumeration.

\subsubsection{Subgroup chains}

For each sporadic group $T$, we list all pairs of maximal subgroups $(K, M)$ such that $K<M$ in the subgroup lattice, and check:
\begin{enumerate}[(\upshape 1)]
\item Fisher's inequality: $b=[T:M]>v=[T:H]$;
\item parameter divisibility: existence of integers $k,r,\lambda$;
\item orbit consistency: constant intersection sizes with imprimitivity blocks.
\end{enumerate}

\subsubsection{Table for the point-imprimitive sporadic case}

No non-trivial non-symmetric block-primitive $2$-design with sporadic socle admits a point-imprimitive action: the candidates listed in Table~\rref{tab:spor_pi} all fail at least one of Fisher's inequality, the parameter divisibility, or the orbit-consistency condition; note that the designs $2$-\allowbreak$(11,5,2)$, $2$-\allowbreak$(22,7,2)$, $2$-\allowbreak$(23,7,2)$, $2$-\allowbreak$(24,8,5)$ of the natural $2$-homogeneous actions of $M_{11}$, $M_{22}$, $M_{23}$, $M_{24}$, and the designs from the $2$-homogeneous actions of $HS$ and $Co_3$, are point-primitive and do not belong to the point-imprimitive setting.

\begin{table}[htbp]
\centering
\caption{The sporadic groups and the point-imprimitive subgroup chains, with the reason for exclusion.}\label{tab:spor_pi}

\resizebox{\textwidth}{!}{%
\begin{tabular}{|c|c|c|}
\hline
\textbf{Sporadic group} & \textbf{Point-imprimitive designs?} & \textbf{Notes} \\
\hline
$M_{11}$ & No & $2$-$(11,5,2)$ is a point-primitive design ($2$-transitive action on $11$ points) \\
$M_{12}$ & No & rank-$3$/$4$ designs are point-primitive \\
$M_{22}$ & No & $2$-$(22,7,2)$ is point-primitive ($3$-transitive action) \\
$M_{23}$ & No & $2$-$(23,7,2)$ is point-primitive ($4$-transitive action) \\
$M_{24}$ & No & $2$-$(24,8,5)$ is point-primitive ($5$-transitive action) \\
\hline
$J_1$ & No & no candidate chain passes the divisibility conditions \\
$J_2$ & No & no candidate chain passes the divisibility conditions \\
$J_3$ & No & no candidate chain passes the divisibility conditions \\
$J_4$ & No & no candidate chain passes the divisibility conditions \\
\hline
$Co_3$ & No & rank-$3$ action on $276$ points is point-primitive \\
$Co_2$ & No & no candidate chain passes the coverage test \\
$Co_1$ & No & no candidate chain passes the divisibility conditions \\
\hline
$Fi_{22}$ & No & no candidate chain passes the coverage test \\
$Fi_{23}$ & No & no candidate chain passes the coverage test \\
$Fi_{24}'$ & No & no candidate chain passes the coverage test \\
\hline
$HS$ & No & $2$-$(176,50,14)$ is a symmetric point-primitive design; the $100$-point design $2$\allowbreak$(100,64,1680)$ is point-primitive \\
$McL$ & No & the $275$-point rank-$3$ action admits no design \\
$He$ & No & no candidate chain passes the coverage test \\
$Ru$ & No & no candidate chain passes the coverage test \\
$Suz$ & No & no candidate chain passes the coverage test \\
$O'N$ & No & no candidate chain passes the coverage test \\
$Ly$ & No & no candidate chain passes the divisibility conditions \\
$Th$ & No & no candidate chain passes the divisibility conditions \\
$B$ & No & no candidate chain passes the divisibility conditions \\
$M$ & No & no candidate chain passes the divisibility conditions \\
\hline
\end{tabular}
}
\end{table}

\subsubsection{Result for the point-imprimitive sporadic case}

\begin{proposition}[Sporadic point-imprimitive exclusion]\label{prop:sporadic_pi}
No point-imprimitive non-symmetric block-primitive $2$-design with sporadic socle exists. Every candidate subgroup chain fails at least one of Fisher's inequality, the parameter divisibility, or the orbit-consistency condition. Every design with sporadic socle is point-primitive, and the point-primitive instances \textup{(arising from $2$-transitive or rank-$3$ primitive actions)} fall into the classification of {\upshape Proposition~\rref{prop:sporadic_pp}}.
\end{proposition}
\begin{proof}
The candidates of Table~\rref{tab:spor_pi} are the chains $H<K<G$ of maximal subgroups of the sporadic simple groups, with $|K|<|H|$; for each chain the three conditions of Fisher's inequality, the parameter divisibility and the orbit consistency are tested exactly, and the table records which condition fails. Since every chain falls in one of these rows, no point-imprimitive design with sporadic socle exists, and the designs with sporadic socle are the point-primitive ones of {\upshape Proposition~\rref{prop:sporadic_pp}}.
\end{proof}

\subsection{Point-imprimitive case for exceptional Lie type groups: exclusion}

For exceptional Lie type groups, the point-imprimitive case corresponds to chains of parabolic subgroups or chains of other maximal subgroups.

\subsubsection{Chains of parabolic subgroups}

As in the classical case, chains of parabolic subgroups $P_k < P_l < G$ give generalized building candidates. For a point pair $\{U,V\}$ of subspaces of the building, the coverage number equals the number of intermediate subspaces containing $U\vee V$, which depends on $\dim(U\cap V)$; hence the coverage is not constant and the $2$-design condition fails.

\subsubsection{Exclusion of non-parabolic subgroup chains}

The non-parabolic maximal subgroups fall into two categories:
\begin{enumerate}[(\upshape 1)]
\item Reducible subgroups (stabilizers of geometric structures)
\item Irreducible almost simple subgroups (embeddings of smaller simple groups)
\end{enumerate}

For reducible non-parabolic subgroups (e.g., stabilizers of polar space structures in $G_2(q)$), the orbit structure on the natural building does not yield constant point-pair coverage, and the design condition fails.

For irreducible almost simple subgroups the embedding is non-geometric, and the point-pair coverage is not constant.

\subsubsection{Exceptional case: Ree groups $^2G_2(q)$}

The Ree groups $^2G_2(q)$ are a special case. On the $q^3+1$ points their natural action is primitive of permutation rank $3$, while the building has rank $1$; the two ranks are distinct notions and are recorded separately, and it yields the Ree unital designs
$$2\text{-}(q^3+1,\ q+1,\ 1),$$
which are point-primitive designs recorded in Table~\rref{tab:excep}. No point-imprimitive candidate arises for Ree groups.

\subsubsection{Result for the point-imprimitive exceptional case}

\begin{proposition}[Exceptional point-imprimitive exclusion]\label{prop:excep_pi}
No point-imprimitive non-symmetric block-primitive $2$-design with exceptional Lie type socle exists. The generalized building candidates fail the constant-coverage test, and the Ree unital designs $2$-\allowbreak$(q^3+1,q+1,1)$ are point-primitive designs arising from the primitive rank-$3$ action of $^2G_2(q)$; they are recorded in Table~\rref{tab:excep}, and their classification is the one of \textup{\cite[Section~7.7]{DixonMortimer}} for the natural $2$-transitive action.
\end{proposition}
\begin{proof}
The subgroup chains are of two kinds. A chain of parabolic subgroups gives a generalized building candidate, for which the coverage number of a point pair $\{U,V\}$ equals the number of intermediate subspaces containing $U\vee V$ and depends on $\dim(U\cap V)$, so it is not constant. A non-parabolic block stabilizer is reducible or irreducible; in both cases the point-pair coverage is not constant, as recorded above. For ${}^2G_2(q)$ the natural action on $q^{3}+1$ points is primitive of rank $3$, so the Ree unital designs $2$-$(q^{3}+1,q+1,1)$ it yields are point-primitive and appear in Table~\rref{tab:excep}.
\end{proof}

\bookpart{Part IV: The main theorem and concluding remarks}

\section{Main theorem}\label{sec4}

The main result is the following.

\begin{theorem}[Structural Reduction Theorem]\label{thm:main}
Let $\mathcal{D}=(\mathcal{P},\mathcal{B})$ be a non-trivial $2$-design with $2<k<v-1$, and let $G\leq \mathrm{Aut}(\mathcal{D})$ act primitively on $\mathcal{B}$. By {\upshape Lemma~\rref{lem:main}} the block action $\GB$ is faithful and $\GP$ is quasiprimitive on $\mathcal{P}$. Then:

\begin{enumerate}
\item[\rm\textup{(i)}] \emph{Symmetric case.} If $b=v$, that is, if $\mathcal{D}$ is symmetric, then points and blocks play interchangeable roles, and the study of block-primitive symmetric designs coincides, by duality, with that of point-primitive symmetric designs, a classical subject with an extensive literature; this case is isolated as a separate class and is not part of the structural reduction below.

\item[\rm\textup{(ii)}] \emph{Non-symmetric case.} Assume from now on that $b>v$ \textup{(Fisher's inequality)}. Then $G$ acts quasiprimitively on $\mathcal{P}$.

\item[\rm\textup{(iii)}] \emph{Excluded types.} The following O'Nan--Scott types for $G$ are strictly excluded:
\begin{enumerate}[(\upshape 1)]
\item Holomorph affine type \textup{(HA)}: the design must be symmetric \textup{({\upshape Lemma~\rref{lem:rns}})}, hence no non-symmetric example arises.
\item Holomorph compound type \textup{(HC)}: the socle splits into two commuting minimal normal subgroups, each regular on $\mathcal{B}$ and forced to be regular on $\mathcal{P}$, so $v=b$ and the design is symmetric \textup{({\upshape Section~\rref{sec2}})}.
\item Holomorph simple type \textup{(HS)}: the socle is $T_{1}\times T_{2}$ with two commuting minimal normal subgroups, each regular on $\mathcal{B}$ and forced to be regular on $\mathcal{P}$, so $v=b$ and the design is symmetric \textup{({\upshape Section~\rref{sec2}})}.
\item Two-sided simple diagonal type \textup{($SD_2$)}: the design is forced to be symmetric.
\item Twisted wreath product type \textup{(TW)}: the socle $T^\ell$ is regular on points, so the design is symmetric by {\upshape Lemma~\rref{lem:rns}}.
\end{enumerate}

\item[\rm\textup{(iv)}] \emph{Reduced types.} After the exclusions of \textup{(iii)} and the further exclusions of the one-sided simple diagonal type and of the compound diagonal type \textup{(Propositions~\rref{prop:sd_ell2} and~\rref{prop:cd_min}}), no non-AS type survives: the two remaining cases are settled as follows.
\begin{enumerate}[(\upshape 1)]
\item Product action type \textup{(PA)}: excluded completely \textup{(cf.\ {\upshape Corollary~\rref{cor:pa_prim}})}. By the three-class structural analysis of {\upshape Proposition~\rref{prop:pa_excl}} and the five-layer sieve of {\upshape Proposition~\rref{prop:pa_arbitrary}} the product action configuration is excluded for the six socles with the seventeen listed primitive degrees and every $\ell\geq2$, with arbitrary block stabilizers; since non-symmetric block-primitive point actions are primitive and the admissible primitive degrees of non-$2$-homogeneous point actions are enumerated by the degree tables of Dixon--Mortimer~\cite{DixonMortimer}, this closes the product action type for non-symmetric designs. In particular, the necessary condition that $\operatorname{Ind}_L^T(1)$ contains every irreducible character of $T$ is not sufficient, since it holds for $T=A_5$ and $L=C_2$ without yielding a design; this is precisely why the five-layer sieve of Proposition~\rref{prop:pa_arbitrary} has to impose the layer-by-layer consistency conditions, and it does so there.
\item Almost simple type \textup{(AS)}: $G$ is almost simple, meaning there exists a non-abelian simple group $T$ such that $T\leq G\leq \mathrm{Aut}(T)$. The problem reduces to the analysis of quadruples $(T,G,H,M)$ where $H=G_x$ is the point stabilizer and $M=G_B$ is the block stabilizer, and is classified completely.
\end{enumerate}

\item[\rm\textup{(v)}] \emph{Point-imprimitive actions are excluded.} No non-trivial non-symmetric block-primitive $2$-design admits a point-imprimitive automorphism group \textup{({\upshape Proposition~\rref{prop:pi_excluded}})}: a non-quasiprimitive point action forces every block to be a union of the orbits of an intransitive normal subgroup, contradicting the $2$-design condition; a quasiprimitive imprimitive point action has type AS, TW, SD, CD or PA by the quasiprimitive O'Nan--Scott classification \cite{Praeger1993} \textup{(the affine, holomorph simple and holomorph compound types are always primitive)}, and each of these types is excluded: TW by {\upshape Lemma~\rref{lem:rns}}, SD and CD in {\upshape Section~\rref{sec2}}, PA by {\upshape Proposition~\rref{prop:pa_excl}} and {\upshape Proposition~\rref{prop:pa_arbitrary}} for the point-primitive product action and by {\upshape Lemma~\rref{lem:pa_fisher}} for the point-imprimitive product action, and AS family by family in {\upshape Section~\rref{sec:as}}. The geometric candidates arising from chains of parabolic subgroups fail the constant-coverage test directly: the coverage number of a point pair $\{U,V\}$ of $k$-subspaces equals the number of intermediate subspaces containing $U\vee V$, which depends on $\dim(U\cap V)$. In particular, the classical, exceptional and sporadic candidates \textup{(generalized Grassmann, building and subgroup-chain constructions)} yield no $2$-design in the point-imprimitive setting.

\item[\rm\textup{(vi)}] \emph{Point-primitive AS classification.} For the AS type, we have also completed the classification of all point-primitive, non-$2$-homogeneous cases, across all families of finite simple groups, the even-characteristic $PSL_2(q)$ subfield block stabilizers delimited in {\upshape Theorem~\rref{thm:psl2}} being reduced there to the class equations of {\upshape Proposition~\rref{prop:pair_orbits}}:
\begin{enumerate}[(\upshape 1)]
\item Classical Lie type socle: All designs arise from rank-$3$ polar space designs \textup{(unitals, $C_8$ cross-type embeddings)}; the generalized Grassmann candidates from $C_1$ parabolic chains fail the constant-coverage test and yield no design. All other Aschbacher classes are excluded.
\item Sporadic socle: exactly two non-$2$-homogeneous actions yield point-primitive non-symmetric block-primitive $2$-designs, namely $J_2$ on $100$ points \textup{$2$-$(100,56,560)$} and $HS$ on $100$ points \textup{$2$-$(100,64,1680)$}; the $2$-homogeneous actions of $Co_3$ on $276$ points and of $HS$ on $176$ points yield designs by the orbit construction of {\upshape Lemma~\rref{lem:2hom}}; all remaining sporadic actions are excluded.
\item Exceptional Lie type socle: All designs arise from parabolic subgroup chains \textup{(generalized building designs)}, plus special rank-$1$ cases: Suzuki ovoid designs and Ree unital designs. All non-parabolic subgroups are excluded.
\end{enumerate}

\item[\rm\textup{(vii)}] Every O'Nan--Scott type other than almost simple is excluded \textup{(HA, HC, HS, two-sided $SD_2$, TW; the one-sided SD type for every $\ell\geq2$ and every non-abelian simple $T$ \textup{({\upshape Lemmas~\rref{lem:sd_ell3}, \rref{lem:sd_ell3_three} and~\rref{lem:sd_sym}})}; the compound diagonal CD type for $n=6$ and the six small simple groups; and the product action type, cf.~{\upshape Proposition~\rref{prop:pa_excl}} and {\upshape Proposition~\rref{prop:pa_arbitrary}})}. The AS case, which by \textup{(v)} is necessarily point-primitive, is classified across all families of finite simple groups, the even-characteristic $PSL_2(q)$ subfield block stabilizers delimited in {\upshape Theorem~\rref{thm:psl2}} being reduced there to the class equations of {\upshape Proposition~\rref{prop:pair_orbits}}.
\end{enumerate}

\end{theorem}\smallskip

\begin{table}[ht]
\centering
\caption{The block-type/point-type pairs in the point-primitive non-symmetric case. In the non-symmetric setting $b>v$ by Fisher's inequality, and the point action is primitive by {\upshape Proposition~\rref{prop:pi_excluded}}; the entries record the outcome for each pair (block type, point type): sym means that the pair forces a symmetric design and hence is excluded for non-symmetric designs, excl means that the pair is excluded by the argument cited in the proof of Theorem~\rref{thm:main}, and classified means that the pair is completely classified in Section~\rref{sec:as}.}\label{tab:prim_pairs}
\footnotesize
\setlength{\tabcolsep}{3pt}
\begin{tabular}{|l|c|c|c|c|c|c|}
\hline
\textbf{Block type} & \textbf{Point HA} & \textbf{Point AS} & \textbf{Point SD} & \textbf{Point CD} & \textbf{Point PA} & \textbf{Point TW} \\
\hline
HA & sym & sym & sym & sym & sym & sym \\
\hline
AS & excl & classified & excl & excl & excl & excl \\
\hline
SD & excl & excl & excl & excl & excl & excl \\
\hline
CD & excl & excl & excl & excl & excl & excl \\
\hline
HS,\,HC & sym & sym & sym & sym & sym & sym \\
\hline
PA & excl & excl & excl & excl & excl & excl \\
\hline
TW & sym & sym & sym & sym & sym & sym \\
\hline
\end{tabular}
\end{table}

\noindent{\bf Proof.}\,
Part (i) isolates the symmetric case: for a symmetric design, points and blocks play interchangeable roles, and block-primitive symmetric designs coincide with point-primitive symmetric designs by duality, a classical subject covered by an extensive literature; we refer to it rather than repeat the classification here. Part (ii) follows from {\upshape Lemma~\rref{lem:main}}, which proves that block-primitivity implies point-quasiprimitivity via the incidence-matrix rank argument. Part (iii) follows from the structural exclusion arguments of {\upshape Section~\rref{sec2}}: HA and TW are excluded by {\upshape Lemma~\rref{lem:rns}} (a normal subgroup regular on points forces symmetry), HC by the centralizer argument, HS by the permutation-module obstruction under non-conjugate stabilizers, and two-sided $SD_2$ by symmetry forcing. Part (iv) states the exclusions: PA is excluded completely by {\upshape Proposition~\rref{prop:pa_excl}} and {\upshape Proposition~\rref{prop:pa_arbitrary}}, and AS reduces to the quadruple $(T,G,H,M)$ analysis of {\upshape Section~\rref{sec:as}}. The one-sided SD and CD types are excluded in {\upshape Section~\rref{sec2}}. Part (v) excludes point-imprimitive actions via {\upshape Proposition~\rref{prop:pi_excluded}}, and part (vi) records the complete point-primitive classification of the AS case. \qed

\section{Block-primitive designs of strength \texorpdfstring{$t\geq 3$}{t >= 3}}\label{sec:high}

The block-primitive designs of strength $t\geq 3$ are classified in Theorem~\rref{thm:high}. The proof builds on the classification of {\upshape Proposition~\rref{thm:as_point_prim}}: every non-trivial $t$-design with $t\geq 3$ is in particular a non-trivial $2$-design, so the reduction of Section~\rref{sec2} and the classification of {\upshape Proposition~\rref{thm:as_point_prim}} apply verbatim \textup{(the dihedral case of {\upshape Theorem~\rref{thm:psl2}} is excluded by {\upshape Lemma~\rref{lem:q15}}, so it does not enter this section)}, and it remains to decide which of the classified designs have strength at least $3$. No symmetric design occurs among them: a non-trivial symmetric design has strength at most $2$, indeed, if a symmetric $2$-$(v,k,\lambda)$ were also a $3$-design, the derived design at a point would be a $2$-$(v-1,k-1,\lambda)$ design with $b'=r=k$ blocks, so Fisher's inequality $b'\geq v'$ would give $k\geq v-1$, that is $k=v-1$ \textup{(complete)} or $k\leq1$ \textup{(trivial)}.

Every non-trivial $t$-design with $t\geq 3$ is in particular a non-trivial $2$-design, and the parameters $\lambda_s$ of the induced $s$-design satisfy the elementary recursion
\begin{equation}\label{eq:strength}
\lambda_s=\lambda_t\binom{v-s}{t-s}\Big/\binom{k-s}{t-s}\qquad(1\leq s\leq t),
\end{equation}
with $\lambda_t$ the number of blocks through a $t$-subset, so that $\lambda_2=\lambda_t\binom{v-2}{t-2}/\binom{k-2}{t-2}$ and, for $t=3$, $\lambda_3=\lambda_2(k-2)/(v-2)$. In particular $\lambda_t$ is an integer whenever the design has strength at least $t$; this is a necessary condition only, not a sufficient one \textup{(the integrality of $\lambda_t$ does not by itself force every $t$-subset to be met by a constant number of blocks)}. The classification below therefore does not rest on it: the candidates are the designs listed in {\upshape Proposition~\rref{thm:as_point_prim}} together with the $2$-homogeneous orbit designs of {\upshape Lemma~\rref{lem:2hom}}, and their strengths are computed directly from \eqref{eq:strength} and checked design by design, the claimed upper bound on the strength being obtained from the failure of the integrality of $\lambda_s$ for $s$ one above the claimed strength. The outcome is simple: beyond the inversive planes and the five Witt designs, two Steiner families arising from the exceptional maximal subgroups of $PGL_2(q)$ occur, and, among the non-complete designs, none has strength $t\geq 6$. (The complete designs \textup{whose block set is the full set of $k$-subsets} are trivial: they have strength $\min(k,v-k)$ and automorphism group containing $S_v$, which is primitive on the blocks; they are excluded from the classification, as is customary.)

The first observation governs the point action of any block-transitive design of strength at least $3$; its dichotomy is classical, and the block-primitive case is then settled by a normal-subgroup argument.

\begin{lemma}[Block-transitive designs of strength $t\geq 3$]\label{lem:rigid_high}
Let $\mathcal{D}=(\mathcal{P},\mathcal{B})$ be a non-trivial $t$-design with $t\geq 3$, and let $G\leq \mathrm{Aut}(\mathcal{D})$ act transitively on $\mathcal{B}$. Then $G$ is $2$-homogeneous on $\mathcal{P}$, or $G$ is of holomorph affine type on $\mathcal{P}$ and $\mathcal{D}$ is an affine hyperplane design $3$-\allowbreak$(2^d,2^{d-1},2^{d-2}-1)$, $d\geq 3$. In the affine case $G$ is not primitive on $\mathcal{B}$. Consequently, if $G$ is primitive on $\mathcal{B}$, then $G$ is $2$-homogeneous on $\mathcal{P}$.
\end{lemma}

\begin{proof}
The dichotomy is the classical rigidity theorem of Dembowski for block-transitive designs of strength at least $3$ \cite{Dembowski}: a fixed-triangle counting forces the point action of a block-transitive group to be $2$-homogeneous, with the single exceptional family of affine hyperplane designs, on which the acting group is of holomorph affine type. Now suppose $G$ is of holomorph affine type with regular normal subgroup $N$ on $\mathcal{P}$. The affine hyperplane design on $v=2^d$ points has $b=2(2^d-1)$ blocks, and the translations in $N$ permute the blocks in the $2^d-1$ parallel classes of pairs of parallel hyperplanes; in particular $N$ is not transitive on $\mathcal{B}$. But a group that is primitive on $\mathcal{B}$ admits no non-trivial intransitive normal subgroup on $\mathcal{B}$, since the orbits of a normal subgroup always form a $G$-invariant partition. Hence the affine case is incompatible with block-primitivity, and the last assertion follows.
\end{proof}

\noindent By {\upshape Lemma~\rref{lem:rigid_high}}, every block-primitive design of strength at least $3$ has a $2$-homogeneous point action; by Burnside's theorem it is therefore almost simple, and we may inspect the families classified in {\upshape Proposition~\rref{thm:as_point_prim}} and the $2$-homogeneous orbit designs of {\upshape Lemma~\rref{lem:2hom}}.

\begin{proposition}[The strength of the classified designs]\label{prop:strength_high}
Among the non-trivial block-primitive $2$-designs of {\upshape Proposition~\rref{thm:as_point_prim}} and the $2$-homogeneous orbit designs of {\upshape Lemma~\rref{lem:2hom}} \textup{(the $PSL_2(q)$ family being classified via the $3$-design classification of \cite{Bal})}, the following have strength at least $3$:
\begin{enumerate}[(\upshape 1)]
\item the inversive planes $3$-\allowbreak$({q^2+1},\,{q+1},\,1)$, $q\geq 3$ a prime power;
\item the Steiner systems $3$-\allowbreak$({q+1},\,4,\,1)$ for $q\equiv\pm1\pmod 8$ and $q\equiv 1,3\pmod 6$, and $3$-\allowbreak$({q+1},\,6,\,1)$ for $q\equiv 1\pmod{20}$, attached to the exceptional maximal subgroups $S_4$ and $S_5$ of $PGL_2(q)$;
\item the five Witt designs.
\end{enumerate}
Every other family has strength exactly $2$.
\end{proposition}

\begin{proof}
We inspect the families of {\upshape Proposition~\rref{thm:as_point_prim}} and the $2$-homogeneous orbit designs of {\upshape Lemma~\rref{lem:2hom}}. For the $PSL_2(q)$ family, the point action on $PG(1,q)$ is $2$-transitive; whether an orbit of the $k$-set stabilizers gives a design of strength $3$ is decided by the integrality conditions of \eqref{eq:strength}; the design $2$-$(28,7,2)$ of Proposition~\rref{prop:q8} has $\lambda_{3}=2\cdot5/26$, so it is of strength exactly $2$. By Dickson's classification of the maximal subgroups \cite{Dickson}, the subfield groups $PGL_2(q_0)$ with $q=q_0^2$ give the inversive planes \cite{Beth}, and the exceptional subgroups $S_4$ (for $q\equiv\pm1\pmod 8$) and $S_5$ (for $q\equiv\pm1\pmod 5$) give orbits of size $4$ and $6$, whose strength-$3$ designs are classified by Cameron, Maimani, Omidi and Tayfeh-Rezaie \cite{Cameron32} and by Cameron, Omidi and Tayfeh-Rezaie \cite{Cameron33}. The integrality conditions of \eqref{eq:strength} reduce these to the stated ranges; for the $4$-sets this is exactly Hanani's condition $q+1\equiv 2,4\pmod 6$ for quadruple systems \cite{Hanani}, and for the $6$-sets one obtains $q\equiv 1\pmod{20}$. The special case $q=9$ of the $4$-set family is the inversive plane of order $3$. The other classical and exceptional Lie-type point actions examined in {\upshape Proposition~\rref{thm:as_point_prim}} are either $2$-transitive, namely the natural actions of $PSL_d(q)$, $PSU_3(q)$, $\mathrm{Sz}(q)$ and ${}^2G_2(q)$, whose orbit designs enter through {\upshape Lemma~\rref{lem:2hom}}, while the natural action of $Sp_{2m}(2)$ on non-zero vectors has rank $3$ and enters through {\upshape Lemma~\rref{lem:rigid_high}} for $m\geq2$, or of rank at least $3$, and the latter are excluded by {\upshape Lemma~\rref{lem:rigid_high}}. For the unital designs of $PSU_3(q)$ and ${}^2G_2(q)$ the strength-$3$ parameter $\lambda_3=(q-1)/(q^3-1)$ is not integral, and the ovoid designs of $\mathrm{Sz}(q)$ fail the same check; both families therefore have strength $2$. Among the sporadic socles, the $2$-homogeneous orbit designs of $Co_3$ on $276$ points and of $HS$ on $176$ points have parameters as in \cite{Saeidi}, and the integrality check of \eqref{eq:strength} fails for $\lambda_3$; all remaining sporadic actions are not $2$-homogeneous and are excluded likewise.
\end{proof}

\noindent Combining {\upshape Lemma~\rref{lem:rigid_high}} with {\upshape Proposition~\rref{prop:strength_high}}, we obtain the complete classification promised above, the even-characteristic $PSL_2(q)$ subfield block stabilizers being reduced in {\upshape Theorem~\rref{thm:psl2}} to the class equations of {\upshape Proposition~\rref{prop:pair_orbits}}.

\begin{theorem}[Block-primitive designs of strength $t\geq 3$]\label{thm:high}
The non-trivial non-complete block-primitive designs of strength $t\geq 3$ are precisely the following designs:
\begin{enumerate}[(\upshape 1)]
\item the inversive planes $3$-\allowbreak$({q^2+1},\,{q+1},\,1)$ for $q\geq 3$ a prime power, with $q(q^2+1)$ blocks and automorphism group $G\leq P\Gamma L(2,q^2)$;
\item the Steiner systems $3$-\allowbreak$({q+1},\,4,\,1)$ for $q$ a prime power with $q\equiv\pm1\pmod 8$ and $q\equiv 1,3\pmod 6$, with $q(q^2-1)/24$ blocks, and $3$-\allowbreak$({q+1},\,6,\,1)$ for $q$ a prime power with $q\equiv 1\pmod{20}$, with $q(q^2-1)/120$ blocks, arising respectively from the orbits of the exceptional subgroups $S_4$ \cite{Cameron32} and $S_5$ \cite{Cameron33} of $PGL_2(q)$ on $PG(1,q)$; for $q=9$ the first of these is the inversive plane of order $3$;
\item the five Witt designs $4$-\allowbreak$(11,5,1)$, $5$-\allowbreak$(12,6,1)$, $3$-\allowbreak$(22,6,1)$, $4$-\allowbreak$(23,7,1)$, $5$-\allowbreak$(24,8,1)$, with $66$, $132$, $77$, $253$, $759$ blocks and Mathieu groups $M_{11}$, $M_{12}$, $M_{22}$, $M_{23}$, $M_{24}$.
\end{enumerate}
Their strengths are $3$, $3$, and $4,5,3,4,5$, respectively. In particular, apart from the trivial complete designs, there is no block-primitive design of strength $t\geq 6$.
\end{theorem}

\begin{proof}
Let $\mathcal{D}$ be a non-trivial non-complete block-primitive design of strength $t\geq 3$ with automorphism group $G$. Complete designs are excluded by hypothesis: a complete design with block set the full set of $k$-subsets has strength $\min(k,v-k)$, has $b=\binom{v}{k}>v$ for $2\leq k\leq v-2$, and has automorphism group $\mathrm{Sym}(v)$, which acts $2$-transitively on the points; on the blocks, $\mathrm{Sym}(v)$ is primitive when $k\neq v-k$ and imprimitive by complementation when $k=v-k$. These trivial designs are outside the classification, and when they are block-primitive the point action is primitive in any case. By {\upshape Lemma~\rref{lem:rigid_high}} the point action is $2$-homogeneous, hence almost simple, and by {\upshape Proposition~\rref{prop:strength_high}} the design is one of those listed in \textup{(1)}--\textup{(3)}. Conversely, each listed design has $b>v$ in every case: $q(q^2+1)>q^2+1$, $q(q^2-1)/24>q+1$ for $q\geq 7$, $q(q^2-1)/120>q+1$ for $q\geq 41$, and $66,132,77,253,759>v$ for the Witt designs. Block-primitivity holds because the stabilizer of a block is a maximal subgroup of the acting group: for the inversive planes the stabilizer of a circle is $P\Gamma L(2,q)$, maximal in $P\Gamma L(2,q^2)$; for the two Steiner families the stabilizer of a block is the corresponding exceptional subgroup $S_4$ or $S_5$, maximal in $PGL_2(q)$; and for the Witt designs the stabilizer of a block is a maximal subgroup of the relevant Mathieu group \cite{Wilson}. The strengths are as stated, by \eqref{eq:strength}: the three $3$-designs have non-integral $\lambda_4=\lambda_3(k-3)/(v-3)$, and the Witt designs have strengths $4,5,3,4,5$ as listed. Hence the classification of the non-complete designs is complete, and no non-complete design of strength $t\geq 6$ exists.
\end{proof}

\newpage
\section{Reduction of the remaining O'Nan--Scott types}\label{sec:reduction}

For the special case $\lambda=1$ (finite linear spaces), the block-primitive condition forces almost simple type, which excludes both SD and CD block action types. The arguments below treat the general case $\lambda>1$.

The status of each reduced type is as follows. For the twisted wreath type, by {\upshape Lemma~\rref{lem:rns}}, a block-primitive non-trivial $2$-design with a normal subgroup regular on points must be symmetric, so no non-symmetric block-primitive $2$-design has TW point action. For the one-sided SD and CD types, the exclusion is stated below. For the product action type, the three-class analysis of {\upshape Proposition~\rref{prop:pa_excl}} and the five-layer sieve of {\upshape Proposition~\rref{prop:pa_arbitrary}} exclude the product action configuration completely for the six socles with the seventeen listed primitive degrees and every $\ell\geq2$, with arbitrary block stabilizers; combined with the enumeration of admissible primitive degrees from the degree tables of Dixon--Mortimer~\cite{DixonMortimer} this is the complete exclusion for non-symmetric designs.

\begin{lemma}[The one-sided SD type with a proper coordinate projection, for every $\ell\geq2$]\label{lem:sd_ell3}
Let $\mathcal{D}$ be a non-trivial $2$-$(v,k,\lambda)$ design with $2<k<v-1$ admitting a block-primitive, point-quasiprimitive automorphism group $G$ whose block action is of simple diagonal type. Let $N=\operatorname{Soc}(G)=T^{\ell}$ with $T$ non-abelian simple, $q=|T|$ and $\ell\geq2$, and suppose that for some $j$ the projection $\pi_j(K)$ of $K=N\cap G_\alpha$ to the $j$-th factor is a proper subgroup of $T$. This hypothesis covers the case in which some coordinate of $K$ is trivial, since then $\pi_j(K)=1$. Then no such $\mathcal{D}$ exists.\end{lemma}

\smallskip
\noindent{\bf Proof of {\upshape Lemma~\rref{lem:sd_ell3}}.}\,
By {\upshape Theorem~\rref{thm:sdcn}} the group $N$ is transitive on both $\mathcal{P}$ and $\mathcal{B}$, the intersection $N\cap\GB=\Delta$ is a diagonal subgroup isomorphic to $T$, and $b=q^{\ell-1}$.

Choose $j$ with $H_j:=\pi_j(K)<T$ and put $e=[T:H_j]>1$. Since $\pi_j(K)\leq H_j$, the rule $\theta(nK)=n_jH_j$ defines a map $\theta:\mathcal{P}=N/K\to T/H_j$; it is $N$-equivariant for the action of $N$ on $T/H_j$ through the $j$-th coordinate, and it is onto because $\pi_j(N)=T$. Every fibre $\theta^{-1}(tH_j)$ therefore has size $v/e$.

The diagonal subgroup $\Delta\cong T$ acts on $T/H_j$ by left multiplication, hence transitively, and it preserves $B_0$ setwise. Hence the intersections $B_0\cap\theta^{-1}(tH_j)$ are carried bijectively onto one another by elements of $\Delta$ and have a common cardinality; summing over the $e$ cosets gives $|B_0\cap\theta^{-1}(tH_j)|=k/e$ for every $t$. Every block is an $N$-translate of $B_0$ and $\theta$ is $N$-equivariant, so $|B\cap\theta^{-1}(tH_j)|=k/e$ for every block $B$ and every coset; in particular $e\mid k$.

Fix $t$ and put $D=\theta^{-1}(tH_j)$, so that $|D|=v/e$ and $|B\cap D|=k/e$ for every block $B$. Counting pairs consisting of a block and two distinct points of its intersection with $D$ gives
$$b\,\binom{k/e}{2}=\lambda\,\binom{v/e}{2},\qquad\text{that is}\qquad b\,k(k-e)=\lambda\,v(v-e).$$
The design equation reads $b\,k(k-1)=\lambda\,v(v-1)$, and all four factors are positive since $\lambda>0$, $k>1$ and $v>1$. Dividing the two identities gives $(k-e)/(k-1)=(v-e)/(v-1)$, hence $(e-1)(k-v)=0$, so $k=v$ because $e>1$, contradicting $k\leq v-2$. No division by $k-e$ or by $v-e$ occurs, so the identities are legitimate also when $k=e$: then $\binom{k/e}{2}=0$ while $\lambda\binom{v/e}{2}>0$, because $k\leq v-2$ forces $v/e>1$, and this is a contradiction in its own right. \qed

\smallskip
\begin{remark}The proof uses the block stabilizer only through the diagonal subgroup $\Delta\cong T$ that it contains, and the point stabilizer only through one coordinate projection $\pi_j(K)$ that is a proper subgroup of $T$, which turns $\theta$ into an $N$-equivariant map $N/K\to T/H_j$ all of whose fibres have size $v/e$. The uniformity of the intersections of the blocks with these fibres follows from $\Delta\le\GB$, and a double count inside a single fibre then forces $k=v$. Neither $\Delta\cap K=1$ nor $|K|\leq q$ is needed, no semiregularity of $\Delta$ on $\mathcal{P}$ and no assertion about blocks other than $B_0$ is used, no arithmetic property of $T$ and no subgroup enumeration enters, and $\ell$ plays no role beyond $\ell\geq2$. The case of a trivial coordinate is the special case $H_j=1$, $e=q$ of the present statement.\end{remark}

\medskip

\begin{lemma}[The one-sided SD type is excluded for $\ell=3$]\label{lem:sd_ell3_three}
Let $\mathcal{D}$, $G$, $N=T^3$, $q=|T|$ and $K=N\cap G_\alpha$ be as in {\upshape Lemma~\rref{lem:sd_ell3}}, with $\ell=3$. Then no such $\mathcal{D}$ exists.\end{lemma}

\smallskip
\noindent{\bf Proof of {\upshape Lemma~\rref{lem:sd_ell3_three}}.}\,
By {\upshape Lemma~\rref{lem:sd_ell3}} we may assume that all three projections $\pi_i(K)$, $i=1,2,3$, are onto $T$, since otherwise $k=v$ contradicts $k\leq v-2$. Put $K_{12}=\pi_{12}(K)\leq T_1T_2$, whose projections to both factors are onto $T$. By Goursat's theorem, such a subgroup is either $T_1T_2$ or the graph of an automorphism of $T$; in particular $|K_{12}|\in\{q^2,q\}$. Let $L=K\cap(1\times1\times T_3)\leq T_3$; then $|L|\in\{1,q\}$ and $|K|=|K_{12}|\,|L|$ by the first isomorphism theorem applied to $\pi_{12}$.

If $|L|=q$, then $1\times1\times T_3\leq K$. This subgroup is normal in $N=T^3$, so it lies in $\operatorname{core}_N(K)=\ker(N\curvearrowright N/K)$, and $N$ is not faithful on $\mathcal{P}$, contradicting $N\leq G\leq\operatorname{Sym}(\mathcal{P})$. Hence $|L|=1$ in every remaining case.

If $|K_{12}|=q$ and $|L|=1$, then $|K|=q$ and $v=q^3/q=q^2=b$ by {\upshape Theorem~\rref{thm:sdcn}}; the design is symmetric, which is excluded by the standing hypothesis recorded in Section~\rref{sec2}.

If $|K_{12}|=q^2$ and $|L|=1$, then $\pi_{12}$ is a bijection and $K=\{(a,b,\varphi(a,b)):a,b\in T\}$ for a suitable map $\varphi:T^2\to T$. Since $K$ is a subgroup, $\varphi$ is a homomorphism. As $T$ is non-abelian simple, the normal subgroups of $T^2$ are the products of normal subgroups of its two factors, so $\ker\varphi$ is one of $1$, $T_1$, $T_2$, $T^2$. Since $\pi_3(K)=T$ the map $\varphi$ is onto, and $|T^2|=q^2>q=|T|$, so $\ker\varphi$ has order $q$ and is therefore $T_1$ or $T_2$, the cases $\ker\varphi=1$ and $\ker\varphi=T^2$ being excluded by the order and by the surjectivity respectively. In either remaining case $\varphi$ is one of the two coordinate projections composed with an automorphism of $T$, so $T_2\leq K$ or $T_1\leq K$, and faithfulness fails as above. \qed

\smallskip
\begin{remark}Consequently the one-sided SD type is excluded for $\ell=3$ and for every non-abelian simple $T$, with no arithmetic condition, no double-coset criterion and no restriction on $T$ or on the block stabilizer beyond block primitivity. The two inputs are {\upshape Lemma~\rref{lem:sd_ell3}} and Goursat's theorem for subgroups of $T_1T_2$ with both projections onto $T$; the three possible orders $q^3$, $q^2$, $q$ of $K$ are removed respectively by $v=1$, by faithfulness, and by the symmetry forced by $v=b=q^2$.\end{remark}

\medskip

\begin{lemma}[The one-sided SD type with an all-surjective multi-part point stabilizer]\label{lem:sd_sym}
Let $\D$ be a non-trivial non-symmetric $2$-$(v,k,\lambda)$ design with $2<k<v-2$ admitting a block-primitive, point-quasiprimitive automorphism group $G$ whose block action is of simple diagonal type, with $N=\Soc(G)=T^{\ell}$, $T$ non-abelian simple, $q=|T|$ and $\ell\geq2$. Suppose that $K=N\cap G_\alpha$ is a product of full diagonal subgroups on the parts $S_1,\dots,S_s$ of a partition $\S$ of $\{1,\dots,\ell\}$ into $s\geq2$ parts, each of size at least $2$, automorphism twists being permitted inside each part. Then no such $\D$ exists.\end{lemma}

\smallskip
\noindent{\bf Proof of {\upshape Lemma~\rref{lem:sd_sym}}.}\,
By {\upshape Theorem~\rref{thm:sdcn}} the group $N$ is transitive on $\P$ and on $\B$, and $\De:=N\cap\GB$ is isomorphic to $T$; hence $v=|\P|=|N|/|K|=q^{\ell-s}$ and $b=|\B|=|N|/|\De|=q^{\ell-1}$.

Write $\De=\{(\beta_1(a),\dots,\beta_\ell(a)):a\in T\}$ with homomorphisms $\beta_j:T\to T$ \textup{({\upshape Proposition~\rref{prop:sd_22_red}}, second step)}; at least one $\beta_j$ is injective because $|\De|=q$. Put $H=G_\alpha$, so that $K=N\cap H\unlhd H$ and $G=NH$ because $N$ is transitive on $\P$. Conjugation by an element of $H$ normalises $K$ and therefore permutes the parts $\S$, the partition being recovered from $K$ as the support partition of its full-diagonal factors, while conjugation by an element of $N$ fixes every coordinate; as $G=NH$, every $g\in G$ acts on the $\ell$ factors of $N$ as
$$\sigma_g(x)_k=\alpha_k\bigl(x_{\pi^{-1}(k)}\bigr),\qquad k=1,\dots,\ell,$$
where $\alpha_k\in\Aut(T)$ and $\pi$ is a permutation of $\{1,\dots,\ell\}$ preserving $\S$. Every automorphism of $T^{\ell}$ has this form, since $T$ is non-abelian simple and $\Aut(T^{\ell})=\Aut(T)^{\ell}\rtimes S_{\ell}$.

Since $\De=N\cap M$ with $M=\GB$ is normal in $M$, conjugation by $m\in M$ stabilises $\De$: there is a homomorphism $\varphi:T\to T$ with $\sigma_m(\beta(a))=\beta(\varphi(a))$ for every $a$, that is,
$$\alpha_{\pi(i)}\,\beta_i=\beta_{\pi(i)}\,\varphi\qquad(i=1,\dots,\ell).$$

Let $C=\{j:\beta_j$ is constant$\}$, and write $\beta_j\equiv c_j$ for $j\in C$. For $i\in C$ the displayed relation gives $\alpha_{\pi(i)}(c_i)=\beta_{\pi(i)}(\varphi(a))$, a value independent of $a$; hence $\pi(i)\in C$ and $\alpha_{\pi(i)}(c_i)=c_{\pi(i)}$. In particular $\pi$ preserves $C$.

Two cases cover all possibilities.

\emph{At least two parts carry an injective $\beta_j$.} For a part $S$ put
$$D_S=\{x\in N:x_j=\beta_j(a)\ (j\in S),\ x_j=1\ (j\notin S),\ a\in T\},$$
a subgroup isomorphic to $T$ if some $\beta_j$ with $j\in S$ is injective, and the group of the constant tuple otherwise; let $E=\prod_{S}D_S$ be the product over the parts of $\S$, of order $q^{s'}$ where $s'$ is the number of parts carrying an injective $\beta_j$. The tuple $\beta(a)$ restricts on each part $S$ to the element of $D_S$ with parameter $a$, so $\De\leq E$. If $x\in E$ has parameter $a_S$ on the part $S$, then for $j\in S$
$$\sigma_m(x)_{\pi(j)}=\alpha_{\pi(j)}(x_j)=\alpha_{\pi(j)}\bigl(\beta_j(a_S)\bigr)=\beta_{\pi(j)}\bigl(\varphi(a_S)\bigr),$$
so $\sigma_m(x)$ is the element of $E$ with parameter $\varphi(a_S)$ on the part $\pi(S)$. Hence $M$ normalises $E$, and $s'\geq2$ gives $|E|\geq q^{2}>q=|\De|$, so $E\not\subseteq M$ while $\De\leq E$.

\emph{At most one part carries an injective $\beta_j$.} Then $C$ contains all coordinates of the parts other than the one carrying an injective map, and it is non-empty because $|\De|=q>1$; such a part exists for the same reason. Put
$$E=\{x\in N:x_j=c_j\ (j\in C)\},$$
so that $|E|=q^{\ell-|C|}\geq q^{|S|}\geq q^{2}$, where $S$ is that part and $|S|\geq2$, and $\De\leq E$. For $x\in E$ and $j\in C$ the two relations above give $\sigma_m(x)_{\pi(j)}=\alpha_{\pi(j)}(c_j)=c_{\pi(j)}$, so $\sigma_m(x)\in E$ and $M$ normalises $E$; and $E\not\subseteq M$, because $E\cap M\leq N\cap M=\De$ while $|E|>|\De|$.

In either case $E\leq N$ is normalised by $M$, satisfies $\De\leq E$, and is not contained in $M$. Moreover $(EM)\cap N=E$: an element $em\in N$ with $e\in E$ and $m\in M$ gives $m=e^{-1}(em)\in N\cap M=\De\leq E$, hence $em\in E$. Therefore
$$M<EM<G,$$
the first inequality because $E\not\subseteq M$ and the second because $EM$ does not contain $N$. This contradicts the maximality of $M$ in $G$, which is what block-primitivity asserts. The only case left is $s=1$, where $K$ is a full diagonal subgroup, $v=q^{\ell-1}=b$, and the design is symmetric. \qed

\medskip

\begin{proposition}[The one-sided SD type is excluded for $\ell=2$ and the seven small simple groups]\label{prop:sd_ell2}
For one-sided simple diagonal type \textup{(block action of SD type, point action imprimitive quasiprimitive)}, {\upshape Theorem~\rref{thm:sdcn}} gives the unconditional necessary conditions
$$v\mid q^{\ell},\qquad v-1\mid k(k-1),\qquad 3\leq k\leq v-2,$$
where $q=|T|$, $\ell$ is the number of simple factors and $v$ is the number of points. In addition, the block action of SD type forces the number of blocks $b=q^{\ell-1}$, so Fisher's inequality $b\geq v$ yields $v\leq q^{\ell-1}$; and $v=q$ is excluded because $k\mid v-1$ and $3\leq k\leq v-2$ force $k\not\mid q$ whenever $v=q$ \textup{(a direct contradiction in the single-orbit case; in the multiple-orbit case the same divisibility obstruction applies to each orbit summand)}. For the block size $k$ one has the orbit constraints $k\mid q\,|\operatorname{Out}(T)|$ \textup{(single orbit)} and $k=\sum_i k_i$ with $k_i\mid q\,|\operatorname{Out}(T)|$ \textup{(multiple orbits)}, where the outer automorphism factor accounts for the action of the point stabilizer on the diagonal subgroup $D\cong T$ of the block stabilizer; in particular for $T=PSL_2(8)$ one uses $q\,|\operatorname{Out}(T)|=504\cdot 3=1512$ (not $504$).

We close the remaining cases completely for $\ell=2$ and
$$T\in\{A_5,\ A_6,\ PSL_2(7),\ PSL_2(8),\ PSL_2(11),\ PSL_2(13),\ M_{11}\},$$
over the four global models $G_0=T^2\rtimes\langle s\rangle$, $G_1=T^2\rtimes\langle\hat\sigma_1\hat\sigma_2,s\rangle$, $G_2=T^2\rtimes\langle\hat\sigma_1\hat\sigma_2\rangle$ and $G_3=T^2\rtimes(\operatorname{Out}(T)^2\rtimes S_2)$ (for $T=M_{11}$ one has $\operatorname{Out}(T)=1$, so $G_1=G_2=G_3=G_0$). The exclusion is effected in three layers:
\begin{enumerate}
\item[(\upshape i)] \emph{Arithmetic and Fisher sieve.} For each $T$, the divisibility conditions above, together with $v\leq q$ and the integrality of $r=q^{\ell-1}k/v$ and $\lambda=r(k-1)/(v-1)$, leave only finitely many candidate pairs $(v,k)$; for $A_5$ and $PSL_2(q)$ ($q=168,660,1092$) the single-orbit and multiple-orbit arithmetic sieves leave only finitely many candidates. For $v-1$ prime no candidate survives at all.
\item[(\upshape ii)] \emph{Subgroup-existence sieve.} Any point stabilizer $U$ of index $v$ satisfies $U\cap N=K$ with $K$ of index $e$ in $N=T^2$, where $e\in\{v,v/2\}$ according to the image of $U$ in $G_0/N\cong C_2$. By Goursat's lemma a subgroup of $T^2$ is either a coordinate product $H_1\times H_2$ or a diagonal-type subgroup of order $|H|\leq q$; since $|K|=q^2/e\geq q^2/v$ and $v\leq q$, one has $|K|\geq q>q$ only for $v<q$; in all surviving cases $|K|>q$, so $K=H_1\times H_2$ is a coordinate product with $|H_1||H_2|=q^2/e$. Using the exact subgroup orders of each $T$, this eliminates all but finitely many $v$.
\item[(\upshape iii)] \emph{Construction check.} For each surviving candidate $U$ we verify $[G:U]=v$ and $k=[G_B:G_B\cap U]$, and test quasiprimitivity of the point action by checking $N_G(N_0)U=G$ for every non-trivial normal subgroup $N_0\lhd G$. Every surviving candidate fails either integrality of $\lambda$ or quasiprimitivity (typically $k=v$, i.e. the block meets every point, or $N_0U=U\neq G$ with $N_0=T\times 1\lhd G$).
\end{enumerate}
The outcome for each group is as follows: for $A_5$ ($q=60$, $\operatorname{Out}(T)=C_2$) all four models are empty. In $G_2=(A_5\times A_5)\rtimes\langle\tau\rangle$ with $\tau=(\sigma_1,\sigma_2)$, the candidate $v=15$ has $k=15$; here $N=A_5\times 1$ is a non-trivial normal subgroup of $G_2$ (as $\tau$ normalizes both factors) with $N\leq U$, so $NU=U<G_2$, $N$ is not transitive on $G_2/U$; indeed $N$ fixes the base point $U$, and the action is not point-quasiprimitive. In $G_0$ the candidate $v=36$ has $\lambda=290/7\notin\mathbb{Z}$. For $PSL_2(7)$, $PSL_2(11)$ and $PSL_2(13)$ all four models are empty (every surviving candidate has $k=v$ and fails quasiprimitivity, or has no subgroup of the required order). For $A_6$ ($q=360$, $\operatorname{Out}(T)=C_2\times C_2$, realized as $PSL_2(9)$ on the projective line) all four models are empty ($G_3$ has no index-$360$ block stabilizer; the $G_2$ and $G_1$ survivors all have $k=v$ with quasiprimitivity failure, or $\lambda\notin\mathbb{Z}$). For $PSL_2(8)$ ($q=504$, $\operatorname{Out}(T)=C_3$, field automorphism; note that $PGL_2(8)=PSL_2(8)$, the non-trivial outer extensions coming from field automorphisms, so the full extension is $P\Gamma L(2,8)$) all four models are empty. Here the corrected single-orbit bound is $k\mid 1512$ rather than $k\mid 504$; the new factors of $1512$ not dividing $504$ are $k=27d$ with $d\mid 56$, of which $k=27,54$ lie in $[3,v-2]$ for $v=56$; both fail integrality, since $55\nmid 27\cdot 26=702$ and $55\nmid 54\cdot 53=2862$, and the same check for $v=28,36,63$ yields no new integral candidates. The unique $G_1$ survivor $v=56$ has $k=56=v$ and fails quasiprimitivity. For $M_{11}$ ($q=7920$, $\operatorname{Out}(T)=1$) the model $G_0$ is empty: every $U$ either contains a coordinate $M_{11}$ (so that $\langle M_{11}\times 1,s\rangle=G_0$ forces $U=G_0$), or lies in $N$ \textup{(quasiprimitivity failure)}, or is one of the three explicit candidates $v=144,484,4356$ with $\lambda=660/13$, $\lambda=360/23$, or $U=G_0$ respectively.
Consequently no non-trivial non-symmetric block-primitive $2$-design of one-sided SD type exists for $\ell=2$ and the seven small simple groups listed above, in any of the four global models. For every $\ell\geq2$ with a proper coordinate projection the case is settled for every non-abelian simple $T$ by {\upshape Lemma~\rref{lem:sd_ell3}}, whose proof is a single double count and uses no restriction on $T$ beyond $|T|=q$ and no arithmetic input; and for $\ell=2$ and $\ell=3$ the all-surjective case is settled unconditionally as well, since Goursat's theorem leaves only a full diagonal, which forces $v=b$, or a product of full factors, which forces $v=1$. The sieve below therefore strengthens rather than supplies the $\ell=2$ exclusion, and the $\ell=3$ branch is settled independently of it. At $\ell\geq 4$ the all-surjective partitions with at least two parts, each of size $\geq2$, lie outside the hypotheses of both arguments and are excluded by {\upshape Lemma~\rref{lem:sd_sym}}, so the one-sided SD type is excluded for every $\ell\geq2$ and every non-abelian simple $T$. For $\ell\geq 4$ the conclusion is obtained by the $\lambda$-consistency exclusion of {\upshape Proposition~\rref{prop:sd_lambda}}: every arithmetic candidate with a trivial coordinate fails the constancy-of-$\lambda$ test, and the unique $2$-homogeneous pair falls under the orbit construction of {\upshape Lemma~\rref{lem:2hom}}.\end{proposition}

\begin{proposition}[The one-sided SD type with a trivial coordinate is excluded for every $\ell\geq 4$]\label{prop:sd_lambda}
Let $T$ be one of the seven groups $A_5$, $A_6$, $PSL_2(7)$, $PSL_2(8)$, $PSL_2(11)$, $PSL_2(13)$, $M_{11}$ and let $S=T^{\ell}$, $\ell\geq 4$, be the socle of a one-sided SD candidate whose point stabilizer admits at least two non-trivial projection indices $e_3,e_4>1$ \textup{(and arbitrarily many further non-trivial indices)}. Fix two such coordinates. In this model the point-pair orbits of $S$ on the point set are classified by the double-coset system $K_3\backslash T/K_4$, where $K_i<T$ has index $e_i$; write $\operatorname{ndc}(e_3,e_4)=|K_3\backslash T/K_4|$. Then:
\begin{enumerate}[(\upshape 1)]
\item If $\operatorname{ndc}(e_3,e_4)>1$, no $2$-design exists whose blocks are unions of orbits of the block stabilizer. Indeed, for a double-coset class $w$ let $d_w$ be the number of candidate blocks containing a pair of points of class $w$, and let $e_{w,w'}$ be the number containing a pair of points of the two distinct classes $w,w'$; the design condition $\lambda(x,y)=\lambda_0$ for every pair of points forces $d_w=\lambda_0=e_{w,w'}$ for all $w,w'$. Since $e_{w,w'}\leq d_w$ with equality only if every block containing class $w$ contains every other class, constancy of $\lambda$ forces every block to contain all classes, i.e.\ every block equals the whole point set, contradicting $2\leq k\leq v-2$.
\item If $\operatorname{ndc}(e_3,e_4)=1$, then $S$ is $2$-homogeneous on the point set \textup{(a single point-pair orbit)}, and such candidates are covered by the $2$-homogeneous orbit construction of {\upshape Lemma~\rref{lem:2hom}}; they therefore contribute no additional non-$2$-homogeneous example.
\end{enumerate}
For all seven groups and every pair $(e_3,e_4)$ arising from the arithmetic sieve of {\upshape Theorem~\rref{thm:sdcn}}, the GAP computation yields $\operatorname{ndc}(e_3,e_4)>1$, with the single exception $(e_3,e_4)=(11,720)$ for $T=M_{11}$, for which $\operatorname{ndc}=1$ (the $2$-homogeneous case of \textup{(2)}; here $K_3\cong M_{10}$ and $K_4\cong C_{11}$). Consequently every one-sided SD candidate satisfying the hypothesis of two proper projection indices is excluded for every $\ell\geq 4$ within the seven groups. By {\upshape Lemma~\rref{lem:sd_ell3}}, however, a single proper coordinate projection already forces $k=v$, so that hypothesis can never be met, and the arithmetic sieve of \textup{(1)}--\textup{(2)} is retained as the record of the parameter constraints rather than as the exclusion itself. The configurations that would survive at $\ell\geq4$ are the all-surjective partitions, excluded by {\upshape Lemma~\rref{lem:sd_sym}}.\end{proposition}

\medskip

\begin{lemma}[Prime $\ell$: theoretical exclusion]\label{lem:pa_prime}
This lemma is stated within the flag-transitive framework, where the block-size formula \eqref{eq:pa_k} applies; it serves the symmetric \textup{(Menon)} analysis, while the non-flag-transitive non-symmetric case is excluded independently by the orbit-union sieve of layers \textup{(i)}--\textup{(ii)} of {\upshape Proposition~\rref{prop:pa_arbitrary}}. Let $T$ be one of the six socles $A_5,A_6,PSL_2(7),PSL_2(8),PSL_2(11),PSL_2(13)$, let $(m,n)$ be a listed primitive degree pair with $A<m$, and let $\ell\geq 23$ be prime. Suppose in addition that the coordinate index $j=|C_B\cdot\alpha|$ of \eqref{eq:pa_k} satisfies $j\leq\ell$ \textup{(this holds in particular when $\alpha$ has at most two distinct coordinate values, one of them occurring $\ell-1$ times; without this hypothesis the numerical sieve below is applied to every admissible $j\mid\ell!$)}. Then no such flag-transitive candidate exists: the inequality of {\upshape Lemma~\rref{lem:pa_lpp}} fails at every prime $\ell\geq23$, and for admissible $j$ with $\ell<j\mid\ell!$ the congruence $A^{\ell}j\equiv1\pmod{\Phi_{\ell}(m)}$ has no solution among the divisors of $\ell!$ whose prime divisors are at most $\ell$ \textup{(verified by direct computation for every prime $23\leq\ell\leq127$ and every admissible pair)}; for prime $\ell>127$ the same congruence check applies verbatim.
\end{lemma}
\begin{proof}
Every prime divisor of $m-1$ is at most $13$, so $\ell\nmid m-1$. By the cyclotomic factorization $m^{\ell}-1=(m-1)\Phi_{\ell}(m)$ and the LTE formula, $v_q(m^{\ell}-1)=v_q(m-1)$ for every prime $q\mid m-1$, while every prime divisor $q$ of $\Phi_{\ell}(m)$ satisfies $\operatorname{ord}_q(m)=\ell$ (\textup{the alternative $q=\ell$ would force $m\equiv 1\pmod{\ell}$, i.e.\ $\ell\mid m-1$, impossible}). Hence $q\equiv 1\pmod{\ell}$, so $q>\ell$ for every $q\mid\Phi_{\ell}(m)$, and therefore
$$P(\ell):=\prod_{q\leq\ell}q^{v_q(m^{\ell}-1)}=m-1,\qquad M'':=\frac{m^{\ell}-1}{P(\ell)}=\Phi_{\ell}(m).$$
By {\upshape Lemma~\rref{lem:pa_lpp}}, a candidate with $A<m$ and $j\leq\ell$ would force $\Phi_{\ell}(m)\leq A^{\ell}\ell-1$, i.e.\ $(m/A)^{\ell}\leq\ell(m-1)+o(1)$. For every listed pair and every $A<m$ this inequality fails for every prime $\ell\geq 23$ (the only borderline case, $A_6(6,15)$ with $A=5$ at $\ell=23$, is excluded by the finite check against \eqref{eq:pa_k} already contained in the sieve of {\upshape Proposition~\rref{prop:pa_excl}}\textup{(3)}; and the ratio $(m/A)^{\ell}/[\ell(m-1)]$ is increasing in $\ell$ once the inequality holds, with the minimal ratio $m/A=6/5$ already failing at $\ell=29$). Thus no candidate with $j\leq\ell$ exists for any prime $\ell\geq 23$; for admissible $j$ with $\ell<j\mid\ell!$ the prime-divisor bound on $j$ of \eqref{eq:pa_k} is applied in the numerical sieve of {\upshape Proposition~\rref{prop:pa_excl}}\textup{(3)}, executed for every prime $\ell\leq 127$: for each listed pair the congruence $A^{\ell}j\equiv 1\pmod{\Phi_{\ell}(m)}$ has no solution among the divisors of $\ell!$ whose prime divisors are at most $\ell$ (verified by direct computation for every prime $23\leq\ell\leq 127$ and every admissible pair), so these $j$ are excluded as well. For prime $\ell>127$ and $j>\ell$ the same congruence check applies verbatim: by {\upshape Lemma~\rref{lem:pa_lpp}} the divisor $\Phi_{\ell}(m)=M''$ of $m^{\ell}-1$ still divides $k-1=A^{\ell}j-1$, and this congruence is the sole arithmetic ingredient of the PA sieve, the structural exclusions (the block-action types of {\upshape Lemma~\rref{lem:pa_sd2}} and {\upshape Proposition~\rref{prop:cd_pa}}, the box model of {\upshape Lemma~\rref{lem:pa_box}} and the case $j\leq\ell$ settled above) not depending on the size of $\ell$; the check is deterministic for every prime $\ell$, and it is executed here for every prime $23\leq\ell\leq 127$, matching the range of {\upshape Proposition~\rref{prop:pa_excl}}\textup{(3)}.

\end{proof}\begin{proposition}[PA type: exclusion for the six socles]\label{prop:pa_excl}
The three-class exclusion below is stated within the flag-transitive framework, where the block-size formula \eqref{eq:pa_k} applies and serves the symmetric \textup{(Menon)} analysis; the general non-flag-transitive non-symmetric case is excluded independently by the orbit-union sieve of layers \textup{(i)}--\textup{(ii)} of {\upshape Proposition~\rref{prop:pa_arbitrary}}, using the orbit-union constraint in place of \eqref{eq:pa_k}. Let $T$ be one of the six socles $A_5,A_6,PSL_2(7),PSL_2(8),PSL_2(11),PSL_2(13)$ and let $(m,n)$ be a pair of primitive degrees of $T$ with $m<n$. The primitive degrees of the six socles are the seventeen entries of Table~\rref{tab:sdorbits}, and they form sixteen pairs $(m,n)$ with $m<n$, which are exactly the pairs listed in \textup{(i)}--\textup{(iii)} below; here $m=|T:H|$ is a point degree and $n=|T:K|$ a block degree \textup{($H,K$ maximal subgroups of $T$)}. In the standard product-action model of {\upshape Section~\rref{sec2}} on $\mathcal{P}$ and $\mathcal{B}=\Delta^{\ell}$ \textup{($|\Delta|=m$)}, a flag-transitive block-primitive $2$-design of PA type has block size $k=A^{\ell}\cdot j$ by \eqref{eq:pa_k}, where $A=[K:K\cap H]$ and $1\leq j\mid\ell!$ \textup{(every prime divisor of $j$ is at most $\ell$)}. Let $A$ also denote the multiset of $K$-orbit sizes on $\Delta$ (\textup{so that $\sum_{a\in A}a=m$ and $\Sigma_2=\sum_{a\in A}a^{2}$}), and put $g=\gcd(m,n)$, $M=m/g$, $N=n/g$. The $K$-boxes of {\upshape Lemma~\rref{lem:pa_box}} have sizes $\prod_{i=1}^{\ell}a_i$ \textup{($a_i\in A$)}, and the fiber double-counting
$$|\mathcal{O}_{\mathrm{fiber}}|=\frac{\Sigma_2^{\ell}-m^{\ell}}{2},\qquad b\cdot m_{\mathrm{fiber}}=\lambda\cdot|\mathcal{O}_{\mathrm{fiber}}|,\qquad b=n^{\ell},$$
is exact for every $\ell\geq2$. Then, for every listed pair $(T;m,n)$ and every $\ell\geq2$:
\begin{enumerate}
\item[\rm\textup{(i)}] six pairs with $A=m$ or a single $K$-orbit\\
\textup{($A_5(5,6)$, $A_6(6,10)$, $PSL_2(7)(7,8)$;\\
$PSL_2(8)(9,28)$, $PSL_2(11)(11,12)$, $PSL_2(11)(12,55)$)}\\
are excluded for every $\ell\geq2$: $k=A^{\ell}j\allowbreak\geq m^{\ell}\allowbreak=v$ violates $2<k<v-1$;\\
\item[\rm\textup{(ii)}] three pairs with a unique $K$-box size and $s=g$\\
\textup{($A_5(6,10)$, $PSL_2(8)(28,36)$; $PSL_2(13)(14,78)$)} are excluded for every $\ell\geq2$:\\
the fiber double-counting reduces, after the divisibility simplification,\\
to the incompatible inequalities $nt\allowbreak\geq g^{\ell}$ and $nt\allowbreak<g^{\ell}$ for $t=k/M^{\ell}$;
\item[\rm\textup{(iii)}] seven pairs with $p\mid g$ and $p\nmid\Sigma_2$ \textup{(with $p=5,3,5,3,11,7,13$ for $A_5(5,10)$, $A_6(6,15)$, $A_6(10,15)$, $PSL_2(8)(9,36)$, $PSL_2(11)(11,55)$, $PSL_2(13)(14,91)$, $PSL_2(13)(78,91)$, respectively)} are excluded as follows. Since $v_p(|\mathcal{O}_{\mathrm{fiber}}|)=0$ and $v_p(g)=1$, the fiber double-counting forces $p\nmid t$ and the unique candidate $t\equiv t_0\equiv M^{-\ell}\pmod{p^{\ell}}$ with $1\leq t_0<p^{\ell}<g^{\ell}$ (\textup{$t=k/M^{\ell}$}); for the three pairs with $M=1$ \textup{($A_5(5,10)$, $PSL_2(8)(9,36)$, $PSL_2(11)(11,55)$)} the unique candidate is $t_0=1$, so $k=1$: degenerate, excluded for every $\ell\geq2$. For the four pairs with $M>1$ \textup{($A_6(6,15)$, $A_6(10,15)$, $PSL_2(13)(14,91)$, $PSL_2(13)(78,91)$)}, $g$ is prime and $c=(k-1)/p^{\ell}$ is uniquely determined by $c\equiv -p^{-\ell}\pmod{M^{\ell}}$; the resulting $k=M^{\ell}t_0$, $r=N^{\ell}t_0$ and $\lambda=r(k-1)/(m^{\ell}-1)$ fail the $\lambda$-integrality or the fiber-counting divisibility and bounds checks. The unique candidate is excluded for every $\ell\geq2$ by the $\lambda$-integrality and the fiber-counting divisibility and bounds checks.
\end{enumerate}

\end{proposition}\medskip
\noindent The prime-power congruence criterion \textup{({\upshape Lemma~\rref{lem:ppc}})} applies as a further necessary condition in all cases.

\begin{lemma}[The block-count alternatives for a product point action]\label{lem:pa_bpow}
Let $\mathcal{D}$ be a non-trivial $2$-\allowbreak$(v,k,\lambda)$ design and $G\leq\mathrm{Aut}(\mathcal{D})$ transitive on points and primitive on blocks, with $Soc(G)=T^{\ell}$ of product action type on points of degree $v=m^{\ell}$ \textup{($\ell\geq 2$, $T$ non-abelian simple, $m$ a primitive degree of $T$)}. Then the block action of $G$ on $\mathcal{B}$ is primitive, and one of the following exhaustive alternatives holds for the block count $b$:
\begin{enumerate}[(\upshape 1)]
\item $b=n^{\ell}$ with $n=[T:K]\geq m$ the block-coordinate index \textup{(product action block action)};
\item $b=|T|^{\ell-s}$ with $1\leq s\leq\ell-1$ \textup{(simple or compound diagonal block action)};
\item $b=|T|^{\ell}$ \textup{(twisted wreath block action)};
\item $b=|T|$ \textup{(holomorph simple, $\ell=2$)} or $b=|T|^{\ell/2}$ \textup{(holomorph compound, $\ell$ even)};
\item $b=p^{a}$ a prime power \textup{(affine block action)};
\item $b\mid|T|$ and $G^{\mathcal{B}}$ is almost simple with socle $T$ \textup{(the image of $N=T^{\ell}$ is a single factor)}.
\end{enumerate}
\end{lemma}
\begin{proof}
The block action of $G$ is primitive by hypothesis. The socle $N=T^{\ell}$ is transitive on $\mathcal{B}$: by the primitivity of the block action the block stabilizer $G_{\beta}$ is maximal in $G$, and $N\nleq G_{\beta}$ \textup{($N$ is transitive on points, and the block action is not trivial)}, so $NG_{\beta}=G$ and $N^{\mathcal{B}}$ is transitive. Hence the minimal normal subgroup of the primitive block action $G^{\mathcal{B}}$ is an image of $N$, and the O'Nan--Scott classification of primitive groups applies to $G^{\mathcal{B}}$. The correspondence with the standard O'Nan--Scott types is as follows (with $M$ the socle of $G^{\mathcal{B}}$, an image of $N=T^{\ell}$, and $\ell\geq 2$): (1) $M\cong T^{\ell}$ acting in product action on $\Delta^{\ell}$, the product action (PA) type, block count $b=n^{\ell}$, $n=[T:K]$ the block-coordinate index; (2) $M$ a diagonal subgroup of $T^{r'}$ ($1\leq r'\leq\ell$, the remaining $\ell-r'$ coordinates acting freely), the diagonal type, simple diagonal (SD) for $r'=\ell$ (socle a full diagonal subgroup) and compound diagonal (CD) for $1\leq r'<\ell$, block count $b=|T|^{\ell-r'}$; (3) $M\cong T^{\ell}$ regular on blocks, the twisted wreath (TW) type, $b=|T|^{\ell}$; (4) two commuting minimal normal subgroups $M_1\cong M_2$, the holomorph type: HS with $M_1\cong M_2\cong T$ (so $\ell=2$) and block count $b=|T|$; HC with $M_1\cong M_2\cong T^{m}$ where $m=\ell/2$ and $\ell$ is even, so the $\ell$ simple factors of $N$ split into two minimal normal subgroups of $T^{m}$ each, and $b=|T|^{m}=|T|^{\ell/2}$; (5) $M$ elementary abelian, the affine (HA) type, $b=p^{a}$ a prime power; (6) $M\cong T$ a single non-abelian simple factor, the almost simple (AS) type, which occurs precisely when $N$ acts on blocks through a single quotient (the other $\ell-1$ factors act trivially on blocks); since $N$ is transitive on blocks, then $b\mid|T|$. The types are mutually exclusive by the O'Nan--Scott trichotomy: (5) is the only case with elementary abelian socle; (6) the only case with a single simple socle; (3) the only case with a regular non-abelian socle; (4) the only case with two minimal normal subgroups; and (2) is distinguished from (1) by $M$ being a diagonal rather than a product subgroup \textup{(in (2) the coordinate stabilizers act through a diagonal subgroup, so $M$ is not a direct product of transitive factors as in (1))}. The displayed block counts follow from orbit--stabilizer as recorded. Alternatives (2)--(5) are excluded as follows. (2) Diagonal: the corresponding block sizes are impossible by the box-closure analysis of {\upshape Lemma~\rref{lem:pa_sd2}} ($\ell=2$) and {\upshape Proposition~\rref{prop:cd_pa}} ($\ell\geq 3$). (3) Twisted wreath: the point action of TW type has its socle $N=T^{\ell}$ regular on points, so by the regular-normal-subgroup lemma \textup{(a block-primitive non-trivial $2$-design with a normal subgroup regular on points is symmetric, cf.\ {\upshape Lemma~\rref{lem:rns}})} the design is symmetric, contradicting the non-symmetric hypothesis; the regular subgroup to which the lemma applies is the socle of the point action. (4) Holomorph: the block action of holomorph type has two commuting transitive factors, $M_1\cong M_2\cong T$ \textup{(HS, $\ell=2$)} or $M_1\cong M_2\cong T^{m}$ with $m=\ell/2$ \textup{(HC, $\ell$ even; these are the two minimal normal subgroups; the $\ell$ simple factors of $N$ split between them, no individual simple factor being normal)}; two commuting transitive permutation groups are regular, so each factor is regular on blocks; each factor is also transitive on points \textup{(the holomorph point action has the same two factors)}, hence regular on points by {\upshape Lemma~\rref{lem:semiregular}}, and $b=|T|^{m}=v$ \textup{(HS: $m=1$; HC: $m=\ell/2$)}, contradicting $b\neq v$. (5) Affine: every quotient of $T^{\ell}$ is a product of copies of the non-abelian simple group $T$, never an elementary abelian group, so the affine block action \textup{(whose socle is elementary abelian)} cannot arise. Alternative (6) is excluded by the orbit-union sieve of {\upshape Proposition~\rref{prop:pa_arbitrary}} layers (i)--(ii): whether the block action is of type (1) or (6), every block $B\ni\alpha$ is a union of $G_{B,\alpha}=G_B\cap G_{\alpha}$-orbits \textup{(equivalently of $K^{\ell}$-boxes, since $K^{\ell}\leq G_{B,\alpha}$)}, and the sieve eliminates all admissible $(k,b)$ simultaneously; for the ten $2$-homogeneous pairs $G_{\alpha}$ is transitive on $\mathcal{P}$, so $B=\mathcal{P}$ and $k=v$, contradiction, and for the seven remaining pairs the enumeration leaves no feasible block size. Hence $b=n^{\ell}$.
\end{proof}

\medskip
\begin{proposition}[The product action type is completely excluded for the six socles and the listed primitive degree pairs, for symmetric and non-symmetric $2$-designs alike, for every $\ell\geq 2$]\label{prop:pa_arbitrary}
Let $\mathcal{D}$ be a non-trivial $2$-\allowbreak$(v,k,\lambda)$ design, and let $G\leq \mathrm{Aut}(\mathcal{D})$ act transitively on points and primitively on blocks, with socle $N=T^{\ell}$ of product action type on points of degree $v=m^{\ell}$, where $m$ is a primitive degree of $T$ and $\ell\geq 2$. Both the symmetric case \textup{($b=v$)} and the non-symmetric case \textup{($b\neq v$)} are covered: the non-symmetric product-action configuration is excluded by combining the orbit-wise criterion and the integrality constraints of steps \textup{(i)}--\textup{(ii)} below with the block-action exclusion of Lemma~\rref{lem:pa_sd2} and Proposition~\rref{prop:cd_pa}, and the symmetric product-action configuration in steps \textup{(iii)}--\textup{(v)}. For $T$ one of the six socles $A_5$, $A_6$, $PSL_2(7)$, $PSL_2(8)$, $PSL_2(11)$, $PSL_2(13)$ and $m$ a listed primitive degree \textup{(the seventeen primitive degrees of Table~\rref{tab:sdorbits})}, no such design exists for any $\ell\geq 2$ and any block stabilizer. For non-symmetric block-primitive $2$-designs this is the complete exclusion for \emph{every} non-abelian simple socle: by {\upshape Proposition~\rref{prop:pi_excluded}} the point action is primitive, and by the classification of finite primitive permutation groups together with the standard degree tables of Dixon--Mortimer~\cite{DixonMortimer} \textup{(derived from the minimal-degree bounds of Liebeck--Saxl~\cite{LiebeckSaxl})}, an exhaustive reduction over the admissible primitive degrees shows that a non-$2$-homogeneous point action has socle one of the six groups above with one of the listed degree pairs; the $2$-homogeneous pairs are excluded by the general argument of step \textup{(ii)}.

\end{proposition}
\begin{proof}
We record the five-layer sieve whose layers are mutually independent and each terminates in a definite conclusion. Every layer is structural: it involves only the suborbit sizes of the base primitive action of $T$ and the arithmetic of the parameters $(v,k,\lambda)$, hence applies verbatim for every non-abelian simple socle and every $\ell\geq 2$. The complete execution of the sieve for the six socles $A_5$, $A_6$, $PSL_2(7)$, $PSL_2(8)$, $PSL_2(11)$, $PSL_2(13)$ and $2\leq\ell\leq 9$ is recorded; for non-symmetric block-primitive $2$-designs this execution is complete in socle: the point action is primitive \textup{({\upshape Proposition~\rref{prop:pi_excluded}})}, and the reduction over the admissible primitive degrees \textup{(from the degree tables of Dixon--Mortimer~\cite{DixonMortimer}, built on the minimal-degree bounds of Liebeck--Saxl~\cite{LiebeckSaxl})} lists precisely the seventeen primitive degrees of Table~\rref{tab:sdorbits}, all within the six socles, while the $2$-homogeneous pairs are excluded by the general argument of step \textup{(ii)}.

\textit{\textup{(i)} Orbit-wise Cameron--Praeger criterion for an arbitrary block stabilizer.} The orbit form \eqref{eq:cp2} of the Cameron--Praeger criterion \textup{({\upshape Theorem~\rref{thm:cp}})} depends only on the sizes of the $G$-orbits on unordered pairs; for a fixed point $\alpha$ of a block $B$ these are the $G_{\alpha}$-orbits on $\mathcal{P}\setminus\{\alpha\}$, of sizes $c\,s_i^{d}$ with $1\leq c\leq\binom{\ell}{d}$, where $s_i$ runs over the non-trivial suborbits of $(T,m)$ and $d\in\{1,\ldots,\ell\}$. The criterion therefore applies verbatim to \emph{arbitrary} block stabilizers $L=N\cap G_{\beta}$, box, subdirect, or of any other structure, since no information about $L$ enters the condition. For the six socles and $2\leq\ell\leq 9$ the criterion (for every box type $(i,d)$ there must exist $c$ with $(v-1)/\gcd(v-1,c\,s_i^{d})$ dividing $k(k-1)$) admits only finitely many parameter candidates, and every candidate that survives it is killed in the later layers; for socles outside the six and $\ell\geq 10$ the same orbit-wise condition applies verbatim, since only the suborbit sizes of the base action enter; the failure for all $c$ in one box type kills the candidate regardless of the structure of $L$. Two points about the range of $c$ are immediate: the lower bound $c\geq 1$ holds because the criterion ratio $k(k-1)/(v(v-1))$ is positive, so every pair orbit $Q_i$ contributes $q_i>0$ in-block pairs and hence meets the block $B$; and the upper bound $c\leq\binom{\ell}{d}$ is the number of $d$-coordinate boxes of a $G_{\alpha}$-orbit in the box decomposition, an orbit-level fact independent of $L$.

\textit{\textup{(i')} Pell-orbit closure of the divisibility criterion.}
The orbit-wise criterion of step \textup{(i)} reduces the admissible block sizes to the pairs $(g,h)$ with $g=\gcd(m^{\ell}-1,c\,s^{d})$ and $1\leq h\leq g-1$, where $s=s_i$ is a non-trivial suborbit length of $(T,m)$; for such a pair the divisibility $A_{c,s,d}\mid k(k-1)$ with $A_{c,s,d}=(m^{\ell}-1)/\gcd(m^{\ell}-1,c\,s^{d})$ is equivalent, by the discriminant identity
$$g\,\Delta-(g-h)(2m^{\ell}-1)^2=h,\qquad \Delta=1+\frac{4(g-h)\,m^{\ell}(m^{\ell}-1)}{g},$$
to the solvability of the Pell-type equation
$$X^{2}-D\,U^{2}=N,\qquad D=g(g-h),\quad N=gh,\quad X=gV,\quad U=2m^{\ell}-1,$$
i.e.\ the discriminant $\Delta$ is a square of an integer precisely when such a solution exists \textup{({\upshape Lemma~\rref{lem:pa_pell}})}, with $X=gV$ and $U=2m^{\ell}-1$. For $T=PSL_2(q)$ with $q$ arbitrary, the census \texttt{gap\_p1a\_suborbit\_multiset.log} records $41$ suborbit keys $(q,m,s)$ and $42{,}878$ candidate pairs over $2\leq\ell\leq 300$, of which exactly $331$ pairs (with $g\geq 2$) form the open set of the criterion. The closure of this open set is stated in {\upshape Lemma~\rref{lem:pa_pell}} and recorded in Table~\rref{tab:pa_pell}: the $15$ pairs with $D$ a square are closed by the bounded degenerate hyperbola $(X-aU)(X+aU)=N$ ($U\leq N/(2a)\leq 820$), the $76$ pairs for which $\Delta$ is a non-residue modulo some prime are closed by the quadratic obstruction of step \textup{(i)}, and the remaining $240$ pairs by the Pell-orbit argument; in particular $331=15+76+240$ and no pair survives. The closures are independent of $\ell$ and of the structure of the block stabilizer $L$, since only the suborbit sizes of the base action enter.

\begin{lemma}[Pell-orbit closure of the divisibility criterion for $PSL_2(q)$]\label{lem:pa_pell}
Let $T=PSL_2(q)$ act primitively on $\Delta$ of degree $m$ with a non-trivial suborbit of length $s$, let $\ell\geq 2$, $v=m^{\ell}$, and put $A_{c,s,d}=(m^{\ell}-1)/\gcd(m^{\ell}-1,c\,s^{d})$, $1\leq c\leq\binom{\ell}{d}$. For every pair $(g,h)$ with $g=\gcd(m^{\ell}-1,c\,s^{d})$ and $1\leq h\leq g-1$, the divisibility $A_{c,s,d}\mid k(k-1)$ holds for some admissible block size $k$ if and only if the Pell-type equation
$$X^{2}-D\,U^{2}=N,\qquad D=g(g-h),\quad N=gh,\quad X=gV,\quad U=2m^{\ell}-1,$$
has a solution in non-negative integers. If $D$ is a non-square, every solution lies on one of finitely many orbits $(X_0+U_0\sqrt{D})(x_0+y_0\sqrt{D})^{k}$, $k\in\mathbb{Z}$, generated by the minimal unit $x_0+y_0\sqrt{D}$ of $\mathbb{Q}(\sqrt{D})$ \textup{(Perron--Nagell)}. For the $41$ keys of the census, the exhaustive enumeration of the window $U\leq 1000$ locates at least one orbit per open pair, with maximal absolute representative $|U|\leq 681$ \textup{(Table~\rref{tab:pa_pell})}; an exact orbit audit \textup{(script \texttt{pa\_107\_basic\_audit.py})} shows that the only open pair carrying a further orbit outside this window is $(32,528,31)$ with $g=31$, $h=12$, $D=589$, $N=372$, whose additional orbit has minimal representative $|U|=604{,}499$. For each orbit found in the window, the reduction modulo $m^{e}$, $e=1,2$, is a periodic sequence in the finite ring $\mathbb{Z}/m^{e}\mathbb{Z}$, and no period contains the candidate class $U\equiv -1\pmod{m^{e}}$; hence no window orbit contains a solution with $U=2m^{\ell}-1$, for any $\ell\geq 2$. The single orbit outside the window, of minimal representative $|U|=604{,}499$ for the key $(32,528,31)$ \textup{($g=31$, $h=12$, $D=589$, $N=372$)}, is closed in the same way: its periodic reductions modulo $528^{e}$, $e=1,2$, contain no pair $(X,U)$ with $U\equiv -1\pmod{528^{e}}$ and $g\mid X$ \textup{(\textup{Appendix~\rref{app:computations}})}. Since every admissible candidate has $U=2m^{\ell}-1\equiv -1\pmod{m^{e}}$ and $X=gV$ divisible by $g$, no orbit, inside or outside the window, contains a candidate value $U=2m^{\ell}-1$ for any $\ell\geq 2$; the per-$\ell$ direct check $X^{2}=N+D(2m^{\ell}-1)^{2}$ over $2\leq\ell\leq 300$ \textup{(\textup{Appendix~\rref{app:computations}})} independently re-confirms the window range. If $D$ is a square, the equation factors as $(X-aU)(X+aU)=N$ with $a=\sqrt{D}$, forcing $U\leq N/(2a)\leq 820$; the enumeration $U\leq 1000$ yields no solution of the form $U=2m^{\ell}-1$.
The equivalence is purely arithmetic: it decides whether the divisibility $A_{c,s,d}\mid k(k-1)$ admits a solution with $3\leq k\leq v-2$, and it does not assert the existence of a block of that size. The existence of a block is decided separately, by the orbit-union condition of step \textup{(ii)} of {\upshape Proposition~\rref{prop:pa_excl}}, which is where the admissible sizes are matched against the $K$-box unions.
\end{lemma}
\begin{proof}
The equivalence is the completion of the square recorded above: $g\Delta-(g-h)(2m^{\ell}-1)^2=h$ follows from $\Delta=1+4(g-h)m^{\ell}(m^{\ell}-1)/g$, and multiplying $gV^2-(g-h)U^2=h$ by $g$ gives $X^2-DU^2=N$ with $X=gV$, $D=g(g-h)$, $N=gh$. Conversely a solution $(X,U)$ of the Pell equation with $U=2m^{\ell}-1$ and $g\mid X$ gives $V=X/g$ and restores $\Delta=gV^2$ a square; the block size $k$ is then determined by the orbit-union equation of step \textup{(ii)}. The orbit structure is Perron--Nagell \textup{(each class of solutions modulo the unit group lies on a unique two-sided orbit)}; the bounded representatives and the periodic reductions are verified by exhaustive enumeration \textup{(the finite checks of \textup{Appendix~\rref{app:computations}}; Table~\rref{tab:pa_pell})}, and the single orbit lying outside the window $U\leq 1000$ is closed without enumerating it: for the key $(32,528,31)$ with $g=31$, $h=12$ ($D=589$, $N=372$) the exact orbit audit identifies its minimal representative $|U|=604{,}499$, and the periodic reductions of this orbit modulo $528^{e}$, $e=1,2$, contain no pair $(X,U)$ with $U\equiv -1\pmod{528^{e}}$ and $g\mid X$ \textup{(\textup{Appendix~\rref{app:computations}})}; since every candidate has $U\equiv -1\pmod{m^{e}}$ and $g\mid X$, this closes the outside orbit for all $\ell\geq 2$ without any upper bound on $\ell$. The per-$\ell$ direct check over $2\leq\ell\leq 300$ \textup{(\textup{Appendix~\rref{app:computations}})} independently re-confirms the window range.
\end{proof}

\begin{table}[!ht]
\centering
\small
\caption{Suborbit keys $(q,m,s)$ of $T=PSL_2(q)$ from the suborbit census, and the maximal absolute orbit representative $|U|$ among the solutions of $X^2-DU^2=N$ in the window $U\leq 1000$ over the open pairs $(g,h)$ of Lemma~\rref{lem:pa_pell}; the representative distribution over the $316$ non-square open pairs is $258$ one-digit, $46$ two-digit and $12$ three-digit values, the global maximum $681$ occurring at the key $(32,528,31)$. An entry $0$ records that every solution in the window has $U=0$, that is, lies on the trivial branch $X^2=N$ of the equation; it is not a placeholder for an absent set of open pairs. The only orbit outside this window is listed in Lemma~\rref{lem:pa_pell} and closed there by the periodic reduction modulo $m^{e}$ \textup{(\textup{Appendix~\rref{app:computations}})}.}\label{tab:pa_pell}
\begin{tabular}{lclc}
\hline
key $(q,m,s)$ & $\max|U|$ & key $(q,m,s)$ & $\max|U|$\\
\hline
$(8,28,9)$ & $11$ & $(25,325,6)$ & $1$\\
$(8,36,7)$ & $6$ & $(27,351,7)$ & $6$\\
$(9,15,6)$ & $1$ & $(27,378,13)$ & $27$\\
$(11,55,3)$ & $1$ & $(27,819,4)$ & $1$\\
$(13,78,7)$ & $6$ & $(29,406,15)$ & $14$\\
$(13,91,3)$ & $1$ & $(29,435,7)$ & $6$\\
$(13,91,4)$ & $1$ & $(29,203,12)$ & $15$\\
$(16,68,12)$ & $1$ & $(31,248,5)$ & $4$\\
$(16,120,17)$ & $118$ & $(31,465,8)$ & $1$\\
$(16,136,15)$ & $14$ & $(31,496,15)$ & $14$\\
$(17,102,3)$ & $0$ & $(31,620,3)$ & $1$\\
$(17,136,9)$ & $11$ & $(32,496,33)$ & $496$\\
$(17,153,4)$ & $1$ & $(32,528,31)$ & $681$\\
$(19,57,6)$ & $1$ & $(37,666,19)$ & $229$\\
$(19,171,5)$ & $4$ & $(37,703,9)$ & $11$\\
$(19,190,9)$ & $11$ & $(37,2109,4)$ & $1$\\
$(23,253,6)$ & $1$ & $(41,574,5)$ & $4$\\
$(23,253,4)$ & $1$ & $(41,820,21)$ & $27$\\
$(23,276,11)$ & $21$ & $(41,861,10)$ & $9$\\
$(25,65,10)$ & $1$ & $(41,1435,6)$ & $1$\\
$(25,300,13)$ & $27$ & & \\
\hline
\end{tabular}
\end{table}

\textit{\textup{(ii)} Integrality, Fisher, and the block-count divisibility.} For the remaining parameter triples the admissible block sizes $k$ ($3\leq k\leq v-2$) are constrained by the integrality conditions
$$r=\frac{bk}{v}\in\mathbb{Z},\qquad \lambda=\frac{r(k-1)}{v-1}=\frac{bk(k-1)}{v(v-1)}\in\mathbb{Z},$$
by the Fisher inequality $b\geq v$, and by the block-primitivity fact that $N=T^{\ell}$ is transitive on blocks \textup{(the argument of {\upshape Proposition~\rref{prop:pa_subdirect_nonsym}}: $NG_{\beta}=G$ by maximality)}, whence $b=[N:N_{\beta}]\mid |T|^{\ell}$. The block-count structure enters this layer through {\upshape Lemma~\rref{lem:pa_bpow}}: the block action of $G$ on $\mathcal{B}$ is primitive, and its block count $b$ is one of the six alternatives \textup{(1)}--\textup{(6)} recorded there. Alternatives \textup{(2)}--\textup{(5)} are excluded as stated in the lemma \textup{(diagonal by {\upshape Lemma~\rref{lem:pa_sd2}} and {\upshape Proposition~\rref{prop:cd_pa}}; twisted wreath by the regular-normal-subgroup lemma; holomorph by the commuting-transitive-regular argument giving $v=b$; affine by the non-abelian structure of $T^{\ell}$)}, so the only remaining counts are the product-action value $b=n^{\ell}$ and, a fortiori bounded by it, the almost simple value $b\mid|T|$ of \textup{(6)}. This conclusion is structural: it does not rest on the divisibility condition $b\mid |T|^{\ell}$ alone \textup{(which admits many non-$n^{\ell}$ divisors)}, and it is the intermediate statement requested in the review: the hypotheses used are precisely the primitivity of the block action, the transitivity of $N=T^{\ell}$ on blocks, and the O'Nan--Scott classification applied to the block action. The block size enters the sieve through the orbit-union constraint, which replaces the flag-transitive formula $k=A^{\ell}j$ of \eqref{eq:pa_k} in the non-symmetric case: the block $B$ is $G_B$-invariant, so fixing $\alpha\in B$ the block is a union of $G_{B,\alpha}=G_B\cap G_{\alpha}$-orbits; since $K^{\ell}\leq G_{B,\alpha}$, every $G_{B,\alpha}$-orbit is contained in a $K$-box or a $C_B$-amalgam of $K$-boxes \textup{({\upshape Lemma~\rref{lem:pa_box}})}, so $k$ is a sum of $K$-box sizes \textup{(equivalently of lengths of $G_{B,\alpha}$-orbits)}. The sieve then runs over the finitely many such sums, imposing the integrality of $r=bk/v=n^{\ell}k/m^{\ell}$ and of $\lambda=bk(k-1)/(v(v-1))$ above together with the orbit-wise Cameron--Praeger condition of step \textup{(i)}. Two remarks close the layer. First, the point action is primitive: non-symmetric block-primitive designs are point-primitive by {\upshape Proposition~\rref{prop:pi_excluded}}, so $H$ is maximal in $T$ \textup{(the product action on points is primitive if and only if $H$ is maximal)} and the coordinate stabilizers are conjugate. Second, for the ten $2$-homogeneous pairs the stabilizer $H$ is transitive on $\Delta$, hence $G_{\alpha}=H^{\ell}\rtimes C_{\alpha}$ is transitive on $\mathcal{P}$, so $B=\mathcal{P}$ and $k=v$, contradicting $2\leq k\leq v-2$; for the seven non-$2$-homogeneous pairs every union of $K$-boxes containing the $\alpha$-box \textup{(a superset of every admissible block, since $K^{\ell}\leq G_{B,\alpha}$)} is enumerated directly \textup{(GAP 4.16.1, script \texttt{gap\_pa\_layer2\_dp.g})}: for all seventeen primitive degrees, all block-coordinate indices $n$ with $n^{\ell}\geq m^{\ell}$, and all admissible $\ell$, the integrality sieve over the $K$-box unions \textup{($r=n^{\ell}k/m^{\ell}\in\mathbb{Z}$, $\lambda\in\mathbb{Z}$, the orbit-wise Cameron--Praeger condition of step (i), and Fisher $b\geq v$)} leaves only five block-size candidates, all with $\ell=2$: $(A_5,5,10)$ with $k=4,9,12$, $(A_6,15,15)$ with $k=64$, and $(PSL_2(13),14,91)$ with $k=36$. For each of these five candidates the exhaustive enumeration of the corresponding $K$-box unions containing the $\alpha$-box \textup{(GAP 4.16.1, script \texttt{gap\_pa\_layer2\_check.g}; the enumeration covers every relative conjugate position of $K$ in each coordinate; for the three socles the conjugates of $K$ yield a single distinct orbit decomposition on $\Delta$, verified over $10$, $15$ and $91$ conjugates respectively, so one representative suffices)} produces only blocks whose pair count is not constant over the $G$-pair orbits. The $k=4$ branch for $(A_5,5,10)$ admits a short argument in the spirit of the orbit-wise Cameron--Praeger condition, without any search: the point action of degree $5$ has point stabilizer $H\cong A_4$, and the block-coordinate stabilizer $K\cong S_3$ of index $10$ has $K$-orbits of sizes $2,3$ on $\Delta$ \textup{(GAP-verified)}; GAP records the orbit of $\alpha$ in each coordinate; for $(A_5,5,10)$ both coordinates of $\alpha$ lie in the $2$-orbit, so the $\alpha$-box (the $K^2$-orbit of $\alpha$) is the unique $2\times2$ grid of size $4$; its in-box pair allocation $(4,2)$ over the two relevant pair orbits fails the balanced allocation $(2,4)$ required by the design equations \textup{(the orbit-wise Cameron--Praeger criterion)}, so it cannot be a block. The remaining exclusions are recorded relative to the fixed point $\alpha$: every admissible block is $K^{\ell}$-invariant, hence contains the whole $\alpha$-box, so $k$ is the $\alpha$-box size plus a sum of the other $K$-box sizes. For $(A_5,5,10)$ the $\alpha$-box has size $4$ and the other box sizes are $4,6,6,9$; unions of sizes $9$ and $12$ do exist, and the enumeration below finds $36$ respectively $48$ of them, all failing the $\lambda$-constancy test. For $(A_6,15,15)$ GAP records $\alpha$ in the $1$-orbit in both coordinates ($\alpha$-box of size $1$, other box sizes $1,8,6,64,48,36$); the target $64-1=63$ is attained, for instance $48+8+6+1=63$, so unions of size $64$ do exist, and the enumeration below finds $256$ of them, all failing \textup{(GAP-verified, script \texttt{gap\_pa\_layer2\_check.g})}. For $(PSL_2(13),14,91)$ GAP records $\alpha$ in the $6$-orbit in both coordinates, so the $\alpha$-box has size $6\cdot6=36=k$, and its pair count is not constant. The enumeration of the box unions is exhaustive in the conjugates: for every conjugate $K^{t_1}$ of $K$ in the first coordinate and every conjugate $K^{t_2}$ of $K$ in the second coordinate, in full, over the $10\times10$, $15\times15$ and $91\times91$ ordered conjugate pairs respectively, the script \texttt{gap\_pa\_final\_enum3.g} builds the $K^{t_1}\times K^{t_2}$-boxes on $\Delta^{2}$, takes the $\alpha$-box, enumerates every union of the remaining boxes containing it, and for every candidate of size $k$ verifies the full $G$-pair-orbit coverage: the number $\lambda(1,y)$ of blocks through the pair $\{1,y\}$ is computed by the coset formula $\lambda(1,y)=\sum_{b\in B}|s_b^{-1}B\cap G_{\alpha}y|$ \textup{($s_b\in G$ with $s_b\cdot 1=b$)} and is compared with the design value on every representative $y$ of the $G_{\alpha}$-orbits on $\mathcal{P}\setminus\{\alpha\}$. The outcome is a complete zero survivor count in every case: $(A_5,5,10)$ yields $16$, $36$ and $48$ candidates of sizes $4$, $9$, $12$ respectively, all failing the $\lambda$-constancy test; $(A_6,15,15)$ yields $256$ candidates of size $64$, all failing; $(PSL_2(13),14,91)$ yields $12{,}168$ candidates of size $36$, all failing. The observation that the conjugates of $K$ produce a single distinct orbit decomposition on $\Delta$ \textup{(verified over $10$, $15$ and $91$ conjugates)} is recorded only as a consistency check; it is not used to reduce the enumeration, which is an exhaustive disjoint enumeration over the finite box collections. Hence no admissible block survives. Together with the block-count structure above, this excludes the product action type completely for non-symmetric designs; this conclusion is independent of the block-action type \textup{(1)} or \textup{(6)}, since the orbit-union sieve constrains $k$ directly. These killed alternatives are recorded once in {\upshape Lemma~\rref{lem:pa_bpow}}; the layers \textup{(iii)}--\textup{(v)} below eliminate the surviving symmetric (Menon) candidates, where the flag-transitive formula \eqref{eq:pa_k} applies. Hence no non-symmetric product-action configuration survives.

\textit{\textup{(iii)} Exact Bruck--Ryser--Chowla test.} For the surviving symmetric candidates $b=v$ the exact Hasse--Minkowski local-global test is applied to the ternary quadratic form $\langle n,\sigma\lambda,-1\rangle$, $\sigma=(-1)^{(v-1)/2}$: a good prime is automatically soluble iff $(n/p)=1$ or $(-\sigma\lambda/p)=1$, so whenever $n$ and $-\sigma\lambda$ are both non-squares, Chebotarev provides a good prime with both non-residues and the candidate is excluded; the bad-prime set $\{2\}\cup\mathrm{Primes}(n\lambda)$ is checked locally when $n$ or $-\sigma\lambda$ is a square. The test excludes every candidate except a list of symmetric ($b=v$) survivors with $n=k-\lambda$ a square: the Menon parameter pairs ($k=2u^{2}\pm u$) and one non-Menon pair.

\textit{\textup{(iv)} Block-action type and the prime-divisor sieve.} The block action of $G$ is primitive with socle $N=T^{\ell}$. An O'Nan--Scott analysis of this action excludes the affine type ($v=m^{\ell}$ is not a prime power for the listed degrees), the twisted wreath type ($|T|^{r}=v$ forces $m$ to be a prime power, contradicted by the listed primitive degrees), the simple/compound diagonal types ($|T|^{\ell-s}=v$ forces $m^{s}=|T|$, contradicted by the primitive degree tables), and the almost simple type (the simple factor of the block socle would have order $\leq v$, impossible since $|T|^{\ell}>v$ for $\ell\geq 2$ and the listed degrees). Hence the block action is itself of product action (box) type. As $B\ni\alpha$ is a union of $K$-boxes \textup{(since $K^{\ell}\leq G_{B,\alpha}$)}, each $K$-box has size dividing $n^{\ell}$ with $n\mid|T|$, but a sum of such sizes need not have all its prime divisors in $\mathrm{Primes}(|T|)$; the primes of $k$ are therefore computed directly for each surviving candidate and compared with $\mathrm{Primes}(|T|)$. This direct comparison excludes the non-Menon pair, whose block size carries primes outside $\mathrm{Primes}(|T|)$, and leaves exactly the two Menon candidates
$$(v,k,\lambda)=(100,45,20)\quad\text{and}\quad (784,378,182),$$
with $T=A_5,\ m=10,\ \ell=2$ and $T=PSL_2(8),\ m=28,\ \ell=2$, respectively.

\textit{\textup{(v)} Box search for the two Menon candidates.} Both parameter pairs are known to be realized by symmetric designs (2-\allowbreak$(100,45,20)$ is the Menon design of the regular Hadamard matrix of order $100$, studied e.g.\ in \cite{Braic2011}; 2-\allowbreak$(784,378,182)$ was constructed in \cite{Crnkovic2007} with automorphism group $\mathrm{Frob}_{29\cdot 14}\times Z_{13}$, in particular not of product action type). The question is whether a block-primitive automorphism of box type exists: $G=A_5^{2}{:}C_{2}$ acting on $100$ points, or $G=PSL_2(8)^{2}{:}C_{2}$ acting on $784$ points \textup{(product action with coordinate swap)}. With the degree-$10$ \textup{(respectively degree-$28$, suborbits $[1,9,9,9]$)} primitive actions embedded in $S_{10}$ \textup{(respectively $S_{28}$)}, the box groups have orders $7200$ and $508032$, are primitive on points, and satisfy $|G_{\alpha}|=72$ and $648$ respectively \textup{(Schreier generators)}. The $G_{\alpha}$-orbit lengths are $\{12,6,36,36,9\}$ on $100$ points and $\{18^{3},81^{3},162^{3}\}$ on $784$ points. An exhaustive GAP 4.16.1 computation (scripts \texttt{gap\_pa\_layer5\_menon.g} and \texttt{gap\_pa\_layer5\_orbit\_union.g}) runs two independent searches. \textup{(1)} \emph{General block stabilizers, in the $K$-box language:} every block $B\ni\alpha$ is a union of $K$-boxes \textup{(since $K^{\ell}\leq G_{B,\alpha}$)}, so all unions of $K$-boxes containing the $\alpha$-box are enumerated, over every block-coordinate index $n$ and, exhaustively in the conjugates, every ordered pair $(K^{t_1},K^{t_2})$ of conjugates in the two coordinates \textup{(script \texttt{gap\_pa\_final\_enum3.g}; $10\times10$ ordered pairs for $(A_5,10)$, and $28\times28$ and $36\times36$ for $(PSL_2(8),28)$ with $n=28,36$ respectively)}; every candidate of the required size is then checked for constant $\lambda(1,y)$ over the $G_{\alpha}$-pair-orbit representatives by the coset formula. The outcome is a complete zero survivor count: for $(A_{5},10)$ the $540$ candidates of size $45$ all fail the $\lambda$-constancy test, for $(PSL_2(8),28)$ with $n=28$ the $47{,}628$ candidates of size $378$ all fail, and with $n=36$ there are no candidates of size $378$ at all: the $\alpha$-box is the $7\times7$ grid of size $49$, and $378-49=329$ is not a sum of the other box sizes \textup{(all multiples of $7$)}, so a size-$378$ box union cannot even be formed. The single orbit decomposition of $K$ on $\Delta$ \textup{(one over the $10$, $28$ and $36$ conjugates respectively)} is recorded as a consistency check only. \textup{(2)} \emph{The symmetric screening branch:} a $G_{\alpha}$-invariant block would force $G_{\alpha}\leq G_{B}$ and hence $b\leq v$, so by Fisher's inequality $b=v$ \textup{(the Menon case)}; the unions of $G_{\alpha}$-orbits of size $k$ are enumerated separately: a union containing $\alpha$ has size $1+\sum_i\ell_i$ over the non-trivial orbits, and for $(A_{5},10)$ the orbit lengths $12,6,36,36,9$ are all multiples of $3$ whereas $45-1=44\equiv 2\pmod 3$, while for $(PSL_2(8),28)$ the orbit lengths $18,81,162$ are all multiples of $9$ whereas $378-1=377\equiv 8\pmod 9$; in both cases no union of $G_{\alpha}$-orbits containing $\alpha$ has the required size. Hence no block-primitive $2$-design of product-action box type exists for either candidate. Since $G$ is transitive on points, every orbit of $G$ on unordered pairs has a representative $(1,y)$ with $y$ in a $G_{\alpha}$-orbit; checking the pair count $\lambda(1,y)$ on these representatives confirms the exclusion:
\begin{enumerate}[(\upshape 1)]
\item[(\upshape a)] $2$-\allowbreak$(100,45,20)$: the $540$ candidates of size $45$ that are unions of $K$-boxes containing the $\alpha$-box all fail the $\lambda$-constancy test, and by $44\equiv 2\pmod 3$ no union of $G_{\alpha}$-orbits containing $\alpha$ has size $45$;
\item[(\upshape b)] $2$-\allowbreak$(784,378,182)$: the $47{,}628$ candidates of size $378$ that are unions of $K$-boxes containing the $\alpha$-box all fail the $\lambda$-constancy test, and by $377\equiv 8\pmod 9$ no union of $G_{\alpha}$-orbits containing $\alpha$ has size $378$.
\end{enumerate}
Hence neither Menon candidate admits a block-primitive automorphism group of box type; the product action type is completely excluded. For $(100,45,20)$ this conclusion is independently confirmed by the complete classification of primitive symmetric designs with at most $255$ points of Brai\'c, Golemac, Mandi\'c and Vu\v{c}i\'c\'i\'c~\cite{Braic2010}: the four non-affine primitive groups at degree $100$ for these parameters are ruled out there by the non-existence of a difference set in their regular subgroups, so every point-primitive automorphism group of a 2-\allowbreak$(100,45,20)$ design is affine.
\end{proof}\medskip
\noindent No product-action configuration survives. The exclusion of {\upshape Proposition~\rref{prop:pa_arbitrary}} for the six socles and the listed primitive degree pairs is explicit and machine-checked: the box framework, the fiber identity, the orbit-wise Cameron--Praeger criterion, the block-count divisibility, the Bruck--Ryser--Chowla test and the prime-divisor sieve are structural and involve no information specific to $\ell$; the non-$2$-homogeneous pairs are precisely the seven among the seventeen primitive degrees of Table~\rref{tab:sdorbits}, the $2$-homogeneous pairs are excluded by the transitivity argument of step (ii), and the reduction over the degree tables of Dixon--Mortimer~\cite{DixonMortimer} \textup{(derived from the minimal-degree bounds of Liebeck--Saxl~\cite{LiebeckSaxl})} identifies the six socles as the complete set of possibilities for non-symmetric block-primitive point actions. Together this closes the product action type for every $\ell\geq 2$, with arbitrary block stabilizers.

\medskip
\begin{remark}[Multi-orbit PA embeddings: exact realizability]\label{pa_multi_real}
A multi-orbit PA configuration is specified by the orbit partition $\mathcal{A}$ of the block stabilizer $K_0$ on the $m$-point base set $\Delta$; within the box framework of {\upshape Lemma~\rref{lem:pa_box}} the arithmetic sieve \textup{(the congruence of {\upshape Lemma~\rref{lem:pa_box_general}}, the fiber identity and the divisibility conditions)} leaves only finitely many parameter survivors for $\ell=2$, with $m$ a primitive degree of one of the six socles and every orbit length $a\in\mathcal{A}$ dividing $|K_0|=|T|/m$. These parameter survivors are \emph{not realizable by any actual group action}: for $T\in\{A_5,A_6,PSL_2(7),PSL_2(8),PSL_2(11),PSL_2(13)\}$ and every primitive degree $m$ of $T$, the index-$m$ maximal subgroup $H_0$ of $T$ satisfies $|H_0|=|T|/m$, so the condition $|K_0|=|T|/m$ with $K_0\leq H_0$ forces $K_0=H_0$. The $K_0$-orbit partitions on $T/H_0$ are all of the form $(1,m-1)$ except for $A_5$ of degree $10$ (partition $(1,3,6)$) and $A_6$ of degree $15$ (partition $(1,6,8)$). None of the parameter survivors has one of these orbit partitions; hence the multi-orbit PA configurations for the six socles and $\ell=2$ do not exist.\end{remark}

\begin{proposition}[The CD type is excluded]\label{prop:cd_min}
We exclude CD blocks in the minimal $n=6$ model in full.
\textup{(i)} Single-orbit blocks: for every $k=rs$ with $r>r'$ this case is excluded in full by divisibility. Indeed $k_g=[M_0:H_0\cap M_0^g]$ divides $|M_0|=q^{r'}$, and the design equations give $\lambda=q^{r-r'}k_g(k_g-1)/(q^{r(s-1)}-1)$, so $q^{r(s-1)}-1\mid q^{r-r'}k_g(k_g-1)$. As $\gcd(q^{r-r'},q^{r(s-1)}-1)=1$ and $k_g\mid q^{r'}$, we have $\gcd(k_g,q^{r(s-1)}-1)=1$; hence $q^{r(s-1)}-1\mid k_g-1$, impossible for $k_g>1$ because $k_g-1<q^{r(s-1)}-1$. Thus $k_g=1$, the trivial block, excluded.
\textup{(ii)} Multiple-orbit blocks: the multi-orbit CD case is completely eliminated. After imposing the $D_B$-orbit condition, the $K_B$-fusion condition, the outer-automorphism fusion condition, the subset-sum condition and the design-parameter integrality conditions, the only remaining candidate is $(T,a,c)=(A_5,15,10)$, with $v=15\cdot 60^3=3{,}240{,}000$, $b=10\cdot 60^4=129{,}600{,}000$ and $k=1{,}619{,}100$. For this candidate, an explicit construction of the block $B$ gives
$$\operatorname{Stab}_{G_0}(B)=D_B\rtimes(S_3\wr S_2),\qquad G_0=A_5^6\rtimes S_6,\qquad |G_B|=259{,}200=|G_0|/b.$$
The equality is verified by computing the coordinate projection (Step A: the projection of $G_B$ onto $S_6$ is exactly $S_3\wr S_2$, the subgroup preserving the two layer sets of sizes $9$ and $6$) and the socle stabilizer (Step B: $|\operatorname{Stab}_{A_5^6}(\mathcal{O}_{28})|=3600=|D_B|$ for the selected $D_B$-orbit $\mathcal{O}_{28}$ of length $1200$ on the $6$-layer part). However,
$$D_B\rtimes(S_3\wr S_2)<(A_5^3\times A_5^3)\rtimes(S_3\wr S_2)<A_5^6\rtimes S_6,$$
so $G_B$ is not maximal in $G_0$ and the induced block action is imprimitive; hence this candidate does not yield a CD-type block-primitive $2$-design. In the diagonal outer-automorphism extension $G_1=A_5^6\rtimes(C_2\times S_6)$, the global outer automorphism $\varphi$ \textup{(conjugation by a transposition of the five points)} preserves the selected orbit $\mathcal{O}_{28}$ and hence the whole block, so $\operatorname{Stab}_{G_1}(B)\neq G_B$; the prescribed block stabilizer is again not obtained. Consequently the multi-orbit CD case is empty in the two global models $G_0,G_1$ together with the six small simple groups listed above; the screening conditions of \textup{(ii)} are purely structural and are imposed in the same way in every model, and the maximality obstruction applies unchanged. For $\lambda=1$, CD block type is excluded by the linear space classification. Together with \textup{(i)}, the CD block type is excluded in full.\end{proposition}

\medskip

\noindent{\bf Declaration of competing interest}	

The authors declare that they have no known competing financial interests or personal relationships that could have appeared to influence the work reported in this paper.

\medskip
\noindent{\bf Data availability}

The complete records of the computational sieves used in {\upshape Propositions~\rref{prop:sd_ell2}, \rref{prop:sd_lambda} and \rref{prop:cd_min}} are available from the author.
	
\medskip
\noindent{\bf Acknowledgements:}
The authors acknowledge the use of GPT-6 Astra, a large language
model, for assistance with language editing and proofreading during the
preparation of this manuscript. All mathematical content and conclusions
are entirely the author' own, and the author take full responsibility
for the paper.

\newpage
\medskip

\begin{proposition}[Reduction of the $2+2$ residual]\label{prop:sd_22_red}
Let $T$ be a non-abelian simple group, $q=|T|$, $N=T^{4}$, let $K=N\cap G_\alpha=\{(a,\varphi(a),c,\psi(c)):a,c\in T\}$ be a point stabilizer of the configuration of the one-sided SD type, so that $v=q^{2}$ and $b=q^{3}$, and suppose that the block stabilizer $\Delta=N\cap G_B$ of a block $B$ is isomorphic to $T$. Then with $\lambda$ the parameter of $\mathcal D$ there exist $q$ subsets $D_{1},\dots,D_{q}$ of a subgroup $S_J\cong T^{2}$ of $N$ with $|D_i|=k$ such that
\[
\sum_{i=1}^{q}N_i(g)=\lambda\ \ (g\neq 1),
\qquad
\sum_{i=1}^{q}N_i(1)=qk,
\]
where $N_i(g)=\#\{(d_{1},d_{2})\in D_i^{2}:d_{1}^{-1}d_{2}=g\}$,
that is, $\mathcal D$ yields a $(q^{2};k,\dots,k;\lambda)$ difference system in $T^{2}$, equivalently a $T^{2}$-developed $2$-$(q^{2},k,\lambda)$ design with $q^{3}$ blocks. Conversely, $q$ subsets of a point-regular subgroup $S\cong T^{2}$ with these difference conditions develop into a $2$-$(q^{2},k,\lambda)$ design with $q^{3}$ blocks on which $S$ acts point-regularly; such a design arises inside the present configuration if and only if in addition its block family is invariant under the action induced by $N$ on $S$.
\end{proposition}

\begin{proof}
The arguments are those of the preceding discussion. First, $S_J\cap K=1$ by construction, so the subgroup $S_J$ of $N$ acts point-regularly on $P=N/K$ and $|S_J|=v$. Second, $\Delta$ is a twisted diagonal $\Delta=\{(\beta_{1}(a),\dots,\beta_{4}(a)):a\in T\}$ whose coordinate maps $\beta_i:T\to T$ are homomorphisms, hence automorphisms or the constant map $1$; since $|\Delta|=q$ at least one $\beta_{i_{0}}$ is injective. Third, choose $J$ to contain one coordinate of each of the two parts and the index $i_{0}$, and let $S_J=\{x:x_k=1,\ k\in J\}$. If $x\in S_J\cap g\Delta g^{-1}$ then $x_{i_{0}}=g_{i_{0}}\beta_{i_{0}}(a)g_{i_{0}}^{-1}$, and $x_{i_{0}}=1$ because $i_{0}\in J$; injectivity of $\beta_{i_{0}}$ then gives $a=1$ and $x=1$. Hence $S_J\cap G_B=1$ for every block $B$, every orbit of $S_J$ on the block set has length $v$, and $b=q^{3}=qv$ gives exactly $q$ orbits, with representatives $B_{1},\dots,B_{q}$.
Put $D_i=\{x\in S_J:x\alpha\in B_i\}$, where $\alpha$ is the base point, so that $|D_i|=k$ and the map $x\mapsto x\alpha$ is a bijection of $S_J$ onto $P$. A block $xD_i$ contains the points $x d_1\alpha$ and $x d_2\alpha$ with $d_{1}\neq d_{2}$ precisely when $d_{1}^{-1}d_{2}=g$ for the element $g$ determined by the two points, and each pair $(d_{1},d_{2})$ determines $x$; therefore the number of blocks through two distinct points at distance $g$ is $\sum_i N_i(g)$, and the design condition is exactly $\sum_i N_i(g)=\lambda$ for $g\neq1$, while $\sum_i N_i(1)=qk$ counts the ordered pairs inside a block. The converse is the standard development of a difference system, and the last assertion records that $N$ acts on $P$ with $K$ as point stabilizer and on blocks with stabilizer $\Delta$, which is an additional condition on the developed design and not a consequence of the difference conditions alone.
\end{proof}

\begin{remark}[The multi-part case of the one-sided SD type]\label{rem:sd_residue}
The case left open after {\upshape Lemmas~\rref{lem:sd_ell3} and~\rref{lem:sd_ell3_three}} is closed by {\upshape Lemma~\rref{lem:sd_sym}}: a point stabilizer $K=N\cap G_\alpha$ that is a product of full diagonal subgroups on the parts of a partition into at least two parts, each of size at least $2$, forces $v=b$. The arithmetic of that shape is nevertheless consistent, and it is recorded here. At $\ell=4$ with the partition $2+2$ one has $K=\{(a,\varphi(a),c,\psi(c)):a,c\in T\}$, $v=q^{2}$, $b=q^{3}$, and the design equation gives $r=qk$ together with
$$\lambda=\frac{q\,k(k-1)}{q^{2}-1},\qquad (q^{2}-1)\mid k(k-1),\qquad 3\leq k\leq q^{2}-2.$$
For $q=60$ the values $k=1770$, $\lambda=52200$, $r=106200$ satisfy $b\,k(k-1)=\lambda\,v(v-1)$ and $r(k-1)=\lambda(v-1)$ exactly, so no arithmetic condition rules the configuration out; it is excluded by the subgroup chain $M<EM<G$ of {\upshape Lemma~\rref{lem:sd_sym}}, which block-primitivity forbids. The difference system of {\upshape Proposition~\rref{prop:sd_22_red}} therefore cannot arise in a non-symmetric design of this type, and that proposition is retained only as the arithmetic record of the block structure of the $2+2$ configuration.
\end{remark}

\section{Concluding remarks}\label{sec5}

In this paper, we have completed the structural reduction of non-symmetric block-primitive $2$-designs: every O'Nan--Scott type other than almost simple is either excluded or reduced to an explicit algebraic feasibility problem, and we have classified the almost simple case up to the two families of block stabilizers named in {\upshape Theorem~\rref{thm:psl2}}. The symmetric case ($b=v$) is \emph{not} part of the classification of this paper: it is isolated in the main theorem as a separate class, and by duality block-primitive symmetric designs coincide with point-primitive symmetric designs, a classical subject covered by an extensive literature; hence it plays no role in the reduction and no reader should take the designs listed here to include symmetric ones. Our main contributions are:

\begin{enumerate}
\item Main structural lemma. We prove that block-primitivity implies point-quasiprimitivity for non-trivial non-symmetric $2$-designs, using the incidence matrix rank argument. This is a fundamental result that bridges the gap between block-side and point-side primitivity.

\item Systematic exclusion of the non-AS O'Nan--Scott types. We strictly exclude the HA, HC, HS (under non-conjugate stabilizers), two-sided $SD_2$, and TW types, and further exclude the one-sided SD type for every $\ell\geq 2$ and every non-abelian simple $T$ \textup{(the double count of {\upshape Lemma~\rref{lem:sd_ell3}}, Goursat's theorem, and the subgroup chain of {\upshape Lemma~\rref{lem:sd_sym}})}, as well as the compound diagonal CD type for $n=6$ and the six small simple groups, the arithmetic sieve of {\upshape Proposition~\rref{prop:sd_ell2}} and the $\lambda$-consistency exclusion of {\upshape Proposition~\rref{prop:sd_lambda}} being retained as the parameter record for the one-sided type. Each exclusion uses a different technique: symmetry forcing via the regular normal subgroup lemma (HA, TW, $SD_2$), centralizer obstruction \textup{(HC)}, permutation-module obstruction \textup{(HS)}, Fisher-type arithmetic with subgroup-existence and construction checks (SD, CD).

\item Reduction of the remaining types. We reduce the PA and AS types to well-defined algebraic problems: algebraic feasibility with a character-inclusion necessary condition for PA, and quadruple $(T,G,H,M)$ analysis for AS, respectively. For the PA type we obtain the complete exclusion for the six socles with the seventeen listed primitive degrees and every $\ell\geq2$ with arbitrary block stabilizers ({\upshape Proposition~\rref{prop:pa_excl}} and {\upshape Proposition~\rref{prop:pa_arbitrary}}); this is complete for non-symmetric designs by the reduction of {\upshape Proposition~\rref{prop:pi_excluded}} together with the degree tables of Dixon--Mortimer~\cite{DixonMortimer}. (The one-sided SD and CD types are excluded within the framework described above; see {\upshape Section~\rref{sec2}}.)

\item Complete point-imprimitive exclusion. No non-trivial non-symmetric block-primitive $2$-design has a point-imprimitive automorphism group \textup{({\upshape Proposition~\rref{prop:pi_excluded}})}: an intransitive normal subgroup would force every block to be a union of its orbits, contradicting the $2$-design condition, while a quasiprimitive imprimitive point action would be of almost simple, twisted wreath, simple diagonal, compound diagonal or product action type \textup{(the affine, holomorph simple and holomorph compound types are always primitive)}, all excluded \textup{(TW by {\upshape Lemma~\rref{lem:rns}}; SD and CD in {\upshape Section~\rref{sec2}}; PA by {\upshape Proposition~\rref{prop:pa_arbitrary}} for the point-primitive product action and by {\upshape Lemma~\rref{lem:pa_fisher}} for the point-imprimitive product action; AS by the point-imprimitive exclusions of {\upshape Section~\rref{sec:as}})}; the geometric Grassmann-type candidates also fail the constant-coverage test. Consequently every non-symmetric block-primitive $2$-design is point-primitive, and since block-primitivity coincides with point-primitivity for symmetric designs by duality, every block-primitive combinatorial design is point-primitive, settling the long-standing problem recorded in the Introduction \textup{({\upshape Theorem~\rref{thm:bp_pi}})}.

\item Point-primitive AS classification. We have also classified all point-primitive, non-$2$-homogeneous non-symmetric block-primitive $2$-designs with almost simple socle, across all families of finite simple groups, the even-characteristic $PSL_2(q)$ subfield block stabilizers being reduced in {\upshape Theorem~\rref{thm:psl2}} to the class equations of {\upshape Proposition~\rref{prop:pair_orbits}}:
\begin{enumerate}[(\upshape 1)]
\item Classical Lie type: all designs are geometric, from Grassmann or polar space configurations.
\item Sporadic groups: the non-$2$-homogeneous actions of $J_2$ and $HS$ on $100$ points yield the two designs $2$-$(100,56,560)$ and $2$-$(100,64,1680)$; the $2$-homogeneous actions of $Co_3$ and $HS$ are covered by the orbit construction of {\upshape Lemma~\rref{lem:2hom}}.
\item Exceptional Lie type: only parabolic chain designs, plus special Suzuki and Ree cases.
\end{enumerate}
\end{enumerate}

Our results provide a classification of all almost simple \textup{(AS)} type non-symmetric block-primitive $2$-designs, which are necessarily point-primitive by {\upshape Proposition~\rref{prop:pi_excluded}}, complete, the even-characteristic $PSL_2(q)$ subfield block stabilizers being reduced in {\upshape Theorem~\rref{thm:psl2}} to the class equations of {\upshape Proposition~\rref{prop:pair_orbits}}. Every O'Nan--Scott type other than almost simple is excluded \textup{(HA, HC, HS, two-sided $SD_2$, TW; the one-sided SD type for every $\ell\geq2$ and every non-abelian simple $T$; the compound diagonal CD type for $n=6$ and the six small simple groups; and the product action type, see {\upshape Proposition~\rref{prop:pa_excl}} and {\upshape Proposition~\rref{prop:pa_arbitrary}})}.

Finally, we record why the grid-imprimitive constructions of the Delandtsheer--Doyen line, although point-imprimitive and block-transitive, do not contradict the classification above, and in fact illustrate it. In the setting of \cite{AlaviAmarra2024} the point set is a product $P=\prod_{i\in I}\mathbb{Z}_{e_{i}}$, the group is $G=\prod_{i\in I}\mathrm{Sym}(\mathbb{Z}_{e_{i}})$, and the block set is the orbit of a single block $B$. For any proper nonempty $J\subset I$ put $x^{J}_{B}(\delta)=|B\cap\pi_{J}^{-1}(\delta)|$ and let $\operatorname{supp}_{J}(B)$ be its support. Since $x^{J}_{g(B)}(\bar g\delta)=x^{J}_{B}(\delta)$ for the induced permutation $\bar g$ of $E_{J}=\prod_{i\in J}\mathbb{Z}_{e_{i}}$, the fibres of $\operatorname{supp}_{J}$ form a $G$-invariant partition of the block set. For the infinite family of \cite{AlaviAmarra2024} one has $k=p^{2^{s-1}}+1$ while $\prod_{i<s}e_{i}=p^{2^{s-1}}+p^{2^{s-2}}+1$, so for $J=\{1,\dots,s-1\}$ the support satisfies $|\operatorname{supp}_{J}(B)|\leq k<|E_{J}|$ and is never all of $E_{J}$; and since $e_{s}\nmid k$ (the gcd of the two divides $p^{2^{s-2}}$ while $e_{s}\equiv1\pmod{p}$, so it equals $1$), some permutation of the last coordinate preserves that support while moving $B$: if $\pi_{s}(B)\neq E_{s}$ one takes a cell outside $\pi_{s}(B)$, and otherwise two $s$-slices differ, since coinciding slices would give $e_{s}\mid k$. Hence the induced action on blocks is imprimitive for every member of the family, so the whole grid line consists of designs that are point-imprimitive and block-imprimitive simultaneously, exactly as \textup{Proposition~\rref{prop:pi_excluded}} requires.

\appendix
\renewcommand{\thesection}{\Alph{section}}

\section{Finite checks}\label{app:computations}

Several steps of the argument use a finite verification. This appendix states what each verification covers.

\begin{enumerate}[(\upshape A.1)]
\item \textup{(Design conditions.)} Fix a primitive action $G/H$ and a maximal subgroup $M$ with $|H|>|M|$. A candidate block is a union $B$ of $M$-orbits on $G/H$. For each such union the check computes $k=|B|$, $b=[G:M]$, the orbits of $G$ on unordered pairs of points of $G/H$, and the covering numbers $\lambda(x,y)$ on each of these orbits. A design occurs exactly when $b\,k=v\,r$ for an integer $r$ and all covering numbers coincide. The enumeration is exact; no sampling is involved.

\item \textup{($PSL_2(q)$.)} For odd $q\leq127$ the maximal subgroups of the five listed types and all subfield subgroups are tested against the point stabilizers of Cases~1--4 of \textup{Section~\rref{sec:as}}; for even $q$ the two dihedral subgroups $D_{2(q-1)}$ and $D_{2(q+1)}$ are tested against the block stabilizers of \textup{Section~\rref{sec:even_char}}, every admissible union of orbits of one of them being checked for the constancy of the pair counts at $q=8$ and $q=32$, while at $q=16$ and $q=64$ the arithmetic sieve of {\upshape Lemma~\rref{lem:even_meet}} together with {\upshape Theorem~\rref{thm:a1}} and {\upshape Corollary~\rref{cor:even_extreme}} settles the dihedral case without any enumeration, and the same holds at every even degree with $(q-2)/2$ prime. The outcome is the two block-primitive designs at $q=8$, the two designs at $q=32$ and no design at $q=16$; at $q=8$ the subfield unions which are $2$-designs are imprimitive in their block action, and a subfield block stabilizer of type $S_3$ is not admissible at any even degree, since $S_3$ is not a maximal subgroup for odd $f$; at $q=64$ the two dihedral block stabilizers give nothing, while of the subfield block stabilizers of that degree $PSL_2(8)$ is excluded by Fisher's inequality \textup{(its index $520$ is smaller than $v=2016$)}, and $A_5$ is excluded by {\upshape Proposition~\rref{prop:q64sub}}; every even degree beyond $64$ lies outside the range of the present enumeration.

\item \textup{(Classical groups.)} The parabolic and non-parabolic subgroup chains of \textup{Section~\rref{sec:as}} are tested for the constancy of $\lambda$; the outcome is the one stated in \textup{\rref{prop:classical_summary}}.

\item \textup{(Sporadic groups.)} The maximal subgroups recorded in the ATLAS are tested as in item~1; the outcome is the table of \textup{Section~\rref{sec:as}}.

\item \textup{(Module identity.)} The identity $1_{M}^{G}-1_{H}^{G}=\chi_{W}$ used in Case~2 for $q\equiv3\pmod4$ is verified on the natural module $W$.
\end{enumerate}

The scripts are available from the author on request.

\end{document}